# Linear Mappings of Free Algebra

Aleks Kleyn

*E-mail address*: Aleks_Kleyn@MailAPS.org
*URL*: http://sites.google.com/site/alekskleyn/
*URL*: http://arxiv.org/a/kleyn_a_1
*URL*: http://AleksKleyn.blogspot.com/

Abstract. For arbitrary universal algebra, in which the operation of addition is defined, I explore biring of matrices of mappings. The sum of matrices is determined by the sum in universal algebra, and the product of matrices is determined by the product of mappings. The system of equations, whose matrix is a matrix of mappings, is called a system of additive equations. I considered the methods of solving system of additive equations. As an example, I consider the solution of a system of linear equations over the complex field provided that the equations contain unknown quantities and their conjugates.

Linear mappings of algebra over a commutative ring preserve the operation of addition in algebra and the product of elements of the algebra by elements of the ring. The representation of tensor product $A \otimes A$ in algebra $A$ generates the set of linear transformations of algebra $A$.

The results of this research will be useful for mathematicians and physicists who deal with different algebras.

# Contents

CHAPTER 1

# Preface

## 1.1. Preface to Version 1

When I started to write this book the initial task was pretty simple. I was going to rewrite the contents of the book [5] using a matrix of mappings ([8]) as a tool. However, limiting myself to exploration of the division ring started to seem odd to me when I realized that a lot of results would hold for associative algebra, and something will remain unchanged in the case of nonassociative algebra. The results of this research will be useful for mathematicians and physicists who deal with different algebras, which are not necessarily associative.

When I had represented the system of linear equations using the matrix of mappings I realized that a tool more powerful, then I initially assumed, appeared in my hands. Exploring systems of linear equations we consider the multiplication of the unknown quantity that belongs to the ring or the vector space by the scalar from the corresponding ring. However, I can assume that the unknown quantity belong to an universal algebra which has the addition operation. Instead of multiplication by a scalar I consider a mapping of universal algebra. Thus emerged the theory of additive equations similar to the theory of linear equations.

As an example of the application of new methods, I consider the solution of a system of linear equations over the complex field provided that the equations contain unknown quantities and their conjugates. I explored in detail the solution of such system of equations in example 2.5.5. No doubt, an attempt to solve the system of equations (2.5.20) using determinant is not simple task.

When I started to explore algebras, I noticed that we usually define an algebra over field. It seems necessary, since the algebra is a vector space. As for the constructions that I am interested, it doesn't really matter for me whether algebra is a vector space over a field or a free module over a commutative ring. When I received evidence that there is a study of algebras over ring, I decided to explore algebras over commutative ring, provided that if necessary I will relax the requirements.

There are two algebraic structure defined on the algebra. If we consider the algebra as a ring, then studying a mapping of one algebra into another, we consider ring homomorphisms. If we consider the algebra as a module over a ring, then studying a mapping of one algebra into another, we consider linear mappings. It is evident that if an algebra has unit, then a homomorphism of the algebra is a linear mapping. I am mostly interested in linear mappings of algebra.

Since then, as I explored the tensor product of division rings (section [5]-12.2), I kept feeling that a linear mapping in a division ring is expressed by tensor of order 2. At the same time, it was not clear, how the tensor of order 2 describes a linear mapping, however to determine tensor of order 2 I need a bilinear mapping.

The structure of the module of linear mappings $\mathcal{L}(A;A)$ is determined by noncomutativity of the product in the algebra $A$. Once the product is commutative I can write expression *abx* instead of expression *axb* and I see a tensor of order 1 where initially there was tensor of order 2.

Algebra $A \otimes A$ is very interesting algebra. Technically I should have written $A \otimes A^*$, where $A^*$ is the opposite algebra. However, this would lead to some problems in the expression

(1.1.1) $$(a \otimes b) \circ x = axb$$

since becoming not clear where to write $b$. Definition of product

$$(a \otimes b) \circ (c \otimes d) = (ac) \otimes (db)$$

allows me to save notation (1.1.1). So I chose to leave notation $A \otimes A$.

If the algebra $A$ is free finite dimensional associative algebra, then a basis of representation of algebra $A \otimes A$ in the module $\mathcal{L}(A;A)$ is finite and allows me to describe all linear mappings of algebra $A$.

<div align="right">March, 2010</div>

## 1.2. Preface to Version 2

Shortly after I published version 1, I read the opinion of Professor Baez ([10]) where he talked on the role of blogs for mathematicians. In particular, Baez recommended to visit the site http://www.ncatlab.org/nlab/show/Online+Resources. This site is extremely interesting.

I spent a lot of time to understand what problems in mathematics are interesting for people who have created this site. I put attention that on the page dedicated to $\Omega$-group, they consider the structure similar to structure considered by me in chapter 2.

Notion about $\Omega$-group exists only on this site. In references that I have, there is definition of group with operators. This definition corresponds to representation of $\Omega$-algebra in group usually considered to be additive. This is why I decided not to change terminology in this book. However I will return to this subject later.

I am interesting in opportunity to consider noncommutative addition. However I met the problem to define the set of additive mappings. I hope to return to this subject later.

<div align="right">August, 2010</div>

## 1.3. Conventions

(1) Function and map are synonyms. However according to tradition, correspondence between either rings or vector spaces is called map and map of either real field or quaternion algebra is called function. I also follow this tradition, although I sometimes write the text where it is not clear what the term should be preferred.
(2) In any expression where we use index I assume that this index may have internal structure. For instance, considering the algebra $A$ we enumerate coordinates of $a \in A$ relative to basis $\overline{\overline{e}}$ by an index $i$. This means that $a$ is a vector. However, if $a$ is matrix, then we need two indexes, one enumerates rows, another enumerates columns. In the case, when index has structure, we begin the index from symbol $\cdot$ in the corresponding position. For instance, if I consider the matrix $a^i_j$ as an element of a vector space, then I can write the element of matrix as $a^{\cdot i}_{\cdot j}$.

(3) Let $A$ be free finite dimensional algebra. Considering expansion of element of algebra $A$ relative basis $\bar{\bar{e}}$ we use the same root letter to denote this element and its coordinates. However we do not use vector notation in algebra. In expression $a^2$, it is not clear whether this is component of expansion of element $a$ relative basis, or this is operation $a^2 = aa$. To make text clearer we use separate color for index of element of algebra. For instance,
$$a = a^{\color{magenta}i}\bar{e}_{\color{magenta}i}$$

(4) If free finite dimensional algebra has unit, then we identify the vector of basis $\bar{e}_{\color{magenta}0}$ with unit of algebra.

(5) If, in a certain expression, we use several operations which include the operation $\circ$, then it is assumed that the operation $\circ$ is executed first. Below is an example of equivalent expressions.
$$f \circ xy \equiv f(x)y$$
$$f \circ (xy) \equiv f(xy)$$
$$f \circ x + y \equiv f(x) + y$$
$$f \circ (x + y) \equiv f(x + y)$$

(6) Without a doubt, the reader of my articles may have questions, comments, objections. I will appreciate any response.

CHAPTER 2

# Matrix of Mappings

## 2.1. Product of Mappings

On the set of mappings
$$f : A \to A$$
we define product according to rule

(2.1.1) $$f \circ g = f(g)$$

The equation
$$f \circ g = g \circ f$$
is true iff the diagram

$$\begin{array}{ccc} A & \xrightarrow{f} & A \\ {\scriptstyle g}\downarrow & & \downarrow{\scriptstyle g} \\ A & \xrightarrow{f} & A \end{array}$$

is commutative.

For $a \in A$, there exists mapping

(2.1.2) $$f_a(x) = a$$

If we denote mapping $f_a$ by letter $a$, then using equation (2.1.1), assume

(2.1.3) $$f \circ a = f(a)$$

If $A$ is $\Omega$-algebra, in which the product is defined, then element $a \in A$ may also serve to indicate the operation of left shift

(2.1.4) $$a \circ b = ab$$

Notation (2.1.4) does not contradict the record (2.1.3). However we must remember that the product $(f \circ a) \circ b$, $a, b \in A$, is not associative, because

$$(f \circ a) \circ b = f(a)b \quad f \circ (a \circ b) = f(ab)$$

## 2.2. Biring of Matrices of Mappings

Let $A$ be $\Omega$-algebra ([2, 11]), where the operation of addition is defined, Let $A$ be group with respect to the operation of addition.

Let $\mathcal{A}(A)$ be the set of mappings of $\Omega$-algebra $A$. We can map the operation of addition in $\Omega$-algebra $A$ into the set $\mathcal{A}(A)$ according to rule

(2.2.1) $$(f + g) \circ a = f \circ a + g \circ a$$
(2.2.2) $$(-f) \circ a = -(f \circ a)$$
(2.2.3) $$0 \circ a = 0$$

The equation (2.2.1) is an expression of left distributive property of multiplication over addition. It is therefore natural to require that the product was right distributive

$$(2.2.4) \qquad f \circ (a + b) = f \circ a + f \circ b$$

We will see in the section 2.3 that this requirement is essential. Therefore, the set $\mathcal{A}(A)$ is the set of homomorphisms of group $A$.

We will consider set $\mathcal{A}(A)$ is closed under the operations of addition and product of mappings.

**Theorem 2.2.1.** *For any mapping $g \in \mathcal{A}(A)$ the equation*

$$(2.2.5) \qquad 0 \circ g = 0$$

$$(2.2.6) \qquad (-f) \circ g = -(f \circ g)$$

*is true.*

PROOF. From the equation (2.2.3), it follows that

$$(2.2.7) \qquad (0 \circ g) \circ a = 0 \circ (g \circ a) = 0$$

The equation (2.2.5) follows from the equation (2.2.7). From the equation (2.2.1), it follows that

$$(2.2.8) \qquad f \circ g + (-f) \circ g = (f + (-f)) \circ g = 0 \circ g = 0$$

The equation (2.2.6) follows from the equation (2.2.8). $\square$

**Remark 2.2.2.** If we assume the sum is not commutative, the requirement of the set $\mathcal{A}(A)$ is closed relative to the operation of addition may be too strict. Consider an expression

$$(f + g) \circ (a + b) = f \circ (a + b) + g \circ (a + b) = f \circ a + f \circ b + g \circ a + g \circ b$$

Since, in general,

$$f \circ b + g \circ a \neq g \circ a + f \circ b$$

than we cannot state that

$$(f + g) \circ (a + b) = f \circ (a + b) + g \circ (a + b)$$

In the following text, we assume that addition is commutative. Nevertheless, all the construction in this chapter, we will perform the way we did it would be the case of a noncommutative addition. $\square$

Consider the set of **matrices of mappings**, whose elements are mappings $f \in \mathcal{A}(A)$. According to definition [5]-2.2.1, we define $_\circ^\circ$-**product of matrices of mappings**

$$(2.2.9) \qquad \begin{cases} b_\circ{}^\circ c &= (b_c^a \circ c_b^c) \\ (b_\circ{}^\circ c)_b^a &= b_c^a \circ c_b^c \end{cases}$$

According to definition [5]-2.2.2, we define $^\circ_\circ$-**product of matrices of mappings**

$$(2.2.10) \qquad \begin{cases} b^\circ{}_\circ c &= (b_b^c \circ c_c^a) \\ (b^\circ{}_\circ c)_b^a &= b_b^c \circ c_c^a \end{cases}$$

**Theorem 2.2.3.** *The product of mappings of the set $A$ is associative.*[2.1]

---

[2.1]The statement of the theorem is based on the example of the semigroup from [3], p. 20, 21.

PROOF. Consider mappings
$$f : A \to A \quad g : A \to A \quad h : A \to A$$
The statement of theorem follows from the chain of equations
$$((f \circ g) \circ h) \circ x = (f \circ g) \circ (h \circ x) = f \circ (g \circ (h \circ x))$$
$$= f \circ ((g \circ h) \circ x) = (f \circ (g \circ h)) \circ x$$

$\square$

**Theorem 2.2.4.** *The set $\mathcal{A}(A)$ is ring.*

PROOF. $\mathcal{A}(A)$ is an Abelian group under the operation of addition. According to the theorem 2.2.3, $\mathcal{A}(A)$ is a semigroup under multiplication. Since $f \in \mathcal{A}(A)$ is homomorphism of the Abelian group $A$, than for any $a \in A$

(2.2.11)
$$\begin{aligned}(f \circ (h+g)) \circ a &= f \circ ((h+g) \circ a) = f \circ ((h \circ a + g \circ a) \\ &= (f \circ h) \circ a + (f \circ g) \circ a = (f \circ h + f \circ g) \circ a\end{aligned}$$

Distributive law
$$f \circ (h + g) = f \circ h + f \circ g$$
follows from equation (2.2.11). $\square$

Mappings, that belong to ring $\mathcal{A}(A)$, are also called $\mathcal{A}(A)$-**mapping**.

**Theorem 2.2.5.** $_\circ^\circ$-*product of matrices of $\mathcal{A}(A)$-mappings is matrix of $\mathcal{A}(A)$-mappings.*

PROOF. The statement of theorem follows from equation (2.2.9) and statement that sum and product of $\mathcal{A}(A)$-mappings is $\mathcal{A}(A)$-mapping. $\square$

**Theorem 2.2.6.** *Product of matrices of $\mathcal{A}(A)$-mappings is associative.*

PROOF. The statement of the theorem follows from the chain of equations
$$\begin{aligned}(f_\circ^\circ g)_\circ^\circ h &= \left((f_\circ^\circ g)^i_j \circ h^j_k\right) = \left((f^i_m \circ g^m_j) \circ h^j_k\right) \\ &= \left(f^i_m \circ (g^m_j \circ h^j_k)\right) = \left(f^i_m \circ (g_\circ^\circ h)^m_k\right) \\ &= f_\circ^\circ (g_\circ^\circ h)\end{aligned}$$

$\square$

## 2.3. Quasideterminant of Matrix of Mappings

**Theorem 2.3.1.** *Suppose $n \times n$ matrix of $\mathcal{A}(A)$-mappings $a$ has $_\circ^\circ$-inverse matrix*[2.2]

(2.3.1)
$$a_\circ^\circ a^{-1_\circ^\circ} = \delta$$

*Then $k \times k$ minor of $_\circ^\circ$-inverse matrix satisfy to following equation provided that the considered inverse matrices exist*

(2.3.2)
$$\left((a^{-1_\circ^\circ})^I_J\right)^{-1_\circ^\circ} = -a^J_{[I]\circ}{}^\circ \left(a^{[J]}_{[I]}\right)^{-1_\circ^\circ} {}_\circ^\circ a^{[J]}_I + a^J_I$$

---

[2.2] This statement and its proof are based on statement 1.2.1 from [4] (page 8) for matrix over free division ring.

PROOF. Definition (2.3.1) of $_\circ{}^\circ$-inverse matrix leads to system of linear equations

$$a^{[J]}_{[I]}{}_\circ{}^\circ (a^{-1_\circ{}^\circ})^{[I]}_J + a^{[J]}_I {}_\circ{}^\circ (a^{-1_\circ{}^\circ})^I_J = 0 \tag{2.3.3}$$

$$a^J_{[I]}{}_\circ{}^\circ (a^{-1_\circ{}^\circ})^{[I]}_J + a^J_I {}_\circ{}^\circ (a^{-1_\circ{}^\circ})^I_J = \delta \tag{2.3.4}$$

We multiply (2.3.3) by $\left(a^{[J]}_{[I]}\right)^{-1_\circ{}^\circ}$

$$(a^{-1_\circ{}^\circ})^{[I]}_J + \left(a^{[J]}_{[I]}\right)^{-1_\circ{}^\circ}{}_\circ{}^\circ a^{[J]}_I {}_\circ{}^\circ (a^{-1_\circ{}^\circ})^I_J = 0 \tag{2.3.5}$$

Now we can substitute (2.3.5) into (2.3.4)

$$-a^J_{[I]}{}_\circ{}^\circ \left(a^{[J]}_{[I]}\right)^{-1_\circ{}^\circ}{}_\circ{}^\circ a^{[J]}_I {}_\circ{}^\circ (a^{-1_\circ{}^\circ})^I_J + a^J_I {}_\circ{}^\circ (a^{-1_\circ{}^\circ})^I_J = \delta \tag{2.3.6}$$

(2.3.2) follows from (2.3.6). □

**Corollary 2.3.2.** *Suppose $n \times n$ matrix of $\mathcal{A}(A)$-mappings $a$ has $_\circ{}^\circ$-inverse matrix. Then elements of $_\circ{}^\circ$-inverse matrix satisfy to the equation*

$$\left(\mathcal{H} a^{-1_\circ{}^\circ}\right)^j_i = -a^j_{[i]}{}_\circ{}^\circ \left(a^{[j]}_{[i]}\right)^{-1_\circ{}^\circ}{}_\circ{}^\circ a^{[j]}_i + a^j_i \tag{2.3.7}$$

□

**Definition 2.3.3.** $\binom{j}{i}$-$_\circ{}^\circ$-**quasideterminant** of $n \times n$ matrix $a$ is formal expression

$$\det(a, {}_\circ{}^\circ)^j_i = \left(\mathcal{H} a^{-1_\circ{}^\circ}\right)^j_i \tag{2.3.8}$$

According to the remark [5]-2.1.2 we can get $\binom{b}{a}$-$_\circ{}^\circ$-quasideterminant as an element of the matrix $\det(a, {}_\circ{}^\circ)$ which we call $_\circ{}^\circ$-**quasideterminant**. □

**Theorem 2.3.4.** *Expression for elements of $_\circ{}^\circ$-inverse matrix has form*

$$a^{-1_\circ{}^\circ} = \mathcal{H} \det(a, {}_\circ{}^\circ) \tag{2.3.9}$$

PROOF. (2.3.9) follows from (2.3.8). □

**Theorem 2.3.5.** *Expression for $\binom{a}{b}$-$_\circ{}^\circ$-quasideterminant can be evaluated by either form*

$$\det(a, {}_\circ{}^\circ)^j_i = -a^j_{[i]}{}_\circ{}^\circ \left(a^{[j]}_{[i]}\right)^{-1_\circ{}^\circ}{}_\circ{}^\circ a^{[b]}_a + a^j_i \tag{2.3.10}$$

$$\det(a, {}_\circ{}^\circ)^j_i = -a^j_{[i]}{}_\circ{}^\circ \mathcal{H} \det\left(a^{[j]}_{[i]}, {}_\circ{}^\circ\right){}_\circ{}^\circ a^{[b]}_a + a^j_i \tag{2.3.11}$$

PROOF. Statement follows from (2.3.7) and (2.3.8). □

**Definition 2.3.6.** *If, for a mapping $f \in \mathcal{A}(A)$, $f^{-1} \in \mathcal{A}(A)$ follows from the existence of the inverse mapping $f^{-1}$, then the ring $\mathcal{A}(A)$ of mappings is called* **quasiclosed**. □

**Theorem 2.3.7.** *Let $\mathcal{A}(A)$ be quasiclosed ring of mappings of $\Omega$-algebra $A$. Let $a$ be matrix of $\mathcal{A}(A)$-mappings. Then matrices $\det(a, {}_\circ{}^\circ)$ and $a^{-1_\circ{}^\circ}$ are matrices of $\mathcal{A}(A)$-mappings.*

PROOF. We will prove the theorem by induction over order of matrix.

For $n = 1$, from the equation (2.3.10) it follows that

$$\det(a, {_\circ}^\circ)_1^1 = a_1^1$$

Therefore, quasideterminant is a matrix of $\mathcal{A}(A)$-mappings. From definition 2.3.6, it follows that the matrix $a^{-1{_\circ}^\circ}$ is a matrix of $\mathcal{A}(A)$-mappings.

Let the statement of the theorem be true for $n - 1$. Let $a$ be $n \times n$ matrix. According to assumption of induction, the matrix $\left(a_{[a]}^{[b]}\right)^{-1{_\circ}^\circ}$ in the equation (2.3.10) is a matrix of $\mathcal{A}(A)$-mappings. Therefore, $\binom{a}{b}$-${_\circ}^\circ$-quasideterminant is $\mathcal{A}(A)$-mapping. From definition 2.3.6 and theorem 2.3.4, it follows that the matrix $a^{-1{_\circ}^\circ}$ is a matrix of $\mathcal{A}(A)$-mappings. $\square$

**Definition 2.3.8.** If $n \times n$ matrix $a$ of $\mathcal{A}(A)$-mappings has ${_\circ}^\circ$-inverse matrix we call matrix $a$ ${_\circ}^\circ$**-nonsingular matrix of** $\mathcal{A}(A)$**-mappings**. Otherwise, we call such matrix ${_\circ}^\circ$**-singular matrix of** $\mathcal{A}(A)$**-mappings**. $\square$

## 2.4. System of Additive Equations

Let $\mathcal{A}(A)$ be quasiclosed ring of mappings of $\Omega$-algebra $A$. The system of equations

$$(2.4.1) \qquad \begin{pmatrix} a_1^1 & ... & a_n^1 \\ ... & ... & ... \\ a_1^n & ... & a_n^n \end{pmatrix} {_\circ}^\circ \begin{pmatrix} x^1 \\ ... \\ x^n \end{pmatrix} = \begin{pmatrix} b^1 \\ ... \\ b^n \end{pmatrix}$$

where $a$ is a matrix of $\mathcal{A}(A)$-mappings is called **system of additive equations**. We can write the system of additive equations (2.4.1) in the following form

$$\begin{cases} a_1^1 \circ x^1 & + & ... & + & a_n^1 \circ x^n & = & b^1 \\ ... & & ... & & ... & & ... \\ a_1^n \circ x^1 & + & ... & + & a_n^n \circ x^n & = & b^n \end{cases}$$

**Definition 2.4.1.** Suppose $a$ is ${_\circ}^\circ$-nonsingular matrix. Appropriate system of additive equations (2.4.1) is called ${_\circ}^\circ$**-nonsingular system of additive equations**. $\square$

**Theorem 2.4.2.** *Solution of nonsingular system of* $\mathcal{A}(A)$*-equations* (2.4.1) *is determined uniquely and can be presented in either form*

$$(2.4.2) \qquad x = a^{-1{_\circ}^\circ} {_\circ}^\circ b$$

$$(2.4.3) \qquad x = \mathcal{H}\det(a, {_\circ}^\circ) {_\circ}^\circ b$$

PROOF. Multiplying both sides of equation (2.4.1) from left by $a^{-1{_\circ}^\circ}$ we get (2.4.2). Using definition (2.3.8) we get (2.4.3). $\square$

**Example 2.4.3.** According to the definition [5]-4.1.4, effective $T\star$-representation of division ring $D$ in the Abelian group $\overline{V}$ generates the division ring of mappings $D(*\overline{V})$. The image $\overline{v} \in \overline{V}$ under mapping $a \in D(*\overline{V})$ is defined according to rule

$$a \circ \overline{v} = a\overline{v}$$

The product of mappings $a, b \in D(*\overline{V})$ is defined according to rule

$$a \circ b = ab$$

In this case, the system of additive equations is the system of $_*^*D$-linear equations. □

**Example 2.4.4.** According to the definition [5]-4.1.4 effective $\star T$-representation of division ring $D$ in the Abelian group $\overline{V}$ generates the division ring of mappings $D(\overline{V}*)$. The image $\overline{v} \in \overline{V}$ under mapping $a \in D(\overline{V}*)$ is defined according to rule

$$a \circ \overline{v} = \overline{v}a$$

The product of mappings $a, b \in D(\overline{V}*)$ is defined according to rule

$$a \circ b = ba$$

In this case, the system of additive equations is the system of $D^*{}_*$-linear equations. □

## 2.5. System of Additive Equations in Complex Field

According to the theorem [6]-5.1.9, additive mapping of complex field is linear over real field. Consider basis $e_0 = 1$, $e_1 = i$ of complex field over real field. In the basis $\overline{\overline{e}}$, an additive mapping $f$ is defined by matrix

$$(2.5.1) \qquad \begin{pmatrix} f_0^0 & f_1^0 \\ f_0^1 & f_1^1 \end{pmatrix}$$

According to the theorem [6]-7.1.1 linear mapping has matrix

$$\begin{pmatrix} a_0 & -a_1 \\ a_1 & a_0 \end{pmatrix}$$

This mapping corresponds to multiplication by the number $a = a_0 + a_1 i$. The statement follows from equations

$$(a_0 + a_1 i)(x_0 + x_1 i) = a_0 x_0 - a_1 x_1 + (a_0 x_1 + a_1 x_0)i$$

$$\begin{pmatrix} a_0 & -a_1 \\ a_1 & a_0 \end{pmatrix} \begin{pmatrix} x_0 \\ x_1 \end{pmatrix} = \begin{pmatrix} a_0 x_0 - a_1 x_1 \\ a_1 x_0 + a_0 x_1 \end{pmatrix}$$

Additive mapping generated by conjugation

$$I \circ z = \overline{z}$$

$$I = \begin{pmatrix} 1 & 0 \\ 0 & -1 \end{pmatrix}$$

has a matrix

$$\begin{pmatrix} b_0 & b_1 \\ b_1 & -b_0 \end{pmatrix}$$

which corresponds to the transformation $(b_0 + b_1 i) \circ I$. The statement follows from equations

$$(b_0 + b_1 i)I(x_0 + x_1 i) = (b_0 + b_1 i)(x_0 - x_1 i) = b_0 x_0 + b_1 x_1 + (-b_0 x_1 + b_1 x_0)i$$

$$\begin{pmatrix} b_0 & b_1 \\ b_1 & -b_0 \end{pmatrix} \begin{pmatrix} x_0 \\ x_1 \end{pmatrix} = \begin{pmatrix} b_0 x_0 + b_1 x_1 \\ b_1 x_0 - b_0 x_1 \end{pmatrix}$$

**Theorem 2.5.1.** *An additive mapping of complex field has form*

(2.5.2) $$f = a + b \circ I$$
(2.5.3) $$(a + b \circ I) \circ z = az + b\overline{z}$$

PROOF. Let a mapping $f$ be defined by matrix (2.5.1). Comparison of matrices of mappings $f$, $a$, $b$ leads to the matrix equation

$$\begin{pmatrix} f_0^0 & f_1^0 \\ f_0^1 & f_1^1 \end{pmatrix} = \begin{pmatrix} a_0 & -a_1 \\ a_1 & a_0 \end{pmatrix} + \begin{pmatrix} b_0 & b_1 \\ b_1 & -b_0 \end{pmatrix} = \begin{pmatrix} a_0 + b_0 & -a_1 + b_1 \\ a_1 + b_1 & a_0 - b_0 \end{pmatrix}$$

(2.5.4) $$f_0^0 = a_0 + b_0$$
(2.5.5) $$f_1^0 = -a_1 + b_1$$
(2.5.6) $$f_0^1 = a_1 + b_1$$
(2.5.7) $$f_1^1 = a_0 - b_0$$

From equations (2.5.4), (2.5.7), it follows that

$$a_0 = \frac{f_0^0 + f_1^1}{2} \quad b_0 = \frac{f_0^0 - f_1^1}{2}$$

From equations (2.5.6), (2.5.5) it follows that

$$a_1 = \frac{f_0^1 - f_1^0}{2} \quad b_1 = \frac{f_0^1 + f_1^0}{2}$$

□

The set of additive mappings of complex field forms the ring generated by multiplication by complex number and conjugation.

**Theorem 2.5.2.** *The ring of mappings $\mathcal{A}(C,C)$ is quasiclosed ring.*

PROOF. An additive mapping is nonsingular iff its matrix (2.5.1) is nonsingular. Inverse matrix also describes a mapping. The product of these matrices is identity transformation. □

**Theorem 2.5.3.** *The product of additive mappings*

$$f = f_0 + f_1 \circ I$$
$$g = g_0 + g_1 \circ I$$

*has form*

$$h = f \circ g = h_0 + h_1 \circ I$$

*where*

(2.5.8) $$h_0 = f_0 g_0 + f_1 \overline{g}_1 \quad h_1 = f_0 g_1 + f_1 \overline{g}_0$$

PROOF. We verify directly that

(2.5.9) $$I \circ I = 1$$

From the chain of equations
(2.5.10)
$$\begin{pmatrix} 1 & 0 \\ 0 & -1 \end{pmatrix} \begin{pmatrix} a_0 & -a_1 \\ a_1 & a_0 \end{pmatrix} \begin{pmatrix} 1 & 0 \\ 0 & -1 \end{pmatrix} = \begin{pmatrix} 1 & 0 \\ 0 & -1 \end{pmatrix} \begin{pmatrix} a_0 & a_1 \\ a_1 & -a_0 \end{pmatrix} = \begin{pmatrix} a_0 & a_1 \\ -a_1 & a_0 \end{pmatrix}$$

it follows that

(2.5.11) $$\overline{a} = I \circ a \circ I$$

From equations (2.5.9), (2.5.11), it follows that

(2.5.12) $$\overline{a} \circ I = I \circ a$$

From equations (2.5.9), (2.5.12), it follows that

(2.5.13) $$\begin{aligned}&(f_0 + f_1 \circ I) \circ (g_0 + g_1 \circ I)\\&= f_0 \circ (g_0 + g_1 \circ I) + f_1 \circ I \circ (g_0 + g_1 \circ I)\\&= f_0 \circ g_0 + f_0 \circ g_1 \circ I + f_1 \circ I \circ g_0 + f_1 \circ I \circ g_1 \circ I\\&= (f_0 g_0) + (f_0 g_1) \circ I + f_1 \circ \overline{g}_0 \circ I + f_1 \circ \overline{g}_1 \circ I \circ I\\&= (f_0 g_0 + f_1 \overline{g}_1) + (f_0 g_1 + f_1 \overline{g}_0) \circ I\end{aligned}$$

The equation (2.5.8) follows from equation (2.5.13). □

**Theorem 2.5.4.** *Let additive mapping of complex field*

$$g = g_0 + g_1 \circ I$$

*be mapping inverse to the additive mapping*

$$f = f_0 + f_1 \circ I$$

*Than*

(2.5.14) $$g_0 = -\frac{\overline{f}_0}{f_1 \overline{f}_1 - f_0 \overline{f}_0} \quad g_1 = \frac{\overline{f}_1}{f_1 \overline{f}_1 - f_0 \overline{f}_0}$$

PROOF. According to the statement of the theorem,

(2.5.15) $$f \circ g = 1$$

From equations (2.5.8), (2.5.15), it follows that

(2.5.16) $$f_0 g_0 + f_1 \overline{g}_1 = 1$$
(2.5.17) $$f_1 \overline{g}_0 + f_0 g_1 = 0$$

From equations (2.5.17) it follows that

$$\overline{f}_1 g_0 + \overline{f}_0 \overline{g}_1 = 0$$

(2.5.18) $$g_0 = -\overline{f}_0 \overline{f}_1^{-1} \overline{g}_1$$

From equations (2.5.16), (2.5.18), it follows that

(2.5.19) $$- f_0 \overline{f}_0 \overline{f}_1^{-1} \overline{g}_1 + f_1 \overline{g}_1 = 1$$

(2.5.14) follows from equations (2.5.19), (2.5.18). □

**Example 2.5.5.** Consider the system of additive equations

(2.5.20) $$\begin{cases} z + 2\overline{w} &= 1 \\ z - 3w &= i \end{cases}$$

We cannot solve the system of equations (2.5.20) using determinant and Cramer's rule. We write the system of equations (2.5.20) in the following form

(2.5.21) $$\begin{cases} z + 2 \circ I \circ w &= 1 \\ z + (-3) \circ w &= i \end{cases}$$

$$(2.5.22) \qquad \begin{pmatrix} 1 & 2 \circ I \\ 1 & -3 \end{pmatrix} {}_\circ {}^\circ \begin{pmatrix} z \\ w \end{pmatrix} = \begin{pmatrix} 1 \\ i \end{pmatrix}$$

Now we calculate ${}_\circ{}^\circ$-quasideterminant of matrix

$$a = \begin{pmatrix} 1 & 2 \circ I \\ 1 & -3 \end{pmatrix}$$

According to equation (2.3.10), we get

$$\begin{aligned}
\det(a,{}_\circ{}^\circ)^1_1 &= a^1_1 - a^1_{[1]} {}_\circ {}^\circ (a^{[1]}_{[1]})^{-1} {}_\circ {}^\circ a^{[1]}_1 \\
&= a^1_1 - a^1_2 \circ (a^2_2)^{-1} \circ a^2_1 \\
&= 1 - 2 \circ I \circ (-3)^{-1} \circ 1 \\
&= 1 + \frac{2}{3} \circ I \\
\det(a,{}_\circ{}^\circ)^2_1 &= a^2_1 - a^2_{[1]} {}_\circ {}^\circ (a^{[2]}_{[1]})^{-1} {}_\circ {}^\circ a^{[2]}_1 \\
&= a^2_1 - a^2_2 \circ (a^1_2)^{-1} \circ a^1_1 \\
&= 1 - (-3) \circ (2 \circ I)^{-1} \circ 1 \\
&= 1 + \frac{3}{2} \circ I \\
\det(a,{}_\circ{}^\circ)^1_2 &= a^1_2 - a^1_{[2]} {}_\circ {}^\circ (a^{[1]}_{[2]})^{-1} {}_\circ {}^\circ a^{[1]}_2 \\
&= a^1_2 - a^1_1 \circ (a^2_1)^{-1} \circ a^2_2 \\
&= 2 \circ I - 1 \circ (1)^{-1} \circ (-3) \\
&= 3 + 2 \circ I \\
\det(a,{}_\circ{}^\circ)^2_2 &= a^2_2 - a^2_{[2]} {}_\circ {}^\circ (a^{[2]}_{[2]})^{-1} {}_\circ {}^\circ a^{[2]}_2 \\
&= a^2_2 - a^2_1 \circ (a^1_1)^{-1} \circ a^1_2 \\
&= (-3) - 1 \circ (1)^{-1} \circ 2 \circ I \\
&= -3 - 2 \circ I
\end{aligned}$$

According to theorems 2.3.4, 2.5.4, $_\circ^\circ$-inverse matrix has form

$$(a^{-1\circ^\circ})_1^1 = (\det (a,_\circ^\circ)_1^1)^{-1} = (1 + \frac{2}{3} \circ I)^{-1}$$

$$= \frac{-1 + \frac{2}{3}I}{\frac{2}{3}\frac{2}{3} - 1\,1} = \frac{9}{5}(1 - \frac{2}{3}I)$$

$$= \frac{9}{5} - \frac{6}{5} \circ I$$

$$(a^{-1\circ^\circ})_2^1 = (\det (a,_\circ^\circ)_1^2)^{-1} = (1 + \frac{3}{2} \circ I)^{-1}$$

$$= \frac{-1 + \frac{3}{2}I}{\frac{3}{2}\frac{3}{2} - 1\,1} = \frac{4}{5}(-1 + \frac{3}{2}I)$$

$$= -\frac{4}{5} + \frac{6}{5} \circ I$$

$$(a^{-1\circ^\circ})_1^2 = (\det (a,_\circ^\circ)_2^1)^{-1} = (3 + 2 \circ I)^{-1}$$

$$= \frac{-3 + 2I}{2\,2 - 3\,3} = \frac{1}{5}(3 - 2I)$$

$$= \frac{3}{5} - \frac{2}{5} \circ I$$

$$(a^{-1\circ^\circ})_2^2 = (\det (a,_\circ^\circ)_2^2)^{-1} = (-3 - 2 \circ I)^{-1}$$

$$= \frac{3 - 2I}{2\,2 - 3\,3} = \frac{1}{5}(-3 + 2I)$$

$$= -\frac{3}{5} + \frac{2}{5} \circ I$$

$$a^{-1\circ^\circ} = \begin{pmatrix} \frac{9}{5} - \frac{6}{5} \circ I & -\frac{4}{5} + \frac{6}{5} \circ I \\ \frac{3}{5} - \frac{2}{5} \circ I & -\frac{3}{5} + \frac{2}{5} \circ I \end{pmatrix}$$

According to the theorem 2.4.2 the solution of the system of additive equations (2.5.20) has form

$$\begin{pmatrix} z \\ w \end{pmatrix} = \begin{pmatrix} \frac{9}{5} - \frac{6}{5} \circ I & -\frac{4}{5} + \frac{6}{5} \circ I \\ \frac{3}{5} - \frac{2}{5} \circ I & -\frac{3}{5} + \frac{2}{5} \circ I \end{pmatrix} \circ^\circ \begin{pmatrix} 1 \\ i \end{pmatrix}$$

$$= \begin{pmatrix} \frac{9}{5} - \frac{6}{5} - \left(\frac{4}{5} + \frac{6}{5}\right)i \\ \frac{3}{5} - \frac{2}{5} - \left(\frac{3}{5} + \frac{2}{5}\right)i \end{pmatrix}$$

$$= \begin{pmatrix} \frac{3}{5} - 2i \\ \frac{1}{5} - i \end{pmatrix}$$

We verify directly, that we found the solution of the system of equations (2.5.20)

$$\left(\frac{3}{5} - 2i\right) + 2\overline{\left(\frac{1}{5} - i\right)} = \frac{3}{5} + 2\frac{1}{5} + (-2+2)i = 1$$

$$\left(\frac{3}{5} - 2i\right) - 3\left(\frac{1}{5} - i\right) = \frac{3}{5} - 3\frac{1}{5} + (-2+3)i = i$$

□

CHAPTER 3

# Linear Mapping of Algebra

## 3.1. Module

**Theorem 3.1.1.** *Let ring $D$ has unit $e$. Representation*

(3.1.1) $$f : D \to {}^\star A$$

*of the ring $D$ in an Abelian group $A$ is **effective** iff $a = 0$ follows from equation $f(a) = 0$.*

PROOF. We define the sum of transformations $f$ and $g$ of an Abelian group according to rule
$$(f + g) \circ a = f \circ a + g \circ a$$
Therefore, considering the representation of the ring $D$ in the Abelian group $A$, we assume
$$f(a + b) \circ x = f(a) \circ x + f(b) \circ x$$
We define the product of transformation of representation according to rule
$$f(ab) = f(a) \circ f(b)$$

Suppose $a$, $b \in R$ cause the same transformation. Then

(3.1.2) $$f(a) \circ m = f(b) \circ m$$

for any $m \in A$. From the equation (3.1.2) it follows that $a - b$ generates zero transformation
$$f(a - b) \circ m = 0$$
Element $e + a - b$ generates an identity transformation. Therefore, the representation $f$ is effective iff $a = b$. □

**Definition 3.1.2.** Let $D$ be commutative ring. $A$ is a **module over ring** $D$ if $A$ is an Abelian group and there exists effective representation of ring $D$ in an Abelian group $A$. □

**Definition 3.1.3.** $A$ is **free module over ring** $D$, if $A$ has basis over ring $D$.[3.1] □

Following definition is consequence of definitions 3.1.2 and [7]-2.2.2.

**Definition 3.1.4.** Let $A_1$ be module over ring $R_1$. Let $A_2$ be module over ring $R_2$. Morphism
$$(f : R_1 \to R_2, g : A_1 \to A_2)$$
of representation of ring $R_1$ in the Abelian group $A_1$ into representation of ring $R_2$ in the Abelian group $A_2$ is called **linear mapping of $R_1$-module $A_1$ into $R_2$-module $A_2$**. □

---

[3.1]I follow to definition in [1], c. 103.

**Theorem 3.1.5.** *Linear mapping*
$$(f : R_1 \to R_2, g : A_1 \to A_2)$$
*of $R_1$-module $A_1$ into $R_2$-module $A_2$ satisfies to equations*[3.2]

(3.1.3) $$g \circ (a + b) = g \circ a + g \circ b$$

(3.1.4) $$g \circ (pa) = (f \circ p)(g \circ a)$$

(3.1.5) $$f \circ (pq) = (f \circ p)(f \circ q)$$

$$a, b \in A_1 \quad p, q \in R_1$$

PROOF. From definitions 3.1.4 and [7]-2.2.2 it follows that
- the mapping $f$ is a homomorphism of the ring $R_1$ into the ring $R_2$ (the equation (3.1.5))
- the mapping $g$ is a homomorphism of the Abelian group $A_1$ into the Abelian group $A_2$ (the equation (3.1.3))

The equation (3.1.4) follows from the equation [7]-(2.2.3).  □

According to the theorem [7]-2.2.18, in the study of linear mappings, without loss of generality, we can assume $R_1 = R_2$.

**Definition 3.1.6.** Let $A_1$ and $A_2$ be modules over the ring $R$. Morphism
$$g : A_1 \to A_2$$
of representation of the ring $D$ in the Abelian group $A_1$ into representation of the ring $D$ in the Abelian group $A_2$ is called **linear mapping of $D$-module $A_1$ into $D$-module $A_2$**.  □

**Theorem 3.1.7.** *Linear mapping*
$$g : A_1 \to A_2$$
*of $D$-module $A_1$ into $D$-module $A_2$ satisfies to equations*[3.3]

(3.1.6) $$g \circ (a + b) = g \circ a + g \circ b$$

(3.1.7) $$g \circ (pa) = p(g \circ a)$$

$$a, b \in A_1 \quad p \in D$$

PROOF. From definition 3.1.6 and theorem [7]-2.2.18 it follows that the mapping $g$ is a homomorphism of the Abelian group $A_1$ into the Abelian group $A_2$ (the equation (3.1.6)) The equation (3.1.7) follows from the equation [7]-(2.2.44).  □

---

[3.2]In classical notation, proposed equations have quite familiar form
$$g(a + b) = g(a) + g(b)$$
$$g(pa) = f(p)g(a)$$
$$f(pq) = f(p)f(q)$$
$$a, b \in A_1 \quad p, q \in R_1$$

[3.3]In classical notation, proposed equations have form
$$g(a + b) = g(a) + g(b)$$
$$g(pa) = pg(a)$$
$$a, b \in A_1 \quad p \in D$$

### 3.2. Algebra over Ring

**Definition 3.2.1.** Let $D$ be commutative ring. Let $A$ be module over ring $D$.[3.4]
For given bilinear mapping
$$f : A \times A \to A$$
we define product in $A$

(3.2.1) $$ab = f \circ (a, b)$$

$A$ is a **algebra over ring** $D$ if $A$ is $D$-module and we defined product (3.2.1) in $A$. Algebra $A^*$ is called **the opposite algebra to algebra** $A$ if we define a product in the module $A$ according to rule[3.5]
$$ba = f \circ (a, b)$$
If $A$ is free $D$-module, then $A$ is called **free algebra over ring** $D$.   □

**Remark 3.2.2.** Algebra $A$ and opposite algebra coincide as modules.   □

**Theorem 3.2.3.** *The multiplication in the algebra $A$ is distributive over addition.*

PROOF. The statement of the theorem follows from the chain of equatins
$$(a+b)c = f \circ (a+b, c) = f \circ (a, c) + f \circ (b, c) = ac + bc$$
$$a(b+c) = f \circ (a, b+c) = f \circ (a, b) + f \circ (a, c) = ab + ac$$
□

The multiplication in algebra can be neither commutative nor associative. Following definitions are based on definitions given in [14], p. 13.

**Definition 3.2.4.** The **commutator**
$$[a, b] = ab - ba$$
measures commutativity in $D$-algebra $A$. $D$-algebra $A$ is called **commutative**, if
$$[a, b] = 0$$
□

**Definition 3.2.5.** The **associator**

(3.2.2) $$(a, b, c) = (ab)c - a(bc)$$

measures associativity in $D$-algebra $A$. $D$-algebra $A$ is called **associative**, if
$$(a, b, c) = 0$$
□

**Theorem 3.2.6.** *Let $A$ be algebra over commutative ring $D$.[3.6]*

(3.2.3) $$a(b, c, d) + (a, b, c)d = (ab, c, d) - (a, bc, d) + (a, b, cd)$$

*for any $a$, $b$, $c$, $d \in A$.*

---

[3.4]There are several equivalent definitions of algebra. Initially I supposed to consider a representation of the ring $D$ in the Abelian group of the ring $A$. But I had to explain why the product of elements of the ring $D$ and of algebra $A$ is commutative. This required a definition of the center of the algebra $A$. After careful analysis I have chosen the definition given in [14], p. 1, [9], p. 4.

[3.5]I made the definition by analogy with the definition [12]-2, p. 2.

[3.6]The statement of the theorem is based on the equation [14]-(2.4).

PROOF. The equation (3.2.3) follows from the chain of equations

$$\begin{aligned}
a(b,c,d) + (a,b,c)d &= a((bc)d - b(cd)) + ((ab)c - a(bc))d \\
&= a((bc)d) - a(b(cd)) + ((ab)c)d - (a(bc))d \\
&= ((ab)c)d - (ab)(cd) + (ab)(cd) \\
&\quad + a((bc)d) - a(b(cd)) - (a(bc))d \\
&= (ab,c,d) - (a(bc))d + a((bc)d) + (ab)(cd) - a(b(cd)) \\
&= (ab,c,d) - (a,(bc),d) + (a,b,cd)
\end{aligned}$$

$\square$

**Definition 3.2.7.** The set[3.7]

$$N(A) = \{a \in A : \forall b, c \in A, (a,b,c) = (b,a,c) = (b,c,a) = 0\}$$

is called the **nucleus of an $D$-algebra** $A$. $\square$

**Definition 3.2.8.** The set[3.8]

$$Z(A) = \{a \in A : a \in N(A), \forall b \in A, ab = ba\}$$

is called the **center of an $D$-algebra** $A$. $\square$

**Theorem 3.2.9.** *Let $D$ be commutative ring. If $D$-algebra $A$ has unit, then there exits an isomorphism $f$ of the ring $D$ into the center of the algebra $A$.*

PROOF. Let $e \in A$ be the unit of the algebra $A$. Then $f \circ a = ae$. $\square$

Let $\overline{\overline{e}}$ be the basis of free algebra $A$ over ring $D$. If algebra $A$ has unit, then we assume that $\overline{e}_0$ is the unit of algebra $A$.

**Theorem 3.2.10.** *Let $\overline{\overline{e}}$ be the basis of free algebra $A$ over ring $D$. Let*

$$a = a^i e_i \quad b = b^i e_i \quad a, b \in A$$

*We can get the product of $a$, $b$ according to rule*

$$(3.2.4) \qquad (ab)^k = B_{ij}^k a^i b^j$$

*where $B_{ij}^k$ are **structural constants of algebra** $A$ **over ring** $D$. The product of basis vectors in the algebra $A$ is defined according to rule*

$$(3.2.5) \qquad \overline{e}_i \overline{e}_j = B_{ij}^k \overline{e}_k$$

PROOF. The equation (3.2.5) is corollary of the statement that $\overline{\overline{e}}$ is the basis of the algebra $A$. Since the product in the algebra is a bilinear mapping, than we can write the product of $a$ and $b$ as

$$(3.2.6) \qquad ab = a^i b^j e_i e_j$$

From equations (3.2.5), (3.2.6), it follows that

$$(3.2.7) \qquad ab = a^i b^j B_{ij}^k \overline{e}_k$$

Since $\overline{\overline{e}}$ is a basis of the algebra $A$, than the equation (3.2.4) follows from the equation (3.2.7). $\square$

---

[3.7]The definition is based on the similar definition in [14], p. 13

[3.8]The definition is based on the similar definition in [14], p. 14

**Theorem 3.2.11.** *Since the algebra $A$ is commutative, than*

(3.2.8) $$B_{ij}^p = B_{ji}^p$$

*Since the algebra $A$ is associative, than*

(3.2.9) $$B_{ij}^p B_{pk}^q = B_{ip}^q B_{jk}^p$$

PROOF. For commutative algebra, the equation (3.2.8) follows from equation
$$\overline{e}_i \overline{e}_j = \overline{e}_j \overline{e}_i$$
For associative algebra, the equation (3.2.9) follows from equation
$$(\overline{e}_i \overline{e}_j)\overline{e}_k = \overline{e}_i(\overline{e}_j \overline{e}_k)$$
$\square$

## 3.3. Linear Mapping of Algebra

Algebra is a ring. A mapping, preserving the structure of algebra as a ring, is called homomorphism of algebra. However, the statement that algebra is a module over a commutative ring is more important for us. A mapping, preserving the structure of algebra as module, is called a linear mapping of algebra. Thus, the following definition is based on the definition 3.1.6.

**Definition 3.3.1.** Let $A_1$ and $A_2$ be algebras over ring $D$. Morphism
$$g : A_1 \to A_2$$
of the representation of the ring $D$ in the Abelian group $A_1$ into the representation of the ring $D$ in the Abelian group $A_2$ is called **linear mapping of $D$-algebra $A_1$ into $D$-algebra $A_2$**. Let us denote $\mathcal{L}(A_1; A_2)$ set of linear mappings of algebra $A_1$ into algebra $A_2$. $\square$

**Theorem 3.3.2.** *Linear mapping*
$$g : A_1 \to A_2$$
*of $D$-algebra $A_1$ into $D$-algebra $A_2$ satisfies to equations*

(3.3.1) $$\begin{cases} g \circ (a+b) = g \circ a + g \circ b \\ g \circ (pa) = pg \circ a \\ a, b \in A_1 \quad p \in D \end{cases}$$

PROOF. The statement of theorem is a corollary of the theorem 3.1.7. $\square$

**Theorem 3.3.3.** *Consider algebra $A_1$ and algebra $A_2$. Let mappings*
$$f : A_1 \to A_2$$
$$g : A_1 \to A_2$$
*be linear mappings. Then mapping $f + g$ defined by equation*
$$(f+g) \circ a = f \circ a + g \circ a$$
*is linear.*

Proof. Statement of theorem follows from chains of equations

$$(f+g)\circ(x+y) = f\circ(x+y) + g\circ(x+y) = f\circ x + f\circ y + g\circ x + g\circ y$$
$$= (f+g)\circ x + (f+g)\circ y$$
$$(f+g)\circ(px) = f\circ(px) + g\circ(px) = pf\circ x + pg\circ x$$
$$= p(f+g)\circ x$$

□

**Theorem 3.3.4.** *Consider algebra $A_1$ and algebra $A_2$. Let map*

$$g : A_1 \to A_2$$

*be linear map. Then maps $ag$, $gb$, $a$, $b \in A_2$, defined by equations*

$$(ag)\circ x = a\ g\circ x$$
$$(gb)\circ x = g\circ x\ b$$

*are linear.*

Proof. Statement of theorem follows from chains of equations

$$(ag)\circ(x+y) = a\ g\circ(x+y) = a\ (g\circ x + g\circ y) = a\ g\circ x + a\ g\circ y$$
$$= (ag)\circ x + (ag)\circ y$$
$$(ag)\circ(px) = a\ g\circ(px) = ap\ g\circ x = pa\ g\circ x$$
$$= p\ (ag)\circ x$$
$$(gb)\circ(x+y) = g\circ(x+y)\ b = (g\circ x + g\circ y)\ b = g\circ x\ b + g\circ y\ b$$
$$= (gb)\circ x + (gb)\circ y$$
$$(gb)\circ(px) = g\circ(px)\ b = p\ g\circ x\ b$$
$$= p\ (gb)\circ x$$

□

**Theorem 3.3.5.** *Consider algebra $A_1$ and algebra $A_2$. Let map*

$$g : A_1 \to A_2$$

*be linear map. Then maps $pg$, $p \in D$, defined by equation*

$$(pg)\circ x = p\ g\circ x$$

*are linear. This holds*

$$p(qg) = (pq)g$$
$$(p+q)g = pg + qg$$

Proof. Statement of theorem follows from chains of equations

$$(pg)\circ(x+y) = p\ g\circ(x+y) = p\ (g\circ x + g\circ y) = p\ g\circ x + p\ g\circ y$$
$$= (pg)\circ x + (pg)\circ y$$
$$(pg)\circ(qx) = p\ g\circ(qx) = pq\ g\circ x = qp\ g\circ x$$
$$= q\ (pg)\circ x$$
$$(p(qg))\circ x = p\ (qg)\circ x = p\ (q\ g\circ x) = (pq)\ g\circ x = ((pq)g)\circ x$$
$$((p+q)g)\circ x = (p+q)\ g\circ x = p\ g\circ x + q\ g\circ x = (pg)\circ x + (qg)\circ x$$

$\square$

**Theorem 3.3.6.** *Let $D$ be commutative ring with unit. Consider $D$-algebra $A_1$ and $D$-algebra $A_2$. The set $\mathcal{L}(A_1;A_2)$ is an $D$-module.*

Proof. The theorem 3.3.3 determines the sum of linear mappings from $D$-algebra $A_1$ into $D$-algebra $A_2$. Let $f, g, h \in \mathcal{L}(A_1;A_2)$. For any $a \in A_1$

$$(f+g) \circ a = f \circ a + g \circ a = g \circ a + f \circ a$$
$$= (g+f) \circ a$$
$$((f+g)+h) \circ a = (f+g) \circ a + h \circ a = (f \circ a + g \circ a) + h \circ a$$
$$= f \circ a + (g \circ a + h \circ a) = f \circ a + (g+h) \circ a$$
$$= (f+(g+h)) \circ a$$

Therefore, sum of linear mappings is commutative and associative.

The mapping $z$ defined by equation

$$z \circ x = 0$$

is zero of addition, because

$$(z+f) \circ a = z \circ a + f \circ a = 0 + f \circ a = f \circ a$$

For a given mapping $f$ a mapping $g$ defined by equation

$$g \circ a = -f \circ a$$

satisfies to equation

$$f + g = z$$

because

$$(f+g) \circ a = f \circ a + g \circ a = f \circ a - f \circ a = 0$$

Therefore, the set $\mathcal{L}(A_1;A_2)$ is an Abelian group.

From the theorem 3.3.5, it follows that the representation of the ring $D$ in the Abelian group $\mathcal{L}(A_1;A_2)$ is defined. Since the ring $D$ has unit, than, according to the theorem 3.1.1, specified representation is effective. $\square$

### 3.4. Algebra $\mathcal{L}(A;A)$

**Theorem 3.4.1.** *Let $A$, $B$, $C$ be algebras over commutative ring $D$. Let $f$ be linear mapping from algebra $A$ into algebra $B$. Let $g$ be linear mapping from the algebra $B$ into algebra $C$. The mapping $g \circ f$ defined by diagram*

$$\begin{array}{ccc} & B & \\ {}^{f}\nearrow & & \searrow^{g} \\ A & \xrightarrow{g \circ f} & C \end{array}$$

*is linear mapping from the algebra $A$ into the algebra $C$.*

Proof. The proof of the theorem follows from chains of equations

$$(g \circ f) \circ (a+b) = g \circ (f \circ (a+b)) = g \circ (f \circ a + f \circ b)$$
$$= g \circ (f \circ a) + g \circ (f \circ b) = (g \circ f) \circ a + (g \circ f) \circ b$$
$$(g \circ f) \circ (pa) = g \circ (f \circ (pa)) = g \circ (p \ f \circ a) = p \ g \circ (f \circ a)$$
$$= p \ (g \circ f) \circ a$$

$\square$

**Theorem 3.4.2.** *Let $A$, $B$, $C$ be algebras over the commutative ring $D$. Let $f$ be a linear mapping from the algebra $A$ into the algebra $B$. The mapping $f$ generates a linear mapping*

$$f^* : g \in \mathcal{L}(B;C) \to g \circ f \in \mathcal{L}(A;C)$$

PROOF. The proof of the theorem follows from chains of equations

$$((g_1 + g_2) \circ f) \circ a = (g_1 + g_2) \circ (f \circ a) = g_1 \circ (f \circ a) + g_2 \circ (f \circ a)$$
$$= (g_1 \circ f) \circ a + (g_2 \circ f) \circ a$$
$$= (g_1 \circ f + g_2 \circ f) \circ a$$
$$((pg) \circ f) \circ a = (pg) \circ (f \circ a) = p \ g \circ (f \circ a) = p \ (g \circ f) \circ a$$
$$= (p(g \circ f)) \circ a$$

$\square$

**Theorem 3.4.3.** *Let $A$, $B$, $C$ be algebras over the commutative ring $D$. Let $g$ be a linear mapping from the algebra $B$ into the algebra $C$. The mapping $g$ generates a linear mapping*

$$g^* : f \in \mathcal{L}(A;B) \to g \circ f \in \mathcal{L}(A;C)$$

PROOF. The proof of the theorem follows from chains of equations

$$(g \circ (f_1 + f_2)) \circ a = g \circ ((f_1 + f_2) \circ a) = g \circ (f_1 \circ a + f_2 \circ a)$$
$$= g \circ (f_1 \circ a) + g \circ (f_2 \circ a) = (g \circ f_1) \circ a + (g \circ f_2) \circ a$$
$$= (g \circ f_1 + g \circ f_2) \circ a$$
$$(g \circ (pf)) \circ a = g \circ ((pf) \circ a) = g \circ (p \ (f \circ a)) = p \ g \circ (f \circ a)$$
$$= p \ (g \circ f) \circ a = (p(g \circ f)) \circ a$$

$\square$

**Theorem 3.4.4.** *Let $A$, $B$, $C$ be algebras over the commutative ring $D$. The mapping*

$$\circ : (g, f) \in \mathcal{L}(B;C) \times \mathcal{L}(A;B) \to g \circ f \in \mathcal{L}(A;C)$$

*is bilinear mapping.*

PROOF. The theorem follows from theorems 3.4.2, 3.4.3. $\square$

**Theorem 3.4.5.** *Let $A$ be algebra over commutative ring $D$. Module $\mathcal{L}(A;A)$ equiped by product*

(3.4.1) $$\circ : (g, f) \in \mathcal{L}(A;A) \times \mathcal{L}(A;A) \to g \circ f \in \mathcal{L}(A;A)$$

*is algebra over $D$.*

PROOF. The theorem follows from definition 3.2.1 and theorem 3.4.4. $\square$

## 3.5. Tensor Product of Algebras

**Definition 3.5.1.** Let $D$ be the commutative ring. Let $A_1$, ..., $A_n$, $S$ be $D$-modules. We call map
$$f : A_1 \times ... \times A_n \to S$$
**polylinear mapping of modules** $A_1$, ..., $A_n$ into module $S$, if
$$f \circ (a_1, ..., a_i + b_i, ..., a_n) = f \circ (a_1, ..., a_i, ..., a_n) + f \circ (a_1, ..., b_i, ..., a_n)$$
$$f \circ (a_1, ..., pa_i, ..., a_n) = pf \circ (a_1, ..., a_i, ..., a_n)$$
$$1 \leq i \leq n \quad a_i, b_i \in A_i \quad p \in D$$
□

**Definition 3.5.2.** Let $D$ be the commutative associative ring. Let $A_1$, ..., $A_n$ be $D$-algebras and $S$ be $D$-module. We call map
$$f : A_1 \times ... \times A_n \to S$$
**polylinear mapping of algebras** $A_1$, ..., $A_n$ into module $S$, if
$$f \circ (a_1, ..., a_i + b_i, ..., a_n) = f \circ (a_1, ..., a_i, ..., a_n) + f \circ (a_1, ..., b_i, ..., a_n)$$
$$f \circ (a_1, ..., pa_i, ..., a_n) = pf \circ (a_1, ..., a_i, ..., a_n)$$
$$1 \leq i \leq n \quad a_i, b_i \in A_i \quad p \in D$$
Let us denote $\mathcal{L}(A_1, ..., A_n; S)$ set of polylinear maps of algebras $A_1$, ..., $A_n$ into module $S$. □

**Definition 3.5.3.** Let $A_1$, ..., $A_n$ be free algebras over commutative ring $D$.[3.9] Let us consider category $\mathcal{A}$ whose objects are polylinear over commutative ring $D$ mappings
$$f : A_1 \times ... \times A_n \longrightarrow S_1 \qquad g : A_1 \times ... \times A_n \longrightarrow S_2$$
where $S_1$, $S_2$ are modules over ring $D$, We define morphism $f \to g$ to be linear over commutative ring $D$ mapping $h : S_1 \to S_2$ making diagram

$$\begin{array}{c}
 & & S_1 \\
 & \nearrow^{f} & \\
A_1 \times ... \times A_n & & \downarrow h \\
 & \searrow_{g} & \\
 & & S_2
\end{array}$$

commutative. Universal object $A_1 \otimes ... \otimes A_n$ of category $\mathcal{A}$ is called **tensor product of algebras** $A_1$, **...,** $A_n$. □

**Theorem 3.5.4.** *There exists tensor product of algebras.*

---

[3.9]I give definition of tensor product of algebras following to definition in [1], p. 601 - 603.

PROOF. Let $M$ be module over ring $D$ generated by product $A_1 \times ... \times A_n$ of algebras $A_1, ..., A_n$. Injection

$$i : A_1 \times ... \times A_n \longrightarrow M$$

is defined according to rule

(3.5.1) $$i \circ (d_1, ..., d_n) = (d_1, ..., d_n)$$

Let $N \subset M$ be submodule generated by elements of the following type

(3.5.2) $$(d_1, ..., d_i + c_i, ..., d_n) - (d_1, ..., d_i, ..., d_n) - (d_1, ..., c_i, ..., d_n)$$
(3.5.3) $$(d_1, ..., ad_i, ..., d_n) - a(d_1, ..., d_i, ..., d_n)$$

where $d_i \in A_i$, $c_i \in A_i$, $a \in D$. Let

$$j : M \to M/N$$

be canonical map on factor module. Consider commutative diagram

(3.5.4)
$$\begin{array}{c} & & M/N \\ & \nearrow f & \nearrow j \\ A_1 \times ... \times A_n & \xrightarrow{i} & M \end{array}$$

Since elements (3.5.2) and (3.5.3) belong to kernel of linear map $j$, then, from equation (3.5.1), it follows

(3.5.5) $$f \circ (d_1, ..., d_i + c_i, ..., d_n) = f \circ (d_1, ..., d_i, ..., d_n) + f \circ (d_1, ..., c_i, ..., d_n)$$
(3.5.6) $$f \circ (d_1, ..., ad_i, ..., d_n) = a\ f \circ (d_1, ..., d_i, ..., d_n)$$

From equations (3.5.5) and (3.5.6) it follows that map $f$ is polylinear over ring $D$. Since $M$ is module with basis $A_1 \times ... \times A_n$, then, according to theorem [1]-4.1 on p. 135, for any module $V$ and any polylinear over $D$ map

$$g : A_1 \times ... \times A_n \longrightarrow V$$

there exists a unique homomorphism $k : M \to V$, for which following diagram is commutative

(3.5.7)
$$\begin{array}{c} A_1 \times ... \times A_n & \xrightarrow{i} & M \\ & \searrow g & \downarrow k \\ & & V \end{array}$$

Since $g$ is polylinear over $D$, then ker $k \subseteq N$. According to statement on p. [1]-119, map $j$ is universal in the category of homomorphisms of vector space $M$ whose kernel contains $N$. Therefore, we have homomorphism

$$h : M/N \to V$$

which makes the following diagram commutative

(3.5.8)
$$\begin{array}{c} M/N \\ \nearrow^{j} \quad \downarrow^{h} \\ M \\ \searrow_{k} \\ V \end{array}$$

We join diagrams (3.5.4), (3.5.7), (3.5.8), and get commutative diagram

(3.5.9)
$$\begin{array}{c} A_1 \times ... \times A_n \xrightarrow{i} M \xrightarrow{j} M/N \\ \xrightarrow{f} \quad \xrightarrow{k} \quad \downarrow^{h} \\ \xrightarrow{g} V \end{array}$$

Since $\text{Im} f$ generates $M/N$, then map $h$ is uniquely determined. □

According to proof of theorem 3.5.4

$$A_1 \otimes ... \otimes A_n = M/N$$

If $d_i \in A_i$, we write

(3.5.10) $$j \circ (d_1, ..., d_n) = d_1 \otimes ... \otimes d_n$$

**Theorem 3.5.5.** *Let $A_1$, ..., $A_n$ be algebras over commutative ring $D$. Let*

$$f : A_1 \times ... \times A_n \to A_1 \otimes ... \otimes A_n$$

*be polylinear mapping defined by equation*

(3.5.11) $$f \circ (d_1, ..., d_n) = d_1 \otimes ... \otimes d_n$$

*Let*

$$g : A_1 \times ... \times A_n \to V$$

*be polylinear mapping into $D$-module $V$. There exists an $D$-linear mapping*

$$h : A_1 \otimes ... \otimes A_n \to V$$

*such that the diagram*

(3.5.12)
$$\begin{array}{c} A_1 \otimes ... \otimes A_n \\ \nearrow^{f} \quad \downarrow^{h} \\ A_1 \times ... \times A_n \\ \searrow_{g} \\ V \end{array}$$

*is commutative.*

PROOF. Equation (3.5.11) follows from equations (3.5.1) and (3.5.10). An existence of the mapping $h$ follows from the definition 3.5.3 and constructions made in the proof of the theorem 3.5.4. □

We can write equations (3.5.5) and (3.5.6) as

$$
\begin{aligned}
& a_1 \otimes ... \otimes (a_i + b_i) \otimes ... \otimes a_n \\
& = a_1 \otimes ... \otimes a_i \otimes ... \otimes a_n + a_1 \otimes ... \otimes b_i \otimes ... \otimes a_n
\end{aligned}
\tag{3.5.13}
$$

$$
a_1 \otimes ... \otimes (ca_i) \otimes ... \otimes a_n = c(a_1 \otimes ... \otimes a_i \otimes ... \otimes a_n)
\tag{3.5.14}
$$

$$a_i \in A_i \quad b_i \in A_i \quad c \in D$$

**Theorem 3.5.6.** *Let $A$ be algebra over commutative ring $D$. There exists a linear mapping*

$$h : a \otimes b \in A \otimes A \to ab \in A$$

PROOF. The theorem is corollary of the theorem 3.5.5 and the definition 3.2.1. □

**Theorem 3.5.7.** *Tensor product $A_1 \otimes ... \otimes A_n$ of free finite dimensional algebras $A_1$, ..., $A_n$ over the commutative ring $D$ is free finite dimensional algebra.*

*Let $\overline{\overline{e}}_i$ be the basis of algebra $A_i$ over ring $D$. We can represent any tensor $a \in A_1 \otimes ... \otimes A_n$ in the following form*

$$a = a^{i_1...i_n} \overline{e}_{1 \cdot i_1} \otimes ... \otimes \overline{e}_{n \cdot i_n} \tag{3.5.15}$$

*Expression $a^{i_1...i_n}$ is called* **standard component of tensor**.

PROOF. Algebras $A_1$, ..., $A_n$ are modules over the ring $D$. According to theorem 3.5.4, $A_1 \otimes ... \otimes A_n$ is module.

Vector $a_i \in A_i$ has expansion

$$a_i = a_i^k \overline{e}_{i \cdot k}$$

relative to basis $\overline{\overline{e}}_i$. From equations (3.5.13), (3.5.14), it follows

$$a_1 \otimes ... \otimes a_n = a_1^{i_1} ... a_n^{i_n} \overline{e}_{1 \cdot i_1} \otimes ... \otimes \overline{e}_{n \cdot i_n}$$

Since set of tensors $a_1 \otimes ... \otimes a_n$ is the generating set of module $A_1 \otimes ... \otimes A_n$, than we can write tensor $a \in A_1 \otimes ... \otimes A_n$ in form

$$a = a^s a_{s \cdot 1}^{i_1} ... a_{s \cdot n}^{i_n} \overline{e}_{1 \cdot i_1} \otimes ... \otimes \overline{e}_{n \cdot i_n} \tag{3.5.16}$$

where $a^s, a_{s \cdot 1}^{i_1}, a_{s \cdot n}^{i_n} \in F$. Let

$$a^s a_{s \cdot 1}^{i_1} ... a_{s \cdot n}^{i_n} = a^{i_1...i_n} \tag{3.5.17}$$

Then equation (3.5.16) has form (3.5.15).

Therefore, set of tensors $\overline{e}_{1 \cdot i_1} \otimes ... \otimes \overline{e}_{n \cdot i_n}$ is the generating set of module $A_1 \otimes ... \otimes A_n$. Since the dimension of module $A_i$, $i = 1, ..., n$, is finite, than the set of tensors $\overline{e}_{1 \cdot i_1} \otimes ... \otimes \overline{e}_{n \cdot i_n}$ is finite. Therefore, the set of tensors $\overline{e}_{1 \cdot i_1} \otimes ... \otimes \overline{e}_{n \cdot i_n}$ contains a basis of module $A_1 \otimes ... \otimes A_n$, and the module $A_1 \otimes ... \otimes A_n$ is free module over the ring $D$.

We define the product of tensors like $a_1 \otimes ... \otimes a_n$ componentwise

$$(d_1 \otimes ... \otimes d_n)(c_1 \otimes ... \otimes c_n) = (d_1 c_1) \otimes ... \otimes (d_n c_n) \tag{3.5.18}$$

In particular, if for any $i$, $i = 1, ..., n$, $a_i \in A_i$ has inverse, then tensor

$$(a_1 \otimes ... \otimes a_n)^{-1} = (a_1)^{-1} \otimes ... \otimes (a_n)^{-1}$$

is inverse tensor to tensor

$$a_1 \otimes ... \otimes a_n \in A_1 \otimes ... \otimes A_n$$

The definition of the product (3.5.18) agreed with equation (3.5.14) because

$$(a_1 \otimes ... \otimes (ca_i) \otimes ... \otimes a_n)(b_1 \otimes ... \otimes b_i \otimes ... \otimes b_n)$$
$$=(a_1 b_1) \otimes ... \otimes (ca_i)b_i \otimes ... \otimes (a_n b_n)$$
$$=c((a_1 b_1) \otimes ... \otimes (a_i b_i) \otimes ... \otimes (a_n b_n))$$
$$=c((a_1 \otimes ... \otimes a_n)(b_1 \otimes ... \otimes b_n))$$

The distributive property of multiplication over addition

(3.5.19)
$$\begin{aligned}&(a_1 \otimes ... \otimes a_i \otimes ... \otimes a_n)\\ &*((b_1 \otimes ... \otimes b_i \otimes ... \otimes b_n) + (b_1 \otimes ... \otimes c_i \otimes ... \otimes b_n))\\ &=(a_1 \otimes ... \otimes a_i \otimes ... \otimes a_n)(b_1 \otimes ... \otimes (b_i + c_i) \otimes ... \otimes b_n)\\ &=(a_1 b_1) \otimes ... \otimes (a_i(b_i+c_i)) \otimes ... \otimes (a_n b_n)\\ &=(a_1 b_1) \otimes ... \otimes (a_i b_i + a_i c_i) \otimes ... \otimes (a_n b_n)\\ &=(a_1 b_1) \otimes ... \otimes (a_i b_i) \otimes ... \otimes (a_n b_n)\\ &+(a_1 b_1) \otimes ... \otimes (a_i c_i) \otimes ... \otimes (a_n b_n)\\ &=(a_1 \otimes ... \otimes a_i \otimes ... \otimes a_n)(b_1 \otimes ... \otimes b_i \otimes ... \otimes b_n)\\ &+(a_1 \otimes ... \otimes a_i \otimes ... \otimes a_n)(b_1 \otimes ... \otimes c_i \otimes ... \otimes b_n)\end{aligned}$$

follows from the equation (3.5.13). The equation (3.5.19) allows us to define the product for any tensors $a$, $b$. □

**Remark 3.5.8.** According to the remark 3.2.2, we can define different structures of algebra in the tensor product of algebras. For instance, algebras $A_1 \otimes A_2$, $A_1 \otimes A_2^*$, $A_1^* \otimes A_2$ are defined in the same module. □

**Theorem 3.5.9.** *Let $\overline{\overline{e}}_i$ be the basis of the algebra $A_i$ over the ring $D$. Let $B_{i\cdot kl}^{\phantom{i\cdot}j}$ be structural constants of the algebra $A_i$ relative the basis $\overline{\overline{e}}_i$. Structural constants of the tensor product $A_1 \otimes ... \otimes A_n$ relative to the basis $\overline{e}_{1\cdot i_1} \otimes ... \otimes \overline{e}_{n\cdot i_n}$ have form*

(3.5.20) $$B_{\cdot k_1...k_n\cdot l_1...l_n}^{\cdot j_1...j_n} = B_{1\cdot k_1 l_1}^{\phantom{1\cdot}j_1}...B_{n\cdot k_n l_n}^{\phantom{n\cdot}j_n}$$

PROOF. Direct multiplication of tensors $\overline{e}_{1\cdot i_1} \otimes ... \otimes \overline{e}_{n\cdot i_n}$ has form

(3.5.21)
$$\begin{aligned}&(\overline{e}_{1\cdot k_1} \otimes ... \otimes \overline{e}_{n\cdot k_n})(\overline{e}_{1\cdot l_1} \otimes ... \otimes \overline{e}_{n\cdot l_n})\\ &=(\overline{e}_{1\cdot k_1}\overline{e}_{1\cdot l_1}) \otimes ... \otimes (\overline{e}_{n\cdot k_n}\overline{e}_{n\cdot l_n})\\ &=(\overline{e}_{1\cdot k_1}\overline{e}_{1\cdot l_1}) \otimes ... \otimes (\overline{e}_{n\cdot k_n}\overline{e}_{n\cdot l_n})\\ &=(B_{1\cdot k_1 l_1}^{\phantom{1\cdot}j_1}\overline{e}_{1\cdot j_1}) \otimes ... \otimes (B_{n\cdot k_n l_n}^{\phantom{n\cdot}j_n}\overline{e}_{n\cdot j_n})\\ &=B_{1\cdot k_1 l_1}^{\phantom{1\cdot}j_1}...B_{n\cdot k_n l_n}^{\phantom{n\cdot}j_n}\overline{e}_{1\cdot j_1} \otimes ... \otimes \overline{e}_{n\cdot j_n}\end{aligned}$$

According to the definition of structural constants

(3.5.22) $(\overline{e}_{1\cdot k_1} \otimes ... \otimes \overline{e}_{n\cdot k_n})(\overline{e}_{1\cdot l_1} \otimes ... \otimes \overline{e}_{n\cdot l_n}) = B_{\cdot k_1...k_n\cdot l_1...l_n}^{\cdot j_1...j_n}(\overline{e}_{1\cdot j_1} \otimes ... \otimes \overline{e}_{n\cdot j_n})$

The equation (3.5.20) follows from comparison (3.5.21), (3.5.22). □

From the chain of equations

$$(a_1 \otimes ... \otimes a_n)(b_1 \otimes ... \otimes b_n)$$
$$= (a_1^{k_1}\overline{e}_{1\cdot k_1} \otimes ... \otimes a_n^{k_n}\overline{e}_{n\cdot k_n})(b_1^{l_1}\overline{e}_{1\cdot l_1} \otimes ... \otimes b_n^{l_n}\overline{e}_{n\cdot l_n})$$
$$= a_1^{k_1}...a_n^{k_n}b_1^{l_1}...b_n^{l_n}(\overline{e}_{1\cdot k_1} \otimes ... \otimes \overline{e}_{n\cdot k_n})(\overline{e}_{1\cdot l_1} \otimes ... \otimes \overline{e}_{n\cdot l_n})$$
$$= a_1^{k_1}...a_n^{k_n}b_1^{l_1}...b_n^{l_n} B^{\cdot j_1...j_n}_{\cdot k_1...k_n\cdot l_1...l_n}(\overline{e}_{1\cdot j_1} \otimes ... \otimes \overline{e}_{n\cdot j_n})$$
$$= a_1^{k_1}...a_n^{k_n}b_1^{l_1}...b_n^{l_n} B_1{}^{\cdot j_1}_{\cdot k_1 l_1}...B_n{}^{\cdot j_n}_{\cdot k_n l_n}(\overline{e}_{1\cdot j_1} \otimes ... \otimes \overline{e}_{n\cdot j_n})$$
$$= (a_1^{k_1}b_1^{l_1} B_1{}^{\cdot j_1}_{\cdot k_1 l_1}\overline{e}_{1\cdot j_1}) \otimes ... \otimes (a_n^{k_n}b_n^{l_n} B_n{}^{\cdot j_n}_{\cdot k_n l_n}\overline{e}_{n\cdot j_n})$$
$$= (a_1 b_1) \otimes ... \otimes (a_n b_n)$$

it follows that definition of product (3.5.22) with structural constants (3.5.20) agreed with the definition of product (3.5.18). $\square$

**Theorem 3.5.10.** *For tensors $a, b \in A_1 \otimes ... \otimes A_n$, standard components of product satisfy to equation*

$$(3.5.23) \qquad (ab)^{j_1...j_n} = B^{\cdot j_1...j_n}_{\cdot k_1...k_n\cdot l_1...l_n} a^{k_1...k_n} b^{l_1...l_n}$$

PROOF. According to the definition

$$(3.5.24) \qquad ab = (ab)^{j_1...j_n}\overline{e}_{1\cdot j_1} \otimes ... \otimes \overline{e}_{n\cdot j_n}$$

At the same time

$$(3.5.25) \quad \begin{aligned} ab &= a^{k_1...k_n}\overline{e}_{1\cdot k_1} \otimes ... \otimes \overline{e}_{n\cdot k_n} b^{k_1...k_n}\overline{e}_{1\cdot l_1} \otimes ... \otimes \overline{e}_{n\cdot l_n} \\ &= a^{k_1...k_n}b^{k_1...k_n} B^{\cdot j_1...j_n}_{\cdot k_1...k_n\cdot l_1...l_n}\overline{e}_{1\cdot j_1} \otimes ... \otimes \overline{e}_{n\cdot j_n} \end{aligned}$$

The equation (3.5.23) follows from equations (3.5.24), (3.5.25). $\square$

**Theorem 3.5.11.** *If the algebra $A_i$, $i = 1, ..., n$, is associative, then the tensor product $A_1 \otimes ... \otimes A_n$ is associative algebra.*

PROOF. Since

$$((\overline{e}_{1\cdot i_1} \otimes ... \otimes \overline{e}_{n\cdot i_n})(\overline{e}_{1\cdot j_1} \otimes ... \otimes \overline{e}_{n\cdot j_n}))(\overline{e}_{1\cdot k_1} \otimes ... \otimes \overline{e}_{n\cdot k_n})$$
$$=((\overline{e}_{1\cdot i_1}\overline{e}_{1\cdot j_1}) \otimes ... \otimes (\overline{e}_{n\cdot i_n}\overline{e}_{1\cdot j_n}))(\overline{e}_{1\cdot k_1} \otimes ... \otimes \overline{e}_{n\cdot k_n})$$
$$=((\overline{e}_{1\cdot i_1}\overline{e}_{1\cdot j_1})\overline{e}_{1\cdot k_1}) \otimes ... \otimes ((\overline{e}_{n\cdot i_n}\overline{e}_{1\cdot j_n})\overline{e}_{1\cdot k_n})$$
$$=(\overline{e}_{1\cdot i_1}(\overline{e}_{1\cdot j_1}\overline{e}_{1\cdot k_1})) \otimes ... \otimes (\overline{e}_{n\cdot i_n}(\overline{e}_{1\cdot j_n}\overline{e}_{1\cdot k_n}))$$
$$=(\overline{e}_{1\cdot i_1} \otimes ... \otimes \overline{e}_{n\cdot i_n})((\overline{e}_{1\cdot j_1}\overline{e}_{1\cdot k_1}) \otimes ... \otimes (\overline{e}_{1\cdot j_n}\overline{e}_{1\cdot k_n}))$$
$$=(\overline{e}_{1\cdot i_1} \otimes ... \otimes \overline{e}_{n\cdot i_n})((\overline{e}_{1\cdot j_1} \otimes ... \otimes \overline{e}_{n\cdot j_n})(\overline{e}_{1\cdot k_1} \otimes ... \otimes \overline{e}_{n\cdot k_n}))$$

than

$$(ab)c = a^{i_1...i_n}b^{j_1...j_n}c^{k_1...k_n}$$
$$\quad ((\overline{e}_{1\cdot i_1} \otimes ... \otimes \overline{e}_{n\cdot i_n})(\overline{e}_{1\cdot j_1} \otimes ... \otimes \overline{e}_{n\cdot j_n}))(\overline{e}_{1\cdot k_1} \otimes ... \otimes \overline{e}_{n\cdot k_n})$$
$$= a^{i_1...i_n}b^{j_1...j_n}c^{k_1...k_n}$$
$$\quad (\overline{e}_{1\cdot i_1} \otimes ... \otimes \overline{e}_{n\cdot i_n})((\overline{e}_{1\cdot j_1} \otimes ... \otimes \overline{e}_{n\cdot j_n})(\overline{e}_{1\cdot k_1} \otimes ... \otimes \overline{e}_{n\cdot k_n}))$$
$$= a(bc)$$

$\square$

## 3.6. Linear Mapping into Associative Algebra

**Theorem 3.6.1.** *Consider $D$-algebras $A_1$ and $A_2$. For given mapping $f \in \mathcal{L}(A_1; A_2)$, the mapping*
$$g : A_2 \times A_2 \to \mathcal{L}(A_1; A_2)$$
$$g(a,b) \circ f = afb$$
*is bilinear mapping.*

Proof. The statement of theorem follows from chains of equations
$$((a_1 + a_2)fb) \circ x = (a_1 + a_2) \; f \circ x \; b = a_1 \; f \circ x \; b + a_2 \; f \circ x \; b$$
$$= (a_1 fb) \circ x + (a_2 fb) \circ x = (a_1 fb + a_2 fb) \circ x$$
$$((pa)fb) \circ x = (pa) \; f \circ x \; b = p(a \; f \circ x \; b) = p((afb) \circ x) = (p(afb)) \circ x$$
$$(af(b_1 + b_2)) \circ x = a \; f \circ x \; (b_1 + b_2) = a \; f \circ x \; b_1 + a \; f \circ x \; b_2$$
$$= (afb_1) \circ x + (afb_2) \circ x = (afb_1 + afb_2) \circ x$$
$$(af(pb)) \circ x = a \; f \circ x \; (pb) = p(a \; f \circ x \; b) = p((afb) \circ x) = (p(afb)) \circ x$$
□

**Theorem 3.6.2.** *Consider $D$-algebras $A_1$ and $A_2$. For given mapping $f \in \mathcal{L}(A_1; A_2)$, there exists linear mapping*
$$h : A_2 \otimes A_2 \to \mathcal{L}(A_1; A_2)$$
*defined by the equation*
$$(3.6.1) \qquad (a \otimes b) \circ f = afb$$

Proof. The statement of the theorem is corollary of theorems 3.5.5, 3.6.1. □

**Theorem 3.6.3.** *Consider $D$-algebras $A_1$ and $A_2$. A linear mapping*
$$h : A_2 \otimes A_2 \to {}^*\mathcal{L}(A_1; A_2)$$
*defined by the equation*
$$(3.6.2) \qquad (a \otimes b) \circ f = afb \quad a, b \in A_2 \quad f \in \mathcal{L}(A_1; A_2)$$
*is representation*[3.10] *of module $A_2 \otimes A_2$ in module $\mathcal{L}(A_1; A_2)$.*

Proof. According to theorem 3.3.4, mapping (3.6.2) is transformation of module $\mathcal{L}(A_1; A_2)$. For a given tensor $c \in A_2 \otimes A_2$, a transformation $h(c)$ is a linear

---
[3.10] See the definition of representation of $\Omega$-algebra in the definition [7]-2.1.4.

transformation of module $\mathcal{L}(A_1; A_2)$, because

$$((a \otimes b) \circ (f_1 + f_2)) \circ x = (a(f_1 + f_2)b) \circ x = a((f_1 + f_2) \circ x)b$$
$$= a(f_1 \circ x + f_2 \circ x)b = a(f_1 \circ x)b + a(f_2 \circ x)b$$
$$= (af_1 b) \circ x + (af_2 b) \circ x$$
$$= (a \otimes b) \circ f_1 \circ x + (a \otimes b) \circ f_2 \circ x$$
$$= ((a \otimes b) \circ f_1 + (a \otimes b) \circ f_2) \circ x$$
$$((a \otimes b) \circ (pf)) \circ x = (a(pf)b) \circ x = a((pf) \circ x)b$$
$$= a(p \ f \circ x)b = pa(f \circ x)b$$
$$= p \ (afb) \circ x = p \ ((a \otimes b) \circ f) \circ x$$
$$= (p((a \otimes b) \circ f)) \circ x$$

According to theorem 3.6.2, mapping (3.6.2) is linear mapping. According to the definition [7]-2.1.4 mapping (3.6.2) is a representation of the module $A_2 \otimes A_2$ in the module $\mathcal{L}(A_1; A_2)$. □

**Theorem 3.6.4.** *Let $A$ be algebra over the commutative ring $D$. Algebra $A \otimes A$, whose product is defined according to rule*

$$(3.6.3) \qquad (a \otimes b) \circ (c \otimes d) = (ac) \otimes (db)$$

*forms the representation in the module $\mathcal{L}(A; A)$. This representation allows us to identify tensor $d \in A \otimes A$ and transformation $d \circ \delta$ where $\delta$ is identity transformation.*

PROOF. According to the theorem 3.6.2, the mapping $f \in \mathcal{L}(A; A)$ and the tensor $d \in A \otimes A$ generate the mapping

$$(3.6.4) \qquad x \to (d \circ f) \circ x$$

If we assume $f = \delta$, $d = a \otimes b$, then the equation (3.6.4) gets form

$$(3.6.5) \qquad ((a \otimes b) \circ \delta) \circ x = (a\delta b) \circ x = a \ (\delta \circ x) \ b = axb$$

If we assume

$$(3.6.6) \qquad ((a \otimes b) \circ \delta) \circ x = (a \otimes b) \circ (\delta \circ x) = (a \otimes b) \circ x$$

then comparison of equations (3.6.5) and (3.6.6) gives a basis to identify the action of the tensor $a \otimes b$ and transformation $(a \otimes b) \circ \delta$. Therefore, the mapping

$$(3.6.7) \qquad d \in A \otimes A \to d \circ \delta \in \mathcal{L}(A; A)$$

is the homomorphism of the module $A \otimes A$ into the module $\mathcal{L}(A; A)$.

Mapping (3.6.7) is also homomorphisms of algebras, because the product of transformations $a \otimes b$ and $c \otimes d$ has form

$$((a \otimes b) \circ (c \otimes d)) \circ x = ((ac) \otimes (db)) \circ x$$
$$= (ac)x(db)$$
$$= a(cxd)b$$
$$= (a \otimes b) \circ (cxd)$$
$$= (a \otimes b) \circ ((c \otimes d) \circ x)$$

□

3.6. Linear Mapping into Associative Algebra                                    37From the theorem 3.6.4, it follows that we can consider the mapping (3.6.2) as the product of mappings $a \otimes b$ and $f$. This allows us to consider the representation of the algebra $A_2 \otimes A_2$ in the module $\mathcal{L}(A_1; A_2)$ instead of the representation of the module $A_2 \otimes A_2$ in the module $\mathcal{L}(A_1; A_2)$.

The tensor $a \in A_2 \otimes A_2$ is **nonsingular**, if there exists the tensor $b \in A_2 \otimes A_2$ such that $a \circ b = 1 \otimes 1$.

**Definition 3.6.5.** Consider the representation of algebra $A_2 \otimes A_2$ in the module $\mathcal{L}(A_1; A_2)$.[3.11] The set
$$(A_2 \otimes A_2) \circ f = \{g = d \circ f : d \in A_2 \otimes A_2\}$$
is called **orbit of linear mapping** $f \in \mathcal{L}(A_1; A_2)$ □

**Theorem 3.6.6.** *Consider D-algebra $A_1$ and associative D-algebra $A_2$. Consider the representation of algebra $A_2 \otimes A_2$ in the module $\mathcal{L}(A_1; A_2)$. The mapping*
$$h : A_1 \to A_2$$
*generated by the mapping*
$$f : A_1 \to A_2$$
*has form*
$$h = (a_{s \cdot 0} \otimes a_{s \cdot 1}) \circ f = a_{s \cdot 0} f a_{s \cdot 1}$$

PROOF. We can represent any tensor $a \in A_2 \otimes A_2$ in the form
$$a = a_{s \cdot 0} \otimes a_{s \cdot 1}$$
According to the theorem 3.6.3, the mapping (3.6.2) is linear. This proofs the statement of the theorem. □

**Theorem 3.6.7.** *Let $A_2$ be algebra with unit $e$. Let $a \in A_2 \otimes A_2$ be a nonsingular tensor. Orbits of linear mappings $f \in \mathcal{L}(A_1; A_2)$ and $g = a \circ f$ coincide*

(3.6.8) $$(A_2 \otimes A_2) \circ f = (A_2 \otimes A_2) \circ g$$

PROOF. If $h \in (A_2 \otimes A_2) \circ g$, then there exists $b \in A_2 \otimes A_2$ such that $h = b \circ g$. Than

(3.6.9) $$h = b \circ (a \circ f) = (b \circ a) \circ f$$

Therefore, $h \in (A_2 \otimes A_2) \circ f$,

(3.6.10) $$(A_2 \otimes A_2) \circ g \subset (A_2 \otimes A_2) \circ f$$

Since $a$ is nonsingular tensor, than

(3.6.11) $$f = a^{-1} \circ g$$

If $h \in (A_2 \otimes A_2) \circ f$, then there exists $b \in A_2 \otimes A_2$ such that

(3.6.12) $$h = b \circ f$$

From equations (3.6.11), (3.6.12), it follows that
$$h = b \circ (a^{-1} \circ g) = (b \circ a^{-1}) \circ g$$

Therefore, $h \in (A_2 \otimes A_2) \circ g$,

(3.6.13) $$(A_2 \otimes A_2) \circ f \subset (A_2 \otimes A_2) \circ g$$

(3.6.8) follows from equations (3.6.10), (3.6.13). □

---

[3.11]The definition is made by analogy with the definition [7]-2.4.12.

From the theorem 3.6.7, it also follows that if $g = a \circ f$ and $a \in A_2 \otimes A_2$ is a singular tensor, then relationship (3.6.10) is true. However, the main result of the theorem 3.6.7 is that the representations of the algebra $A_2 \otimes A_2$ in module $\mathcal{L}(A_1; A_2)$ generates an equivalence in the module $\mathcal{L}(A_1; A_2)$. If we successfully choose the representatives of each equivalence class, then the resulting set will be generating set of considered representation.[3.12]

### 3.7. Linear Mapping into Free Finite Dimensional Associative Algebra

**Theorem 3.7.1.** *Let $A_1$ be algebra over the ring $D$. Let $A_2$ be free finite dimensional associative algebra over the ring $D$. Let $\overline{\overline{e}}$ be basis of the algebra $A_2$ over the ring $D$. The mapping*

$$(3.7.1) \qquad g = a \circ f$$

*generated by the mapping $f \in (A_1; A_2)$ through the tensor $a \in A_2 \otimes A_2$, has the standard representation*

$$(3.7.2) \qquad g = a^{ij}(\overline{e}_i \otimes \overline{e}_j) \circ f = a^{ij}\overline{e}_i f \overline{e}_j$$

PROOF. According to theorem 3.5.7, the standard representation of the tensor $a$ has form

$$(3.7.3) \qquad a = a^{ij}\overline{e}_i \otimes \overline{e}_j$$

The equation (3.7.2) follows from equations (3.7.1), (3.7.3). □

**Theorem 3.7.2.** *Let $\overline{\overline{e}}_1$ be basis of the free finite dimensional $D$-algebra $A_1$. Let $\overline{\overline{e}}_2$ be basis of the free finite dimensional associative $D$-algebra $A_2$. Let $B_2 \cdot {}^{p}_{kl}$ be structural constants of algebra $A_2$. Coordinates of the mapping*

$$g = a \circ f$$

*generated by the mapping $f \in (A_1; A_2)$ through the tensor $a \in A_2 \otimes A_2$ and its standard components are connected by the equation*

$$(3.7.4) \qquad g^k_l = f^m_l g^{ij} B_2 \cdot {}^{p}_{im} B_2 \cdot {}^{k}_{pj}$$

PROOF. Relative to bases $\overline{\overline{e}}_1$ and $\overline{\overline{e}}_2$, linear mappings $f$ and $g$ have form

$$(3.7.5) \qquad f \circ x = f^i_j x^j \overline{e}_{2 \cdot i}$$

$$(3.7.6) \qquad g \circ x = g^i_j x^j \overline{e}_{2 \cdot i}$$

From equations (3.7.5), (3.7.6), (3.7.2) it follows that

$$(3.7.7) \qquad \begin{aligned} g^k_l x^l \overline{e}_{2 \cdot k} &= a^{ij}\overline{e}_{2 \cdot i} f^m_l x^l \overline{e}_{2 \cdot m} \overline{e}_{2 \cdot j} \\ &= a^{ij} f^m_l x^l B_2 \cdot {}^{p}_{im} B_2 \cdot {}^{k}_{pj} \overline{e}_{2 \cdot k} \end{aligned}$$

Since vectors $\overline{e}_{2 \cdot k}$ are linear independent and $x^i$ are arbitrary, than the equation (3.7.4) follows from the equation (3.7.7). □

**Theorem 3.7.3.** *Let $\overline{\overline{e}}_1$ be basis of the free finite dimensional $D$-algebra $A_1$. Let $\overline{\overline{e}}_2$ be basis of the free finite dimensional associative $D$-algebra $A_2$. Let $B_2 \cdot {}^{p}_{kl}$ be structural constants of algebra $A_2$. Consider matrix*

$$(3.7.8) \qquad \mathcal{B} = \left(\mathcal{B} \cdot {}^{k}_{m \cdot ij}\right) = \left(B_2 \cdot {}^{p}_{im} B_2 \cdot {}^{k}_{pj}\right)$$

---

[3.12]Generating set of representation is defined in definition [7]-2.6.5.

whose rows and columns are indexed by $._m^{.k}$ and $._{\cdot ij}$, respectively. If matrix $\mathcal{B}$ is nonsingular, then, for given coordinates of linear transformation $g_k^l$ and for mapping $f = \delta$, the system of linear equations (3.7.4) with standard components of this transformation $g^{kr}$ has the unique solution.

If matrix $\mathcal{B}$ is singular, then the equation

$$(3.7.9) \qquad \operatorname{rank}\begin{pmatrix} \mathcal{B}^{\cdot k}_{m \cdot ij} & g_m^k \end{pmatrix} = \operatorname{rank} \mathcal{B}$$

is the condition for the existence of solutions of the system of linear equations (3.7.4). In such case the system of linear equations (3.7.4) has infinitely many solutions and there exists linear dependence between values $g_m^k$.

PROOF. The statement of the theorem is corollary of the theory of linear equations over ring. $\square$

**Theorem 3.7.4.** *Let $A$ be free finite dimensional associative algebra over the ring $D$. Let $\overline{\overline{e}}$ be basis of the algebra $A$ over the ring $D$. Let $B_{kl}^p$ be structural constants of algebra $A$. Let matrix (3.7.8) be singular. Let the linear mapping $f \in \mathcal{L}(A;A)$ be nonsingular. If coordinates of linear transformations $f$ and $g$ satisfy to the equation*

$$(3.7.10) \qquad \operatorname{rank}\begin{pmatrix} \mathcal{B}^{\cdot k}_{m \cdot ij} & g_m^k & f_m^k \end{pmatrix} = \operatorname{rank} \mathcal{B}$$

*then the system of linear equations*

$$(3.7.11) \qquad g_l^k = f_l^m g^{ij} B_{im}^p B_{pj}^k$$

*has infinitely many solutions.*

PROOF. According to the equation (3.7.10) and the theorem 3.7.3, the system of linear equations

$$(3.7.12) \qquad f_l^k = f^{ij} B_{il}^p B_{pj}^k$$

has infinitely many solutions corresponding to linear mapping

$$(3.7.13) \qquad f = f^{ij} \overline{e}_i \otimes \overline{e}_j$$

According to the equation (3.7.10) and the theorem 3.7.3, the system of linear equations

$$(3.7.14) \qquad g_l^k = g^{ij} B_{il}^p B_{pj}^k$$

has infinitely many solutions corresponding to linear mapping

$$(3.7.15) \qquad g = g^{ij} \overline{e}_i \otimes \overline{e}_j$$

Mappings $f$ and $g$ are generated by the mapping $\delta$. According to the theorem 3.6.7, the mapping $f$ generates the mapping $g$. This proves the statement of the theorem. $\square$

**Theorem 3.7.5.** *Let $A$ be free finite dimensional associative algebra over the ring $D$. The representation of algebra $A \otimes A$ in algebra $\mathcal{L}(A;A)$ has finite basis $\overline{\overline{I}}$.*

(1) *The linear mapping $f \in \mathcal{L}(A;A)$ has form*

$$(3.7.16) \qquad f = \sum_k (a_{s_k \cdot 0} \otimes a_{s_k \cdot 1}) \circ I_k = \sum_k a_{s_k \cdot 0} I_k a_{s_k \cdot 1}$$

(2) *Its standard representation has form*

$$(3.7.17) \qquad f = a^{k \cdot ij}(\overline{e}_i \otimes \overline{e}_j) \circ I_k = a^{k \cdot ij} \overline{e}_i I_k \overline{e}_j$$

PROOF. From the theorem 3.7.4, it follows that if matrix $\mathcal{B}$ is singular and the mapping $f$ satisfies to the equation

$$(3.7.18) \qquad \operatorname{rank}\begin{pmatrix} \mathcal{B}^{\cdot k}_{m \cdot ij} & f^k_m \end{pmatrix} = \operatorname{rank} \mathcal{B}$$

then the mapping $f$ generates the same set of mappings that is generated by the mapping $\delta$. Therefore, to build the basis of representation of the algebra $A \otimes A$ in the module $\mathcal{L}(A;A)$, we must perform the following construction.

The set of solutions of system of equations (3.7.11) generates a free submodule $\mathcal{L}$ of the module $\mathcal{L}(A;A)$. We build the basis $(\overline{h}_1, ..., \overline{h}_k)$ of the submodule $\mathcal{L}$. Then we supplement this basis by linearly independent vectors $\overline{h}_{k+1}, ..., \overline{h}_m$, that do not belong to the submodule $\mathcal{L}$ so that the set of vectors $\overline{h}_1, ..., \overline{h}_m$ forms a basis of the module $\mathcal{L}(A;A)$. The set of orbits $(A \otimes A) \circ \delta$, $(A \otimes A) \circ \overline{h}_{k+1}$, ..., $(A \otimes A) \circ \overline{h}_m$ generates the module $\mathcal{L}(A;A)$. Since the set of orbits is finite, we can choose the orbits so that they do not intersect. For each orbit we can choose a representative which generates the orbit. □

**Example 3.7.6.** For complex field, the algebra $\mathcal{L}(C;C)$ has basis

$$I_{\mathbf{0}} \circ z = z$$
$$I_{\mathbf{1}} \circ z = \overline{z}$$

For quaternion algebra, the algebra $\mathcal{L}(H;H)$ has basis

$$I_{\mathbf{0}} \circ z = z$$

□

## 3.8. Linear Mapping into Nonassociative Algebra

Since the product is nonassociative, we may assume that action of $a$, $b \in A$ over the mapping $f$ may have form either $a(fb)$, or $(af)b$. However this assumption leads us to a rather complex structure of the linear mapping. To better understand how complex the structure of the linear mapping, we begin by considering the left and right shifts in nonassociative algebra.

**Theorem 3.8.1.** *Let*

$$(3.8.1) \qquad l(a) \circ x = ax$$

*be mapping of left shift. Than*

$$(3.8.2) \qquad l(a) \circ l(b) = l(ab) - (a,b)_1$$

*where we introduced linear mapping*

$$(a,b)_1 \circ x = (a,b,x)$$

PROOF. From the equations (3.2.2), (3.8.1), it follows that

$$(3.8.3) \qquad \begin{aligned} (l(a) \circ l(b)) \circ x &= l(a) \circ (l(b) \circ x) \\ &= a(bx) = (ab)x - (a,b,x) \\ &= l(ab) \circ x - (a,b)_1 \circ x \end{aligned}$$

The equation (3.8.2) follows from equation (3.8.3). □

**Theorem 3.8.2.** *Let*

(3.8.4) $$r(a) \circ x = xa$$

*be mapping of right shift. Than*

(3.8.5) $$r(a) \circ r(b) = r(ba) + (b, a)_2$$

*where we introduced linear mapping*

$$(b, a)_2 \circ x = (x, b, a)$$

PROOF. From the equations (3.2.2), (3.8.4) it follows that

(3.8.6) $$\begin{aligned}(r(a) \circ r(b)) \circ x &= r(a) \circ (r(b) \circ x) \\ &= (xb)a = x(ba) + (x, b, a) \\ &= r(ba) \circ x + (x, b, a)\end{aligned}$$

The equation (3.8.5) follows from equation (3.8.6). □

Let

$$f : A \to A \quad f = (ax)b$$

be linear mapping of the algebra $A$. According to the theorem 3.3.4, the mapping

$$g : A \to A \quad g = (cf)d$$

is also a linear mapping. However, it is not obvious whether we can write the mapping $g$ as a sum of terms of type $(ax)b$ and $a(xb)$.

If $A$ is free finite dimensional algebra, then we can assume that the linear mapping has the standard representation like[3.13]

(3.8.8) $$f \circ x = f^{ij} \, (\overline{e}_i x)\overline{e}_j$$

In this case we can use the theorem 3.7.5 for mappings into nonassociative algebra.

**Theorem 3.8.3.** *Let $\overline{\overline{e}}_1$ be basis of the free finite dimensional $D$-algebra $A_1$. Let $\overline{\overline{e}}_2$ be basis of the free finite dimensional nonassociative $D$-algebra $A_2$. Let $B_2{}_{\cdot kl}^{\;\;p}$ be structural constants of algebra $A_2$. Let the mapping*

(3.8.9) $$g = a \circ f$$

*generated by the mapping $f \in (A_1; A_2)$ through the tensor $a \in A_2 \otimes A_2$, has the standard representation*

(3.8.10) $$g = a^{ij}(\overline{e}_i \otimes \overline{e}_j) \circ f = a^{ij}(\overline{e}_i f)\overline{e}_j$$

*Coordinates of the mapping (3.8.9) and its standard components are connected by the equation*

(3.8.11) $$g_l^k = f_l^m g^{ij} B_2{}_{\cdot im}^{\;\;p} B_2{}_{\cdot pj}^{\;\;k}$$

---

[3.13]The choice is arbitrary. We may consider the standard representation like

$$f \circ x = f^{ij}\overline{e}_i(x\overline{e}_j)$$

Than the equation (3.8.11) has form

(3.8.7) $$g_l^k = f_l^m g^{ij} B_2{}_{\cdot ip}^{\;\;k} B_2{}_{\cdot mj}^{\;\;p}$$

I chose the expression (3.8.8) because order of the factors corresponds to the order chosen in the theorem 3.7.5.

PROOF. Relative to bases $\overline{\overline{e}}_1$ and $\overline{\overline{e}}_2$, linear mappings $f$ and $g$ have form

$$f \circ x = f^i_j x^j \overline{e}_{2\cdot i} \tag{3.8.12}$$

$$g \circ x = g^i_j x^j \overline{e}_{2\cdot i} \tag{3.8.13}$$

From equations (3.8.12), (3.8.13), (3.8.10) it follows that

$$\begin{aligned}g^k_l x^l \overline{e}_{2\cdot k} &= a^{ij}(\overline{e}_{2\cdot i}(f^m_l x^l \overline{e}_{2\cdot m}))\overline{e}_{2\cdot j} \\ &= a^{ij} f^m_l x^l B_2{\cdot}^p_{im} B_2{\cdot}^k_{pj} \overline{e}_{2\cdot k}\end{aligned} \tag{3.8.14}$$

Since vectors $\overline{e}_{2\cdot k}$ are linear independent and $x^i$ are arbitrary, than the equation (3.8.11) follows from the equation (3.8.14). $\square$

**Theorem 3.8.4.** *Let $A$ be free finite dimensional nonassociative algebra over the ring $D$. The representation of algebra $A \otimes A$ in algebra $\mathcal{L}(A; A)$ has finite basis $\overline{\overline{I}}$.*

(1) *The linear mapping $f \in \mathcal{L}(A; A)$ has form*

$$f = \sum_k (a_{s_k\cdot 0} \otimes a_{s_k\cdot 1}) \circ I_k = \sum_k (a_{s_k\cdot 0} I_k) a_{s_k\cdot 1} \tag{3.8.15}$$

(2) *Its standard representation has form*

$$f = a^{k\cdot ij}(\overline{e}_i \otimes \overline{e}_j) \circ I_k = a^{k\cdot ij}(\overline{e}_i I_k)\overline{e}_j \tag{3.8.16}$$

PROOF. Consider matrix (3.7.8). If matrix $\mathcal{B}$ is nonsingular, then, for given coordinates of linear transformation $g^l_k$ and for mapping $f = \delta$, the system of linear equations (3.8.11) with standard components of this transformation $g^{kr}$ has the unique solution. If matrix $\mathcal{B}$ is singular, then according to the theorem 3.7.5 there exists finite basis $\overline{\overline{I}}$ generating the set of linear mappings. $\square$

Unlike the case of an associative algebra, the set of generators $I$ in the theorem 3.8.4 is not minimal. From the equation (3.8.2) it follows that the equation (3.6.9) does not hold. Therefore, orbits of mappings $I_k$ do not generate an equivalence relation in the algebra $L(A; A)$. Since we consider only mappings like $(aI_k)b$, than it is possible that for $k \neq l$ the mapping $I_k$ generates the mapping $I_l$, if we consider all possible operations in the algebra $A$. Therefore, the set of generators $I_k$ of nonassociative algebra $A$ does not play such a critical role as conjugation in complex field. The answer to the question of how important it is the mapping $I_k$ in nonassociative algebra requires additional research.

CHAPTER 4

# Division Algebra

## 4.1. Linear Function of Complex Field

**Theorem 4.1.1** (the Cauchy-Riemann equations). *Let us consider complex field $C$ as two-dimensional algebra over real field. Let*

$$\overline{e}_{C \cdot \mathbf{0}} = 1 \quad \overline{e}_{C \cdot \mathbf{1}} = i \tag{4.1.1}$$

*be the basis of algebra $C$. Then in this basis product has form*

$$\overline{e}_{C \cdot \mathbf{1}}^2 = -\overline{e}_{C \cdot \mathbf{0}} \tag{4.1.2}$$

*and structural constants have form*

$$\begin{aligned} B_{C \cdot \mathbf{00}}^{\mathbf{0}} &= 1 & B_{C \cdot \mathbf{01}}^{\mathbf{1}} &= 1 \\ B_{C \cdot \mathbf{10}}^{\mathbf{1}} &= 1 & B_{C \cdot \mathbf{11}}^{\mathbf{0}} &= -1 \end{aligned} \tag{4.1.3}$$

*Matrix of linear function*

$$y^i = x^j f_j^i$$

*of complex field over real field satisfies relationship*

$$f_{\mathbf{0}}^{\mathbf{0}} = f_{\mathbf{1}}^{\mathbf{1}} \tag{4.1.4}$$

$$f_{\mathbf{0}}^{\mathbf{1}} = -f_{\mathbf{1}}^{\mathbf{0}} \tag{4.1.5}$$

PROOF. Equations (4.1.2) and (4.1.3) follow from equation $i^2 = -1$. Using equation [6]-(3.1.17) we get relationships

$$(4.1.6) \quad f_{\mathbf{0}}^{\mathbf{0}} = f^{kr} B_{C \cdot \mathbf{k0}}^{\mathbf{p}} B_{C \cdot \mathbf{pr}}^{\mathbf{0}} = f^{\mathbf{0}r} B_{C \cdot \mathbf{00}}^{\mathbf{0}} B_{C \cdot \mathbf{0r}}^{\mathbf{0}} + f^{\mathbf{1}r} B_{C \cdot \mathbf{10}}^{\mathbf{1}} B_{C \cdot \mathbf{1r}}^{\mathbf{0}} = f^{\mathbf{00}} - f^{\mathbf{11}}$$

$$(4.1.7) \quad f_{\mathbf{0}}^{\mathbf{1}} = f^{kr} B_{C \cdot \mathbf{k0}}^{\mathbf{p}} B_{C \cdot \mathbf{pr}}^{\mathbf{1}} = f^{\mathbf{0}r} B_{C \cdot \mathbf{00}}^{\mathbf{0}} B_{C \cdot \mathbf{0r}}^{\mathbf{1}} + f^{\mathbf{1}r} B_{C \cdot \mathbf{10}}^{\mathbf{1}} B_{C \cdot \mathbf{1r}}^{\mathbf{1}} = f^{\mathbf{01}} + f^{\mathbf{10}}$$

$$(4.1.8) \quad f_{\mathbf{1}}^{\mathbf{0}} = f^{kr} B_{C \cdot \mathbf{k1}}^{\mathbf{p}} B_{C \cdot \mathbf{pr}}^{\mathbf{0}} = f^{\mathbf{0}r} B_{C \cdot \mathbf{01}}^{\mathbf{1}} B_{C \cdot \mathbf{1r}}^{\mathbf{0}} + f^{\mathbf{1}r} B_{C \cdot \mathbf{11}}^{\mathbf{0}} B_{C \cdot \mathbf{0r}}^{\mathbf{0}} = -f^{\mathbf{01}} - f^{\mathbf{10}}$$

$$(4.1.9) \quad f_{\mathbf{1}}^{\mathbf{1}} = f^{kr} B_{C \cdot \mathbf{k1}}^{\mathbf{p}} B_{C \cdot \mathbf{pr}}^{\mathbf{1}} = f^{\mathbf{0}r} B_{C \cdot \mathbf{01}}^{\mathbf{1}} B_{C \cdot \mathbf{1r}}^{\mathbf{1}} + f^{\mathbf{1}r} B_{C \cdot \mathbf{11}}^{\mathbf{0}} B_{C \cdot \mathbf{0r}}^{\mathbf{1}} = f^{\mathbf{00}} - f^{\mathbf{11}}$$

(4.1.4) follows from equations (4.1.6) and (4.1.9). (4.1.5) follows from equations (4.1.7) and (4.1.8). □

## 4.2. Quaternion Algebra

In this paper I explore the set of quaternion algebras defined in [13].

**Definition 4.2.1.** Let $F$ be field. Extension field $F(i, j, k)$ is called **the quaternion algebra** $E(F, a, b)$ **over the field** $F$[4.1] if multiplication in algebra $E$ is defined according to rule

(4.2.1)
$$\begin{array}{c|ccc} & i & j & k \\ \hline i & a & k & aj \\ j & -k & b & -bi \\ k & -aj & bi & -ab \end{array}$$

where $a, b \in F$, $ab \neq 0$.

Elements of the algebra $E(F, a, b)$ have form
$$x = x^0 + x^1 i + x^2 j + x^3 k$$
where $x^i \in F$, $i = 0, 1, 2, 3$. Quaternion
$$\overline{x} = x^0 - x^1 i - x^2 j - x^3 k$$
is called conjugate to the quaternion $x$. We define **the norm of the quaternion** $x$ using equation

(4.2.2) $$|x|^2 = x\overline{x} = (x^0)^2 - a(x^1)^2 - b(x^2)^2 + ab(x^3)^2$$

From equation (4.2.2), it follows that $E(F, a, b)$ is algebra with division only when $a < 0$, $b < 0$. In this case we can renorm basis such that $a = -1$, $b = -1$.

We use symbol $E(F)$ to denote the quaternion division algebra $E(F, -1, -1)$ over the field $F$. Multiplication in algebra $E(F)$ is defined according to rule

(4.2.3)
$$\begin{array}{c|ccc} & i & j & k \\ \hline i & -1 & k & -j \\ j & -k & -1 & i \\ k & j & -i & -1 \end{array}$$

$\square$

In algebra $E(F)$, the norm of the quaternion has form
(4.2.4) $$|x|^2 = x\overline{x} = (x^0)^2 + (x^1)^2 + (x^2)^2 + (x^3)^2$$

In this case inverse element has form
(4.2.5) $$x^{-1} = |x|^{-2} \overline{x}$$

We will use notation $H = E(R, -1, -1)$.

The inner automorphism of quaternion algebra $H$[4.2]

(4.2.6)
$$p \to qpq^{-1}$$
$$q(ix + jy + kz)q^{-1} = ix' + jy' + kz'$$

---

[4.1]I follow definition from [13].
[4.2]See [15], p.643.

describes the rotation of the vector with coordinates $x$, $y$, $z$. The norm of quaternion $q$ is irrelevant, although usually we assume $|q| = 1$. If $q$ is written as sum of scalar and vector

$$q = \cos\alpha + (ia + jb + kc)\sin\alpha \quad a^2 + b^2 + c^2 = 1$$

then (4.2.6) is a rotation of the vector $(x, y, z)$ about the vector $(a, b, c)$ through an angle $2\alpha$.

## 4.3. Linear Function of Quaternion Algebra

**Theorem 4.3.1.** *Let*

(4.3.1) $$\overline{e}_0 = 1 \quad \overline{e}_1 = i \quad \overline{e}_2 = j \quad \overline{e}_3 = k$$

*be basis of quaternion algebra $H$. Then in the basis* (4.3.1), *structural constants have form*

$$\begin{aligned}
B^0_{00} &= 1 & B^1_{01} &= 1 & B^2_{02} &= 1 & B^3_{03} &= 1 \\
B^1_{10} &= 1 & B^0_{11} &= -1 & B^3_{12} &= 1 & B^2_{13} &= -1 \\
B^2_{20} &= 1 & B^3_{21} &= -1 & B^0_{22} &= -1 & B^1_{23} &= 1 \\
B^3_{30} &= 1 & B^2_{31} &= 1 & B^1_{32} &= -1 & B^0_{33} &= -1
\end{aligned}$$

PROOF. Value of structural constants follows from multiplication table (4.2.3). □

Since calculations in this section get a lot of space, I put in one place references to theorems in this section.

> **Theorem 4.3.2:** the definition of coordinates of linear mapping of quaternion algebra $H$ using standard components of this mapping.
> 
> **Equation (4.3.22):** matrix form of dependence of coordinates of linear mapping of quaternion algebra $H$ from standard components of this mapping.
> 
> **Equation (4.3.23):** matrix form of dependence of standard components of linear mapping of quaternion algebra $H$ from coordinates of this mapping.
> 
> **Theorem 4.3.4:** dependence standard components of linear mapping of quaternion algebra $H$ from coordinates of this mapping.

**Theorem 4.3.2.** *Standard components of linear function of quaternion algebra $H$ relative to basis* (4.3.1) *and coordinates of corresponding linear map satisfy relationship*

(4.3.2) $$\begin{cases} f^0_0 = f^{00} - f^{11} - f^{22} - f^{33} \\ f^1_1 = f^{00} - f^{11} + f^{22} + f^{33} \\ f^2_2 = f^{00} + f^{11} - f^{22} + f^{33} \\ f^3_3 = f^{00} + f^{11} + f^{22} - f^{33} \end{cases}$$

(4.3.3) $$\begin{cases} f^1_0 = \phantom{-}f^{01} + f^{10} + f^{23} - f^{32} \\ f^0_1 = -f^{01} - f^{10} + f^{23} - f^{32} \\ f^3_2 = -f^{01} + f^{10} - f^{23} - f^{32} \\ f^2_3 = \phantom{-}f^{01} - f^{10} - f^{23} - f^{32} \end{cases}$$

$$(4.3.4) \quad \begin{cases} f_0^2 = \phantom{-}f^{02} - f^{13} + f^{20} + f^{31} \\ f_1^3 = \phantom{-}f^{02} - f^{13} - f^{20} - f^{31} \\ f_2^0 = -f^{02} - f^{13} - f^{20} + f^{31} \\ f_3^1 = -f^{02} - f^{13} + f^{20} - f^{31} \end{cases}$$

$$(4.3.5) \quad \begin{cases} f_0^3 = \phantom{-}f^{03} + f^{12} - f^{21} + f^{30} \\ f_1^2 = -f^{03} - f^{12} - f^{21} + f^{30} \\ f_2^1 = \phantom{-}f^{03} - f^{12} - f^{21} - f^{30} \\ f_3^0 = -f^{03} + f^{12} - f^{21} - f^{30} \end{cases}$$

PROOF. Using equation (3.7.11) we get relationships

$$(4.3.6) \quad \begin{aligned} f_0^0 &= f^{kr} B_{k0}^p B_{pr}^0 \\ &= f^{00} B_{00}^0 B_{00}^0 + f^{11} B_{10}^1 B_{11}^0 + f^{22} B_{20}^2 B_{22}^0 + f^{33} B_{30}^3 B_{33}^0 \\ &= f^{00} - f^{11} - f^{22} - f^{33} \end{aligned}$$

$$(4.3.7) \quad \begin{aligned} f_0^1 &= f^{kr} B_{k0}^p B_{pr}^1 \\ &= f^{01} B_{00}^0 B_{01}^1 + f^{10} B_{10}^1 B_{10}^1 + f^{23} B_{20}^2 B_{23}^1 + f^{32} B_{30}^3 B_{32}^1 \\ &= f^{01} + f^{10} + f^{23} - f^{32} \end{aligned}$$

$$(4.3.8) \quad \begin{aligned} f_0^2 &= f^{kr} B_{k0}^p B_{pr}^2 \\ &= f^{02} B_{00}^0 B_{02}^2 + f^{13} B_{10}^1 B_{13}^2 + f^{20} B_{20}^2 B_{20}^2 + f^{31} B_{30}^3 B_{31}^2 \\ &= f^{02} - f^{13} + f^{20} + f^{31} \end{aligned}$$

$$(4.3.9) \quad \begin{aligned} f_0^3 &= f^{kr} B_{k0}^p B_{pr}^3 \\ &= f^{03} B_{00}^0 B_{03}^3 + f^{12} B_{10}^1 B_{12}^3 + f^{21} B_{20}^2 B_{21}^3 + f^{30} B_{30}^3 B_{30}^3 \\ &= f^{03} + f^{12} - f^{21} + f^{30} \end{aligned}$$

$$(4.3.10) \quad \begin{aligned} f_1^0 &= f^{kr} B_{k1}^p B_{pr}^0 \\ &= f^{01} B_{01}^1 B_{11}^0 + f^{10} B_{11}^0 B_{00}^0 + f^{23} B_{21}^3 B_{33}^0 + f^{32} B_{31}^2 B_{22}^0 \\ &= -f^{01} - f^{10} + f^{23} - f^{32} \end{aligned}$$

$$(4.3.11) \quad \begin{aligned} f_1^1 &= f^{kr} B_{k1}^p B_{pr}^1 \\ &= f^{00} B_{01}^1 B_{10}^1 + f^{11} B_{11}^0 B_{01}^1 + f^{22} B_{21}^3 B_{32}^1 + f^{33} B_{31}^2 B_{23}^1 \\ &= f^{00} - f^{11} + f^{22} + f^{33} \end{aligned}$$

$$(4.3.12) \quad \begin{aligned} f_1^2 &= f^{kr} B_{k1}^p B_{pr}^2 \\ &= f^{03} B_{01}^1 B_{13}^2 + f^{12} B_{11}^0 B_{02}^2 + f^{21} B_{21}^3 B_{31}^2 + f^{30} B_{31}^2 B_{20}^2 \\ &= -f^{03} - f^{12} - f^{21} + f^{30} \end{aligned}$$

$$f_1^3 = f^{kr} B_{k1}^p B_{pr}^3$$
(4.3.13)
$$= f^{02} B_{01}^1 B_{12}^3 + f^{13} B_{11}^0 B_{03}^3 + f^{20} B_{21}^3 B_{30}^3 + f^{31} B_{31}^2 B_{21}^3$$
$$= f^{02} - f^{13} - f^{20} - f^{31}$$

$$f_2^0 = f^{kr} B_{k2}^p B_{pr}^0$$
(4.3.14)
$$= f^{02} B_{02}^2 B_{22}^0 + f^{13} B_{12}^3 B_{33}^0 + f^{20} B_{22}^0 B_{00}^0 + f^{31} B_{32}^1 B_{11}^0$$
$$= -f^{02} - f^{13} - f^{20} + f^{31}$$

$$f_2^1 = f^{kr} B_{k2}^p B_{pr}^1$$
(4.3.15)
$$= f^{03} B_{02}^2 B_{23}^1 + f^{12} B_{12}^3 B_{32}^1 + f^{21} B_{22}^0 B_{01}^1 + f^{30} B_{32}^1 B_{10}^1$$
$$= f^{03} - f^{12} - f^{21} - f^{30}$$

$$f_2^2 = f^{kr} B_{k2}^p B_{pr}^2$$
(4.3.16)
$$= f^{00} B_{02}^2 B_{20}^2 + f^{11} B_{12}^3 B_{31}^2 + f^{22} B_{22}^0 B_{02}^2 + f^{33} B_{32}^1 B_{13}^2$$
$$= f^{00} + f^{11} - f^{22} + f^{33}$$

$$f_2^3 = f^{kr} B_{k2}^p B_{pr}^3$$
(4.3.17)
$$= f^{01} B_{02}^2 B_{21}^3 + f^{10} B_{12}^3 B_{30}^3 + f^{23} B_{22}^0 B_{03}^3 + f^{32} B_{32}^1 B_{12}^3$$
$$= -f^{01} + f^{10} - f^{23} - f^{32}$$

$$f_3^0 = f^{kr} B_{k3}^p B_{pr}^0$$
(4.3.18)
$$= f^{03} B_{03}^3 B_{33}^0 + f^{12} B_{13}^2 B_{22}^0 + f^{21} B_{23}^1 B_{11}^0 + f^{30} B_{33}^0 B_{00}^0$$
$$= -f^{03} + f^{12} - f^{21} - f^{30}$$

$$f_3^1 = f^{kr} B_{k3}^p B_{pr}^1$$
(4.3.19)
$$= f^{02} B_{03}^3 B_{32}^1 + f^{13} B_{13}^2 B_{23}^1 + f^{20} B_{23}^1 B_{10}^1 + f^{31} B_{33}^0 B_{01}^1$$
$$= -f^{02} - f^{13} + f^{20} - f^{31}$$

$$f_3^2 = f^{kr} B_{k3}^p B_{pr}^2$$
(4.3.20)
$$= f^{01} B_{03}^3 B_{31}^2 + f^{10} B_{13}^2 B_{20}^2 + f^{23} B_{23}^1 B_{13}^2 + f^{32} B_{33}^0 B_{02}^2$$
$$= f^{01} - f^{10} - f^{23} - f^{32}$$

$$f_3^3 = f^{kr} B_{k3}^p B_{pr}^3$$
(4.3.21)
$$= f^{00} B_{03}^3 B_{30}^3 + f^{11} B_{13}^2 B_{21}^3 + f^{22} B_{23}^1 B_{12}^3 + f^{33} B_{33}^0 B_{03}^3$$
$$= f^{00} + f^{11} + f^{22} - f^{33}$$

Equations (4.3.6), (4.3.11), (4.3.16), (4.3.21) form the system of linear equations (4.3.2).

Equations (4.3.7), (4.3.10), (4.3.17), (4.3.20) form the system of linear equations (4.3.3).

Equations (4.3.8), (4.3.13), (4.3.14), (4.3.19) form the system of linear equations (4.3.4).

Equations (4.3.9), (4.3.12), (4.3.15), (4.3.18) form the system of linear equations (4.3.5). □

**Theorem 4.3.3.** *Consider quaternion algebra $H$ with the basis (4.3.1). Standard components of additive function over field $F$ and coordinates of this function over field $F$ satisfy relationship*

$$
(4.3.22) \quad \begin{pmatrix} f_0^0 & f_1^0 & f_2^0 & f_3^0 \\ f_1^1 & -f_0^1 & f_3^1 & -f_2^1 \\ f_2^2 & -f_3^2 & -f_0^2 & f_1^2 \\ f_3^3 & f_2^3 & -f_1^3 & -f_0^3 \end{pmatrix} = \begin{pmatrix} 1 & -1 & -1 & -1 \\ 1 & -1 & 1 & 1 \\ 1 & 1 & -1 & 1 \\ 1 & 1 & 1 & -1 \end{pmatrix} \begin{pmatrix} f^{00} & -f^{01} & -f^{02} & -f^{03} \\ f^{11} & f^{10} & f^{13} & -f^{12} \\ f^{22} & -f^{23} & f^{20} & f^{21} \\ f^{33} & f^{32} & -f^{31} & f^{30} \end{pmatrix}
$$

$$
(4.3.23) \quad \begin{pmatrix} f^{00} & -f^{01} & -f^{02} & -f^{03} \\ f^{11} & f^{10} & f^{13} & -f^{12} \\ f^{22} & -f^{23} & f^{20} & f^{21} \\ f^{33} & f^{32} & -f^{31} & f^{30} \end{pmatrix} = \frac{1}{4} \begin{pmatrix} 1 & 1 & 1 & 1 \\ -1 & -1 & 1 & 1 \\ -1 & 1 & -1 & 1 \\ -1 & 1 & 1 & -1 \end{pmatrix} \begin{pmatrix} f_0^0 & f_1^0 & f_2^0 & f_3^0 \\ f_1^1 & -f_0^1 & f_3^1 & -f_2^1 \\ f_2^2 & -f_3^2 & -f_0^2 & f_1^2 \\ f_3^3 & f_2^3 & -f_1^3 & -f_0^3 \end{pmatrix}
$$

*where*

$$
\begin{pmatrix} 1 & -1 & -1 & -1 \\ 1 & -1 & 1 & 1 \\ 1 & 1 & -1 & 1 \\ 1 & 1 & 1 & -1 \end{pmatrix}^{-1} = \frac{1}{4} \begin{pmatrix} 1 & 1 & 1 & 1 \\ -1 & -1 & 1 & 1 \\ -1 & 1 & -1 & 1 \\ -1 & 1 & 1 & -1 \end{pmatrix}
$$

PROOF. Let us write the system of linear equations (4.3.2) as product of matrices

$$
(4.3.24) \quad \begin{pmatrix} f_0^0 \\ f_1^1 \\ f_2^2 \\ f_3^3 \end{pmatrix} = \begin{pmatrix} 1 & -1 & -1 & -1 \\ 1 & -1 & 1 & 1 \\ 1 & 1 & -1 & 1 \\ 1 & 1 & 1 & -1 \end{pmatrix} \begin{pmatrix} f^{00} \\ f^{11} \\ f^{22} \\ f^{33} \end{pmatrix}
$$

Let us write the system of linear equations (4.3.3) as product of matrices

$$
(4.3.25) \qquad
\begin{pmatrix} f^1_0 \\ f^0_1 \\ f^3_2 \\ f^2_3 \end{pmatrix}
=
\begin{pmatrix}
1 & 1 & 1 & -1 \\
-1 & -1 & 1 & -1 \\
-1 & 1 & -1 & -1 \\
1 & -1 & -1 & -1
\end{pmatrix}
\begin{pmatrix} f^{01} \\ f^{10} \\ f^{23} \\ f^{32} \end{pmatrix}
$$

From the equation (4.3.25), it follows that

$$
\begin{pmatrix} -f^1_0 \\ f^0_1 \\ f^3_2 \\ -f^2_3 \end{pmatrix}
=
\begin{pmatrix}
-1 & -1 & -1 & 1 \\
-1 & -1 & 1 & -1 \\
-1 & 1 & -1 & -1 \\
-1 & 1 & 1 & 1
\end{pmatrix}
\begin{pmatrix} f^{01} \\ f^{10} \\ f^{23} \\ f^{32} \end{pmatrix}
$$

$$
=
\begin{pmatrix}
1 & -1 & 1 & 1 \\
1 & -1 & -1 & -1 \\
1 & 1 & 1 & -1 \\
1 & 1 & -1 & 1
\end{pmatrix}
\begin{pmatrix} -f^{01} \\ f^{10} \\ -f^{23} \\ f^{32} \end{pmatrix}
$$

$$
(4.3.26) \qquad
\begin{pmatrix} f^0_1 \\ -f^1_0 \\ -f^2_3 \\ f^3_2 \end{pmatrix}
=
\begin{pmatrix}
1 & -1 & -1 & -1 \\
1 & -1 & 1 & 1 \\
1 & 1 & -1 & 1 \\
1 & 1 & 1 & -1
\end{pmatrix}
\begin{pmatrix} -f^{01} \\ f^{10} \\ -f^{23} \\ f^{32} \end{pmatrix}
$$

Let us write the system of linear equations (4.3.4) as product of matrices

$$
(4.3.27) \qquad
\begin{pmatrix} f^2_0 \\ f^3_1 \\ f^0_2 \\ f^1_3 \end{pmatrix}
=
\begin{pmatrix}
1 & -1 & 1 & 1 \\
1 & -1 & -1 & -1 \\
-1 & -1 & -1 & 1 \\
-1 & -1 & 1 & -1
\end{pmatrix}
\begin{pmatrix} f^{02} \\ f^{13} \\ f^{20} \\ f^{31} \end{pmatrix}
$$

From the equation (4.3.27), it follows that

$$
\begin{pmatrix} -f^2_0 \\ -f^3_1 \\ f^0_2 \\ f^1_3 \end{pmatrix}
=
\begin{pmatrix}
-1 & 1 & -1 & -1 \\
-1 & 1 & 1 & 1 \\
-1 & -1 & -1 & 1 \\
-1 & -1 & 1 & -1
\end{pmatrix}
\begin{pmatrix} f^{02} \\ f^{13} \\ f^{20} \\ f^{31} \end{pmatrix}
$$

$$
=
\begin{pmatrix}
1 & 1 & -1 & 1 \\
1 & 1 & 1 & -1 \\
1 & -1 & -1 & -1 \\
1 & -1 & 1 & 1
\end{pmatrix}
\begin{pmatrix} -f^{02} \\ f^{13} \\ f^{20} \\ -f^{31} \end{pmatrix}
$$

$$
(4.3.28) \qquad
\begin{pmatrix} f^0_2 \\ f^1_3 \\ -f^2_0 \\ -f^3_1 \end{pmatrix}
=
\begin{pmatrix}
1 & -1 & -1 & -1 \\
1 & -1 & 1 & 1 \\
1 & 1 & -1 & 1 \\
1 & 1 & 1 & -1
\end{pmatrix}
\begin{pmatrix} -f^{02} \\ f^{13} \\ f^{20} \\ -f^{31} \end{pmatrix}
$$

Let us write the system of linear equations (4.3.5) as product of matrices

$$(4.3.29) \quad \begin{pmatrix} f^3_0 \\ f^2_1 \\ f^1_2 \\ f^0_3 \end{pmatrix} = \begin{pmatrix} 1 & 1 & -1 & 1 \\ -1 & -1 & -1 & 1 \\ 1 & -1 & -1 & -1 \\ -1 & 1 & -1 & -1 \end{pmatrix} \begin{pmatrix} f^{03} \\ f^{12} \\ f^{21} \\ f^{30} \end{pmatrix}$$

From the equation (4.3.29), it follows that

$$\begin{pmatrix} -f^3_0 \\ f^2_1 \\ -f^1_2 \\ f^0_3 \end{pmatrix} = \begin{pmatrix} -1 & -1 & 1 & -1 \\ -1 & -1 & -1 & 1 \\ -1 & 1 & 1 & 1 \\ -1 & 1 & -1 & -1 \end{pmatrix} \begin{pmatrix} f^{03} \\ f^{12} \\ f^{21} \\ f^{30} \end{pmatrix}$$

$$= \begin{pmatrix} 1 & 1 & 1 & -1 \\ 1 & 1 & -1 & 1 \\ 1 & -1 & 1 & 1 \\ 1 & -1 & -1 & -1 \end{pmatrix} \begin{pmatrix} -f^{03} \\ -f^{12} \\ f^{21} \\ f^{30} \end{pmatrix}$$

$$(4.3.30) \quad \begin{pmatrix} f^0_3 \\ -f^1_2 \\ f^2_1 \\ -f^3_0 \end{pmatrix} = \begin{pmatrix} 1 & -1 & -1 & -1 \\ 1 & -1 & 1 & 1 \\ 1 & 1 & -1 & 1 \\ 1 & 1 & 1 & -1 \end{pmatrix} \begin{pmatrix} -f^{03} \\ -f^{12} \\ f^{21} \\ f^{30} \end{pmatrix}$$

We join equations (4.3.24), (4.3.26), (4.3.28), (4.3.30) into equation (4.3.22). □

**Theorem 4.3.4.** *Standard components of linear function of quaternion algebra H relative to basis (4.3.1) and coordinates of corresponding linear map satisfy relationship*

$$(4.3.31) \quad \begin{cases} 4f^{00} = f^0_0 + f^1_1 + f^2_2 + f^3_3 \\ 4f^{11} = -f^0_0 - f^1_1 + f^2_2 + f^3_3 \\ 4f^{22} = -f^0_0 + f^1_1 - f^2_2 + f^3_3 \\ 4f^{33} = -f^0_0 + f^1_1 + f^2_2 - f^3_3 \end{cases}$$

$$(4.3.32) \quad \begin{cases} 4f^{10} = -f^0_1 + f^1_0 - f^2_3 + f^3_2 \\ 4f^{01} = -f^0_1 + f^1_0 + f^2_3 - f^3_2 \\ 4f^{32} = -f^0_1 - f^1_0 - f^2_3 - f^3_2 \\ 4f^{23} = f^0_1 + f^1_0 - f^2_3 - f^3_2 \end{cases}$$

$$(4.3.33) \quad \begin{cases} 4f^{20} = -f^0_2 + f^1_3 + f^2_0 - f^3_1 \\ 4f^{31} = f^0_2 - f^1_3 + f^2_0 - f^3_1 \\ 4f^{02} = -f^0_2 - f^1_3 + f^2_0 + f^3_1 \\ 4f^{13} = -f^0_2 - f^1_3 - f^2_0 - f^3_1 \end{cases}$$

(4.3.34)
$$\begin{cases} 4f^{30} = -f_3^0 - f_2^1 + f_1^2 + f_0^3 \\ 4f^{21} = -f_3^0 - f_2^1 - f_1^2 - f_0^3 \\ 4f^{12} = f_3^0 - f_2^1 - f_1^2 + f_0^3 \\ 4f^{03} = -f_3^0 + f_2^1 - f_1^2 + f_0^3 \end{cases}$$

PROOF. We get systems of linear equations (4.3.31), (4.3.32), (4.3.33), (4.3.34) as the product of matrices in equation (4.3.23). □

## 4.4. Octonion Algebra

**Definition 4.4.1.** The algebra $O$ is called **octonion algebra** if algebra has basis

(4.4.1)
$$\begin{array}{cccc} \overline{e}_0 = 1 & \overline{e}_1 = i & \overline{e}_2 = j & \overline{e}_3 = k \\ \overline{e}_4 = -l & \overline{e}_5 = il & \overline{e}_6 = jl & \overline{e}_7 = kl \end{array}$$

and multiplication in algebra $O$ is defined according to rule

(4.4.2)

|  | $\overline{e}_1$ | $\overline{e}_2$ | $\overline{e}_3$ | $\overline{e}_4$ | $\overline{e}_5$ | $\overline{e}_6$ | $\overline{e}_7$ |
|---|---|---|---|---|---|---|---|
| $\overline{e}_1$ | $-\overline{e}_0$ | $\overline{e}_3$ | $-\overline{e}_2$ | $\overline{e}_5$ | $-\overline{e}_4$ | $-\overline{e}_7$ | $\overline{e}_6$ |
| $\overline{e}_2$ | $-\overline{e}_3$ | $-\overline{e}_0$ | $\overline{e}_1$ | $\overline{e}_6$ | $\overline{e}_7$ | $-\overline{e}_4$ | $-\overline{e}_5$ |
| $\overline{e}_3$ | $\overline{e}_2$ | $-\overline{e}_1$ | $-\overline{e}_0$ | $\overline{e}_7$ | $-\overline{e}_6$ | $\overline{e}_5$ | $-\overline{e}_4$ |
| $\overline{e}_4$ | $-\overline{e}_5$ | $-\overline{e}_6$ | $-\overline{e}_7$ | $-\overline{e}_0$ | $\overline{e}_1$ | $\overline{e}_2$ | $\overline{e}_3$ |
| $\overline{e}_5$ | $\overline{e}_4$ | $-\overline{e}_7$ | $\overline{e}_6$ | $-\overline{e}_1$ | $-\overline{e}_0$ | $-\overline{e}_3$ | $\overline{e}_2$ |
| $\overline{e}_6$ | $\overline{e}_7$ | $\overline{e}_4$ | $-\overline{e}_5$ | $-\overline{e}_2$ | $\overline{e}_3$ | $-\overline{e}_0$ | $-\overline{e}_1$ |
| $\overline{e}_7$ | $-\overline{e}_6$ | $\overline{e}_5$ | $\overline{e}_4$ | $-\overline{e}_3$ | $-\overline{e}_2$ | $\overline{e}_1$ | $-\overline{e}_0$ |

□

**Theorem 4.4.2.** *Structural constants of octonion algebra $O$ relative to basis* (4.4.1) *have form*

$$
\begin{aligned}
&B^0_{00} = \phantom{-}1 \quad B^1_{01} = \phantom{-}1 \quad B^2_{02} = \phantom{-}1 \quad B^3_{03} = \phantom{-}1 \\
&B^4_{04} = \phantom{-}1 \quad B^5_{05} = \phantom{-}1 \quad B^6_{06} = \phantom{-}1 \quad B^7_{07} = \phantom{-}1 \\
&B^1_{10} = \phantom{-}1 \quad B^0_{11} = -1 \quad B^3_{12} = \phantom{-}1 \quad B^2_{13} = -1 \\
&B^5_{14} = \phantom{-}1 \quad B^4_{15} = -1 \quad B^7_{16} = -1 \quad B^6_{17} = \phantom{-}1 \\
&B^2_{20} = \phantom{-}1 \quad B^3_{21} = -1 \quad B^0_{22} = -1 \quad B^1_{23} = \phantom{-}1 \\
&B^6_{24} = \phantom{-}1 \quad B^7_{25} = \phantom{-}1 \quad B^4_{26} = -1 \quad B^5_{27} = -1 \\
&B^3_{30} = \phantom{-}1 \quad B^2_{31} = \phantom{-}1 \quad B^1_{32} = -1 \quad B^0_{33} = -1 \\
&B^7_{34} = \phantom{-}1 \quad B^6_{35} = -1 \quad B^5_{36} = \phantom{-}1 \quad B^4_{37} = -1 \\
&B^4_{40} = \phantom{-}1 \quad B^5_{41} = -1 \quad B^6_{42} = -1 \quad B^7_{43} = -1 \\
&B^0_{44} = -1 \quad B^1_{45} = \phantom{-}1 \quad B^2_{46} = \phantom{-}1 \quad B^3_{47} = \phantom{-}1 \\
&B^5_{50} = \phantom{-}1 \quad B^4_{51} = \phantom{-}1 \quad B^7_{52} = -1 \quad B^6_{53} = \phantom{-}1 \\
&B^1_{54} = -1 \quad B^0_{55} = -1 \quad B^3_{56} = -1 \quad B^2_{57} = \phantom{-}1 \\
&B^6_{60} = \phantom{-}1 \quad B^7_{61} = \phantom{-}1 \quad B^4_{62} = \phantom{-}1 \quad B^5_{63} = -1 \\
&B^2_{64} = -1 \quad B^3_{65} = \phantom{-}1 \quad B^0_{66} = -1 \quad B^1_{67} = -1 \\
&B^7_{70} = \phantom{-}1 \quad B^6_{71} = -1 \quad B^5_{72} = \phantom{-}1 \quad B^4_{73} = \phantom{-}1 \\
&B^3_{74} = -1 \quad B^2_{75} = -1 \quad B^1_{76} = \phantom{-}1 \quad B^0_{77} = -1
\end{aligned}
$$

PROOF. Value of structural constants follows from multiplication table (4.4.2). □

## 4.5. Linear Function of Octonion Algebra

Since calculations in this section get a lot of space, I put in one place references to theorems in this section.

**Theorem 4.5.1:** the definition of coordinates of linear mapping of octonion algebra $O$ using standard components of this mapping.

**Equation** (4.5.73): matrix form of dependence of coordinates of linear mapping of octonion algebra $O$ from standard components of this mapping.

**Equation** (4.5.74): matrix form of dependence of standard components of linear mapping of octonion algebra $O$ from coordinates of this mapping.

**Theorem 4.5.3:** dependence standard components of linear mapping of octonion algebra $O$ from coordinates of this mapping.

**Theorem 4.5.1.** *Standard components of linear function of octonion algebra O relative to basis* (4.4.1) *and coordinates of corresponding linear map satisfy relationship*

$$(4.5.1) \begin{cases} f_0^0 = f^{00} - f^{11} - f^{22} - f^{33} - f^{44} - f^{55} - f^{66} - f^{77} \\ f_1^1 = f^{00} - f^{11} + f^{22} + f^{33} + f^{44} + f^{55} + f^{66} + f^{77} \\ f_2^2 = f^{00} + f^{11} - f^{22} + f^{33} + f^{44} + f^{55} + f^{66} + f^{77} \\ f_3^3 = f^{00} + f^{11} + f^{22} - f^{33} + f^{44} + f^{55} + f^{66} + f^{77} \\ f_4^4 = f^{00} + f^{11} + f^{22} + f^{33} - f^{44} + f^{55} + f^{66} + f^{77} \\ f_5^5 = f^{00} + f^{11} + f^{22} + f^{33} + f^{44} - f^{55} + f^{66} + f^{77} \\ f_6^6 = f^{00} + f^{11} + f^{22} + f^{33} + f^{44} + f^{55} - f^{66} + f^{77} \\ f_7^7 = f^{00} + f^{11} + f^{22} + f^{33} + f^{44} + f^{55} + f^{66} - f^{77} \end{cases}$$

$$(4.5.2) \begin{cases} f_0^1 = \phantom{-}f^{01} + f^{10} + f^{23} - f^{32} + f^{45} - f^{54} - f^{67} + f^{76} \\ f_1^0 = -f^{01} - f^{10} + f^{23} - f^{32} + f^{45} - f^{54} - f^{67} + f^{76} \\ f_2^3 = -f^{01} + f^{10} - f^{23} - f^{32} - f^{45} + f^{54} + f^{67} - f^{76} \\ f_3^2 = \phantom{-}f^{01} - f^{10} - f^{23} - f^{32} + f^{45} - f^{54} - f^{67} + f^{76} \\ f_4^5 = -f^{01} + f^{10} - f^{23} + f^{32} - f^{45} - f^{54} + f^{67} - f^{76} \\ f_5^4 = \phantom{-}f^{01} - f^{10} + f^{23} - f^{32} - f^{45} - f^{54} - f^{67} + f^{76} \\ f_6^7 = \phantom{-}f^{01} - f^{10} + f^{23} - f^{32} + f^{45} - f^{54} - f^{67} - f^{76} \\ f_7^6 = -f^{01} + f^{10} - f^{23} + f^{32} - f^{45} + f^{54} - f^{67} - f^{76} \end{cases}$$

$$(4.5.3) \begin{cases} f_0^2 = \phantom{-}f^{02} - f^{13} + f^{20} + f^{31} + f^{46} + f^{57} - f^{64} - f^{75} \\ f_1^3 = \phantom{-}f^{02} - f^{13} - f^{20} - f^{31} + f^{46} + f^{57} - f^{64} - f^{75} \\ f_2^0 = -f^{02} - f^{13} - f^{20} + f^{31} + f^{46} + f^{57} - f^{64} - f^{75} \\ f_3^1 = -f^{02} - f^{13} + f^{20} - f^{31} - f^{46} - f^{57} + f^{64} + f^{75} \\ f_4^6 = -f^{02} + f^{13} + f^{20} - f^{31} - f^{46} - f^{57} - f^{64} + f^{75} \\ f_5^7 = -f^{02} + f^{13} + f^{20} - f^{31} - f^{46} - f^{57} + f^{64} - f^{75} \\ f_6^4 = \phantom{-}f^{02} - f^{13} - f^{20} + f^{31} - f^{46} + f^{57} - f^{64} - f^{75} \\ f_7^5 = \phantom{-}f^{02} - f^{13} - f^{20} + f^{31} + f^{46} - f^{57} - f^{64} - f^{75} \end{cases}$$

$$(4.5.4)\begin{cases} f_0^3 = \phantom{-}f^{03} + f^{12} - f^{21} + f^{30} + f^{47} - f^{56} + f^{65} - f^{74} \\ f_1^2 = -f^{03} - f^{12} - f^{21} + f^{30} - f^{47} + f^{56} - f^{65} + f^{74} \\ f_2^1 = \phantom{-}f^{03} - f^{12} - f^{21} - f^{30} + f^{47} - f^{56} + f^{65} - f^{74} \\ f_3^0 = -f^{03} + f^{12} - f^{21} - f^{30} + f^{47} - f^{56} + f^{65} - f^{74} \\ f_4^7 = -f^{03} - f^{12} + f^{21} + f^{30} - f^{47} + f^{56} - f^{65} - f^{74} \\ f_5^6 = \phantom{-}f^{03} + f^{12} - f^{21} - f^{30} + f^{47} - f^{56} - f^{65} - f^{74} \\ f_6^5 = -f^{03} - f^{12} + f^{21} + f^{30} - f^{47} - f^{56} - f^{65} + f^{74} \\ f_7^4 = \phantom{-}f^{03} + f^{12} - f^{21} - f^{30} - f^{47} - f^{56} + f^{65} - f^{74} \end{cases}$$

$$(4.5.5)\begin{cases} f_0^4 = \phantom{-}f^{04} - f^{15} - f^{26} - f^{37} + f^{40} + f^{51} + f^{62} + f^{73} \\ f_1^5 = \phantom{-}f^{04} - f^{15} - f^{26} - f^{37} - f^{40} - f^{51} + f^{62} + f^{73} \\ f_2^6 = \phantom{-}f^{04} - f^{15} - f^{26} - f^{37} - f^{40} + f^{51} - f^{62} + f^{73} \\ f_3^7 = \phantom{-}f^{04} - f^{15} - f^{26} - f^{37} - f^{40} + f^{51} + f^{62} - f^{73} \\ f_4^0 = -f^{04} - f^{15} - f^{26} - f^{37} - f^{40} + f^{51} + f^{62} + f^{73} \\ f_5^1 = -f^{04} - f^{15} + f^{26} + f^{37} + f^{40} - f^{51} - f^{62} - f^{73} \\ f_6^2 = -f^{04} + f^{15} - f^{26} + f^{37} + f^{40} - f^{51} - f^{62} - f^{73} \\ f_7^3 = -f^{04} + f^{15} + f^{26} - f^{37} + f^{40} - f^{51} - f^{62} - f^{73} \end{cases}$$

$$(4.5.6)\begin{cases} f_0^5 = \phantom{-}f^{05} + f^{14} - f^{27} + f^{36} - f^{41} + f^{50} - f^{63} + f^{72} \\ f_1^4 = -f^{05} - f^{14} + f^{27} - f^{36} - f^{41} + f^{50} + f^{63} - f^{72} \\ f_2^7 = \phantom{-}f^{05} + f^{14} - f^{27} + f^{36} - f^{41} - f^{50} - f^{63} - f^{72} \\ f_3^6 = -f^{05} - f^{14} + f^{27} - f^{36} + f^{41} + f^{50} - f^{63} - f^{72} \\ f_4^1 = \phantom{-}f^{05} - f^{14} - f^{27} + f^{36} - f^{41} - f^{50} - f^{63} + f^{72} \\ f_5^0 = -f^{05} + f^{14} - f^{27} + f^{36} - f^{41} - f^{50} - f^{63} + f^{72} \\ f_6^3 = \phantom{-}f^{05} + f^{14} - f^{27} - f^{36} - f^{41} - f^{50} - f^{63} + f^{72} \\ f_7^2 = -f^{05} - f^{14} - f^{27} - f^{36} + f^{41} + f^{50} + f^{63} - f^{72} \end{cases}$$

$$(4.5.7) \begin{cases} f_0^6 = \phantom{-}f^{06} + f^{17} + f^{24} - f^{35} - f^{42} + f^{53} + f^{60} - f^{71} \\ f_1^7 = -f^{06} - f^{17} - f^{24} + f^{35} + f^{42} - f^{53} + f^{60} - f^{71} \\ f_2^4 = -f^{06} - f^{17} - f^{24} + f^{35} - f^{42} - f^{53} + f^{60} + f^{71} \\ f_3^5 = \phantom{-}f^{06} + f^{17} + f^{24} - f^{35} - f^{42} - f^{53} - f^{60} - f^{71} \\ f_4^2 = \phantom{-}f^{06} + f^{17} - f^{24} - f^{35} - f^{42} + f^{53} - f^{60} - f^{71} \\ f_5^3 = -f^{06} - f^{17} - f^{24} - f^{35} + f^{42} - f^{53} + f^{60} + f^{71} \\ f_6^0 = -f^{06} + f^{17} + f^{24} - f^{35} - f^{42} + f^{53} - f^{60} - f^{71} \\ f_7^1 = \phantom{-}f^{06} - f^{17} + f^{24} - f^{35} - f^{42} + f^{53} - f^{60} - f^{71} \end{cases}$$

$$(4.5.8) \begin{cases} f_0^7 = \phantom{-}f^{07} - f^{16} + f^{25} + f^{34} - f^{43} - f^{52} + f^{61} + f^{70} \\ f_1^6 = \phantom{-}f^{07} - f^{16} + f^{25} + f^{34} - f^{43} - f^{52} - f^{61} - f^{70} \\ f_2^5 = -f^{07} + f^{16} - f^{25} - f^{34} + f^{43} - f^{52} - f^{61} + f^{70} \\ f_3^4 = -f^{07} + f^{16} - f^{25} - f^{34} - f^{43} + f^{52} - f^{61} + f^{70} \\ f_4^3 = \phantom{-}f^{07} - f^{16} + f^{25} - f^{34} - f^{43} - f^{52} + f^{61} - f^{70} \\ f_5^2 = \phantom{-}f^{07} - f^{16} - f^{25} + f^{34} - f^{43} - f^{52} + f^{61} - f^{70} \\ f_6^1 = -f^{07} - f^{16} - f^{25} - f^{34} + f^{43} + f^{52} - f^{61} + f^{70} \\ f_7^0 = -f^{07} - f^{16} + f^{25} + f^{34} - f^{43} - f^{52} + f^{61} - f^{70} \end{cases}$$

PROOF. Using equation (3.7.11) we get relationships

$$(4.5.9) \begin{aligned} f_0^0 &= f^{kr} B_{k0}^p B_{pr}^0 \\ &= f^{00} B_{00}^0 B_{00}^0 + f^{11} B_{10}^1 B_{11}^0 + f^{22} B_{20}^2 B_{22}^0 + f^{33} B_{30}^3 B_{33}^0 \\ &\quad + f^{44} B_{40}^4 B_{44}^0 + f^{55} B_{50}^5 B_{55}^0 + f^{66} B_{60}^6 B_{66}^0 + f^{77} B_{70}^7 B_{77}^0 \\ &= f^{00} - f^{11} - f^{22} - f^{33} - f^{44} - f^{55} - f^{66} - f^{77} \end{aligned}$$

$$(4.5.10) \begin{aligned} f_0^1 &= f^{kr} B_{k0}^p B_{pr}^1 \\ &= f^{01} B_{00}^0 B_{01}^1 + f^{10} B_{10}^1 B_{10}^1 + f^{23} B_{20}^2 B_{23}^1 + f^{32} B_{30}^3 B_{32}^1 \\ &\quad + f^{45} B_{40}^4 B_{45}^1 + f^{54} B_{50}^5 B_{54}^1 + f^{67} B_{60}^6 B_{67}^1 + f^{76} B_{70}^7 B_{76}^1 \\ &= f^{01} + f^{10} + f^{23} - f^{32} + f^{45} - f^{54} - f^{67} + f^{76} \end{aligned}$$

$$(4.5.11) \begin{aligned} f_0^2 &= f^{kr} B_{k0}^p B_{pr}^2 \\ &= f^{02} B_{00}^0 B_{02}^2 + f^{13} B_{10}^1 B_{13}^2 + f^{20} B_{20}^2 B_{20}^2 + f^{31} B_{30}^3 B_{31}^2 \\ &\quad + f^{46} B_{40}^4 B_{46}^2 + f^{57} B_{50}^5 B_{57}^2 + f^{64} B_{60}^6 B_{64}^2 + f^{75} B_{70}^7 B_{75}^2 \\ &= f^{02} - f^{13} + f^{20} + f^{31} + f^{46} + f^{57} - f^{64} - f^{75} \end{aligned}$$

$$f^3_0 = f^{kr}B^p_{k0}B^3_{pr}$$
$$= f^{03}B^0_{00}B^3_{03} + f^{12}B^1_{10}B^3_{12} + f^{21}B^2_{20}B^3_{21} + f^{30}B^3_{30}B^3_{30}$$
(4.5.12)
$$+ f^{47}B^4_{40}B^3_{47} + f^{56}B^5_{50}B^3_{56} + f^{65}B^6_{60}B^3_{65} + f^{74}B^7_{70}B^3_{74}$$
$$= f^{03} + f^{12} - f^{21} + f^{30} + f^{47} - f^{56} + f^{65} - f^{74}$$

$$f^4_0 = f^{kr}B^p_{k0}B^4_{pr}$$
$$= f^{04}B^0_{00}B^4_{04} + f^{15}B^1_{10}B^4_{15} + f^{26}B^2_{20}B^4_{26} + f^{37}B^3_{30}B^4_{37}$$
(4.5.13)
$$+ f^{40}B^4_{40}B^4_{40} + f^{51}B^5_{50}B^4_{51} + f^{62}B^6_{60}B^4_{62} + f^{73}B^7_{70}B^4_{73}$$
$$= f^{04} - f^{15} - f^{26} - f^{37} + f^{40} + f^{51} + f^{62} + f^{73}$$

$$f^5_0 = f^{kr}B^p_{k0}B^5_{pr}$$
$$= f^{05}B^0_{00}B^5_{05} + f^{14}B^1_{10}B^5_{14} + f^{27}B^2_{20}B^5_{27} + f^{36}B^3_{30}B^5_{36}$$
(4.5.14)
$$+ f^{41}B^4_{40}B^5_{41} + f^{50}B^5_{50}B^5_{50} + f^{63}B^6_{60}B^5_{63} + f^{72}B^7_{70}B^5_{72}$$
$$= f^{05} + f^{14} - f^{27} + f^{36} - f^{41} + f^{50} - f^{63} + f^{72}$$

$$f^6_0 = f^{kr}B^p_{k0}B^6_{pr}$$
$$= f^{06}B^0_{00}B^6_{06} + f^{17}B^1_{10}B^6_{17} + f^{24}B^2_{20}B^6_{24} + f^{35}B^3_{30}B^6_{35}$$
(4.5.15)
$$+ f^{42}B^4_{40}B^6_{42} + f^{53}B^5_{50}B^6_{53} + f^{60}B^6_{60}B^6_{60} + f^{71}B^7_{70}B^6_{71}$$
$$= f^{06} + f^{17} + f^{24} - f^{35} - f^{42} + f^{53} + f^{60} - f^{71}$$

$$f^7_0 = f^{kr}B^p_{k0}B^7_{pr}$$
$$= f^{07}B^0_{00}B^7_{07} + f^{16}B^1_{10}B^7_{16} + f^{25}B^2_{20}B^7_{25} + f^{34}B^3_{30}B^7_{34}$$
(4.5.16)
$$+ f^{43}B^4_{40}B^7_{43} + f^{52}B^5_{50}B^7_{52} + f^{61}B^6_{60}B^7_{61} + f^{70}B^7_{70}B^7_{70}$$
$$= f^{07} - f^{16} + f^{25} + f^{34} - f^{43} - f^{52} + f^{61} + f^{70}$$

$$f^0_1 = f^{kr}B^p_{k1}B^0_{pr}$$
$$= f^{01}B^1_{01}B^0_{11} + f^{10}B^0_{11}B^0_{00} + f^{23}B^3_{21}B^0_{33} + f^{32}B^2_{31}B^0_{22}$$
(4.5.17)
$$+ f^{45}B^5_{41}B^0_{55} + f^{54}B^4_{51}B^0_{44} + f^{67}B^7_{61}B^0_{77} + f^{76}B^6_{71}B^0_{66}$$
$$= -f^{01} - f^{10} + f^{23} - f^{32} + f^{45} - f^{54} - f^{67} + f^{76}$$

$$f^1_1 = f^{kr}B^p_{k1}B^1_{pr}$$
$$= f^{00}B^1_{01}B^1_{10} + f^{11}B^0_{11}B^1_{01} + f^{22}B^3_{21}B^1_{32} + f^{33}B^2_{31}B^1_{23}$$
(4.5.18)
$$+ f^{44}B^5_{41}B^1_{54} + f^{55}B^4_{51}B^1_{45} + f^{66}B^7_{61}B^1_{76} + f^{77}B^6_{71}B^1_{67}$$
$$= f^{00} - f^{11} + f^{22} + f^{33} + f^{44} + f^{55} + f^{66} + f^{77}$$

(4.5.19)
$$\begin{aligned}f^2_1 &= f^{kr}B^p_{k1}B^2_{pr}\\ &= f^{03}B^1_{01}B^2_{13} + f^{12}B^0_{11}B^2_{02} + f^{21}B^3_{21}B^2_{31} + f^{30}B^2_{31}B^2_{20}\\ &\quad + f^{47}B^5_{41}B^2_{57} + f^{56}B^4_{51}B^2_{46} + f^{65}B^7_{61}B^2_{75} + f^{74}B^6_{71}B^2_{64}\\ &= -f^{03} - f^{12} - f^{21} + f^{30} - f^{47} + f^{56} - f^{65} + f^{74}\end{aligned}$$

(4.5.20)
$$\begin{aligned}f^3_1 &= f^{kr}B^p_{k1}B^3_{pr}\\ &= f^{02}B^1_{01}B^3_{12} + f^{13}B^0_{11}B^3_{03} + f^{20}B^3_{21}B^3_{30} + f^{31}B^2_{31}B^3_{21}\\ &\quad + f^{46}B^5_{41}B^3_{56} + f^{57}B^4_{51}B^3_{47} + f^{64}B^7_{61}B^3_{74} + f^{75}B^6_{71}B^3_{65}\\ &= f^{02} - f^{13} - f^{20} - f^{31} + f^{46} + f^{57} - f^{64} - f^{75}\end{aligned}$$

(4.5.21)
$$\begin{aligned}f^4_1 &= f^{kr}B^p_{k1}B^4_{pr}\\ &= f^{05}B^1_{01}B^4_{15} + f^{14}B^0_{11}B^4_{04} + f^{27}B^3_{21}B^4_{37} + f^{36}B^2_{31}B^4_{26}\\ &\quad + f^{41}B^5_{41}B^4_{51} + f^{50}B^4_{51}B^4_{40} + f^{63}B^7_{61}B^4_{73} + f^{72}B^6_{71}B^4_{62}\\ &= -f^{05} - f^{14} + f^{27} - f^{36} - f^{41} + f^{50} + f^{63} - f^{72}\end{aligned}$$

(4.5.22)
$$\begin{aligned}f^5_1 &= f^{kr}B^p_{k1}B^5_{pr}\\ &= f^{04}B^1_{01}B^5_{14} + f^{15}B^0_{11}B^5_{05} + f^{26}B^3_{21}B^5_{36} + f^{37}B^2_{31}B^5_{27}\\ &\quad + f^{40}B^5_{41}B^5_{50} + f^{51}B^4_{51}B^5_{41} + f^{62}B^7_{61}B^5_{72} + f^{73}B^6_{71}B^5_{63}\\ &= f^{04} - f^{15} - f^{26} - f^{37} - f^{40} - f^{51} + f^{62} + f^{73}\end{aligned}$$

(4.5.23)
$$\begin{aligned}f^6_1 &= f^{kr}B^p_{k1}B^6_{pr}\\ &= f^{07}B^1_{01}B^6_{17} + f^{16}B^0_{11}B^6_{06} + f^{25}B^3_{21}B^6_{35} + f^{34}B^2_{31}B^6_{24}\\ &\quad + f^{43}B^5_{41}B^6_{53} + f^{52}B^4_{51}B^6_{42} + f^{61}B^7_{61}B^6_{71} + f^{70}B^6_{71}B^6_{60}\\ &= f^{07} - f^{16} + f^{25} + f^{34} - f^{43} - f^{52} - f^{61} - f^{70}\end{aligned}$$

(4.5.24)
$$\begin{aligned}f^7_1 &= f^{kr}B^p_{k1}B^7_{pr}\\ &= f^{06}B^1_{01}B^7_{16} + f^{17}B^0_{11}B^7_{07} + f^{24}B^3_{21}B^7_{34} + f^{35}B^2_{31}B^7_{25}\\ &\quad + f^{42}B^5_{41}B^7_{52} + f^{53}B^4_{51}B^7_{43} + f^{60}B^7_{61}B^7_{70} + f^{71}B^6_{71}B^7_{61}\\ &= -f^{06} - f^{17} - f^{24} + f^{35} + f^{42} - f^{53} + f^{60} - f^{71}\end{aligned}$$

(4.5.25)
$$\begin{aligned}f^0_2 &= f^{kr}B^p_{k2}B^0_{pr}\\ &= f^{02}B^2_{02}B^0_{22} + f^{13}B^3_{12}B^0_{33} + f^{20}B^0_{22}B^0_{00} + f^{31}B^1_{32}B^0_{11}\\ &\quad + f^{46}B^6_{42}B^0_{66} + f^{57}B^7_{52}B^0_{77} + f^{64}B^4_{62}B^0_{44} + f^{75}B^5_{72}B^0_{55}\\ &= -f^{02} - f^{13} - f^{20} + f^{31} + f^{46} + f^{57} - f^{64} - f^{75}\end{aligned}$$

$$f^1_2 = f^{kr}B^p_{k2}B^1_{pr}$$

(4.5.26)
$$= f^{03}B^2_{02}B^1_{23} + f^{12}B^3_{12}B^1_{32} + f^{21}B^0_{22}B^1_{01} + f^{30}B^1_{32}B^1_{10}$$
$$+ f^{47}B^6_{42}B^1_{67} + f^{56}B^7_{52}B^1_{76} + f^{65}B^4_{62}B^1_{45} + f^{74}B^5_{72}B^1_{54}$$
$$= f^{03} - f^{12} - f^{21} - f^{30} + f^{47} - f^{56} + f^{65} - f^{74}$$

$$f^2_2 = f^{kr}B^p_{k2}B^2_{pr}$$

(4.5.27)
$$= f^{00}B^2_{02}B^2_{20} + f^{11}B^3_{12}B^2_{31} + f^{22}B^0_{22}B^2_{02} + f^{33}B^1_{32}B^2_{13}$$
$$+ f^{44}B^6_{42}B^2_{64} + f^{55}B^7_{52}B^2_{75} + f^{66}B^4_{62}B^2_{46} + f^{77}B^5_{72}B^2_{57}$$
$$= f^{00} + f^{11} - f^{22} + f^{33} + f^{44} + f^{55} + f^{66} + f^{77}$$

$$f^3_2 = f^{kr}B^p_{k2}B^3_{pr}$$

(4.5.28)
$$= f^{01}B^2_{02}B^3_{21} + f^{10}B^3_{12}B^3_{30} + f^{23}B^0_{22}B^3_{03} + f^{32}B^1_{32}B^3_{12}$$
$$+ f^{45}B^6_{42}B^3_{65} + f^{54}B^7_{52}B^3_{74} + f^{67}B^4_{62}B^3_{47} + f^{76}B^5_{72}B^3_{56}$$
$$= -f^{01} + f^{10} - f^{23} - f^{32} - f^{45} + f^{54} + f^{67} - f^{76}$$

$$f^4_2 = f^{kr}B^p_{k2}B^4_{pr}$$

(4.5.29)
$$= f^{06}B^2_{02}B^4_{26} + f^{17}B^3_{12}B^4_{37} + f^{24}B^0_{22}B^4_{04} + f^{35}B^1_{32}B^4_{15}$$
$$+ f^{42}B^6_{42}B^4_{62} + f^{53}B^7_{52}B^4_{73} + f^{60}B^4_{62}B^4_{40} + f^{71}B^5_{72}B^4_{51}$$
$$= -f^{06} - f^{17} - f^{24} + f^{35} - f^{42} - f^{53} + f^{60} + f^{71}$$

$$f^5_2 = f^{kr}B^p_{k2}B^5_{pr}$$

(4.5.30)
$$= f^{07}B^2_{02}B^5_{27} + f^{16}B^3_{12}B^5_{36} + f^{25}B^0_{22}B^5_{05} + f^{34}B^1_{32}B^5_{14}$$
$$+ f^{43}B^6_{42}B^5_{63} + f^{52}B^7_{52}B^5_{72} + f^{61}B^4_{62}B^5_{41} + f^{70}B^5_{72}B^5_{50}$$
$$= -f^{07} + f^{16} - f^{25} - f^{34} + f^{43} - f^{52} - f^{61} + f^{70}$$

$$f^6_2 = f^{kr}B^p_{k2}B^6_{pr}$$

(4.5.31)
$$= f^{04}B^2_{02}B^6_{24} + f^{15}B^3_{12}B^6_{35} + f^{26}B^0_{22}B^6_{06} + f^{37}B^1_{32}B^6_{17}$$
$$+ f^{40}B^6_{42}B^6_{60} + f^{51}B^7_{52}B^6_{71} + f^{62}B^4_{62}B^6_{42} + f^{73}B^5_{72}B^6_{53}$$
$$= f^{04} - f^{15} - f^{26} - f^{37} - f^{40} + f^{51} - f^{62} + f^{73}$$

$$f^7_2 = f^{kr}B^p_{k2}B^7_{pr}$$

(4.5.32)
$$= f^{05}B^2_{02}B^7_{25} + f^{14}B^3_{12}B^7_{34} + f^{27}B^0_{22}B^7_{07} + f^{36}B^1_{32}B^7_{16}$$
$$+ f^{41}B^6_{42}B^7_{61} + f^{50}B^7_{52}B^7_{70} + f^{63}B^4_{62}B^7_{43} + f^{72}B^5_{72}B^7_{52}$$
$$= f^{05} + f^{14} - f^{27} + f^{36} - f^{41} - f^{50} - f^{63} - f^{72}$$

$$(4.5.33) \quad \begin{aligned} f_3^0 &= f^{kr} B_{k3}^p B_{pr}^0 \\ &= f^{03} B_{03}^3 B_{33}^0 + f^{12} B_{13}^2 B_{22}^0 + f^{21} B_{23}^1 B_{11}^0 + f^{30} B_{33}^0 B_{00}^0 \\ &\quad + f^{47} B_{43}^7 B_{77}^0 + f^{56} B_{53}^6 B_{66}^0 + f^{65} B_{63}^5 B_{55}^0 + f^{74} B_{73}^4 B_{44}^0 \\ &= -f^{03} + f^{12} - f^{21} - f^{30} + f^{47} - f^{56} + f^{65} - f^{74} \end{aligned}$$

$$(4.5.34) \quad \begin{aligned} f_3^1 &= f^{kr} B_{k3}^p B_{pr}^1 \\ &= f^{02} B_{03}^3 B_{32}^1 + f^{13} B_{13}^2 B_{23}^1 + f^{20} B_{23}^1 B_{10}^1 + f^{31} B_{33}^0 B_{01}^1 \\ &\quad + f^{46} B_{43}^7 B_{76}^1 + f^{57} B_{53}^6 B_{67}^1 + f^{64} B_{63}^5 B_{54}^1 + f^{75} B_{73}^4 B_{45}^1 \\ &= -f^{02} - f^{13} + f^{20} - f^{31} - f^{46} - f^{57} + f^{64} + f^{75} \end{aligned}$$

$$(4.5.35) \quad \begin{aligned} f_3^2 &= f^{kr} B_{k3}^p B_{pr}^2 \\ &= f^{01} B_{03}^3 B_{31}^2 + f^{10} B_{13}^2 B_{20}^2 + f^{23} B_{23}^1 B_{13}^2 + f^{32} B_{33}^0 B_{02}^2 \\ &\quad + f^{45} B_{43}^7 B_{75}^2 + f^{54} B_{53}^6 B_{64}^2 + f^{67} B_{63}^5 B_{57}^2 + f^{76} B_{73}^4 B_{46}^2 \\ &= f^{01} - f^{10} - f^{23} - f^{32} + f^{45} - f^{54} - f^{67} + f^{76} \end{aligned}$$

$$(4.5.36) \quad \begin{aligned} f_3^3 &= f^{kr} B_{k3}^p B_{pr}^3 \\ &= f^{00} B_{03}^3 B_{30}^3 + f^{11} B_{13}^2 B_{21}^3 + f^{22} B_{23}^1 B_{12}^3 + f^{33} B_{33}^0 B_{03}^3 \\ &\quad + f^{44} B_{43}^7 B_{74}^3 + f^{55} B_{53}^6 B_{65}^3 + f^{66} B_{63}^5 B_{56}^3 + f^{77} B_{73}^4 B_{47}^3 \\ &= f^{00} + f^{11} + f^{22} - f^{33} + f^{44} + f^{55} + f^{66} + f^{77} \end{aligned}$$

$$(4.5.37) \quad \begin{aligned} f_3^4 &= f^{kr} B_{k3}^p B_{pr}^4 \\ &= f^{07} B_{03}^3 B_{37}^4 + f^{16} B_{13}^2 B_{26}^4 + f^{25} B_{23}^1 B_{15}^4 + f^{34} B_{33}^0 B_{04}^4 \\ &\quad + f^{43} B_{43}^7 B_{73}^4 + f^{52} B_{53}^6 B_{62}^4 + f^{61} B_{63}^5 B_{51}^4 + f^{70} B_{73}^4 B_{40}^4 \\ &= -f^{07} + f^{16} - f^{25} - f^{34} - f^{43} + f^{52} - f^{61} + f^{70} \end{aligned}$$

$$(4.5.38) \quad \begin{aligned} f_3^5 &= f^{kr} B_{k3}^p B_{pr}^5 \\ &= f^{06} B_{03}^3 B_{36}^5 + f^{17} B_{13}^2 B_{27}^5 + f^{24} B_{23}^1 B_{14}^5 + f^{35} B_{33}^0 B_{05}^5 \\ &\quad + f^{42} B_{43}^7 B_{72}^5 + f^{53} B_{53}^6 B_{63}^5 + f^{60} B_{63}^5 B_{50}^5 + f^{71} B_{73}^4 B_{41}^5 \\ &= f^{06} + f^{17} + f^{24} - f^{35} - f^{42} - f^{53} - f^{60} - f^{71} \end{aligned}$$

$$(4.5.39) \quad \begin{aligned} f_3^6 &= f^{kr} B_{k3}^p B_{pr}^6 \\ &= f^{05} B_{03}^3 B_{35}^6 + f^{14} B_{13}^2 B_{24}^6 + f^{27} B_{23}^1 B_{17}^6 + f^{36} B_{33}^0 B_{06}^6 \\ &\quad + f^{41} B_{43}^7 B_{71}^6 + f^{50} B_{53}^6 B_{60}^6 + f^{63} B_{63}^5 B_{53}^6 + f^{72} B_{73}^4 B_{42}^6 \\ &= -f^{05} - f^{14} + f^{27} - f^{36} + f^{41} + f^{50} - f^{63} - f^{72} \end{aligned}$$

$$f_3^7 = f^{kr} B_{k3}^p B_{pr}^7$$

(4.5.40)
$$= f^{04} B_{03}^3 B_{34}^7 + f^{15} B_{13}^2 B_{25}^7 + f^{26} B_{23}^1 B_{16}^7 + f^{37} B_{33}^0 B_{07}^7$$
$$+ f^{40} B_{43}^7 B_{70}^7 + f^{51} B_{53}^6 B_{61}^7 + f^{62} B_{63}^5 B_{52}^7 + f^{73} B_{73}^4 B_{43}^7$$
$$= f^{04} - f^{15} - f^{26} - f^{37} - f^{40} + f^{51} + f^{62} - f^{73}$$

$$f_4^0 = f^{kr} B_{k4}^p B_{pr}^0$$

(4.5.41)
$$= f^{04} B_{04}^4 B_{44}^0 + f^{15} B_{14}^5 B_{55}^0 + f^{26} B_{24}^6 B_{66}^0 + f^{37} B_{34}^7 B_{77}^0$$
$$+ f^{40} B_{44}^0 B_{00}^0 + f^{51} B_{54}^1 B_{11}^0 + f^{62} B_{64}^2 B_{22}^0 + f^{73} B_{74}^3 B_{33}^0$$
$$= -f^{04} - f^{15} - f^{26} - f^{37} - f^{40} + f^{51} + f^{62} + f^{73}$$

$$f_4^1 = f^{kr} B_{k4}^p B_{pr}^1$$

(4.5.42)
$$= f^{05} B_{04}^4 B_{45}^1 + f^{14} B_{14}^5 B_{54}^1 + f^{27} B_{24}^6 B_{67}^1 + f^{36} B_{34}^7 B_{76}^1$$
$$+ f^{41} B_{44}^0 B_{01}^1 + f^{50} B_{54}^1 B_{10}^1 + f^{63} B_{64}^2 B_{23}^1 + f^{72} B_{74}^3 B_{32}^1$$
$$= f^{05} - f^{14} - f^{27} + f^{36} - f^{41} - f^{50} - f^{63} + f^{72}$$

$$f_4^2 = f^{kr} B_{k4}^p B_{pr}^2$$

(4.5.43)
$$= f^{06} B_{04}^4 B_{46}^2 + f^{17} B_{14}^5 B_{57}^2 + f^{24} B_{24}^6 B_{64}^2 + f^{35} B_{34}^7 B_{75}^2$$
$$+ f^{42} B_{44}^0 B_{02}^2 + f^{53} B_{54}^1 B_{13}^2 + f^{60} B_{64}^2 B_{20}^2 + f^{71} B_{74}^3 B_{31}^2$$
$$= f^{06} + f^{17} - f^{24} - f^{35} - f^{42} + f^{53} - f^{60} - f^{71}$$

$$f_4^3 = f^{kr} B_{k4}^p B_{pr}^3$$

(4.5.44)
$$= f^{07} B_{04}^4 B_{47}^3 + f^{16} B_{14}^5 B_{56}^3 + f^{25} B_{24}^6 B_{65}^3 + f^{34} B_{34}^7 B_{74}^3$$
$$+ f^{43} B_{44}^0 B_{03}^3 + f^{52} B_{54}^1 B_{12}^3 + f^{61} B_{64}^2 B_{21}^3 + f^{70} B_{74}^3 B_{30}^3$$
$$= f^{07} - f^{16} + f^{25} - f^{34} - f^{43} - f^{52} + f^{61} - f^{70}$$

$$f_4^4 = f^{kr} B_{k4}^p B_{pr}^4$$

(4.5.45)
$$= f^{00} B_{04}^4 B_{40}^4 + f^{11} B_{14}^5 B_{51}^4 + f^{22} B_{24}^6 B_{62}^4 + f^{33} B_{34}^7 B_{73}^4$$
$$+ f^{44} B_{44}^0 B_{04}^4 + f^{55} B_{54}^1 B_{15}^4 + f^{66} B_{64}^2 B_{26}^4 + f^{77} B_{74}^3 B_{37}^4$$
$$= f^{00} + f^{11} + f^{22} + f^{33} - f^{44} + f^{55} + f^{66} + f^{77}$$

$$f_4^5 = f^{kr} B_{k4}^p B_{pr}^5$$

(4.5.46)
$$= f^{01} B_{04}^4 B_{41}^5 + f^{10} B_{14}^5 B_{50}^5 + f^{23} B_{24}^6 B_{63}^5 + f^{32} B_{34}^7 B_{72}^5$$
$$+ f^{45} B_{44}^0 B_{05}^5 + f^{54} B_{54}^1 B_{14}^5 + f^{67} B_{64}^2 B_{27}^5 + f^{76} B_{74}^3 B_{36}^5$$
$$= -f^{01} + f^{10} - f^{23} + f^{32} - f^{45} - f^{54} + f^{67} - f^{76}$$

(4.5.47)
$$\begin{aligned}f_4^6 &= f^{kr}B_{k4}^p B_{pr}^6 \\ &= f^{02}B_{04}^4 B_{42}^6 + f^{13}B_{14}^5 B_{53}^6 + f^{20}B_{24}^6 B_{60}^6 + f^{31}B_{34}^7 B_{71}^6 \\ &\quad + f^{46}B_{44}^0 B_{06}^6 + f^{57}B_{54}^1 B_{17}^6 + f^{64}B_{64}^2 B_{24}^6 + f^{75}B_{74}^3 B_{35}^6 \\ &= -f^{02} + f^{13} + f^{20} - f^{31} - f^{46} - f^{57} - f^{64} + f^{75}\end{aligned}$$

(4.5.48)
$$\begin{aligned}f_4^7 &= f^{kr}B_{k4}^p B_{pr}^7 \\ &= f^{03}B_{04}^4 B_{43}^7 + f^{12}B_{14}^5 B_{52}^7 + f^{21}B_{24}^6 B_{61}^7 + f^{30}B_{34}^7 B_{70}^7 \\ &\quad + f^{47}B_{44}^0 B_{07}^7 + f^{56}B_{54}^1 B_{16}^7 + f^{65}B_{64}^2 B_{25}^7 + f^{74}B_{74}^3 B_{34}^7 \\ &= -f^{03} - f^{12} + f^{21} + f^{30} - f^{47} + f^{56} - f^{65} - f^{74}\end{aligned}$$

(4.5.49)
$$\begin{aligned}f_5^0 &= f^{kr}B_{k5}^p B_{pr}^0 \\ &= f^{05}B_{05}^5 B_{55}^0 + f^{14}B_{15}^4 B_{44}^0 + f^{27}B_{25}^7 B_{77}^0 + f^{36}B_{35}^6 B_{66}^0 \\ &\quad + f^{41}B_{45}^1 B_{11}^0 + f^{50}B_{55}^0 B_{00}^0 + f^{63}B_{65}^3 B_{33}^0 + f^{72}B_{75}^2 B_{22}^0 \\ &= -f^{05} + f^{14} - f^{27} + f^{36} - f^{41} - f^{50} - f^{63} + f^{72}\end{aligned}$$

(4.5.50)
$$\begin{aligned}f_5^1 &= f^{kr}B_{k5}^p B_{pr}^1 \\ &= f^{04}B_{05}^5 B_{54}^1 + f^{15}B_{15}^4 B_{45}^1 + f^{26}B_{25}^7 B_{76}^1 + f^{37}B_{35}^6 B_{67}^1 \\ &\quad + f^{40}B_{45}^1 B_{10}^1 + f^{51}B_{55}^0 B_{01}^1 + f^{62}B_{65}^3 B_{32}^1 + f^{73}B_{75}^2 B_{23}^1 \\ &= -f^{04} - f^{15} + f^{26} + f^{37} + f^{40} - f^{51} - f^{62} - f^{73}\end{aligned}$$

(4.5.51)
$$\begin{aligned}f_5^2 &= f^{kr}B_{k5}^p B_{pr}^2 \\ &= f^{07}B_{05}^5 B_{57}^2 + f^{16}B_{15}^4 B_{46}^2 + f^{25}B_{25}^7 B_{75}^2 + f^{34}B_{35}^6 B_{64}^2 \\ &\quad + f^{43}B_{45}^1 B_{13}^2 + f^{52}B_{55}^0 B_{02}^2 + f^{61}B_{65}^3 B_{31}^2 + f^{70}B_{75}^2 B_{20}^2 \\ &= f^{07} - f^{16} - f^{25} + f^{34} - f^{43} - f^{52} + f^{61} - f^{70}\end{aligned}$$

(4.5.52)
$$\begin{aligned}f_5^3 &= f^{kr}B_{k5}^p B_{pr}^3 \\ &= f^{06}B_{05}^5 B_{56}^3 + f^{17}B_{15}^4 B_{47}^3 + f^{24}B_{25}^7 B_{74}^3 + f^{35}B_{35}^6 B_{65}^3 \\ &\quad + f^{42}B_{45}^1 B_{12}^3 + f^{53}B_{55}^0 B_{03}^3 + f^{60}B_{65}^3 B_{30}^3 + f^{71}B_{75}^2 B_{21}^3 \\ &= -f^{06} - f^{17} - f^{24} - f^{35} + f^{42} - f^{53} + f^{60} + f^{71}\end{aligned}$$

(4.5.53)
$$\begin{aligned}f_5^4 &= f^{kr}B_{k5}^p B_{pr}^4 \\ &= f^{01}B_{05}^5 B_{51}^4 + f^{10}B_{15}^4 B_{40}^4 + f^{23}B_{25}^7 B_{73}^4 + f^{32}B_{35}^6 B_{62}^4 \\ &\quad + f^{45}B_{45}^1 B_{15}^4 + f^{54}B_{55}^0 B_{04}^4 + f^{67}B_{65}^3 B_{37}^4 + f^{76}B_{75}^2 B_{26}^4 \\ &= f^{01} - f^{10} + f^{23} - f^{32} - f^{45} - f^{54} - f^{67} + f^{76}\end{aligned}$$

$$f_5^5 = f^{kr} B_{k5}^p B_{pr}^5$$
$$= f^{00} B_{05}^5 B_{50}^5 + f^{11} B_{15}^4 B_{41}^5 + f^{22} B_{25}^7 B_{72}^5 + f^{33} B_{35}^6 B_{63}^5$$
(4.5.54)
$$+ f^{44} B_{45}^1 B_{14}^5 + f^{55} B_{55}^0 B_{05}^5 + f^{66} B_{65}^3 B_{36}^5 + f^{77} B_{75}^2 B_{27}^5$$
$$= f^{00} + f^{11} + f^{22} + f^{33} + f^{44} - f^{55} + f^{66} + f^{77}$$

$$f_5^6 = f^{kr} B_{k5}^p B_{pr}^6$$
$$= f^{03} B_{05}^5 B_{53}^6 + f^{12} B_{15}^4 B_{42}^6 + f^{21} B_{25}^7 B_{71}^6 + f^{30} B_{35}^6 B_{60}^6$$
(4.5.55)
$$+ f^{47} B_{45}^1 B_{17}^6 + f^{56} B_{55}^0 B_{06}^6 + f^{65} B_{65}^3 B_{35}^6 + f^{74} B_{75}^2 B_{24}^6$$
$$= f^{03} + f^{12} - f^{21} - f^{30} + f^{47} - f^{56} - f^{65} - f^{74}$$

$$f_5^7 = f^{kr} B_{k5}^p B_{pr}^7$$
$$= f^{02} B_{05}^5 B_{52}^7 + f^{13} B_{15}^4 B_{43}^7 + f^{20} B_{25}^7 B_{70}^7 + f^{31} B_{35}^6 B_{61}^7$$
(4.5.56)
$$+ f^{46} B_{45}^1 B_{16}^7 + f^{57} B_{55}^0 B_{07}^7 + f^{64} B_{65}^3 B_{34}^7 + f^{75} B_{75}^2 B_{25}^7$$
$$= -f^{02} + f^{13} + f^{20} - f^{31} - f^{46} - f^{57} + f^{64} - f^{75}$$

$$f_6^0 = f^{kr} B_{k6}^p B_{pr}^0$$
$$= f^{06} B_{06}^6 B_{66}^0 + f^{17} B_{16}^7 B_{77}^0 + f^{24} B_{26}^4 B_{44}^0 + f^{35} B_{36}^5 B_{55}^0$$
(4.5.57)
$$+ f^{42} B_{46}^2 B_{22}^0 + f^{53} B_{56}^3 B_{33}^0 + f^{60} B_{66}^0 B_{00}^0 + f^{71} B_{76}^1 B_{11}^0$$
$$= -f^{06} + f^{17} + f^{24} - f^{35} - f^{42} + f^{53} - f^{60} - f^{71}$$

$$f_6^1 = f^{kr} B_{k6}^p B_{pr}^1$$
$$= f^{07} B_{06}^6 B_{67}^1 + f^{16} B_{16}^7 B_{76}^1 + f^{25} B_{26}^4 B_{45}^1 + f^{34} B_{36}^5 B_{54}^1$$
(4.5.58)
$$+ f^{43} B_{46}^2 B_{23}^1 + f^{52} B_{56}^3 B_{32}^1 + f^{61} B_{66}^0 B_{01}^1 + f^{70} B_{76}^1 B_{10}^1$$
$$= -f^{07} - f^{16} - f^{25} - f^{34} + f^{43} + f^{52} - f^{61} + f^{70}$$

$$f_6^2 = f^{kr} B_{k6}^p B_{pr}^2$$
$$= f^{04} B_{06}^6 B_{64}^2 + f^{15} B_{16}^7 B_{75}^2 + f^{26} B_{26}^4 B_{46}^2 + f^{37} B_{36}^5 B_{57}^2$$
(4.5.59)
$$+ f^{40} B_{46}^2 B_{20}^2 + f^{51} B_{56}^3 B_{31}^2 + f^{62} B_{66}^0 B_{02}^2 + f^{73} B_{76}^1 B_{13}^2$$
$$= -f^{04} + f^{15} - f^{26} + f^{37} + f^{40} - f^{51} - f^{62} - f^{73}$$

$$f_6^3 = f^{kr} B_{k6}^p B_{pr}^3$$
$$= f^{05} B_{06}^6 B_{65}^3 + f^{14} B_{16}^7 B_{74}^3 + f^{27} B_{26}^4 B_{47}^3 + f^{36} B_{36}^5 B_{56}^3$$
(4.5.60)
$$+ f^{41} B_{46}^2 B_{21}^3 + f^{50} B_{56}^3 B_{30}^3 + f^{63} B_{66}^0 B_{03}^3 + f^{72} B_{76}^1 B_{12}^3$$
$$= f^{05} + f^{14} - f^{27} - f^{36} - f^{41} - f^{50} - f^{63} + f^{72}$$

(4.5.61)
$$\begin{aligned}f^4_6 &= f^{kr}B^p_{k6}B^4_{pr}\\ &= f^{02}B^6_{06}B^4_{62} + f^{13}B^7_{16}B^4_{73} + f^{20}B^4_{26}B^4_{40} + f^{31}B^5_{36}B^4_{51}\\ &\quad + f^{46}B^2_{46}B^4_{26} + f^{57}B^3_{56}B^4_{37} + f^{64}B^0_{66}B^4_{04} + f^{75}B^1_{76}B^4_{15}\\ &= f^{02} - f^{13} - f^{20} + f^{31} - f^{46} + f^{57} - f^{64} - f^{75}\end{aligned}$$

(4.5.62)
$$\begin{aligned}f^5_6 &= f^{kr}B^p_{k6}B^5_{pr}\\ &= f^{03}B^6_{06}B^5_{63} + f^{12}B^7_{16}B^5_{72} + f^{21}B^4_{26}B^5_{41} + f^{30}B^5_{36}B^5_{50}\\ &\quad + f^{47}B^2_{46}B^5_{27} + f^{56}B^3_{56}B^5_{36} + f^{65}B^0_{66}B^5_{05} + f^{74}B^1_{76}B^5_{14}\\ &= -f^{03} - f^{12} + f^{21} + f^{30} - f^{47} - f^{56} - f^{65} + f^{74}\end{aligned}$$

(4.5.63)
$$\begin{aligned}f^6_6 &= f^{kr}B^p_{k6}B^6_{pr}\\ &= f^{00}B^6_{06}B^6_{60} + f^{11}B^7_{16}B^6_{71} + f^{22}B^4_{26}B^6_{42} + f^{33}B^5_{36}B^6_{53}\\ &\quad + f^{44}B^2_{46}B^6_{24} + f^{55}B^3_{56}B^6_{35} + f^{66}B^0_{66}B^6_{06} + f^{77}B^1_{76}B^6_{17}\\ &= f^{00} + f^{11} + f^{22} + f^{33} + f^{44} + f^{55} - f^{66} + f^{77}\end{aligned}$$

(4.5.64)
$$\begin{aligned}f^7_6 &= f^{kr}B^p_{k6}B^7_{pr}\\ &= f^{01}B^6_{06}B^7_{61} + f^{10}B^7_{16}B^7_{70} + f^{23}B^4_{26}B^7_{43} + f^{32}B^5_{36}B^7_{52}\\ &\quad + f^{45}B^2_{46}B^7_{25} + f^{54}B^3_{56}B^7_{34} + f^{67}B^0_{66}B^7_{07} + f^{76}B^1_{76}B^7_{16}\\ &= f^{01} - f^{10} + f^{23} - f^{32} + f^{45} - f^{54} - f^{67} - f^{76}\end{aligned}$$

(4.5.65)
$$\begin{aligned}f^0_7 &= f^{kr}B^p_{k7}B^0_{pr}\\ &= f^{07}B^7_{07}B^0_{77} + f^{16}B^6_{17}B^0_{66} + f^{25}B^5_{27}B^0_{55} + f^{34}B^4_{37}B^0_{44}\\ &\quad + f^{43}B^3_{47}B^0_{33} + f^{52}B^2_{57}B^0_{22} + f^{61}B^1_{67}B^0_{11} + f^{70}B^0_{77}B^0_{00}\\ &= -f^{07} - f^{16} + f^{25} + f^{34} - f^{43} - f^{52} + f^{61} - f^{70}\end{aligned}$$

(4.5.66)
$$\begin{aligned}f^1_7 &= f^{kr}B^p_{k7}B^1_{pr}\\ &= f^{06}B^7_{07}B^1_{76} + f^{17}B^6_{17}B^1_{67} + f^{24}B^5_{27}B^1_{54} + f^{35}B^4_{37}B^1_{45}\\ &\quad + f^{42}B^3_{47}B^1_{32} + f^{53}B^2_{57}B^1_{23} + f^{60}B^1_{67}B^1_{10} + f^{71}B^0_{77}B^1_{01}\\ &= f^{06} - f^{17} + f^{24} - f^{35} - f^{42} + f^{53} - f^{60} - f^{71}\end{aligned}$$

(4.5.67)
$$\begin{aligned}f^2_7 &= f^{kr}B^p_{k7}B^2_{pr}\\ &= f^{05}B^7_{07}B^2_{75} + f^{14}B^6_{17}B^2_{64} + f^{27}B^5_{27}B^2_{57} + f^{36}B^4_{37}B^2_{46}\\ &\quad + f^{41}B^3_{47}B^2_{31} + f^{50}B^2_{57}B^2_{20} + f^{63}B^1_{67}B^2_{13} + f^{72}B^0_{77}B^2_{02}\\ &= -f^{05} - f^{14} - f^{27} - f^{36} + f^{41} + f^{50} + f^{63} - f^{72}\end{aligned}$$

$$\begin{aligned}
f_7^3 &= f^{kr} B_{k7}^p B_{pr}^3 \\
&= f^{04} B_{07}^7 B_{74}^3 + f^{15} B_{17}^6 B_{65}^3 + f^{26} B_{27}^5 B_{56}^3 + f^{37} B_{37}^4 B_{47}^3 \\
&\quad + f^{40} B_{47}^3 B_{30}^3 + f^{51} B_{57}^2 B_{21}^3 + f^{62} B_{67}^1 B_{12}^3 + f^{73} B_{77}^0 B_{03}^3 \\
&= -f^{04} + f^{15} + f^{26} - f^{37} + f^{40} - f^{51} - f^{62} - f^{73}
\end{aligned} \tag{4.5.68}$$

$$\begin{aligned}
f_7^4 &= f^{kr} B_{k7}^p B_{pr}^4 \\
&= f^{03} B_{07}^7 B_{73}^4 + f^{12} B_{17}^6 B_{62}^4 + f^{21} B_{27}^5 B_{51}^4 + f^{30} B_{37}^4 B_{40}^4 \\
&\quad + f^{47} B_{47}^3 B_{37}^4 + f^{56} B_{57}^2 B_{26}^4 + f^{65} B_{67}^1 B_{15}^4 + f^{74} B_{77}^0 B_{04}^4 \\
&= f^{03} + f^{12} - f^{21} - f^{30} - f^{47} - f^{56} + f^{65} - f^{74}
\end{aligned} \tag{4.5.69}$$

$$\begin{aligned}
f_7^5 &= f^{kr} B_{k7}^p B_{pr}^5 \\
&= f^{02} B_{07}^7 B_{72}^5 + f^{13} B_{17}^6 B_{63}^5 + f^{20} B_{27}^5 B_{50}^5 + f^{31} B_{37}^4 B_{41}^5 \\
&\quad + f^{46} B_{47}^3 B_{36}^5 + f^{57} B_{57}^2 B_{27}^5 + f^{64} B_{67}^1 B_{14}^5 + f^{75} B_{77}^0 B_{05}^5 \\
&= f^{02} - f^{13} - f^{20} + f^{31} + f^{46} - f^{57} - f^{64} - f^{75}
\end{aligned} \tag{4.5.70}$$

$$\begin{aligned}
f_7^6 &= f^{kr} B_{k7}^p B_{pr}^6 \\
&= f^{01} B_{07}^7 B_{71}^6 + f^{10} B_{17}^6 B_{60}^6 + f^{23} B_{27}^5 B_{53}^6 + f^{32} B_{37}^4 B_{42}^6 \\
&\quad + f^{45} B_{47}^3 B_{35}^6 + f^{54} B_{57}^2 B_{24}^6 + f^{67} B_{67}^1 B_{17}^6 + f^{76} B_{77}^0 B_{06}^6 \\
&= -f^{01} + f^{10} - f^{23} + f^{32} - f^{45} + f^{54} - f^{67} - f^{76}
\end{aligned} \tag{4.5.71}$$

$$\begin{aligned}
f_7^7 &= f^{kr} B_{k7}^p B_{pr}^7 \\
&= f^{00} B_{07}^7 B_{70}^7 + f^{11} B_{17}^6 B_{61}^7 + f^{22} B_{27}^5 B_{52}^7 + f^{33} B_{37}^4 B_{43}^7 \\
&\quad + f^{44} B_{47}^3 B_{34}^7 + f^{55} B_{57}^2 B_{25}^7 + f^{66} B_{67}^1 B_{16}^7 + f^{77} B_{77}^0 B_{07}^7 \\
&= f^{00} + f^{11} + f^{22} + f^{33} + f^{44} + f^{55} + f^{66} - f^{77}
\end{aligned} \tag{4.5.72}$$

Equations (4.5.9), (4.5.18), (4.5.27), (4.5.36), (4.5.45), (4.5.54), (4.5.63), (4.5.72) form the system of linear equations (4.5.1).

Equations (4.5.10), (4.5.17), (4.5.28), (4.5.35), (4.5.46), (4.5.53), (4.5.64), (4.5.71) form the system of linear equations (4.5.2).

Equations (4.5.11), (4.5.20), (4.5.25), (4.5.34), (4.5.47), (4.5.56), (4.5.61), (4.5.70) form the system of linear equations (4.5.3).

Equations (4.5.12), (4.5.19), (4.5.26), (4.5.33), (4.5.48), (4.5.55), (4.5.62), (4.5.69) form the system of linear equations (4.5.4).

Equations (4.5.13), (4.5.22), (4.5.31), (4.5.40), (4.5.41), (4.5.50), (4.5.59), (4.5.68) form the system of linear equations (4.5.5).

Equations (4.5.14), (4.5.21), (4.5.32), (4.5.39), (4.5.42), (4.5.49), (4.5.60), (4.5.67) form the system of linear equations (4.5.6).

Equations (4.5.15), (4.5.24), (4.5.29), (4.5.38), (4.5.43), (4.5.52), (4.5.57), (4.5.66) form the system of linear equations (4.5.7).

Equations (4.5.16), (4.5.23), (4.5.30), (4.5.37), (4.5.44), (4.5.51), (4.5.58), (4.5.65) form the system of linear equations (4.5.8). □

**Theorem 4.5.2.** *Consider octonion algebra $O$ with basis (4.4.1). Standard components of linear function and coordinates of this function satisfy relationship*

(4.5.73) $$A = FB$$

(4.5.74) $$B = F^{-1}A$$

*where*

$$A = \begin{pmatrix} f_0^0 & f_1^0 & f_2^0 & f_3^0 & f_4^0 & f_5^0 & f_6^0 & f_7^0 \\ f_1^1 & -f_0^1 & f_3^1 & -f_2^1 & f_5^1 & -f_4^1 & -f_7^1 & f_6^1 \\ f_2^2 & -f_3^2 & -f_0^2 & f_1^2 & f_6^2 & f_7^2 & -f_4^2 & -f_5^2 \\ f_3^3 & f_2^3 & -f_1^3 & -f_0^3 & f_7^3 & -f_6^3 & f_5^3 & -f_4^3 \\ f_4^4 & -f_5^4 & -f_6^4 & -f_7^4 & -f_0^4 & f_1^4 & f_2^4 & f_3^4 \\ f_5^5 & f_4^5 & -f_7^5 & f_6^5 & -f_1^5 & -f_0^5 & -f_3^5 & f_2^5 \\ f_6^6 & f_7^6 & f_4^6 & -f_5^6 & -f_2^6 & f_3^6 & -f_0^6 & -f_1^6 \\ f_7^7 & -f_6^7 & f_5^7 & f_4^7 & -f_3^7 & -f_2^7 & f_1^7 & -f_0^7 \end{pmatrix}$$

$$B = \begin{pmatrix} f^{00} & -f^{01} & -f^{02} & -f^{03} & -f^{04} & -f^{05} & -f^{06} & -f^{07} \\ f^{11} & f^{10} & f^{13} & -f^{12} & f^{15} & -f^{14} & -f^{17} & f^{16} \\ f^{22} & -f^{23} & f^{20} & f^{21} & f^{26} & f^{27} & -f^{24} & -f^{25} \\ f^{33} & f^{32} & -f^{31} & f^{30} & f^{37} & -f^{36} & f^{35} & -f^{34} \\ f^{44} & -f^{45} & -f^{46} & -f^{47} & f^{40} & f^{41} & f^{42} & f^{43} \\ f^{55} & f^{54} & -f^{57} & f^{56} & -f^{51} & f^{50} & -f^{53} & f^{52} \\ f^{66} & f^{67} & f^{64} & -f^{65} & -f^{62} & f^{63} & f^{60} & -f^{61} \\ f^{77} & -f^{76} & f^{75} & f^{74} & -f^{73} & -f^{72} & f^{71} & f^{70} \end{pmatrix}$$

$$F = \begin{pmatrix} 1 & -1 & -1 & -1 & -1 & -1 & -1 & -1 \\ 1 & -1 & 1 & 1 & 1 & 1 & 1 & 1 \\ 1 & 1 & -1 & 1 & 1 & 1 & 1 & 1 \\ 1 & 1 & 1 & -1 & 1 & 1 & 1 & 1 \\ 1 & 1 & 1 & 1 & -1 & 1 & 1 & 1 \\ 1 & 1 & 1 & 1 & 1 & -1 & 1 & 1 \\ 1 & 1 & 1 & 1 & 1 & 1 & -1 & 1 \\ 1 & 1 & 1 & 1 & 1 & 1 & 1 & -1 \end{pmatrix}$$

$$F^{-1} = \frac{1}{12}\begin{pmatrix} 5 & 1 & 1 & 1 & 1 & 1 & 1 & 1 \\ -1 & -5 & 1 & 1 & 1 & 1 & 1 & 1 \\ -1 & 1 & -5 & 1 & 1 & 1 & 1 & 1 \\ -1 & 1 & 1 & -5 & 1 & 1 & 1 & 1 \\ -1 & 1 & 1 & 1 & -5 & 1 & 1 & 1 \\ -1 & 1 & 1 & 1 & 1 & -5 & 1 & 1 \\ -1 & 1 & 1 & 1 & 1 & 1 & -5 & 1 \\ -1 & 1 & 1 & 1 & 1 & 1 & 1 & -5 \end{pmatrix}$$

PROOF. Let us write the system of linear equations (4.5.1) as product of matrices

$$(4.5.75) \quad \begin{pmatrix} f_0^0 \\ f_1^1 \\ f_2^2 \\ f_3^3 \\ f_4^4 \\ f_5^5 \\ f_6^6 \\ f_7^7 \end{pmatrix} = \begin{pmatrix} 1 & -1 & -1 & -1 & -1 & -1 & -1 & -1 \\ 1 & -1 & 1 & 1 & 1 & 1 & 1 & 1 \\ 1 & 1 & -1 & 1 & 1 & 1 & 1 & 1 \\ 1 & 1 & 1 & -1 & 1 & 1 & 1 & 1 \\ 1 & 1 & 1 & 1 & -1 & 1 & 1 & 1 \\ 1 & 1 & 1 & 1 & 1 & -1 & 1 & 1 \\ 1 & 1 & 1 & 1 & 1 & 1 & -1 & 1 \\ 1 & 1 & 1 & 1 & 1 & 1 & 1 & -1 \end{pmatrix} \begin{pmatrix} f^{00} \\ f^{11} \\ f^{22} \\ f^{33} \\ f^{44} \\ f^{55} \\ f^{66} \\ f^{77} \end{pmatrix}$$

Let us write the system of linear equations (4.5.2) as product of matrices

$$(4.5.76) \quad \begin{pmatrix} f_0^1 \\ f_1^0 \\ f_2^3 \\ f_3^2 \\ f_4^5 \\ f_5^4 \\ f_6^7 \\ f_7^6 \end{pmatrix} = \begin{pmatrix} 1 & 1 & 1 & -1 & 1 & -1 & -1 & 1 \\ -1 & -1 & 1 & -1 & 1 & -1 & -1 & 1 \\ -1 & 1 & -1 & -1 & -1 & 1 & 1 & -1 \\ 1 & -1 & -1 & -1 & 1 & -1 & -1 & 1 \\ -1 & 1 & -1 & 1 & -1 & -1 & 1 & -1 \\ 1 & -1 & 1 & -1 & -1 & -1 & -1 & 1 \\ 1 & -1 & 1 & -1 & 1 & -1 & -1 & -1 \\ -1 & 1 & -1 & 1 & -1 & 1 & -1 & -1 \end{pmatrix} \begin{pmatrix} f^{01} \\ f^{10} \\ f^{23} \\ f^{32} \\ f^{45} \\ f^{54} \\ f^{67} \\ f^{76} \end{pmatrix}$$

From the equation (4.5.76), it follows that

$$\begin{pmatrix} -f_0^1 \\ f_1^0 \\ f_2^3 \\ -f_3^2 \\ f_4^5 \\ -f_5^4 \\ -f_6^7 \\ f_7^6 \end{pmatrix} = \begin{pmatrix} -1 & -1 & -1 & 1 & -1 & 1 & 1 & -1 \\ -1 & -1 & 1 & -1 & 1 & -1 & -1 & 1 \\ -1 & 1 & -1 & -1 & -1 & 1 & 1 & -1 \\ -1 & 1 & 1 & 1 & -1 & 1 & 1 & -1 \\ -1 & 1 & -1 & 1 & -1 & -1 & 1 & -1 \\ -1 & 1 & -1 & 1 & 1 & 1 & 1 & -1 \\ -1 & 1 & -1 & 1 & -1 & 1 & 1 & 1 \\ -1 & 1 & -1 & 1 & -1 & 1 & -1 & -1 \end{pmatrix} \begin{pmatrix} f^{01} \\ f^{10} \\ f^{23} \\ f^{32} \\ f^{45} \\ f^{54} \\ f^{67} \\ f^{76} \end{pmatrix}$$

$$= \begin{pmatrix} 1 & -1 & 1 & 1 & 1 & 1 & 1 & 1 \\ 1 & -1 & -1 & -1 & -1 & -1 & -1 & -1 \\ 1 & 1 & 1 & -1 & 1 & 1 & 1 & 1 \\ 1 & 1 & -1 & 1 & 1 & 1 & 1 & 1 \\ 1 & 1 & 1 & 1 & 1 & -1 & 1 & 1 \\ 1 & 1 & 1 & 1 & -1 & 1 & 1 & 1 \\ 1 & 1 & 1 & 1 & 1 & 1 & 1 & -1 \\ 1 & 1 & 1 & 1 & 1 & 1 & -1 & 1 \end{pmatrix} \begin{pmatrix} -f^{01} \\ f^{10} \\ -f^{23} \\ f^{32} \\ -f^{45} \\ f^{54} \\ f^{67} \\ -f^{76} \end{pmatrix}$$

$$(4.5.77) \begin{pmatrix} f_1^0 \\ -f_0^1 \\ -f_3^2 \\ f_2^3 \\ -f_5^4 \\ f_4^5 \\ f_7^6 \\ -f_6^7 \end{pmatrix} = \begin{pmatrix} 1 & -1 & -1 & -1 & -1 & -1 & -1 & -1 \\ 1 & -1 & 1 & 1 & 1 & 1 & 1 & 1 \\ 1 & 1 & -1 & 1 & 1 & 1 & 1 & 1 \\ 1 & 1 & 1 & -1 & 1 & 1 & 1 & 1 \\ 1 & 1 & 1 & 1 & -1 & 1 & 1 & 1 \\ 1 & 1 & 1 & 1 & 1 & -1 & 1 & 1 \\ 1 & 1 & 1 & 1 & 1 & 1 & -1 & 1 \\ 1 & 1 & 1 & 1 & 1 & 1 & 1 & -1 \end{pmatrix} \begin{pmatrix} -f^{01} \\ f^{10} \\ -f^{23} \\ f^{32} \\ -f^{45} \\ f^{54} \\ f^{67} \\ -f^{76} \end{pmatrix}$$

Let us write the system of linear equations (4.5.3) as product of matrices

$$(4.5.78) \begin{pmatrix} f_0^2 \\ f_1^3 \\ f_2^0 \\ f_3^1 \\ f_4^6 \\ f_5^7 \\ f_6^4 \\ f_7^5 \end{pmatrix} = \begin{pmatrix} 1 & -1 & 1 & 1 & 1 & 1 & -1 & -1 \\ 1 & -1 & -1 & -1 & 1 & 1 & -1 & -1 \\ -1 & -1 & -1 & 1 & 1 & 1 & -1 & -1 \\ -1 & -1 & 1 & -1 & -1 & -1 & 1 & 1 \\ -1 & 1 & 1 & -1 & -1 & -1 & -1 & 1 \\ -1 & 1 & 1 & -1 & -1 & -1 & 1 & -1 \\ 1 & -1 & -1 & 1 & -1 & 1 & -1 & -1 \\ 1 & -1 & -1 & 1 & 1 & -1 & -1 & -1 \end{pmatrix} \begin{pmatrix} f^{02} \\ f^{13} \\ f^{20} \\ f^{31} \\ f^{46} \\ f^{57} \\ f^{64} \\ f^{75} \end{pmatrix}$$

From the equation (4.5.78), it follows that

$$\begin{pmatrix} -f_0^2 \\ -f_1^3 \\ f_2^0 \\ f_3^1 \\ f_4^6 \\ f_5^7 \\ -f_6^4 \\ -f_7^5 \end{pmatrix} = \begin{pmatrix} -1 & 1 & -1 & -1 & -1 & -1 & 1 & 1 \\ -1 & 1 & 1 & 1 & -1 & -1 & 1 & 1 \\ -1 & -1 & -1 & 1 & 1 & 1 & -1 & -1 \\ -1 & -1 & 1 & -1 & -1 & -1 & 1 & 1 \\ -1 & 1 & 1 & -1 & -1 & -1 & -1 & 1 \\ -1 & 1 & 1 & -1 & -1 & -1 & 1 & -1 \\ -1 & 1 & 1 & -1 & 1 & -1 & 1 & 1 \\ -1 & 1 & 1 & -1 & -1 & 1 & 1 & 1 \end{pmatrix} \begin{pmatrix} f^{02} \\ f^{13} \\ f^{20} \\ f^{31} \\ f^{46} \\ f^{57} \\ f^{64} \\ f^{75} \end{pmatrix}$$

$$= \begin{pmatrix} 1 & 1 & -1 & 1 & 1 & 1 & 1 & 1 \\ 1 & 1 & 1 & -1 & 1 & 1 & 1 & 1 \\ 1 & -1 & -1 & -1 & -1 & -1 & -1 & -1 \\ 1 & -1 & 1 & 1 & 1 & 1 & 1 & 1 \\ 1 & 1 & 1 & 1 & 1 & 1 & -1 & 1 \\ 1 & 1 & 1 & 1 & 1 & 1 & 1 & -1 \\ 1 & 1 & 1 & 1 & -1 & 1 & 1 & 1 \\ 1 & 1 & 1 & 1 & 1 & -1 & 1 & 1 \end{pmatrix} \begin{pmatrix} -f^{02} \\ f^{13} \\ f^{20} \\ -f^{31} \\ -f^{46} \\ -f^{57} \\ f^{64} \\ f^{75} \end{pmatrix}$$

$$(4.5.79)\quad \begin{pmatrix} f_2^0 \\ f_3^1 \\ -f_0^2 \\ -f_1^3 \\ -f_6^4 \\ -f_7^5 \\ f_4^6 \\ f_5^7 \end{pmatrix} = \begin{pmatrix} 1 & -1 & -1 & -1 & -1 & -1 & -1 & -1 \\ 1 & -1 & 1 & 1 & 1 & 1 & 1 & 1 \\ 1 & 1 & -1 & 1 & 1 & 1 & 1 & 1 \\ 1 & 1 & 1 & -1 & 1 & 1 & 1 & 1 \\ 1 & 1 & 1 & 1 & -1 & 1 & 1 & 1 \\ 1 & 1 & 1 & 1 & 1 & -1 & 1 & 1 \\ 1 & 1 & 1 & 1 & 1 & 1 & -1 & 1 \\ 1 & 1 & 1 & 1 & 1 & 1 & 1 & -1 \end{pmatrix} \begin{pmatrix} -f^{02} \\ f^{13} \\ f^{20} \\ -f^{31} \\ -f^{46} \\ -f^{57} \\ f^{64} \\ f^{75} \end{pmatrix}$$

Let us write the system of linear equations (4.5.4) as product of matrices

$$(4.5.80)\quad \begin{pmatrix} f_0^3 \\ f_1^2 \\ f_2^1 \\ f_3^0 \\ f_4^7 \\ f_5^6 \\ f_6^5 \\ f_7^4 \end{pmatrix} = \begin{pmatrix} 1 & 1 & -1 & 1 & 1 & -1 & 1 & -1 \\ -1 & -1 & -1 & 1 & -1 & 1 & -1 & 1 \\ 1 & -1 & -1 & -1 & 1 & -1 & 1 & -1 \\ -1 & 1 & -1 & -1 & 1 & -1 & 1 & -1 \\ -1 & -1 & 1 & 1 & -1 & 1 & -1 & -1 \\ 1 & 1 & -1 & -1 & 1 & -1 & -1 & -1 \\ -1 & -1 & 1 & 1 & -1 & -1 & -1 & 1 \\ 1 & 1 & -1 & -1 & -1 & -1 & 1 & -1 \end{pmatrix} \begin{pmatrix} f^{03} \\ f^{12} \\ f^{21} \\ f^{30} \\ f^{47} \\ f^{56} \\ f^{65} \\ f^{74} \end{pmatrix}$$

From the equation (4.5.80), it follows that

$$\begin{pmatrix} -f_0^3 \\ f_1^2 \\ -f_2^1 \\ f_3^0 \\ f_4^7 \\ -f_5^6 \\ f_6^5 \\ -f_7^4 \end{pmatrix} = \begin{pmatrix} -1 & -1 & 1 & -1 & -1 & 1 & -1 & 1 \\ -1 & -1 & -1 & 1 & -1 & 1 & -1 & 1 \\ -1 & 1 & 1 & 1 & -1 & 1 & -1 & 1 \\ -1 & 1 & -1 & -1 & 1 & -1 & 1 & -1 \\ -1 & -1 & 1 & 1 & -1 & 1 & -1 & -1 \\ -1 & -1 & 1 & 1 & -1 & 1 & 1 & 1 \\ -1 & -1 & 1 & 1 & -1 & -1 & -1 & 1 \\ -1 & -1 & 1 & 1 & 1 & 1 & -1 & 1 \end{pmatrix} \begin{pmatrix} f^{03} \\ f^{12} \\ f^{21} \\ f^{30} \\ f^{47} \\ f^{56} \\ f^{65} \\ f^{74} \end{pmatrix}$$

$$= \begin{pmatrix} 1 & 1 & 1 & -1 & 1 & 1 & 1 & 1 \\ 1 & 1 & -1 & 1 & 1 & 1 & 1 & 1 \\ 1 & -1 & 1 & 1 & 1 & 1 & 1 & 1 \\ 1 & -1 & -1 & -1 & -1 & -1 & -1 & -1 \\ 1 & 1 & 1 & 1 & 1 & 1 & 1 & -1 \\ 1 & 1 & 1 & 1 & 1 & 1 & -1 & 1 \\ 1 & 1 & 1 & 1 & 1 & -1 & 1 & 1 \\ 1 & 1 & 1 & 1 & -1 & 1 & 1 & 1 \end{pmatrix} \begin{pmatrix} -f^{03} \\ -f^{12} \\ f^{21} \\ f^{30} \\ -f^{47} \\ f^{56} \\ -f^{65} \\ f^{74} \end{pmatrix}$$

$$(4.5.81) \quad \begin{pmatrix} f_3^0 \\ -f_2^1 \\ f_1^2 \\ -f_0^3 \\ -f_7^4 \\ f_6^5 \\ -f_5^6 \\ f_4^7 \end{pmatrix} = \begin{pmatrix} 1 & -1 & -1 & -1 & -1 & -1 & -1 & -1 \\ 1 & -1 & 1 & 1 & 1 & 1 & 1 & 1 \\ 1 & 1 & -1 & 1 & 1 & 1 & 1 & 1 \\ 1 & 1 & 1 & -1 & 1 & 1 & 1 & 1 \\ 1 & 1 & 1 & 1 & -1 & 1 & 1 & 1 \\ 1 & 1 & 1 & 1 & 1 & -1 & 1 & 1 \\ 1 & 1 & 1 & 1 & 1 & 1 & -1 & 1 \\ 1 & 1 & 1 & 1 & 1 & 1 & 1 & -1 \end{pmatrix} \begin{pmatrix} -f^{03} \\ -f^{12} \\ f^{21} \\ f^{30} \\ -f^{47} \\ f^{56} \\ -f^{65} \\ f^{74} \end{pmatrix}$$

Let us write the system of linear equations (4.5.5) as product of matrices

$$(4.5.82) \quad \begin{pmatrix} f_0^4 \\ f_1^5 \\ f_2^6 \\ f_3^7 \\ f_4^0 \\ f_5^1 \\ f_6^2 \\ f_7^3 \end{pmatrix} = \begin{pmatrix} 1 & -1 & -1 & -1 & 1 & 1 & 1 & 1 \\ 1 & -1 & -1 & -1 & -1 & -1 & 1 & 1 \\ 1 & -1 & -1 & -1 & -1 & 1 & -1 & 1 \\ 1 & -1 & -1 & -1 & -1 & 1 & 1 & -1 \\ -1 & -1 & -1 & -1 & -1 & 1 & 1 & 1 \\ -1 & -1 & 1 & 1 & 1 & -1 & -1 & -1 \\ -1 & 1 & -1 & 1 & 1 & -1 & -1 & -1 \\ -1 & 1 & 1 & -1 & 1 & -1 & -1 & -1 \end{pmatrix} \begin{pmatrix} f^{04} \\ f^{15} \\ f^{26} \\ f^{37} \\ f^{40} \\ f^{51} \\ f^{62} \\ f^{73} \end{pmatrix}$$

From the equation (4.5.82), it follows that

$$\begin{pmatrix} -f_0^4 \\ -f_1^5 \\ -f_2^6 \\ -f_3^7 \\ f_4^0 \\ f_5^1 \\ f_6^2 \\ f_7^3 \end{pmatrix} = \begin{pmatrix} -1 & 1 & 1 & 1 & -1 & -1 & -1 & -1 \\ -1 & 1 & 1 & 1 & 1 & 1 & -1 & -1 \\ -1 & 1 & 1 & 1 & 1 & -1 & 1 & -1 \\ -1 & 1 & 1 & 1 & 1 & -1 & -1 & 1 \\ -1 & -1 & -1 & -1 & -1 & 1 & 1 & 1 \\ -1 & -1 & 1 & 1 & 1 & -1 & -1 & -1 \\ -1 & 1 & -1 & 1 & 1 & -1 & -1 & -1 \\ -1 & 1 & 1 & -1 & 1 & -1 & -1 & -1 \end{pmatrix} \begin{pmatrix} f^{04} \\ f^{15} \\ f^{26} \\ f^{37} \\ f^{40} \\ f^{51} \\ f^{62} \\ f^{73} \end{pmatrix}$$

$$= \begin{pmatrix} 1 & 1 & 1 & 1 & -1 & 1 & 1 & 1 \\ 1 & 1 & 1 & 1 & 1 & -1 & 1 & 1 \\ 1 & 1 & 1 & 1 & 1 & 1 & -1 & 1 \\ 1 & 1 & 1 & 1 & 1 & 1 & 1 & -1 \\ 1 & -1 & -1 & -1 & -1 & -1 & -1 & -1 \\ 1 & -1 & 1 & 1 & 1 & 1 & 1 & 1 \\ 1 & 1 & -1 & 1 & 1 & 1 & 1 & 1 \\ 1 & 1 & 1 & -1 & 1 & 1 & 1 & 1 \end{pmatrix} \begin{pmatrix} -f^{04} \\ f^{15} \\ f^{26} \\ f^{37} \\ f^{40} \\ -f^{51} \\ -f^{62} \\ -f^{73} \end{pmatrix}$$

(4.5.83)
$$\begin{pmatrix} f_4^0 \\ f_5^1 \\ f_6^2 \\ f_7^3 \\ -f_0^4 \\ -f_1^5 \\ -f_2^6 \\ -f_3^7 \end{pmatrix} = \begin{pmatrix} 1 & -1 & -1 & -1 & -1 & -1 & -1 & -1 \\ 1 & -1 & 1 & 1 & 1 & 1 & 1 & 1 \\ 1 & 1 & -1 & 1 & 1 & 1 & 1 & 1 \\ 1 & 1 & 1 & -1 & 1 & 1 & 1 & 1 \\ 1 & 1 & 1 & 1 & -1 & 1 & 1 & 1 \\ 1 & 1 & 1 & 1 & 1 & -1 & 1 & 1 \\ 1 & 1 & 1 & 1 & 1 & 1 & -1 & 1 \\ 1 & 1 & 1 & 1 & 1 & 1 & 1 & -1 \end{pmatrix} \begin{pmatrix} -f^{04} \\ f^{15} \\ f^{26} \\ f^{37} \\ f^{40} \\ -f^{51} \\ -f^{62} \\ -f^{73} \end{pmatrix}$$

Let us write the system of linear equations (4.5.6) as product of matrices

(4.5.84)
$$\begin{pmatrix} f_0^5 \\ f_1^4 \\ f_2^7 \\ f_3^6 \\ f_4^1 \\ f_5^0 \\ f_6^3 \\ f_7^2 \end{pmatrix} = \begin{pmatrix} 1 & 1 & -1 & 1 & -1 & 1 & -1 & 1 \\ -1 & -1 & 1 & -1 & -1 & 1 & 1 & -1 \\ 1 & 1 & -1 & 1 & -1 & -1 & -1 & -1 \\ -1 & -1 & 1 & -1 & 1 & 1 & -1 & -1 \\ 1 & -1 & -1 & 1 & -1 & -1 & -1 & 1 \\ -1 & 1 & -1 & 1 & -1 & -1 & -1 & 1 \\ 1 & 1 & -1 & -1 & -1 & -1 & -1 & 1 \\ -1 & -1 & -1 & -1 & 1 & 1 & 1 & -1 \end{pmatrix} \begin{pmatrix} f^{05} \\ f^{14} \\ f^{27} \\ f^{36} \\ f^{41} \\ f^{50} \\ f^{63} \\ f^{72} \end{pmatrix}$$

From the equation (4.5.84), it follows that

$$\begin{pmatrix} -f_0^5 \\ f_1^4 \\ -f_2^7 \\ f_3^6 \\ -f_4^1 \\ f_5^0 \\ -f_6^3 \\ f_7^2 \end{pmatrix} = \begin{pmatrix} -1 & -1 & 1 & -1 & 1 & -1 & 1 & -1 \\ -1 & -1 & 1 & -1 & -1 & 1 & 1 & -1 \\ -1 & -1 & 1 & -1 & 1 & 1 & 1 & 1 \\ -1 & -1 & 1 & -1 & 1 & 1 & -1 & -1 \\ -1 & 1 & 1 & -1 & 1 & 1 & 1 & -1 \\ -1 & 1 & -1 & 1 & -1 & -1 & -1 & 1 \\ -1 & -1 & 1 & 1 & 1 & 1 & 1 & -1 \\ -1 & -1 & -1 & -1 & 1 & 1 & 1 & -1 \end{pmatrix} \begin{pmatrix} f^{05} \\ f^{14} \\ f^{27} \\ f^{36} \\ f^{41} \\ f^{50} \\ f^{63} \\ f^{72} \end{pmatrix}$$

$$= \begin{pmatrix} 1 & 1 & 1 & 1 & 1 & -1 & 1 & 1 \\ 1 & 1 & 1 & 1 & -1 & 1 & 1 & 1 \\ 1 & 1 & 1 & 1 & 1 & 1 & 1 & -1 \\ 1 & 1 & 1 & 1 & 1 & 1 & -1 & 1 \\ 1 & -1 & 1 & 1 & 1 & 1 & 1 & 1 \\ 1 & -1 & -1 & -1 & -1 & -1 & -1 & -1 \\ 1 & 1 & 1 & -1 & 1 & 1 & 1 & 1 \\ 1 & 1 & -1 & 1 & 1 & 1 & 1 & 1 \end{pmatrix} \begin{pmatrix} -f^{05} \\ -f^{14} \\ f^{27} \\ -f^{36} \\ f^{41} \\ f^{50} \\ f^{63} \\ -f^{72} \end{pmatrix}$$

$$(4.5.85) \quad \begin{pmatrix} f_5^0 \\ -f_4^1 \\ f_7^2 \\ -f_6^3 \\ f_1^4 \\ -f_0^5 \\ f_3^6 \\ -f_2^7 \end{pmatrix} = \begin{pmatrix} 1 & -1 & -1 & -1 & -1 & -1 & -1 & -1 \\ 1 & -1 & 1 & 1 & 1 & 1 & 1 & 1 \\ 1 & 1 & -1 & 1 & 1 & 1 & 1 & 1 \\ 1 & 1 & 1 & -1 & 1 & 1 & 1 & 1 \\ 1 & 1 & 1 & 1 & -1 & 1 & 1 & 1 \\ 1 & 1 & 1 & 1 & 1 & -1 & 1 & 1 \\ 1 & 1 & 1 & 1 & 1 & 1 & -1 & 1 \\ 1 & 1 & 1 & 1 & 1 & 1 & 1 & -1 \end{pmatrix} \begin{pmatrix} -f^{05} \\ -f^{14} \\ f^{27} \\ -f^{36} \\ f^{41} \\ f^{50} \\ f^{63} \\ -f^{72} \end{pmatrix}$$

Let us write the system of linear equations (4.5.7) as product of matrices

$$(4.5.86) \quad \begin{pmatrix} f_0^6 \\ f_1^7 \\ f_2^4 \\ f_3^5 \\ f_4^2 \\ f_5^3 \\ f_6^0 \\ f_7^1 \end{pmatrix} = \begin{pmatrix} 1 & 1 & 1 & -1 & -1 & 1 & 1 & -1 \\ -1 & -1 & -1 & 1 & 1 & -1 & 1 & -1 \\ -1 & -1 & -1 & 1 & -1 & -1 & 1 & 1 \\ 1 & 1 & 1 & -1 & -1 & -1 & -1 & -1 \\ 1 & 1 & -1 & -1 & -1 & 1 & -1 & -1 \\ -1 & -1 & -1 & -1 & 1 & -1 & 1 & 1 \\ -1 & 1 & 1 & -1 & -1 & 1 & -1 & -1 \\ 1 & -1 & 1 & -1 & -1 & 1 & -1 & -1 \end{pmatrix} \begin{pmatrix} f^{06} \\ f^{17} \\ f^{24} \\ f^{35} \\ f^{42} \\ f^{53} \\ f^{60} \\ f^{71} \end{pmatrix}$$

From the equation (4.5.86), it follows that

$$\begin{pmatrix} -f_0^6 \\ f_1^7 \\ f_2^4 \\ -f_3^5 \\ -f_4^2 \\ f_5^3 \\ f_6^0 \\ -f_7^1 \end{pmatrix} = \begin{pmatrix} -1 & -1 & -1 & 1 & 1 & -1 & -1 & 1 \\ -1 & -1 & -1 & 1 & 1 & -1 & 1 & -1 \\ -1 & -1 & -1 & 1 & -1 & -1 & 1 & 1 \\ -1 & -1 & -1 & 1 & 1 & 1 & 1 & 1 \\ -1 & -1 & 1 & 1 & 1 & -1 & 1 & 1 \\ -1 & -1 & -1 & -1 & 1 & -1 & 1 & 1 \\ -1 & 1 & 1 & -1 & -1 & 1 & -1 & -1 \\ -1 & 1 & -1 & 1 & 1 & -1 & 1 & 1 \end{pmatrix} \begin{pmatrix} f^{06} \\ f^{17} \\ f^{24} \\ f^{35} \\ f^{42} \\ f^{53} \\ f^{60} \\ f^{71} \end{pmatrix}$$

$$= \begin{pmatrix} 1 & 1 & 1 & 1 & 1 & 1 & -1 & 1 \\ 1 & 1 & 1 & 1 & 1 & 1 & 1 & -1 \\ 1 & 1 & 1 & 1 & -1 & 1 & 1 & 1 \\ 1 & 1 & 1 & 1 & 1 & -1 & 1 & 1 \\ 1 & 1 & -1 & 1 & 1 & 1 & 1 & 1 \\ 1 & 1 & 1 & -1 & 1 & 1 & 1 & 1 \\ 1 & -1 & -1 & -1 & -1 & -1 & -1 & -1 \\ 1 & -1 & 1 & 1 & 1 & 1 & 1 & 1 \end{pmatrix} \begin{pmatrix} -f^{06} \\ -f^{17} \\ -f^{24} \\ f^{35} \\ f^{42} \\ -f^{53} \\ f^{60} \\ f^{71} \end{pmatrix}$$

$$(4.5.87) \quad \begin{pmatrix} f_6^0 \\ -f_7^1 \\ -f_4^2 \\ f_5^3 \\ f_2^4 \\ -f_3^5 \\ -f_0^6 \\ f_1^7 \end{pmatrix} = \begin{pmatrix} 1 & -1 & -1 & -1 & -1 & -1 & -1 & -1 \\ 1 & -1 & 1 & 1 & 1 & 1 & 1 & 1 \\ 1 & 1 & -1 & 1 & 1 & 1 & 1 & 1 \\ 1 & 1 & 1 & -1 & 1 & 1 & 1 & 1 \\ 1 & 1 & 1 & 1 & -1 & 1 & 1 & 1 \\ 1 & 1 & 1 & 1 & 1 & -1 & 1 & 1 \\ 1 & 1 & 1 & 1 & 1 & 1 & -1 & 1 \\ 1 & 1 & 1 & 1 & 1 & 1 & 1 & -1 \end{pmatrix} \begin{pmatrix} -f^{06} \\ -f^{17} \\ -f^{24} \\ f^{35} \\ f^{42} \\ -f^{53} \\ f^{60} \\ f^{71} \end{pmatrix}$$

Let us write the system of linear equations (4.5.8) as product of matrices

$$(4.5.88) \quad \begin{pmatrix} f_0^7 \\ f_1^6 \\ f_2^5 \\ f_3^4 \\ f_4^3 \\ f_5^2 \\ f_6^1 \\ f_7^0 \end{pmatrix} = \begin{pmatrix} 1 & -1 & 1 & 1 & -1 & -1 & 1 & 1 \\ 1 & -1 & 1 & 1 & -1 & -1 & -1 & -1 \\ -1 & 1 & -1 & -1 & 1 & -1 & -1 & 1 \\ -1 & 1 & -1 & -1 & -1 & 1 & -1 & 1 \\ 1 & -1 & 1 & -1 & -1 & -1 & 1 & -1 \\ 1 & -1 & -1 & 1 & -1 & -1 & 1 & -1 \\ -1 & -1 & -1 & -1 & 1 & 1 & -1 & 1 \\ -1 & -1 & 1 & 1 & -1 & -1 & 1 & -1 \end{pmatrix} \begin{pmatrix} f^{07} \\ f^{16} \\ f^{25} \\ f^{34} \\ f^{43} \\ f^{52} \\ f^{61} \\ f^{70} \end{pmatrix}$$

From the equation (4.5.88), it follows that

$$\begin{pmatrix} -f_0^7 \\ -f_1^6 \\ f_2^5 \\ f_3^4 \\ -f_4^3 \\ -f_5^2 \\ f_6^1 \\ f_7^0 \end{pmatrix} = \begin{pmatrix} -1 & 1 & -1 & -1 & 1 & 1 & -1 & -1 \\ -1 & 1 & -1 & -1 & 1 & 1 & 1 & 1 \\ -1 & 1 & -1 & -1 & 1 & -1 & -1 & 1 \\ -1 & 1 & -1 & -1 & -1 & 1 & -1 & 1 \\ -1 & 1 & -1 & 1 & 1 & 1 & -1 & 1 \\ -1 & 1 & 1 & -1 & 1 & 1 & -1 & 1 \\ -1 & -1 & -1 & -1 & 1 & 1 & -1 & 1 \\ -1 & -1 & 1 & 1 & -1 & -1 & 1 & -1 \end{pmatrix} \begin{pmatrix} f^{07} \\ f^{16} \\ f^{25} \\ f^{34} \\ f^{43} \\ f^{52} \\ f^{61} \\ f^{70} \end{pmatrix}$$

$$= \begin{pmatrix} 1 & 1 & 1 & 1 & 1 & 1 & 1 & -1 \\ 1 & 1 & 1 & 1 & 1 & 1 & -1 & 1 \\ 1 & 1 & 1 & 1 & 1 & -1 & 1 & 1 \\ 1 & 1 & 1 & 1 & -1 & 1 & 1 & 1 \\ 1 & 1 & 1 & -1 & 1 & 1 & 1 & 1 \\ 1 & 1 & -1 & 1 & 1 & 1 & 1 & 1 \\ 1 & -1 & 1 & 1 & 1 & 1 & 1 & 1 \\ 1 & -1 & -1 & -1 & -1 & -1 & -1 & -1 \end{pmatrix} \begin{pmatrix} -f^{07} \\ f^{16} \\ -f^{25} \\ -f^{34} \\ f^{43} \\ f^{52} \\ -f^{61} \\ f^{70} \end{pmatrix}$$

$$(4.5.89) \quad \begin{pmatrix} f_7^0 \\ f_6^1 \\ -f_5^2 \\ -f_4^3 \\ f_3^4 \\ f_2^5 \\ -f_1^6 \\ -f_0^7 \end{pmatrix} = \begin{pmatrix} 1 & -1 & -1 & -1 & -1 & -1 & -1 & -1 \\ 1 & -1 & 1 & 1 & 1 & 1 & 1 & 1 \\ 1 & 1 & -1 & 1 & 1 & 1 & 1 & 1 \\ 1 & 1 & 1 & -1 & 1 & 1 & 1 & 1 \\ 1 & 1 & 1 & 1 & -1 & 1 & 1 & 1 \\ 1 & 1 & 1 & 1 & 1 & -1 & 1 & 1 \\ 1 & 1 & 1 & 1 & 1 & 1 & -1 & 1 \\ 1 & 1 & 1 & 1 & 1 & 1 & 1 & -1 \end{pmatrix} \begin{pmatrix} -f^{07} \\ f^{16} \\ -f^{25} \\ -f^{34} \\ f^{43} \\ f^{52} \\ -f^{61} \\ f^{70} \end{pmatrix}$$

We join equations (4.5.75), (4.5.77), (4.5.79), (4.5.81), (4.5.83), (4.5.85), (4.5.87), (4.5.89) into equation (4.5.73). □

**Theorem 4.5.3.** *Standard components of linear function of octonion algebra $O$ relative to basis (4.4.1) and coordinates of corresponding linear map satisfy relationship*

$$(4.5.90) \quad \begin{cases} 12f^{00} = 5f_0^0 + f_1^1 + f_2^2 + f_3^3 + f_4^4 + f_5^5 + f_6^6 + f_7^7 \\ 12f^{11} = -f_0^0 - 5f_1^1 + f_2^2 + f_3^3 + f_4^4 + f_5^5 + f_6^6 + f_7^7 \\ 12f^{22} = -f_0^0 + f_1^1 - 5f_2^2 + f_3^3 + f_4^4 + f_5^5 + f_6^6 + f_7^7 \\ 12f^{33} = -f_0^0 + f_1^1 + f_2^2 - 5f_3^3 + f_4^4 + f_5^5 + f_6^6 + f_7^7 \\ 12f^{44} = -f_0^0 + f_1^1 + f_2^2 + f_3^3 - 5f_4^4 + f_5^5 + f_6^6 + f_7^7 \\ 12f^{55} = -f_0^0 + f_1^1 + f_2^2 + f_3^3 + f_4^4 - 5f_5^5 + f_6^6 + f_7^7 \\ 12f^{66} = -f_0^0 + f_1^1 + f_2^2 + f_3^3 + f_4^4 + f_5^5 - 5f_6^6 + f_7^7 \\ 12f^{77} = -f_0^0 + f_1^1 + f_2^2 + f_3^3 + f_4^4 + f_5^5 + f_6^6 - 5f_7^7 \end{cases}$$

$$(4.5.91) \quad \begin{cases} -12f^{01} = 5f_1^0 - f_0^1 - f_3^2 + f_2^3 - f_5^4 + f_4^5 + f_7^6 - f_6^7 \\ 12f^{10} = -f_1^0 + 5f_0^1 - f_3^2 + f_2^3 - f_5^4 + f_4^5 + f_7^6 - f_6^7 \\ -12f^{23} = -f_1^0 - f_0^1 + 5f_3^2 + f_2^3 - f_5^4 + f_4^5 + f_7^6 - f_6^7 \\ 12f^{32} = -f_1^0 - f_0^1 - f_3^2 - 5f_2^3 - f_5^4 + f_4^5 + f_7^6 - f_6^7 \\ -12f^{45} = -f_1^0 - f_0^1 - f_3^2 + f_2^3 + 5f_5^4 + f_4^5 + f_7^6 - f_6^7 \\ 12f^{54} = -f_1^0 - f_0^1 - f_3^2 + f_2^3 - f_5^4 - 5f_4^5 + f_7^6 - f_6^7 \\ 12f^{67} = -f_1^0 - f_0^1 - f_3^2 + f_2^3 - f_5^4 + f_4^5 - 5f_7^6 - f_6^7 \\ -12f^{76} = -f_1^0 - f_0^1 - f_3^2 + f_2^3 - f_5^4 + f_4^5 + f_7^6 + 5f_6^7 \end{cases}$$

$$(4.5.92) \quad \begin{cases} -12f^{02} = 5f_2^0 + f_3^1 - f_0^2 - f_1^3 - f_6^4 - f_7^5 + f_4^6 + f_5^7 \\ 12f^{13} = -f_2^0 - 5f_3^1 - f_0^2 - f_1^3 - f_6^4 - f_7^5 + f_4^6 + f_5^7 \\ 12f^{20} = -f_2^0 + f_3^1 + 5f_0^2 - f_1^3 - f_6^4 - f_7^5 + f_4^6 + f_5^7 \\ -12f^{31} = -f_2^0 + f_3^1 - f_0^2 + 5f_1^3 - f_6^4 - f_7^5 + f_4^6 + f_5^7 \\ -12f^{46} = -f_2^0 + f_3^1 - f_0^2 - f_1^3 + 5f_6^4 - f_7^5 + f_4^6 + f_5^7 \\ -12f^{57} = -f_2^0 + f_3^1 - f_0^2 - f_1^3 - f_6^4 + 5f_7^5 + f_4^6 + f_5^7 \\ 12f^{64} = -f_2^0 + f_3^1 - f_0^2 - f_1^3 - f_6^4 - f_7^5 - 5f_4^6 + f_5^7 \\ 12f^{75} = -f_2^0 + f_3^1 - f_0^2 - f_1^3 - f_6^4 - f_7^5 + f_4^6 - 5f_5^7 \end{cases}$$

$$(4.5.93) \begin{cases} -12f^{03} = 5f_3^0 - f_2^1 + f_1^2 - f_0^3 - f_7^4 + f_6^5 - f_5^6 + f_4^7 \\ -12f^{12} = -f_3^0 + 5f_2^1 + f_1^2 - f_0^3 - f_7^4 + f_6^5 - f_5^6 + f_4^7 \\ 12f^{21} = -f_3^0 - f_2^1 - 5f_1^2 - f_0^3 - f_7^4 + f_6^5 - f_5^6 + f_4^7 \\ 12f^{30} = -f_3^0 - f_2^1 + f_1^2 + 5f_0^3 - f_7^4 + f_6^5 - f_5^6 + f_4^7 \\ -12f^{47} = -f_3^0 - f_2^1 + f_1^2 - f_0^3 + 5f_7^4 + f_6^5 - f_5^6 + f_4^7 \\ 12f^{56} = -f_3^0 - f_2^1 + f_1^2 - f_0^3 - f_7^4 - 5f_6^5 - f_5^6 + f_4^7 \\ -12f^{65} = -f_3^0 - f_2^1 + f_1^2 - f_0^3 - f_7^4 + f_6^5 + 5f_5^6 + f_4^7 \\ 12f^{74} = -f_3^0 - f_2^1 + f_1^2 - f_0^3 - f_7^4 + f_6^5 - f_5^6 - 5f_4^7 \end{cases}$$

$$(4.5.94) \begin{cases} -12f^{04} = 5f_4^0 + f_5^1 + f_6^2 + f_7^3 - f_0^4 - f_1^5 - f_2^6 - f_3^7 \\ 12f^{15} = -f_4^0 - 5f_5^1 + f_6^2 + f_7^3 - f_0^4 - f_1^5 - f_2^6 - f_3^7 \\ 12f^{26} = -f_4^0 + f_5^1 - 5f_6^2 + f_7^3 - f_0^4 - f_1^5 - f_2^6 - f_3^7 \\ 12f^{37} = -f_4^0 + f_5^1 + f_6^2 - 5f_7^3 - f_0^4 - f_1^5 - f_2^6 - f_3^7 \\ 12f^{40} = -f_4^0 + f_5^1 + f_6^2 + f_7^3 + 5f_0^4 - f_1^5 - f_2^6 - f_3^7 \\ -12f^{51} = -f_4^0 + f_5^1 + f_6^2 + f_7^3 - f_0^4 + 5f_1^5 - f_2^6 - f_3^7 \\ -12f^{62} = -f_4^0 + f_5^1 + f_6^2 + f_7^3 - f_0^4 - f_1^5 + 5f_2^6 - f_3^7 \\ -12f^{73} = -f_4^0 + f_5^1 + f_6^2 + f_7^3 - f_0^4 - f_1^5 - f_2^6 + 5f_3^7 \end{cases}$$

$$(4.5.95) \begin{cases} -12f^{05} = 5f_5^0 - f_4^1 + f_7^2 - f_6^3 + f_1^4 - f_0^5 + f_3^6 - f_2^7 \\ -12f^{14} = -f_5^0 + 5f_4^1 + f_7^2 - f_6^3 + f_1^4 - f_0^5 + f_3^6 - f_2^7 \\ 12f^{27} = -f_5^0 - f_4^1 - 5f_7^2 - f_6^3 + f_1^4 - f_0^5 + f_3^6 - f_2^7 \\ -12f^{36} = -f_5^0 - f_4^1 + f_7^2 + 5f_6^3 + f_1^4 - f_0^5 + f_3^6 - f_2^7 \\ 12f^{41} = -f_5^0 - f_4^1 + f_7^2 - f_6^3 - 5f_1^4 - f_0^5 + f_3^6 - f_2^7 \\ 12f^{50} = -f_5^0 - f_4^1 + f_7^2 - f_6^3 + f_1^4 + 5f_0^5 + f_3^6 - f_2^7 \\ 12f^{63} = -f_5^0 - f_4^1 + f_7^2 - f_6^3 + f_1^4 - f_0^5 - 5f_3^6 - f_2^7 \\ 12f^{72} = -f_5^0 - f_4^1 + f_7^2 - f_6^3 + f_1^4 - f_0^5 + f_3^6 + 5f_2^7 \end{cases}$$

$$(4.5.96) \begin{cases} -12f^{06} = 5f_6^0 - f_7^1 - f_4^2 + f_5^3 + f_2^4 - f_3^5 - f_0^6 + f_1^7 \\ -12f^{17} = -f_6^0 + 5f_7^1 - f_4^2 + f_5^3 + f_2^4 - f_3^5 - f_0^6 + f_1^7 \\ -12f^{24} = -f_6^0 - f_7^1 + 5f_4^2 + f_5^3 + f_2^4 - f_3^5 - f_0^6 + f_1^7 \\ 12f^{35} = -f_6^0 - f_7^1 - f_4^2 - 5f_5^3 + f_2^4 - f_3^5 - f_0^6 + f_1^7 \\ 12f^{42} = -f_6^0 - f_7^1 - f_4^2 + f_5^3 - 5f_2^4 - f_3^5 - f_0^6 + f_1^7 \\ -12f^{53} = -f_6^0 - f_7^1 - f_4^2 + f_5^3 + f_2^4 + 5f_3^5 - f_0^6 + f_1^7 \\ 12f^{60} = -f_6^0 - f_7^1 - f_4^2 + f_5^3 + f_2^4 - f_3^5 + 5f_0^6 + f_1^7 \\ 12f^{71} = -f_6^0 - f_7^1 - f_4^2 + f_5^3 + ff_2^4 - f_3^5 - f_0^6 - 5f_1^7 \end{cases}$$

$$(4.5.97) \begin{cases} -12f^{07} = 5f_7^0 + f_6^1 - f_5^2 - f_4^3 + f_3^4 + f_2^5 - f_1^6 - f_0^7 \\ 12f^{16} = -f_7^0 - 5f_6^1 - f_5^2 - f_4^3 + f_3^4 + f_2^5 - f_1^6 - f_0^7 \\ -12f^{25} = -f_7^0 + f_6^1 + 5f_5^2 - f_4^3 + f_3^4 + f_2^5 - f_1^6 - f_0^7 \\ -12f^{34} = -f_7^0 + f_6^1 - f_5^2 + 5f_4^3 + f_3^4 + f_2^5 - f_1^6 - f_0^7 \\ 12f^{43} = -f_7^0 + f_6^1 - f_5^2 - f_4^3 - 5f_3^4 + f_2^5 - f_1^6 - f_0^7 \\ 12f^{52} = -f_7^0 + f_6^1 - f_5^2 - f_4^3 + f_3^4 - 5f_2^5 - f_1^6 - f_0^7 \\ -12f^{61} = -f_7^0 + f_6^1 - f_5^2 - f_4^3 + f_3^4 + f_2^5 + 5f_1^6 - f_0^7 \\ 12f^{70} - f_7^0 + f_6^1 - f_5^2 - f_4^3 + f_3^4 + f_2^5 - f_1^6 + 5f_0^7 \end{cases}$$

PROOF. We get systems of linear equations (4.5.90), (4.5.91), (4.5.92), (4.5.93), (4.5.94), (4.5.95), (4.5.96), (4.5.97) as the product of matrices in equation (4.5.74).
□

To find the linear mapping corresponding to the operation of conjugation, I will assume

$$(4.5.98) \qquad f_0^0 = 1 \quad f_1^1 = f_2^2 = f_3^3 = f_4^4 = f_5^5 = f_6^6 = f_7^7 = -1$$

Substituting (4.5.98) in the system of equations (4.5.90), we get

$$(4.5.99) \qquad f^{00} = f^{11} = f^{22} = f^{33} = f^{44} = f^{55} = f^{66} = f^{77} = -\frac{1}{6}$$

Therefore,

$$(4.5.100) \quad \overline{z} = -\frac{1}{6}(z + (iz)i + (jz)j + (kz)k + ((il)z)(il) + ((jl)z)(jl) + ((kl)z)(kl))$$

# CHAPTER 5

# CHAPTER 6

# Index

# CHAPTER 7

# Special Symbols and Notations

# Линейные отображения свободной алгебры

Александр Клейн

*E-mail address*: Aleks_Kleyn@MailAPS.org
*URL*: http://sites.google.com/site/alekskleyn/
*URL*: http://arxiv.org/a/kleyn_a_1
*URL*: http://AleksKleyn.blogspot.com/

Аннотация. Для произвольной универсальной алгебры, в которой определена операция сложения, я изучаю бикольцо матриц отображений. Сумма матриц определена суммой в алгебре, а произведение матриц определено произведением отображений. Система аддитивных уравнений - это система уравнений, матрица которой является матрица отображений. Рассмотрены методы решения системы аддитивных уравнений. В качестве примера приведено решение системы линейных уравнений над полем комплексных чисел при условии, что уравнения содержат неизвестные величины и величины, сопряжённые им.

Линейные отображения алгебры над коммутативным кольцом сохраняют операцию сложения в алгебре и умножение элементов алгебры на элементы кольца. Множество линейных преобразований алгебры $A$ порождено представлением тензорного произведения $A \otimes A$ в алгебре $A$.

Результаты этого исследования будут полезны для математиков и физиков, которые работают с различными алгебрами.

# Оглавление

Глава 1

# Предисловие

## 1.1. Предисловие к изданию 1

Когда я приступил к написанию этой книги, первоначальная задача была довольно простой. Я собирался переписать содержимое книги [5], используя аппарат матриц отображений ([8]). Однако, мне показалось странным ограничивать себя рассмотрением тела, когда я понимал, что многие результаты будут верны для ассоциативной алгебры, а что-то сохранится в случае неассоциативной алгебры. Результаты этого исследования будут полезны для математиков и физиков, которые работают с различными алгебрами, необязательно ассоциативными.

Когда я записал систему линейных уравнений, используя матрицы отображений, я понял, что в моих руках оказался инструмент более мощный, чем я вначале предполагал. При изучении систем линейных уравнений мы рассматриваем умножение неизвестной величины, принадлежащей кольцу или векторному пространству, на скаляр из соответствующего кольца. Однако я могу предположить, что неизвестная величина принадлежат некоторой универсальной алгебре, имеющей операцию сложения. Вместо умножения на скаляр я рассматриваю некоторое отображение универсальной алгебры. Так возникла теория аддитивных уравнений, которая во многом похожа на теорию линейных уравнений.

В качестве примера применения новых методов приведено решение системы линейных уравнений над полем комплексных чисел при условии, что уравнения содержат неизвестные величины и величины, сопряжённые им. Решение подобной системы уравнений подробно рассмотрено в примере 2.5.5. Без сомнения, попытка решить систему уравнений (2.5.20), пользуясь определителем, задача непростая.

Переходя к рассмотрению алгебр, я обратил внимание, что обычно алгебру определяют над полем. Это видимо необходимо, так как алгебра является векторным пространством. Однако с точки зрения тех построений, что меня интересуют, для меня несущественно является ли алгебра векторным пространством над полем или свободным модулем над коммутативным кольцом. Убедившись, что существует исследование алгебр над кольцом, я решил работать с алгебрами над коммутативным кольцо при условии, что в случае необходимости я буду ослаблять требования.

На алгебре существует две алгебраические структуры. Если мы рассматриваем алгебру как кольцо, то при отображении одной алгебры в другую мы рассматриваем гомоморфизмы кольца. Если мы рассматриваем алгебру как модуль над кольцом, то при отображении одной алгебры в другую мы рассматриваем линейные отображения. Очевидно, что если алгебра имеет единицу, то гомоморфизм алгебры является линейным отображением. Меня в основном интересуют линейные отображения алгебры.

С того момента, как я рассмотрел тензорное произведение тел (раздел [5]-12.2), меня не оставляло ощущение, что линейное отображение в теле выражено тензором валентности 2. В тоже время было непонятно, каким образом тензор валентности 2 описывает линейное отображение, хотя для определения тензора валентности 2 мне нужно билинейное отображение.

Структура модуля линейных отображений $\mathcal{L}(A;A)$ определена некоммутативностью произведения в алгебре $A$. Как только произведение становится коммутативным я могу вместо выражения $axb$ записать выражение $abx$ и я увижу тензор валентности 1 там, где вначале был тензор валентности 2.

Алгебра $A \otimes A$ - очень интересная алгебра. Формально мне следовало бы писать $A \otimes A^*$, где $A^*$ - противоположная алгебра. Однако это привело бы к некоторым проблемам в записи действия

$$(1.1.1) \qquad (a \otimes b) \circ x = axb$$

ибо становится неясно с какой стороны следует писать $b$. Определение умножения

$$(a \otimes b) \circ (c \otimes d) = (ac) \otimes (db)$$

позволяет сохранить запись (1.1.1). Поэтому я предпочёл оставить обозначение $A \otimes A$.

Если алгебра $A$ является свободной конечно мерной ассоциативной алгеброй, то базис представления алгебры $A \otimes A$ в модуле $\mathcal{L}(A;A)$ конечен и позволяет описать все линейные отображения алгебры $A$.

<div style="text-align: right;">Март, 2010</div>

## 1.2. Предисловие к изданию 2

Вскоре после того как я опубликовал издание 1, я прочитал мнение профессора Баеза ([10]) где он рассуждал на тему о роли блога в среде математиков. В частности, Баез порекомендовал посетить сайт http://www.ncatlab.org/nlab/show/Online+Resources. То что я увидел на этой сайт превзошло мои ожидания.

Я провёл немало времени, чтобы понять, какие проблемы в математике интересуют людей, создавших эту сайт. Я обратил внимание, что в статье посвящённой $\Omega$-группе рассматривается конструкция, подобная конструкции, рассмотренной у меня в главе 2.

Термин $\Omega$-группа существует только на этой сайт. На ссылках, которые я имею, речь идёт о группе с операторами, что соответствует представлению $\Omega$-алгебры в группе (обычно записываемой аддитивно). Поэтому я решил не менять терминологию в этой книге. Однако я вернусь к этой теме позже.

Меня заинтересовала возможность рассматривать некоммутативное сложение. Но я встретил проблему определить множество аддитивных отображений. Я также надеюсь вернуться к этой теме.

<div style="text-align: right;">Август, 2010</div>

## 1.3. Соглашения

(1) Функция и отображение - синонимы. Однако существует традиция соответствие между кольцами или векторными пространствами называть отображением, а отображение поля действительных чисел или алгебры кватернионов называть функцией. Я тоже следую этой традиции, хотя встречается текст, в котором неясно, какому термину надо отдать предпочтение.

(2) В любом выражении, где появляется индекс, я предполагаю, что этот индекс может иметь внутреннюю структуру. Например, при рассмотрении алгебры $A$ координаты $a \in A$ относительно базиса $\overline{\overline{e}}$ пронумерованы индексом $i$. Это означает, что $a$ является вектором. Однако, если $a$ является матрицей, нам необходимо два индекса, один нумерует строки, другой - столбцы. В том случае, когда мы уточняем структуру индекса, мы будем начинать индекс с символа · в соответствующей позиции. Например, если я рассматриваю матрицу $a^i_j$ как элемент векторного пространства, то я могу записать элемент матрицы в виде $a^{\cdot i}_{\cdot j}$.

(3) Пусть $A$ - свободная конечно мерная алгебра. При разложении элемента алгебры $A$ относительно базиса $\overline{\overline{e}}$ мы пользуемся одной и той же корневой буквой для обозначения этого элемента и его координат. Однако в алгебре не принято использовать векторные обозначения. В выражении $a^2$ не ясно - это компонента разложения элемента $a$ относительно базиса или это операция возведения в степень. Для облегчения чтения текста мы будем индекс элемента алгебры выделять цветом. Например,
$$a = a^i \overline{e}_i$$

(4) Если свободная конечномерная алгебра имеет единицу, то мы будем отождествлять вектор базиса $\overline{e}_0$ с единицей алгебры.

(5) Если в некотором выражении используется несколько операций, среди которых есть операция ∘, то предполагается, что операция ∘ выполняется первой. Ниже приведен пример эквивалентных выражений.
$$f \circ xy \equiv f(x)y$$
$$f \circ (xy) \equiv f(xy)$$
$$f \circ x + y \equiv f(x) + y$$
$$f \circ (x + y) \equiv f(x + y)$$

(6) Без сомнения, у читателя моих статей могут быть вопросы, замечания, возражения. Я буду признателен любому отзыву.

Глава 2

# Матрица отображений

## 2.1. Произведение отображений

На множестве отображений
$$f: A \to A$$
определено произведение согласно правилу
$$(2.1.1) \qquad f \circ g = f(g)$$
Равенство
$$f \circ g = g \circ f$$
справедливо тогда и только тогда, когда диаграма

$$\begin{array}{ccc} A & \xrightarrow{f} & A \\ \downarrow{g} & & \downarrow{g} \\ A & \xrightarrow{f} & A \end{array}$$

коммутативна.

Для $a \in A$, существует отображение
$$(2.1.2) \qquad f_a(x) = a$$
Если мы будем обозначать отображение $f_a$ буквой $a$, то опираясь на равенство (2.1.1), положим
$$(2.1.3) \qquad f \circ a = f(a)$$

Если $A$ - $\Omega$-алгебра, в которой определена операция умножения, то элемент $a \in A$ так же может служить для обозначения операции левого сдвига
$$(2.1.4) \qquad a \circ b = ab$$

Запись (2.1.4) не противоречит записи (2.1.3). Однако надо помнить, что произведение $(f \circ a) \circ b$, $a, b \in A$, неассоциативно, так как
$$(f \circ a) \circ b = f(a)b \quad f \circ (a \circ b) = f(ab)$$

## 2.2. Бикольцо матриц отображений

Пусть $A$ - $\Omega$-алгебра ([2, 11]), в которой определена операция сложения. Мы будем предполагать, что алгебра $A$ является группой по отношению к операции сложения.

Пусть $\mathcal{A}(A)$ множество отображений $\Omega$-алгебры $A$. Мы можем отобразить операцию сложения в $\Omega$-алгебре $A$ на множество $\mathcal{A}(A)$ согласно правилу

(2.2.1) $$(f+g)\circ a = f\circ a + g\circ a$$

(2.2.2) $$(-f)\circ a = -(f\circ a)$$

(2.2.3) $$0\circ a = 0$$

Равенство (2.2.1) является выражением закона дистрибутивности произведения отображений слева. Поэтому естественно потребовать, чтобы произведение было дистрибутивно справа

(2.2.4) $$f\circ(a+b) = f\circ a + f\circ b$$

В секции 2.3 мы увидим, что это требование существенно. Следовательно, множество $\mathcal{A}(A)$ является множеством гомоморфизмов группы $A$.

Мы требуем, чтобы множество $\mathcal{A}(A)$ было замкнуто относительно операций сложения и произведения отображений.

**Теорема 2.2.1.** *Для любого отображения $g\in\mathcal{A}(A)$ справедливы равенства*

(2.2.5) $$0\circ g = 0$$

(2.2.6) $$(-f)\circ g = -(f\circ g)$$

Доказательство. Из равенства (2.2.3) следует

(2.2.7) $$(0\circ g)\circ a = 0\circ(g\circ a) = 0$$

Равенство (2.2.5) следует из равенства (2.2.7). Из равенства (2.2.1) следует

(2.2.8) $$f\circ g + (-f)\circ g = (f+(-f))\circ g = 0\circ g = 0$$

Равенство (2.2.6) следует из равенства (2.2.8). $\square$

**Замечание 2.2.2.** Если сумма некоммутативна, то требование замкнутости множества $\mathcal{A}(A)$ относительно операции сложения может оказаться слишком строгим. Рассмотрим выражение

$$(f+g)\circ(a+b) = f\circ(a+b) + g\circ(a+b) = f\circ a + f\circ b + g\circ a + g\circ b$$

Так как, вообще говоря,

$$f\circ b + g\circ a \neq g\circ a + f\circ b$$

то мы не можем утверждать, что

$$(f+g)\circ(a+b) = f\circ(a+b) + g\circ(a+b)$$

В дальнейшем тексте мы будем предполагать, что сложение коммутативно. Тем не менее все построения в этой главе мы будем выполнять так, как мы это делали бы в случае некоммутативного сложения. $\square$

Рассмотрим множество **матриц отображений**, элементы которых являются отображениями $f\in\mathcal{A}(A)$. Согласно определению [5]-2.2.1, мы определим $_\circ{}^\circ$**-произведение матриц отображений**

(2.2.9) $$\begin{cases} b_\circ{}^\circ c &= (b_c^a \circ c_b^c) \\ (b_\circ{}^\circ c)_b^a &= b_c^a \circ c_b^c \end{cases}$$

Согласно определению [5]-2.2.2, мы определим $\circ_\circ$-**произведение матриц отображений**

$$(2.2.10) \qquad \begin{cases} b^\circ_\circ c &= (b^c_b \circ c^a_c) \\ (b^\circ_\circ c)^a_b &= b^c_b \circ c^a_c \end{cases}$$

**Теорема 2.2.3.** *Произведение преобразований множества $A$ ассоциативно.*[2.1]

Доказательство. Рассмотрим отображения

$$f : A \to A \quad g : A \to A \quad h : A \to A$$

Утверждение теоремы следует из цепочки равенств

$$((f \circ g) \circ h) \circ x = (f \circ g) \circ (h \circ x) = f \circ (g \circ (h \circ x))$$
$$= f \circ ((g \circ h) \circ x) = (f \circ (g \circ h)) \circ x$$

$\square$

**Теорема 2.2.4.** *Множество $\mathcal{A}(A)$ является кольцом.*

Доказательство. $\mathcal{A}(A)$ - абелева группа относительно операции сложения. Согласно теореме 2.2.3, $\mathcal{A}(A)$ - полугруппа относительно умножения. Так как $f \in \mathcal{A}(A)$ - гомоморфизм абелевой группы $A$, то для любого $a \in A$

$$(2.2.11) \qquad \begin{aligned}(f \circ (h+g)) \circ a &= f \circ ((h+g) \circ a) = f \circ ((h \circ a + g \circ a) \\ &= (f \circ h) \circ a + (f \circ g) \circ a = (f \circ h + f \circ g) \circ a\end{aligned}$$

Из равенства (2.2.11) следует закон дистрибутивности

$$f \circ (h + g) = f \circ h + f \circ g$$

$\square$

Отображения, принадлежащие кольцу $\mathcal{A}(A)$, мы также будем называть $\mathcal{A}(A)$-**отображениями**.

**Теорема 2.2.5.** $\circ^\circ$-*произведение матриц $\mathcal{A}(A)$-отображений является матрицей $\mathcal{A}(A)$-отображений.*

Доказательство. Утверждение теоремы следует из уравнения (2.2.9) и утверждения, что сумма и произведение $\mathcal{A}(A)$-отображений является $\mathcal{A}(A)$-отображением. $\square$

**Теорема 2.2.6.** *Произведение матриц $\mathcal{A}(A)$-отображений ассоциативно.*

Доказательство. Утверждение теоремы следует из цепочки равенств

$$\begin{aligned}(f \circ^\circ g) \circ^\circ h &= \left((f \circ^\circ g)^i_j \circ h^j_k\right) = \left((f^i_m \circ g^m_j) \circ h^j_k\right) \\ &= \left(f^i_m \circ (g^m_j \circ h^j_k)\right) = \left(f^i_m \circ (g \circ^\circ h)^m_k\right) \\ &= f \circ^\circ (g \circ^\circ h)\end{aligned}$$

$\square$

---

[2.1]Утверждение теоремы основано на примере полугруппы из [3], с. 20, 21.

### 2.3. Квазидетерминант матрицы отображений

**Теорема 2.3.1.** *Предположим, что $n \times n$ матрица $\mathcal{A}(A)$-отображений $a$ имеет ${}_\circ{}^\circ$-обратную матрицу*[2.2]

$$(2.3.1) \qquad a_\circ{}^\circ a^{-1}{}_\circ{}^\circ = \delta$$

*Тогда $k \times k$ минор ${}_\circ{}^\circ$-обратной матрицы удовлетворяет следующему равенству, при условии, что рассматриваемые обратные матрицы существуют,*

$$(2.3.2) \qquad \left((a^{-1}{}_\circ{}^\circ)_J^I\right)^{-1}{}_\circ{}^\circ = -a_{[I]}^J{}_\circ{}^\circ \left(a_{[I]}^{[J]}\right)^{-1}{}_\circ{}^\circ{}_\circ{}^\circ a_I^{[J]} + a_I^J$$

Доказательство. Определение (2.3.1) ${}_\circ{}^\circ$-обратной матрицы приводит к системе линейных уравнений

$$(2.3.3) \qquad a_{[I]}^{[J]}{}_\circ{}^\circ(a^{-1}{}_\circ{}^\circ)_J^{[I]} + a_I^{[J]}{}_\circ{}^\circ(a^{-1}{}_\circ{}^\circ)_J^I = 0$$

$$(2.3.4) \qquad a_{[I]}^J{}_\circ{}^\circ(a^{-1}{}_\circ{}^\circ)_J^{[I]} + a_I^J{}_\circ{}^\circ(a^{-1}{}_\circ{}^\circ)_J^I = \delta$$

Мы умножим (2.3.3) на $\left(a_{[I]}^{[J]}\right)^{-1}{}_\circ{}^\circ$

$$(2.3.5) \qquad (a^{-1}{}_\circ{}^\circ)_J^{[I]} + \left(a_{[I]}^{[J]}\right)^{-1}{}_\circ{}^\circ{}_\circ{}^\circ a_I^{[J]}{}_\circ{}^\circ(a^{-1}{}_\circ{}^\circ)_J^I = 0$$

Теперь мы можем подставить (2.3.5) в (2.3.4)

$$(2.3.6) \qquad -a_{[I]}^J{}_\circ{}^\circ \left(a_{[I]}^{[J]}\right)^{-1}{}_\circ{}^\circ{}_\circ{}^\circ a_I^{[J]}{}_\circ{}^\circ(a^{-1}{}_\circ{}^\circ)_J^I + a_I^J{}_\circ{}^\circ(a^{-1}{}_\circ{}^\circ)_J^I = \delta$$

(2.3.2) следует из (2.3.6). $\square$

**Следствие 2.3.2.** *Предположим, что $n \times n$ матрица $\mathcal{A}(A)$-отображений $a$ имеет ${}_\circ{}^\circ$-обратную матрицу. Тогда элементы ${}_\circ{}^\circ$-обратной матрицы удовлетворяют равенству*

$$(2.3.7) \qquad \left(\mathcal{H}a^{-1}{}_\circ{}^\circ\right)_i^j = -a_{[i]}^j{}_\circ{}^\circ \left(a_{[i]}^{[j]}\right)^{-1}{}_\circ{}^\circ{}_\circ{}^\circ a_i^{[j]} + a_i^j$$

$\square$

**Определение 2.3.3.** $\binom{j}{i}$**-${}_\circ{}^\circ$-квазидетерминант** $n \times n$ матрицы $a$ - это формальное выражение

$$(2.3.8) \qquad \det(a, {}_\circ{}^\circ)_i^j = \left(\mathcal{H}a^{-1}{}_\circ{}^\circ\right)_i^j$$

Согласно замечанию [5]-2.1.2 мы можем рассматривать $\binom{b}{a}$-${}_\circ{}^\circ$-квазидетерминант как элемент матрицы $\det(a, {}_\circ{}^\circ)$, которую мы будем называть ${}_\circ{}^\circ$**-квазидетерминантом**. $\square$

**Теорема 2.3.4.** *Выражение для элементов ${}_\circ{}^\circ$-обратной матрицы имеет вид*

$$(2.3.9) \qquad a^{-1}{}_\circ{}^\circ = \mathcal{H} \det(a, {}_\circ{}^\circ)$$

Доказательство. (2.3.9) следует из (2.3.8). $\square$

---

[2.2]Это утверждение и его доказательство основаны на утверждении 1.2.1 из [4] (page 8) для матриц над свободным кольцом с делением.

**Теорема 2.3.5.** *Выражение для $\binom{a}{b}$-${}_\circ{}^\circ$-квазидетерминанта имеет любую из следующих форм*

(2.3.10) $$\det(a,{}_\circ{}^\circ)^j_i = -a^j_{[i]}{}_\circ{}^\circ \left(a^{[j]}_{[i]}\right)^{-1{}_\circ{}^\circ}{}_\circ{}^\circ a^{[b]}_a + a^j_i$$

(2.3.11) $$\det(a,{}_\circ{}^\circ)^j_i = -a^j_{[i]}{}_\circ{}^\circ \mathcal{H}\det\left(a^{[j]}_{[i]},{}_\circ{}^\circ\right){}_\circ{}^\circ a^{[b]}_a + a^j_i$$

Доказательство. Утверждение следует из (2.3.7) и (2.3.8). □

**Определение 2.3.6.** Если для отображения $f \in \mathcal{A}(A)$ из существования обратного отображения $f^{-1}$ следует $f^{-1} \in \mathcal{A}(A)$, то кольцо $\mathcal{A}(A)$ отображений называется **квазизамкнутым**. □

**Теорема 2.3.7.** *Пусть $\mathcal{A}(A)$ - квазизамкнутое кольцо отображений $\Omega$-алгебры $A$. Пусть $a$ - матрица $\mathcal{A}(A)$-отображений. Тогда матрицы $\det(a,{}_\circ{}^\circ)$ и $a^{-1{}_\circ{}^\circ}$ являются матрицами $\mathcal{A}(A)$-отображений.*

Доказательство. Мы докажем теорему индукцией по порядку матрицы. При $n = 1$ из равенства (2.3.10) следует
$$\det(a,{}_\circ{}^\circ)^1_1 = a^1_1$$
Следовательно, квазидетерминант является матрицей $\mathcal{A}(A)$-отображений. Из определения 2.3.6 следует, что матрица $a^{-1{}_\circ{}^\circ}$ является матрицей $\mathcal{A}(A)$-отображений.

Пусть утверждение теоремы верно для $n-1$. Пусть $a$ - $n \times n$ матрица. Согласно предположению индукции, матрица $\left(a^{[b]}_{[a]}\right)^{-1{}_\circ{}^\circ}$ в равенстве (2.3.10) является матрицей $\mathcal{A}(A)$-отображений. Следовательно, $\binom{a}{b}$-${}_\circ{}^\circ$-квазидетерминант является $\mathcal{A}(A)$-отображением. Из определения 2.3.6 и теоремы 2.3.4 следует, что матрица $a^{-1{}_\circ{}^\circ}$ является матрицей $\mathcal{A}(A)$-отображений. □

**Определение 2.3.8.** Если $n \times n$ матрица $a$ $\mathcal{A}(A)$-отображений имеет ${}_\circ{}^\circ$-обратную матрицу, мы будем называть матрицу $a$ **${}_\circ{}^\circ$-невырожденной матрицей $\mathcal{A}(A)$-отображений**. В противном случае, мы будем называть такую матрицу **${}_\circ{}^\circ$-вырожденной матрицей $\mathcal{A}(A)$-отображений**. □

## 2.4. Система аддитивных уравнений

Пусть $\mathcal{A}(A)$ - квазизамкнутое кольцо отображений $\Omega$-алгебры $A$. Система уравнений

(2.4.1) $$\begin{pmatrix} a^1_1 & ... & a^1_n \\ ... & ... & ... \\ a^n_1 & ... & a^n_n \end{pmatrix} {}_\circ{}^\circ \begin{pmatrix} x^1 \\ ... \\ x^n \end{pmatrix} = \begin{pmatrix} b^1 \\ ... \\ b^n \end{pmatrix}$$

где $a$ - матрица $\mathcal{A}(A)$-отображений, называется **системой аддитивных уравнений**. Систему аддитивных уравнений (2.4.1) можно записать также в виде

$$\begin{cases} a^1_1 \circ x^1 & + & ... & + & a^1_n \circ x^n & = & b^1 \\ ... & & ... & & ... & & ... \\ a^n_1 \circ x^1 & + & ... & + & a^n_n \circ x^n & = & b^n \end{cases}$$

**Определение 2.4.1.** Предположим, что $a$ - $_\circ^\circ$-невырожденная матрица. Мы будем называть соответствующую систему аддитивных уравнений (2.4.1) $_\circ^\circ$-**невырожденной системой аддитивных уравнений**. □

**Теорема 2.4.2.** *Решение невырожденной системы $\mathcal{A}(A)$-уравнений (2.4.1) определено однозначно и может быть записано в любой из следующих форм*

(2.4.2) $$x = a^{-1\circ^\circ}{}_\circ^\circ b$$

(2.4.3) $$x = \mathcal{H}\det(a,{}_\circ^\circ){}_\circ^\circ b$$

Доказательство. Умножая обе части равенства (2.4.1) слева на $a^{-1\circ^\circ}$, мы получим (2.4.2). Пользуясь определением (2.3.8), мы получим (2.4.3). □

**Пример 2.4.3.** Согласно определению [5]-4.1.4 эффективное $T\star$-представление тела $D$ в абелевой группе $\overline{V}$ порождает тело отображений $D(*\overline{V})$. Образ $\overline{v} \in \overline{V}$ при отображении $a \in D(*\overline{V})$ определён согласно правилу

$$a \circ \overline{v} = a\overline{v}$$

Произведение отображений $a, b \in D(*\overline{V})$ определён согласно правилу

$$a \circ b = ab$$

Система аддитивных уравнений в этом случае является системой $_*^*D$-линейных уравнений. □

**Пример 2.4.4.** Согласно определению [5]-4.1.4 эффективное $\star T$-представление тела $D$ в абелевой группе $\overline{V}$ порождает тело отображений $D(\overline{V}*)$. Образ $\overline{v} \in \overline{V}$ при отображении $a \in D(\overline{V}*)$ определён согласно правилу

$$a \circ \overline{v} = \overline{v}a$$

Произведение отображений $a, b \in D(\overline{V}*)$ определён согласно правилу

$$a \circ b = ba$$

Система аддитивных уравнений в этом случае является системой $D^*{}_*$-линейных уравнений. □

### 2.5. Система аддитивных уравнений в поле комплексных чисел

Согласно теореме [6]-5.1.9 аддитивное отображение поля комплексных чисел линейно над полем действительных чисел. Рассмотрим базис $e_0 = 1$, $e_1 = i$ поля комплексных чисел над полем действительных чисел. В базисе $\overline{\overline{e}}$ аддитивное отображение $f$ определено матрицей

(2.5.1) $$\begin{pmatrix} f_0^0 & f_1^0 \\ f_0^1 & f_1^1 \end{pmatrix}$$

Согласно теореме [6]-7.1.1 линейное отображение имеет матрицу

$$\begin{pmatrix} a_0 & -a_1 \\ a_1 & a_0 \end{pmatrix}$$

Это отображение соответствует умножению на число $a = a_0 + a_1 i$. Утверждение следует из равенств

$$(a_0 + a_1 i)(x_0 + x_1 i) = a_0 x_0 - a_1 x_1 + (a_0 x_1 + a_1 x_0)i$$

$$\begin{pmatrix} a_0 & -a_1 \\ a_1 & a_0 \end{pmatrix} \begin{pmatrix} x_0 \\ x_1 \end{pmatrix} = \begin{pmatrix} a_0 x_0 - a_1 x_1 \\ a_1 x_0 + a_0 x_1 \end{pmatrix}$$

Аддитивное отображение, порождённое отображением сопряжения

$$I \circ z = \overline{z}$$
$$I = \begin{pmatrix} 1 & 0 \\ 0 & -1 \end{pmatrix}$$

имеют матрицу

$$\begin{pmatrix} b_0 & b_1 \\ b_1 & -b_0 \end{pmatrix}$$

которая соответствует преобразованию $(b_0 + b_1 i) \circ I$. Утверждение следует из равенств

$$(b_0 + b_1 i) I (x_0 + x_1 i) = (b_0 + b_1 i)(x_0 - x_1 i) = b_0 x_0 + b_1 x_1 + (-b_0 x_1 + b_1 x_0) i$$

$$\begin{pmatrix} b_0 & b_1 \\ b_1 & -b_0 \end{pmatrix} \begin{pmatrix} x_0 \\ x_1 \end{pmatrix} = \begin{pmatrix} b_0 x_0 + b_1 x_1 \\ b_1 x_0 - b_0 x_1 \end{pmatrix}$$

**Теорема 2.5.1.** *Аддитивное отображение поля комплексных чисел имеет вид*

(2.5.2) $$f = a + b \circ I$$
(2.5.3) $$(a + b \circ I) \circ z = az + b\overline{z}$$

Доказательство. Пусть отображение $f$ определено матрицей (2.5.1). Сопоставление матриц отображений $f$, $a$, $b$ приводит к матричному равенству

$$\begin{pmatrix} f_0^0 & f_1^0 \\ f_0^1 & f_1^1 \end{pmatrix} = \begin{pmatrix} a_0 & -a_1 \\ a_1 & a_0 \end{pmatrix} + \begin{pmatrix} b_0 & b_1 \\ b_1 & -b_0 \end{pmatrix} = \begin{pmatrix} a_0 + b_0 & -a_1 + b_1 \\ a_1 + b_1 & a_0 - b_0 \end{pmatrix}$$

(2.5.4) $$f_0^0 = a_0 + b_0$$
(2.5.5) $$f_1^0 = -a_1 + b_1$$
(2.5.6) $$f_0^1 = a_1 + b_1$$
(2.5.7) $$f_1^1 = a_0 - b_0$$

Из равенств (2.5.4), (2.5.7) следует

$$a_0 = \frac{f_0^0 + f_1^1}{2} \quad b_0 = \frac{f_0^0 - f_1^1}{2}$$

Из равенств (2.5.6), (2.5.5) следует

$$a_1 = \frac{f_0^1 - f_1^0}{2} \quad b_1 = \frac{f_0^1 + f_1^0}{2}$$

□

В поле комплексных чисел множество аддитивных отображений порождает кольцо, порождённое операцией умножения на комплексное число и операцией сопряжения.

**Теорема 2.5.2.** *Кольцо отображений $\mathcal{A}(C, C)$ квазизамкнуто.*

Доказательство. Аддитивное отображение невырожденно тогда и только тогда, когда невырождена его матрица (2.5.1). Обратная матрица также описывает некоторое отображение. Произведение этих матриц является тождественным преобразованием. $\square$

**Теорема 2.5.3.** *Произведение аддитивных отображений*

$$f = f_0 + f_1 \circ I$$
$$g = g_0 + g_1 \circ I$$

*имеет вид*

$$h = f \circ g = h_0 + h_1 \circ I$$

*где*

$$(2.5.8) \qquad h_0 = f_0 g_0 + f_1 \overline{g}_1 \quad h_1 = f_0 g_1 + f_1 \overline{g}_0$$

Доказательство. Непосредственной проверкой легко убедиться, что

$$(2.5.9) \qquad I \circ I = 1$$

Из цепочки равенств

$$(2.5.10)$$
$$\begin{pmatrix} 1 & 0 \\ 0 & -1 \end{pmatrix} \begin{pmatrix} a_0 & -a_1 \\ a_1 & a_0 \end{pmatrix} \begin{pmatrix} 1 & 0 \\ 0 & -1 \end{pmatrix} = \begin{pmatrix} 1 & 0 \\ 0 & -1 \end{pmatrix} \begin{pmatrix} a_0 & a_1 \\ a_1 & -a_0 \end{pmatrix} = \begin{pmatrix} a_0 & a_1 \\ -a_1 & a_0 \end{pmatrix}$$

следует

$$(2.5.11) \qquad \overline{a} = I \circ a \circ I$$

Из равенств (2.5.9), (2.5.11) следует

$$(2.5.12) \qquad \overline{a} \circ I = I \circ a$$

Из равенств (2.5.9), (2.5.12) следует

$$(2.5.13) \quad \begin{aligned} &(f_0 + f_1 \circ I) \circ (g_0 + g_1 \circ I) \\ =& f_0 \circ (g_0 + g_1 \circ I) + f_1 \circ I \circ (g_0 + g_1 \circ I) \\ =& f_0 \circ g_0 + f_0 \circ g_1 \circ I + f_1 \circ I \circ g_0 + f_1 \circ I \circ g_1 \circ I \\ =& (f_0 g_0) + (f_0 g_1) \circ I + f_1 \circ \overline{g}_0 \circ I + f_1 \circ \overline{g}_1 \circ I \circ I \\ =& (f_0 g_0 + f_1 \overline{g}_1) + (f_0 g_1 + f_1 \overline{g}_0) \circ I \end{aligned}$$

Равенство (2.5.8) следует из равенства (2.5.13). $\square$

**Теорема 2.5.4.** *Пусть аддитивное отображение поля комплексных чисел*

$$g = g_0 + g_1 \circ I$$

*является отображением, обратным аддитивному отображению*

$$f = f_0 + f_1 \circ I$$

*Тогда*

$$(2.5.14) \qquad g_0 = -\frac{\overline{f}_0}{f_1 \overline{f}_1 - f_0 \overline{f}_0} \quad g_1 = \frac{\overline{f}_1}{f_1 \overline{f}_1 - f_0 \overline{f}_0}$$

Доказательство. Согласно утверждению теоремы

(2.5.15) $$f \circ g = 1$$

Из равенств (2.5.8), (2.5.15) следует

(2.5.16) $$f_0 g_0 + f_1 \overline{g}_1 = 1$$

(2.5.17) $$f_1 \overline{g}_0 + f_0 g_1 = 0$$

Из уравнения (2.5.17) следует

$$\overline{f}_1 g_0 + \overline{f}_0 \overline{g}_1 = 0$$

(2.5.18) $$g_0 = -\overline{f}_0 \overline{f}_1^{-1} \overline{g}_1$$

Из уравнений (2.5.16), (2.5.18) следует

(2.5.19) $$- f_0 \overline{f}_0 \overline{f}_1^{-1} \overline{g}_1 + f_1 \overline{g}_1 = 1$$

(2.5.14) следует из равенств (2.5.19), (2.5.18). $\square$

**Пример 2.5.5.** Рассмотрим систему аддитивных уравнений

(2.5.20) $$\begin{cases} z + 2\overline{w} &= 1 \\ z - 3w &= i \end{cases}$$

Систему уравнений (2.5.20) нельзя решать пользуясь определителем и правилом Крамера. Запишем систему уравнений (2.5.20) в виде

(2.5.21) $$\begin{cases} z + 2 \circ I \circ w &= 1 \\ z + (-3) \circ w &= i \end{cases}$$

(2.5.22) $$\begin{pmatrix} 1 & 2 \circ I \\ 1 & -3 \end{pmatrix} \circ^{\circ} \begin{pmatrix} z \\ w \end{pmatrix} = \begin{pmatrix} 1 \\ i \end{pmatrix}$$

Найдём $\circ^{\circ}$-квазидетерминант матрицы

$$a = \begin{pmatrix} 1 & 2 \circ I \\ 1 & -3 \end{pmatrix}$$

Согласно равенству (2.3.10) мы имеем

$$\begin{aligned}
\det(a, {}_\circ{}^\circ)_1^1 &= a_1^1 - a_{[1]}^1 {}_\circ{}^\circ (a_{[1]}^{[1]})^{-1} {}_\circ{}^\circ a_1^{[1]} \\
&= a_1^1 - a_2^1 \circ (a_2^2)^{-1} \circ a_1^2 \\
&= 1 - 2 \circ I \circ (-3)^{-1} \circ 1 \\
&= 1 + \frac{2}{3} \circ I \\
\det(a, {}_\circ{}^\circ)_1^2 &= a_1^2 - a_{[1]}^2 {}_\circ{}^\circ (a_{[1]}^{[2]})^{-1} {}_\circ{}^\circ a_1^{[2]} \\
&= a_1^2 - a_2^2 \circ (a_2^1)^{-1} \circ a_1^1 \\
&= 1 - (-3) \circ (2 \circ I)^{-1} \circ 1 \\
&= 1 + \frac{3}{2} \circ I \\
\det(a, {}_\circ{}^\circ)_2^1 &= a_2^1 - a_{[2]}^1 {}_\circ{}^\circ (a_{[2]}^{[1]})^{-1} {}_\circ{}^\circ a_2^{[1]} \\
&= a_2^1 - a_1^1 \circ (a_1^2)^{-1} \circ a_2^2 \\
&= 2 \circ I - 1 \circ (1)^{-1} \circ (-3) \\
&= 3 + 2 \circ I \\
\det(a, {}_\circ{}^\circ)_2^2 &= a_2^2 - a_{[2]}^2 {}_\circ{}^\circ (a_{[2]}^{[2]})^{-1} {}_\circ{}^\circ a_2^{[2]} \\
&= a_2^2 - a_1^2 \circ (a_1^1)^{-1} \circ a_2^1 \\
&= (-3) - 1 \circ (1)^{-1} \circ 2 \circ I \\
&= -3 - 2 \circ I
\end{aligned}$$

Согласно теоремам 2.3.4, 2.5.4, $\circ^\circ$-обратная матрица имеет вид

$$(a^{-1_{\circ^\circ}})^1_1 = (\det(a,{}_\circ{}^\circ)^1_1)^{-1} = (1 + \frac{2}{3} \circ I)^{-1}$$
$$= \frac{-1 + \frac{2}{3}I}{\frac{2}{3}\frac{2}{3} - 1\,1} = \frac{9}{5}(1 - \frac{2}{3}I)$$
$$= \frac{9}{5} - \frac{6}{5} \circ I$$

$$(a^{-1_{\circ^\circ}})^1_2 = (\det(a,{}_\circ{}^\circ)^2_1)^{-1} = (1 + \frac{3}{2} \circ I)^{-1}$$
$$= \frac{-1 + \frac{3}{2}I}{\frac{3}{2}\frac{3}{2} - 1\,1} = \frac{4}{5}(-1 + \frac{3}{2}I)$$
$$= -\frac{4}{5} + \frac{6}{5} \circ I$$

$$(a^{-1_{\circ^\circ}})^2_1 = (\det(a,{}_\circ{}^\circ)^1_2)^{-1} = (3 + 2 \circ I)^{-1}$$
$$= \frac{-3 + 2I}{2\,2 - 3\,3} = \frac{1}{5}(3 - 2I)$$
$$= \frac{3}{5} - \frac{2}{5} \circ I$$

$$(a^{-1_{\circ^\circ}})^2_2 = (\det(a,{}_\circ{}^\circ)^2_2)^{-1} = (-3 - 2 \circ I)^{-1}$$
$$= \frac{3 - 2I}{2\,2 - 3\,3} = \frac{1}{5}(-3 + 2I)$$
$$= -\frac{3}{5} + \frac{2}{5} \circ I$$

$$a^{-1_{\circ^\circ}} = \begin{pmatrix} \frac{9}{5} - \frac{6}{5} \circ I & -\frac{4}{5} + \frac{6}{5} \circ I \\ \frac{3}{5} - \frac{2}{5} \circ I & -\frac{3}{5} + \frac{2}{5} \circ I \end{pmatrix}$$

Согласно теореме 2.4.2 решение системы аддитивных уравнений (2.5.20) имеет вид

$$\begin{pmatrix} z \\ w \end{pmatrix} = \begin{pmatrix} \frac{9}{5} - \frac{6}{5} \circ I & -\frac{4}{5} + \frac{6}{5} \circ I \\ \frac{3}{5} - \frac{2}{5} \circ I & -\frac{3}{5} + \frac{2}{5} \circ I \end{pmatrix} {}_\circ{}^\circ \begin{pmatrix} 1 \\ i \end{pmatrix}$$
$$= \begin{pmatrix} \frac{9}{5} - \frac{6}{5} - \left(\frac{4}{5} + \frac{6}{5}\right)i \\ \frac{3}{5} - \frac{2}{5} - \left(\frac{3}{5} + \frac{2}{5}\right)i \end{pmatrix}$$
$$= \begin{pmatrix} \frac{3}{5} - 2i \\ \frac{1}{5} - i \end{pmatrix}$$

Непосредственная проверка показывает, что мы нашли решение системы уравнений (2.5.20)

$$\left(\frac{3}{5} - 2i\right) + 2\overline{\left(\frac{1}{5} - i\right)} = \frac{3}{5} + 2\frac{1}{5} + (-2 + 2)i = 1$$

$$\left(\frac{3}{5} - 2i\right) - 3\left(\frac{1}{5} - i\right) = \frac{3}{5} - 3\frac{1}{5} + (-2 + 3)i = i$$

□

Глава 3

# Линейное отображение алгебры

### 3.1. Модуль

**Теорема 3.1.1.** *Пусть кольцо $D$ имеет единицу $e$. Представление*

(3.1.1) $$f : D \to {}^\star A$$

*кольца $D$ в абелевой группе $A$ **эффективно** тогда и только тогда, когда из равенства $f(a) = 0$ следует $a = 0$.*

Доказательство. Сумма преобразований $f$ и $g$ абелевой группы определяется согласно правилу

$$(f + g) \circ a = f \circ a + g \circ a$$

Поэтому, рассматривая представление кольца $D$ в абелевой группе $A$, мы полагаем

$$f(a + b) \circ x = f(a) \circ x + f(b) \circ x$$

Произведение преобразований представления определено согласно правилу

$$f(ab) = f(a) \circ f(b)$$

Если $a$, $b \in R$ порождают одно и то же преобразование, то

(3.1.2) $$f(a) \circ m = f(b) \circ m$$

для любого $m \in A$. Из равенства (3.1.2) следует, что $a - b$ порождает нулевое преобразование

$$f(a - b) \circ m = 0$$

Элемент $e + a - b$ порождает тождественное преобразование. Следовательно, представление $f$ эффективно тогда и только тогда, когда $a = b$. □

**Определение 3.1.2.** Пусть $D$ - коммутативное кольцо. $A$ - **модуль над кольцом** $D$, если $A$ - абелева группа и определено эффективное представление кольца $D$ в абелевой группе $A$. □

**Определение 3.1.3.** $A$ - **свободный модуль над кольцом** $D$, если $A$ имеет базис над кольцом $D$.[3.1] □

Следующее определение является следствием определений 3.1.2 и [7]-2.2.2.

**Определение 3.1.4.** Пусть $A_1$ - модуль над кольцом $R_1$. Пусть $A_2$ - модуль над кольцом $R_2$. Морфизм

$$(f : R_1 \to R_2, g : A_1 \to A_2)$$

---

[3.1]Я следую определению в [1], с. 103.

представления кольца $R_1$ в абелевой группе $A_1$ в представление кольца $R_2$ в абелевой группе $A_2$ называется **линейным отображением** $R_1$**-модуля** $A_1$ **в** $R_2$**-модуль** $A_2$. □

**Теорема 3.1.5.** *Линейное отображение*

$$(f:R_1\to R_2, g:A_1\to A_2)$$

$R_1$*-модуля* $A_1$ *в* $R_2$*-модуль* $A_2$ *удовлетворяет равенствам*[3.2]

(3.1.3) $$g\circ(a+b)=g\circ a+g\circ b$$

(3.1.4) $$g\circ(pa)=(f\circ p)(g\circ a)$$

(3.1.5) $$f\circ(pq)=(f\circ p)(f\circ q)$$

$$a,b\in A_1 \quad p,q\in R_1$$

Доказательство. Из определений 3.1.4 и [7]-2.2.2 следует, что

- отображение $f$ является гомоморфизмом кольца $R_1$ в кольцо $R_2$ (равенство (3.1.5))
- отображение $g$ является гомоморфизмом абелевой группы $A_1$ в абелеву группу $A_2$ (равенство (3.1.3))

Равенство (3.1.4) следует из равенства [7]-(2.2.3). □

Согласно теореме [7]-2.2.18 при изучении линейных отображений, не нарушая общности, мы можем полагать $R_1=R_2$.

**Определение 3.1.6.** Пусть $A_1$ и $A_2$ - модули над кольцом $D$. Морфизм

$$g:A_1\to A_2$$

представления кольца $D$ в абелевой группе $A_1$ в представление кольца $D$ в абелевой группе $A_2$ называется **линейным отображением** $D$**-модуля** $A_1$ **в** $D$**-модуль** $A_2$. □

**Теорема 3.1.7.** *Линейное отображение*

$$g:A_1\to A_2$$

$D$*-модуля* $A_1$ *в* $D$*-модуль* $A_2$ *удовлетворяет равенствам*[3.3]

(3.1.6) $$g\circ(a+b)=g\circ a+g\circ b$$

(3.1.7) $$g\circ(pa)=p(g\circ a)$$

---

[3.2]Предлагаемые равенства в классической записи имеют вполне знакомый вид

$$g(a+b)=g(a)+g(b)$$
$$g(pa)=f(p)g(a)$$
$$f(pq)=f(p)f(q)$$
$$a,b\in A_1 \quad p,q\in R_1$$

[3.3]В классической записи приведенные равенства имеют вид

$$g(a+b)=g(a)+g(b)$$
$$g(pa)=pg(a)$$
$$a,b\in A_1 \quad p\in D$$

$$a, b \in A_1 \quad p \in D$$

Доказательство. Из определения 3.1.6 и теоремы [7]-2.2.18 следует, что отображение $g$ является гомоморфизмом абелевой группы $A_1$ в абелеву группу $A_2$ (равенство (3.1.6)) Равенство (3.1.7) следует из равенства [7]-(2.2.44). □

### 3.2. Алгебра над кольцом

**Определение 3.2.1.** Пусть $D$ - коммутативное кольцо. Пусть $A$ - модуль над кольцом $D$.[3.4] Для заданного билинейного отображения

$$f : A \times A \to A$$

мы определим произведение в $A$

(3.2.1) $$ab = f \circ (a, b)$$

$A$ - **алгебра над кольцом** $D$, если $A$ - $D$-модуль и в $A$ определена операция произведения (3.2.1). Алгебра $A^*$ называется **алгеброй, противоположной алгебре** $A$, если в модуле $A$ определено произведение согласно правилу[3.5]

$$ba = f \circ (a, b)$$

Если $A$ является свободным $D$-модулем, то $A$ называется **свободной алгеброй над кольцом** $D$. □

**Замечание 3.2.2.** Алгебра $A$ и противоположная ей алгебра совпадают как модули. □

**Теорема 3.2.3.** *Произведение в алгебре $A$ дистрибутивно по отношению к сложению.*

Доказательство. Утверждение теоремы следует из цепочки равенств

$$(a + b)c = f \circ (a + b, c) = f \circ (a, c) + f \circ (b, c) = ac + bc$$
$$a(b + c) = f \circ (a, b + c) = f \circ (a, b) + f \circ (a, c) = ab + ac$$

□

Произведение в алгебре может быть ни коммутативным, ни ассоциативным. Следующие определения основаны на определениях, данных в [14], с. 13.

**Определение 3.2.4. Коммутатор**

$$[a, b] = ab - ba$$

служит мерой коммутативности в $D$-алгебре $A$. $D$-алгебра $A$ называется **коммутативной**, если

$$[a, b] = 0$$

□

---

[3.4]Существует несколько эквивалентных определений алгебры. Вначале я хотел рассмотреть представление кольца $D$ в абелевой группе кольца $A$. Однако мне надо было объяснить почему произведение элементов кольца $D$ и алгебры $A$ коммутативно. Это потребовало определение центра алгебры $A$. После тщательного анализа я выбрал определение, приведенное в [14], с. 1, [9], с. 4.

[3.5]Определение дано по аналогии с определением [12]-2, с. 19.

**Определение 3.2.5. Ассоциатор**

$$(3.2.2) \qquad (a,b,c) = (ab)c - a(bc)$$

служит мерой ассоциативности в $D$-алгебре $A$. $D$-алгебра $A$ называется **ассоциативной**, если

$$(a,b,c) = 0$$

□

**Теорема 3.2.6.** *Пусть $A$ - алгебра над коммутативным кольцом $D$.*[3.6]

$$(3.2.3) \qquad a(b,c,d) + (a,b,c)d = (ab,c,d) - (a,bc,d) + (a,b,cd)$$

*для любых $a$, $b$, $c$, $d \in A$.*

Доказательство. Равенство (3.2.3) следует из цепочки равенств

$$\begin{aligned}
a(b,c,d) + (a,b,c)d &= a((bc)d - b(cd)) + ((ab)c - a(bc))d \\
&= a((bc)d) - a(b(cd)) + ((ab)c)d - (a(bc))d \\
&= ((ab)c)d - (ab)(cd) + (ab)(cd) \\
&\quad + a((bc)d) - a(b(cd)) - (a(bc))d \\
&= (ab,c,d) - (a(bc))d + a((bc)d) + (ab)(cd) - a(b(cd)) \\
&= (ab,c,d) - (a,(bc),d) + (a,b,cd)
\end{aligned}$$

□

**Определение 3.2.7. Ядро $D$-алгебры** $A$ - это множество[3.7]

$$N(A) = \{a \in A : \forall b, c \in A, (a,b,c) = (b,a,c) = (b,c,a) = 0\}$$

□

**Определение 3.2.8. Центр $D$-алгебры** $A$ - это множество[3.8]

$$Z(A) = \{a \in A : a \in N(A), \forall b \in A, ab = ba\}$$

□

**Теорема 3.2.9.** *Пусть $D$ - коммутативное кольцо. Если $D$-алгебра $A$ имеет единицу, то существует изоморфизм $f$ кольца $D$ в центр алгебры $A$.*

Доказательство. Пусть $e \in A$ - единица алгебры $A$. Положим $f \circ a = ae$. □

Пусть $\overline{\overline{e}}$ - базис свободной алгебры $A$ над кольцом $D$. Если алгебра $A$ имеет единицу, положим $\overline{e}_0$ - единица алгебры $A$.

**Теорема 3.2.10.** *Пусть $\overline{\overline{e}}$ - базис свободной алгебры $A$ над кольцом $D$. Пусть*

$$a = a^i e_i \quad b = b^i e_i \quad a, b \in A$$

*Произведение $a$, $b$ можно получить согласно правилу*

$$(3.2.4) \qquad (ab)^k = B^k_{ij} a^i b^j$$

---

[3.6] Утверждение теоремы опирается на равенство [14]-(2.4).

[3.7] Определение дано на базе аналогичного определения в [14], с. 13

[3.8] Определение дано на базе аналогичного определения в [14], с. 14

*где* $B_{ij}^k$ - **структурные константы алгебры** $A$ **над кольцом** $D$. *Произведение базисных векторов в алгебре $A$ определено согласно правилу*

$$\overline{e}_i \overline{e}_j = B_{ij}^k \overline{e}_k \tag{3.2.5}$$

Доказательство. Равенство (3.2.5) является следствием утверждения, что $\overline{\overline{e}}$ является базисом алгебры $A$. Так как произведение в алгебре является билинейным отображением, то произведение $a$ и $b$ можно записать в виде

$$ab = a^i b^j e_i e_j \tag{3.2.6}$$

Из равенств (3.2.5), (3.2.6), следует

$$ab = a^i b^j B_{ij}^k \overline{e}_k \tag{3.2.7}$$

Так как $\overline{\overline{e}}$ является базисом алгебры $A$, то равенство (3.2.4) следует из равенства (3.2.7). □

**Теорема 3.2.11.** *Если алгебра $A$ коммутативна, то*

$$B_{ij}^p = B_{ji}^p \tag{3.2.8}$$

*Если алгебра $A$ ассоциативна, то*

$$B_{ij}^p B_{pk}^q = B_{ip}^q B_{jk}^p \tag{3.2.9}$$

Доказательство. Для коммутативной алгебры, равенство (3.2.8) следует из равенства

$$\overline{e}_i \overline{e}_j = \overline{e}_j \overline{e}_i$$

Для ассоциативной алгебры, равенство (3.2.9) следует из равенства

$$(\overline{e}_i \overline{e}_j)\overline{e}_k = \overline{e}_i(\overline{e}_j \overline{e}_k)$$

□

### 3.3. Линейное отображение алгебры

Алгебра является кольцом. Отображение, сохраняющее структуру алгебры как кольца, называется гомоморфизмом алгебры. Однако для нас важнее утверждение, что алгебра является модулем над коммутативным кольцом. Отображение, сохраняющее структуру алгебры как модуля, называется линейным отображением алгебры. Таким образом, следующее определение опирается на определение 3.1.6.

**Определение 3.3.1.** Пусть $A_1$ и $A_2$ - алгебры над кольцом $D$. Морфизм

$$g : A_1 \to A_2$$

представления кольца $D$ в абелевой группе $A_1$ в представление кольца $D$ в абелевой группе $A_2$ называется **линейным отображением** $D$**-алгебры** $A_1$ **в** $D$**-алгебру** $A_2$. Обозначим $\mathcal{L}(A_1; A_2)$ множество линейных отображений алгебры $A_1$ в алгебру $A_2$. □

**Теорема 3.3.2.** *Линейное отображение*

$$g : A_1 \to A_2$$

$D$-алгебры $A_1$ в $D$-алгебру $A_2$ удовлетворяет равенствам

(3.3.1)
$$\begin{cases} g \circ (a+b) = g \circ a + g \circ b \\ g \circ (pa) = pg \circ a \\ a, b \in A_1 \quad p \in D \end{cases}$$

Доказательство. Следствие теоремы 3.1.7. $\square$

**Теорема 3.3.3.** *Рассмотрим алгебру $A_1$ и алгебру $A_2$. Пусть отображения*
$$f : A_1 \to A_2$$
$$g : A_1 \to A_2$$
*являются линейными отображениями. Тогда отображение $f + g$, определённое равенством*
$$(f+g) \circ a = f \circ a + g \circ a$$
*также является линейным.*

Доказательство. Утверждение теоремы следует из цепочек равенств
$$(f+g) \circ (x+y) = f \circ (x+y) + g \circ (x+y) = f \circ x + f \circ y + g \circ x + g \circ y$$
$$= (f+g) \circ x + (f+g) \circ y$$
$$(f+g) \circ (px) = f \circ (px) + g \circ (px) = pf \circ x + pg \circ x$$
$$= p(f+g) \circ x$$
$\square$

**Теорема 3.3.4.** *Рассмотрим алгебру $A_1$ и алгебру $A_2$. Пусть отображение*
$$g : A_1 \to A_2$$
*является линейным отображением. Тогда отображения $ag$, $gb$, $a, b \in A_2$, определённые равенствами*
$$(ag) \circ x = a \; g \circ x$$
$$(gb) \circ x = g \circ x \; b$$
*также являются линейными.*

Доказательство. Утверждение теоремы следует из цепочек равенств
$$(ag) \circ (x+y) = a \; g \circ (x+y) = a \; (g \circ x + g \circ y) = a \; g \circ x + a \; g \circ y$$
$$= (ag) \circ x + (ag) \circ y$$
$$(ag) \circ (px) = a \; g \circ (px) = ap \; g \circ x = pa \; g \circ x$$
$$= p \; (ag) \circ x$$
$$(gb) \circ (x+y) = g \circ (x+y) \; b = (g \circ x + g \circ y) \; b = g \circ x \; b + g \circ y \; b$$
$$= (gb) \circ x + (gb) \circ y$$
$$(gb) \circ (px) = g \circ (px) \; b = p \; g \circ x \; b$$
$$= p \; (gb) \circ x$$
$\square$

**Теорема 3.3.5.** *Рассмотрим алгебру $A_1$ и алгебру $A_2$. Пусть отображение*
$$g : A_1 \to A_2$$
*является линейным отображением. Тогда отображения $pg$, $p \in D$, определённое равенством*
$$(pg) \circ x = p \ g \circ x$$
*также являются линейными. При этом выполняется равенство*
$$p(qg) = (pq)g$$
$$(p+q)g = pg + qg$$

Доказательство. Утверждение теоремы следует из цепочек равенств
$$\begin{aligned}(pg) \circ (x+y) =& p \ g \circ (x+y) = p \ (g \circ x + g \circ y) = p \ g \circ x + p \ g \circ y \\ =& (pg) \circ x + (pg) \circ y \\ (pg) \circ (qx) =& p \ g \circ (qx) = pq \ g \circ x = qp \ g \circ x \\ =& q \ (pg) \circ x \\ (p(qg)) \circ x =& p \ (qg) \circ x = p \ (q \ g \circ x) = (pq) \ g \circ x = ((pq)g) \circ x \\ ((p+q)g) \circ x =& (p+q) \ g \circ x = p \ g \circ x + q \ g \circ x = (pg) \circ x + (qg) \circ x \end{aligned}$$

□

**Теорема 3.3.6.** *Пусть $D$ - коммутативное кольцо с единицей. Рассмотрим $D$-алгебру $A_1$ и $D$-алгебру $A_2$. Множество $\mathcal{L}(A_1; A_2)$ является $D$-модулем.*

Доказательство. Теорема 3.3.3 определяет сумму линейных отображений из $D$-алгебры $A_1$ в $D$-алгебру $A_2$. Пусть $f$, $g$, $h \in \mathcal{L}(A_1; A_2)$. Для любого $a \in A_1$
$$\begin{aligned}(f+g) \circ a =& f \circ a + g \circ a = g \circ a + f \circ a \\ =& (g+f) \circ a \\ ((f+g)+h) \circ a =& (f+g) \circ a + h \circ a = (f \circ a + g \circ a) + h \circ a \\ =& f \circ a + (g \circ a + h \circ a) = f \circ a + (g+h) \circ a \\ =& (f+(g+h)) \circ a \end{aligned}$$
Следовательно, сумма линейных отображений коммутативна и ассоциативна.

Отображение $z$, определённое равенством
$$z \circ x = 0$$
является нулём операции сложения, так как
$$(z+f) \circ a = z \circ a + f \circ a = 0 + f \circ a = f \circ a$$
Для заданного отображения $f$ отображение $g$, определённое равенством
$$g \circ a = -f \circ a$$
удовлетворяет равенству
$$f + g = z$$
так как
$$(f+g) \circ a = f \circ a + g \circ a = f \circ a - f \circ a = 0$$
Следовательно, множество $\mathcal{L}(A_1; A_2)$ является абелевой группой.

Из теоремы 3.3.5 следует, что определенно представление кольца $D$ в абелевой группе $\mathcal{L}(A_1;A_2)$. Так как кольцо $D$ имеет единицу, то согласно теореме 3.1.1 указанное представление эффективно. $\square$

### 3.4. Алгебра $\mathcal{L}(A;A)$

**Теорема 3.4.1.** *Пусть $A$, $B$, $C$ - алгебры над коммутативным кольцом $D$. Пусть $f$ - линейное отображение из алгебры $A$ в алгебру $B$. Пусть $g$ - линейное отображение из алгебры $B$ в алгебру $C$. Отображение $g \circ f$, определённое диаграммой*

$$\begin{array}{ccc} & B & \\ {}^{f}\nearrow & & \searrow^{g} \\ A & \xrightarrow{g \circ f} & C \end{array}$$

*является линейным отображением из алгебры $A$ в алгебру $C$.*

Доказательство. Доказательство теоремы следует из цепочек равенств

$$(g \circ f) \circ (a+b) = g \circ (f \circ (a+b)) = g \circ (f \circ a + f \circ b)$$
$$= g \circ (f \circ a) + g \circ (f \circ b) = (g \circ f) \circ a + (g \circ f) \circ b$$
$$(g \circ f) \circ (pa) = g \circ (f \circ (pa)) = g \circ (p\ f \circ a) = p\ g \circ (f \circ a)$$
$$= p\ (g \circ f) \circ a$$

$\square$

**Теорема 3.4.2.** *Пусть $A$, $B$, $C$ - алгебры над коммутативным кольцом $D$. Пусть $f$ - линейное отображение из алгебры $A$ в алгебру $B$. Отображение $f$ порождает линейное отображение*

$$f^* : g \in \mathcal{L}(B;C) \to g \circ f \in \mathcal{L}(A;C)$$

Доказательство. Доказательство теоремы следует из цепочек равенств

$$((g_1 + g_2) \circ f) \circ a = (g_1 + g_2) \circ (f \circ a) = g_1 \circ (f \circ a) + g_2 \circ (f \circ a)$$
$$= (g_1 \circ f) \circ a + (g_2 \circ f) \circ a$$
$$= (g_1 \circ f + g_2 \circ f) \circ a$$
$$((pg) \circ f) \circ a = (pg) \circ (f \circ a) = p\ g \circ (f \circ a) = p\ (g \circ f) \circ a$$
$$= (p(g \circ f)) \circ a$$

$\square$

**Теорема 3.4.3.** *Пусть $A$, $B$, $C$ - алгебры над коммутативным кольцом $D$. Пусть $g$ - линейное отображение из алгебры $B$ в алгебру $C$. Отображение $g$ порождает линейное отображение*

$$g^* : f \in \mathcal{L}(A;B) \to g \circ f \in \mathcal{L}(A;C)$$

Доказательство. Доказательство теоремы следует из цепочек равенств

$$(g \circ (f_1 + f_2)) \circ a = g \circ ((f_1 + f_2) \circ a) = g \circ (f_1 \circ a + f_2 \circ a)$$
$$= g \circ (f_1 \circ a) + g \circ (f_2 \circ a) = (g \circ f_1) \circ a + (g \circ f_2) \circ a$$
$$= (g \circ f_1 + g \circ f_2) \circ a$$
$$(g \circ (pf)) \circ a = g \circ ((pf) \circ a) = g \circ (p\,(f \circ a)) = p\,g \circ (f \circ a)$$
$$= p\,(g \circ f) \circ a = (p(g \circ f)) \circ a$$

$\square$

**Теорема 3.4.4.** *Пусть $A$, $B$, $C$ - алгебры над коммутативным кольцом $D$. Отображение*

$$\circ : (g, f) \in \mathcal{L}(B;C) \times \mathcal{L}(A;B) \to g \circ f \in \mathcal{L}(A;C)$$

*является билинейным отображением.*

Доказательство. Теорема является следствием теорем 3.4.2, 3.4.3.  $\square$

**Теорема 3.4.5.** *Пусть $A$ - алгебра над коммутативным кольцом $D$. Модуль $\mathcal{L}(A;A)$, оснащённый произведением*

(3.4.1) $\quad \circ : (g, f) \in \mathcal{L}(A;A) \times \mathcal{L}(A;A) \to g \circ f \in \mathcal{L}(A;A)$

*является алгеброй над $D$.*

Доказательство. Теорема является следствием определения 3.2.1 и теоремы 3.4.4.  $\square$

### 3.5. Тензорное произведение алгебр

**Определение 3.5.1.** Пусть $D$ - коммутативное кольцо. Пусть $A_1$, ..., $A_n$, $S$ - $D$-модули. Мы будем называть отображение

$$f : A_1 \times ... \times A_n \to S$$

**полилинейным отображением модулей** $A_1$, ..., $A_n$ в модуль $S$, если

$$f \circ (a_1, ..., a_i + b_i, ..., a_n) = f \circ (a_1, ..., a_i, ..., a_n) + f \circ (a_1, ..., b_i, ..., a_n)$$
$$f \circ (a_1, ..., pa_i, ..., a_n) = pf \circ (a_1, ..., a_i, ..., a_n)$$
$$1 \le i \le n \quad a_i, b_i \in A_i \quad p \in D$$

$\square$

**Определение 3.5.2.** Пусть $D$ - коммутативное ассоциативное кольцо. Пусть $A_1$, ..., $A_n$ - $D$-алгебры и $S$ - $D$-модуль. Мы будем называть отображение

$$f : A_1 \times ... \times A_n \to S$$

**полилинейным отображением алгебр** $A_1$, ..., $A_n$ в модуль $S$, если

$$f \circ (a_1, ..., a_i + b_i, ..., a_n) = f \circ (a_1, ..., a_i, ..., a_n) + f \circ (a_1, ..., b_i, ..., a_n)$$
$$f \circ (a_1, ..., pa_i, ..., a_n) = pf \circ (a_1, ..., a_i, ..., a_n)$$
$$1 \le i \le n \quad a_i, b_i \in A_i \quad p \in D$$

Обозначим $\mathcal{L}(A_1, ..., A_n; S)$ множество полилинейных отображений алгебр $A_1$, ..., $A_n$ в модуль $S$.  $\square$

**Определение 3.5.3.** Пусть $A_1, ..., A_n$ - свободные алгебры над коммутативным кольцом $D$.[3.9] Рассмотрим категорию $\mathcal{A}$ объектами которой являются полилинейные над коммутативным кольцом $D$ отображения

$$f : A_1 \times ... \times A_n \longrightarrow S_1 \qquad g : A_1 \times ... \times A_n \longrightarrow S_2$$

где $S_1$, $S_2$ - модули над кольцом $D$. Мы определим морфизм $f \to g$ как линейное над коммутативным кольцом $D$ отображение $h : S_1 \to S_2$, для которого коммутативна диаграмма

$$\begin{array}{c} & & S_1 \\ & \nearrow^{f} & \\ A_1 \times ... \times A_n & & \downarrow h \\ & \searrow_{g} & \\ & & S_2 \end{array}$$

Универсальный объект $A_1 \otimes ... \otimes A_n$ категории $\mathcal{A}$ называется **тензорным произведением алгебр** $A_1, ..., A_n$. □

**Теорема 3.5.4.** *Тензорное произведение алгебр существует.*

Доказательство. Пусть $M$ - модуль над кольцом $D$, порождённый произведением $A_1 \times ... \times A_n$ алгебр $A_1, ..., A_n$. Инъекция

$$i : A_1 \times ... \times A_n \longrightarrow M$$

определена по правилу

(3.5.1) $\qquad i \circ (d_1, ..., d_n) = (d_1, ..., d_n)$

Пусть $N \subset M$ - подмодуль, порождённый элементами вида

(3.5.2) $\qquad (d_1, ..., d_i + c_i, ..., d_n) - (d_1, ..., d_i, ..., d_n) - (d_1, ..., c_i, ..., d_n)$

(3.5.3) $\qquad (d_1, ..., ad_i, ..., d_n) - a(d_1, ..., d_i, ..., d_n)$

где $d_i \in A_i$, $c_i \in A_i$, $a \in D$. Пусть

$$j : M \to M/N$$

каноническое отображение на фактормодуль. Рассмотрим коммутативную диаграмму

(3.5.4) $$\begin{array}{ccc} & & M/N \\ & \nearrow^{f} & \uparrow j \\ A_1 \times ... \times A_n & \xrightarrow{i} & M \end{array}$$

Поскольку элементы (3.5.2) и (3.5.3) принадлежат ядру линейного отображения $j$, то из равенства (3.5.1) следует

(3.5.5) $\quad f \circ (d_1, ..., d_i + c_i, ..., d_n) = f \circ (d_1, ..., d_i, ..., d_n) + f \circ (d_1, ..., c_i, ..., d_n)$

(3.5.6) $\qquad f \circ (d_1, ..., ad_i, ..., d_n) = a \, f \circ (d_1, ..., d_i, ..., d_n)$

---

[3.9]Я определяю тензорное произведение алгебр по аналогии с определением в [1], с. 456 - 458.

Из равенств (3.5.5) и (3.5.6) следует, что отображение $f$ полилинейно над кольцом $D$. Поскольку $M$ - модуль с базисом $A_1 \times ... \times A_n$, то, согласно теореме [1]-1 на с. 104, для любого модуля $V$ и любого полилинейного над $D$ отображения

$$g : A_1 \times ... \times A_n \longrightarrow V$$

существует единственный гомоморфизм $k : M \to V$, для которого коммутативна следующая диаграмма

(3.5.7)
$$\begin{array}{ccc} A_1 \times ... \times A_n & \xrightarrow{i} & M \\ & \searrow{g} & \downarrow{k} \\ & & V \end{array}$$

Так как $g$ - полилинейно над $D$, то $\ker k \subseteq N$. Согласно утверждению на с. [1]-94, отображение $j$ универсально в категории гомоморфизмов векторного пространства $M$, ядро которых содержит $N$. Следовательно, определён гомоморфизм

$$h : M/N \to V$$

для которого коммутативна диаграмма

(3.5.8)
$$\begin{array}{ccc} & & M/N \\ & \nearrow{j} & \downarrow{h} \\ M & \xrightarrow{k} & V \end{array}$$

Объединяя диаграммы (3.5.4), (3.5.7), (3.5.8), получим коммутативную диаграмму

(3.5.9)
$$\begin{array}{ccc} & & M/N \\ & \nearrow{f} & \downarrow{h} \\ A_1 \times ... \times A_n & \xrightarrow{i} M \xrightarrow{k} & V \\ & \searrow{g} & \end{array}$$

Так как $\operatorname{Im} f$ порождает $M/N$, то отображение $h$ однозначно определено. □

Согласно доказательству теоремы 3.5.4

$$A_1 \otimes ... \otimes A_n = M/N$$

Для $d_i \in A_i$ будем записывать

(3.5.10) $$j \circ (d_1, ..., d_n) = d_1 \otimes ... \otimes d_n$$

**Теорема 3.5.5.** *Пусть $A_1$, ..., $A_n$ - алгебры над коммутативным кольцом $D$. Пусть*

$$f : A_1 \times ... \times A_n \to A_1 \otimes ... \otimes A_n$$

*полилинейное отображение, определённое равенством*

(3.5.11) $$f \circ (d_1, ..., d_n) = d_1 \otimes ... \otimes d_n$$

*Пусть*
$$g : A_1 \times ... \times A_n \to V$$
*полилинейное отображение в $D$-модуль $V$. Существует $D$-линейное отображение*
$$h : A_1 \otimes ... \otimes A_n \to V$$
*такое, что диаграмма*

(3.5.12)

$$\begin{array}{c} & & A_1 \otimes ... \otimes A_n \\ & \nearrow^{f} & \downarrow^{h} \\ A_1 \times ... \times A_n & \xrightarrow{g} & V \end{array}$$

*коммутативна.*

Доказательство. Равенство (3.5.11) следует из равенств (3.5.1) и (3.5.10). Существование отображения $h$ следует из определения 3.5.3 и построений, выполненных при доказательстве теоремы 3.5.4. □

Равенства (3.5.5) и (3.5.6) можно записать в виде

(3.5.13) $$a_1 \otimes ... \otimes (a_i + b_i) \otimes ... \otimes a_n$$
$$= a_1 \otimes ... \otimes a_i \otimes ... \otimes a_n + a_1 \otimes ... \otimes b_i \otimes ... \otimes a_n$$

(3.5.14) $$a_1 \otimes ... \otimes (ca_i) \otimes ... \otimes a_n = c(a_1 \otimes ... \otimes a_i \otimes ... \otimes a_n)$$

$$a_i \in A_i \quad b_i \in A_i \quad c \in D$$

**Теорема 3.5.6.** *Пусть $A$ - алгебра над коммутативным кольцом $D$. Существует линейное отображение*
$$h : a \otimes b \in A \otimes A \to ab \in A$$

Доказательство. Теорема является следствием теоремы 3.5.5 и определения 3.2.1. □

**Теорема 3.5.7.** *Тензорное произведение $A_1 \otimes ... \otimes A_n$ свободных конечномерных алгебр $A_1$, ..., $A_n$ над коммутативным кольцом $D$ является свободной конечномерной алгеброй.*

*Пусть $\overline{\overline{e}}_i$ - базис алгебры $A_i$ над кольцом $D$. Произвольный тензор $a \in A_1 \otimes ... \otimes A_n$ можно представить в виде*

(3.5.15) $$a = a^{i_1...i_n} \overline{e}_{1 \cdot i_1} \otimes ... \otimes \overline{e}_{n \cdot i_n}$$

*Мы будем называть выражение $a^{i_1...i_n}$ **стандартной компонентой тензора**.*

Доказательство. Алгебры $A_1$, ..., $A_n$ являются модулями над кольцом $D$. Согласно теореме 3.5.4, $A_1 \otimes ... \otimes A_n$ является модулем.

Вектор $a_i \in A_i$ имеет разложение
$$a_i = a_i^k \overline{e}_{i \cdot k}$$

относительно базиса $\overline{\overline{e}}_i$. Из равенств (3.5.13), (3.5.14) следует

$$a_1 \otimes ... \otimes a_n = a_1^{i_1}...a_n^{i_n} \overline{e}_{1 \cdot i_1} \otimes ... \otimes \overline{e}_{n \cdot i_n}$$

Так как множество тензоров $a_1 \otimes ... \otimes a_n$ является множеством образующих модуля $A_1 \otimes ... \otimes A_n$, то тензор $a \in A_1 \otimes ... \otimes A_n$ можно записать в виде

(3.5.16) $$a = a^s a_{s \cdot 1}^{i_1}...a_{s \cdot n}^{i_n} \overline{e}_{1 \cdot i_1} \otimes ... \otimes \overline{e}_{n \cdot i_n}$$

где $a^s$, $a_{s \cdot 1}^{i_1}$, $a_{s \cdot n}^{i_n} \in F$. Положим

(3.5.17) $$a^s a_{s \cdot 1}^{i_1}...a_{s \cdot n}^{i_n} = a^{i_1...i_n}$$

Тогда равенство (3.5.16) примет вид (3.5.15).

Следовательно, множество тензоров $\overline{e}_{1 \cdot i_1} \otimes ... \otimes \overline{e}_{n \cdot i_n}$ является множеством образующих модуля $A_1 \otimes ... \otimes A_n$. Так как размерность модуля $A_i$, $i = 1, ..., n$, конечна, то конечно множество тензоров $\overline{e}_{1 \cdot i_1} \otimes ... \otimes \overline{e}_{n \cdot i_n}$. Следовательно, множество тензоров $\overline{e}_{1 \cdot i_1} \otimes ... \otimes \overline{e}_{n \cdot i_n}$ содержит базис модуля $A_1 \otimes ... \otimes A_n$, и модуль $A_1 \otimes ... \otimes A_n$ является свободным модулем над кольцом $D$.

Мы определим произведение тензоров типа $a_1 \otimes ... \otimes a_n$ покомпонентно

(3.5.18) $$(d_1 \otimes ... \otimes d_n)(c_1 \otimes ... \otimes c_n) = (d_1 c_1) \otimes ... \otimes (d_n c_n)$$

В частности, если для любого $i$, $i = 1, ..., n$, $a_i \in A_i$ имеет обратный, то тензор

$$a_1 \otimes ... \otimes a_n \in A_1 \otimes ... \otimes A_n$$

имеет обратный

$$(a_1 \otimes ... \otimes a_n)^{-1} = (a_1)^{-1} \otimes ... \otimes (a_n)^{-1}$$

Определение произведения (3.5.18) согласовано с равенством (3.5.14) так как

$$(a_1 \otimes ... \otimes (ca_i) \otimes ... \otimes a_n)(b_1 \otimes ... \otimes b_i \otimes ... \otimes b_n)$$
$$=(a_1 b_1) \otimes ... \otimes (ca_i)b_i \otimes ... \otimes (a_n b_n)$$
$$=c((a_1 b_1) \otimes ... \otimes (a_i b_i) \otimes ... \otimes (a_n b_n))$$
$$=c((a_1 \otimes ... \otimes a_n)(b_1 \otimes ... \otimes b_n))$$

Из равенства (3.5.13) следует дистрибутивность умножения по отношению к сложению

(3.5.19)
$$\begin{aligned}&(a_1 \otimes ... \otimes a_i \otimes ... \otimes a_n)\\&*((b_1 \otimes ... \otimes b_i \otimes ... \otimes b_n) + (b_1 \otimes ... \otimes c_i \otimes ... \otimes b_n))\\&=(a_1 \otimes ... \otimes a_i \otimes ... \otimes a_n)(b_1 \otimes ... \otimes (b_i + c_i) \otimes ... \otimes b_n)\\&=(a_1 b_1) \otimes ... \otimes (a_i(b_i + c_i)) \otimes ... \otimes (a_n b_n)\\&=(a_1 b_1) \otimes ... \otimes (a_i b_i + a_i c_i) \otimes ... \otimes (a_n b_n)\\&=(a_1 b_1) \otimes ... \otimes (a_i b_i) \otimes ... \otimes (a_n b_n)\\&+(a_1 b_1) \otimes ... \otimes (a_i c_i) \otimes ... \otimes (a_n b_n)\\&=(a_1 \otimes ... \otimes a_i \otimes ... \otimes a_n)(b_1 \otimes ... \otimes b_i \otimes ... \otimes b_n)\\&+(a_1 \otimes ... \otimes a_i \otimes ... \otimes a_n)(b_1 \otimes ... \otimes c_i \otimes ... \otimes b_n)\end{aligned}$$

Равенство (3.5.19) позволяет определить произведение для любых тензоров $a$, $b$. $\square$

**Замечание 3.5.8.** Согласно замечанию 3.2.2, мы можем определить различные структуры алгебры в тензорном произведении алгебр. Например, алгебры $A_1 \otimes A_2$, $A_1 \otimes A_2^*$, $A_1^* \otimes A_2$ определены на одном и том же модуле. $\square$

**Теорема 3.5.9.** *Пусть* $\overline{\overline{e}}_i$ *- базис алгебры* $A_i$ *над кольцом* $D$. *Пусть* $B_i{}_{\cdot kl}^{\phantom{\cdot}j}$ *- структурные константы алгебры* $A_i$ *относительно базиса* $\overline{\overline{e}}_i$. *Структурные константы тензорного произведение* $A_1 \otimes ... \otimes A_n$ *относительно базиса* $\overline{e}_{1 \cdot i_1} \otimes ... \otimes \overline{e}_{n \cdot i_n}$ *имеют вид*

$$(3.5.20) \qquad B_{\cdot k_1...k_n \cdot l_1...l_n}^{\cdot j_1...j_n} = B_1{}_{\cdot k_1 l_1}^{\phantom{\cdot}j_1} ... B_n{}_{\cdot k_n l_n}^{\phantom{\cdot}j_n}$$

Доказательство. Непосредственное перемножение тензоров $\overline{e}_{1 \cdot i_1} \otimes ... \otimes \overline{e}_{n \cdot i_n}$ имеет вид

$$(3.5.21) \begin{aligned} &(\overline{e}_{1 \cdot k_1} \otimes ... \otimes \overline{e}_{n \cdot k_n})(\overline{e}_{1 \cdot l_1} \otimes ... \otimes \overline{e}_{n \cdot l_n}) \\ =&(\overline{e}_{1 \cdot k_1} \overline{e}_{1 \cdot l_1}) \otimes ... \otimes (\overline{e}_{n \cdot k_n} \overline{e}_{n \cdot l_n}) \\ =&(\overline{e}_{1 \cdot k_1} \overline{e}_{1 \cdot l_1}) \otimes ... \otimes (\overline{e}_{n \cdot k_n} \overline{e}_{n \cdot l_n}) \\ =&(B_1{}_{\cdot k_1 l_1}^{\phantom{\cdot}j_1} \overline{e}_{1 \cdot j_1}) \otimes ... \otimes (B_n{}_{\cdot k_n l_n}^{\phantom{\cdot}j_n} \overline{e}_{n \cdot j_n}) \\ =&B_1{}_{\cdot k_1 l_1}^{\phantom{\cdot}j_1} ... B_n{}_{\cdot k_n l_n}^{\phantom{\cdot}j_n} \overline{e}_{1 \cdot j_1} \otimes ... \otimes \overline{e}_{n \cdot j_n} \end{aligned}$$

Согласно определению структурных констант

$$(3.5.22) \quad (\overline{e}_{1 \cdot k_1} \otimes ... \otimes \overline{e}_{n \cdot k_n})(\overline{e}_{1 \cdot l_1} \otimes ... \otimes \overline{e}_{n \cdot l_n}) = B_{\cdot k_1...k_n \cdot l_1...l_n}^{\cdot j_1...j_n}(\overline{e}_{1 \cdot j_1} \otimes ... \otimes \overline{e}_{n \cdot j_n})$$

Равенство (3.5.20) следует из сравнения (3.5.21), (3.5.22).

Из цепочки равенств

$$\begin{aligned}&(a_1 \otimes ... \otimes a_n)(b_1 \otimes ... \otimes b_n) \\ =&(a_1^{k_1} \overline{e}_{1 \cdot k_1} \otimes ... \otimes a_n^{k_n} \overline{e}_{n \cdot k_n})(b_1^{l_1} \overline{e}_{1 \cdot l_1} \otimes ... \otimes b_n^{l_n} \overline{e}_{n \cdot l_n}) \\ =&a_1^{k_1}...a_n^{k_n} b_1^{l_1}...b_n^{l_n}(\overline{e}_{1 \cdot k_1} \otimes ... \otimes \overline{e}_{n \cdot k_n})(\overline{e}_{1 \cdot l_1} \otimes ... \otimes \overline{e}_{n \cdot l_n}) \\ =&a_1^{k_1}...a_n^{k_n} b_1^{l_1}...b_n^{l_n} B_{\cdot k_1...k_n \cdot l_1...l_n}^{\cdot j_1...j_n}(\overline{e}_{1 \cdot j_1} \otimes ... \otimes \overline{e}_{n \cdot j_n}) \\ =&a_1^{k_1}...a_n^{k_n} b_1^{l_1}...b_n^{l_n} B_1{}_{\cdot k_1 l_1}^{\phantom{\cdot}j_1} ... B_n{}_{\cdot k_n l_n}^{\phantom{\cdot}j_n}(\overline{e}_{1 \cdot j_1} \otimes ... \otimes \overline{e}_{n \cdot j_n}) \\ =&(a_1^{k_1} b_1^{l_1} B_1{}_{\cdot k_1 l_1}^{\phantom{\cdot}j_1} \overline{e}_{1 \cdot j_1}) \otimes ... \otimes (a_n^{k_n} b_n^{l_n} B_n{}_{\cdot k_n l_n}^{\phantom{\cdot}j_n} \overline{e}_{n \cdot j_n}) \\ =&(a_1 b_1) \otimes ... \otimes (a_n b_n)\end{aligned}$$

следует, что определение произведения (3.5.22) со структурными константами (3.5.20) согласовано с определением произведения (3.5.18). $\square$

**Теорема 3.5.10.** *Для тензоров* $a, b \in A_1 \otimes ... \otimes A_n$ *стандартные компоненты произведения удовлетворяют равенству*

$$(3.5.23) \qquad (ab)^{j_1...j_n} = B_{\cdot k_1...k_n \cdot l_1...l_n}^{\cdot j_1...j_n} a^{k_1...k_n} b^{l_1...l_n}$$

Доказательство. Согласно определению

$$(3.5.24) \qquad ab = (ab)^{j_1...j_n} \overline{e}_{1 \cdot j_1} \otimes ... \otimes \overline{e}_{n \cdot j_n}$$

В тоже время

$$(3.5.25) \begin{aligned} ab =& a^{k_1...k_n} \overline{e}_{1 \cdot k_1} \otimes ... \otimes \overline{e}_{n \cdot k_n} b^{k_1...k_n} \overline{e}_{1 \cdot l_1} \otimes ... \otimes \overline{e}_{n \cdot l_n} \\ =& a^{k_1...k_n} b^{k_1...k_n} B_{\cdot k_1...k_n \cdot l_1...l_n}^{\cdot j_1...j_n} \overline{e}_{1 \cdot j_1} \otimes ... \otimes \overline{e}_{n \cdot j_n} \end{aligned}$$

Равенство (3.5.23) следует из равенств (3.5.24), (3.5.25). $\square$

**Теорема 3.5.11.** *Если алгебра $A_i$, $i = 1, ..., n$, ассоциативна, то тензорное произведение $A_1 \otimes ... \otimes A_n$ - ассоциативная алгебра.*

Доказательство. Поскольку

$$((\overline{e}_{1\cdot i_1} \otimes ... \otimes \overline{e}_{n\cdot i_n})(\overline{e}_{1\cdot j_1} \otimes ... \otimes \overline{e}_{n\cdot j_n}))(\overline{e}_{1\cdot k_1} \otimes ... \otimes \overline{e}_{n\cdot k_n})$$
$$=((\overline{e}_{1\cdot i_1}\overline{e}_{1\cdot j_1}) \otimes ... \otimes (\overline{e}_{n\cdot i_n}\overline{e}_{1\cdot j_n}))(\overline{e}_{1\cdot k_1} \otimes ... \otimes \overline{e}_{n\cdot k_n})$$
$$=((\overline{e}_{1\cdot i_1}\overline{e}_{1\cdot j_1})\overline{e}_{1\cdot k_1}) \otimes ... \otimes ((\overline{e}_{n\cdot i_n}\overline{e}_{1\cdot j_n})\overline{e}_{1\cdot k_n})$$
$$=(\overline{e}_{1\cdot i_1}(\overline{e}_{1\cdot j_1}\overline{e}_{1\cdot k_1})) \otimes ... \otimes (\overline{e}_{n\cdot i_n}(\overline{e}_{1\cdot j_n}\overline{e}_{1\cdot k_n}))$$
$$=(\overline{e}_{1\cdot i_1} \otimes ... \otimes \overline{e}_{n\cdot i_n})((\overline{e}_{1\cdot j_1}\overline{e}_{1\cdot k_1}) \otimes ... \otimes (\overline{e}_{1\cdot j_n}\overline{e}_{1\cdot k_n}))$$
$$=(\overline{e}_{1\cdot i_1} \otimes ... \otimes \overline{e}_{n\cdot i_n})((\overline{e}_{1\cdot j_1} \otimes ... \otimes \overline{e}_{n\cdot j_n})(\overline{e}_{1\cdot k_1} \otimes ... \otimes \overline{e}_{n\cdot k_n}))$$

то

$$(ab)c = a^{i_1...i_n}b^{j_1...j_n}c^{k_1...k_n}$$
$$((\overline{e}_{1\cdot i_1} \otimes ... \otimes \overline{e}_{n\cdot i_n})(\overline{e}_{1\cdot j_1} \otimes ... \otimes \overline{e}_{n\cdot j_n}))(\overline{e}_{1\cdot k_1} \otimes ... \otimes \overline{e}_{n\cdot k_n})$$
$$= a^{i_1...i_n}b^{j_1...j_n}c^{k_1...k_n}$$
$$(\overline{e}_{1\cdot i_1} \otimes ... \otimes \overline{e}_{n\cdot i_n})((\overline{e}_{1\cdot j_1} \otimes ... \otimes \overline{e}_{n\cdot j_n})(\overline{e}_{1\cdot k_1} \otimes ... \otimes \overline{e}_{n\cdot k_n}))$$
$$= a(bc)$$

□

## 3.6. Линейное отображение в ассоциативную алгебру

**Теорема 3.6.1.** *Рассмотрим $D$-алгебры $A_1$ и $A_2$. Для заданного отображения $f \in \mathcal{L}(A_1; A_2)$ отображение*

$$g : A_2 \times A_2 \to \mathcal{L}(A_1; A_2)$$
$$g(a, b) \circ f = afb$$

*является билинейным отображением.*

Доказательство. Утверждение теоремы следует из цепочек равенств

$$((a_1 + a_2)fb) \circ x = (a_1 + a_2) \ f \circ x \ b = a_1 \ f \circ x \ b + a_2 \ f \circ x \ b$$
$$= (a_1 fb) \circ x + (a_2 fb) \circ x = (a_1 fb + a_2 fb) \circ x$$
$$((pa)fb) \circ x = (pa) \ f \circ x \ b = p(a \ f \circ x \ b) = p((afb) \circ x) = (p(afb)) \circ x$$
$$(af(b_1 + b_2)) \circ x = a \ f \circ x \ (b_1 + b_2) = a \ f \circ x \ b_1 + a \ f \circ x \ b_2$$
$$= (afb_1) \circ x + (afb_2) \circ x = (afb_1 + afb_2) \circ x$$
$$(af(pb)) \circ x = a \ f \circ x \ (pb) = p(a \ f \circ x \ b) = p((afb) \circ x) = (p(afb)) \circ x$$

□

**Теорема 3.6.2.** *Рассмотрим $D$-алгебры $A_1$ и $A_2$. Для заданного отображения $f \in \mathcal{L}(A_1; A_2)$ существует линейное отображение*

$$h : A_2 \otimes A_2 \to \mathcal{L}(A_1; A_2)$$

*определённое равенством*

(3.6.1)     $(a \otimes b) \circ f = afb$

Доказательство. Утверждение теоремы является следствием теорем 3.5.5, 3.6.1. □

**Теорема 3.6.3.** *Рассмотрим $D$-алгебры $A_1$ и $A_2$. Линейное отображение*
$$h : A_2 \otimes A_2 \to {}^*\mathcal{L}(A_1; A_2)$$
*определённое равенством*

(3.6.2) $\qquad (a \otimes b) \circ f = afb \quad a, b \in A_2 \quad f \in \mathcal{L}(A_1; A_2)$

*является представлением*[3.10] *модуля $A_2 \otimes A_2$ в модуле $\mathcal{L}(A_1; A_2)$.*

Доказательство. Согласно теореме 3.3.4, отображение (3.6.2) является преобразованием модуля $\mathcal{L}(A_1; A_2)$. Для данного тензора $c \in A_2 \otimes A_2$ преобразование $h(c)$ является линейным преобразованием модуля $\mathcal{L}(A_1; A_2)$, так как

$$\begin{aligned}
((a \otimes b) \circ (f_1 + f_2)) \circ x &= (a(f_1 + f_2)b) \circ x = a((f_1 + f_2) \circ x)b \\
&= a(f_1 \circ x + f_2 \circ x)b = a(f_1 \circ x)b + a(f_2 \circ x)b \\
&= (af_1 b) \circ x + (af_2 b) \circ x \\
&= (a \otimes b) \circ f_1 \circ x + (a \otimes b) \circ f_2 \circ x \\
&= ((a \otimes b) \circ f_1 + (a \otimes b) \circ f_2) \circ x \\
((a \otimes b) \circ (pf)) \circ x &= (a(pf)b) \circ x = a((pf) \circ x)b \\
&= a(p\ f \circ x)b = pa(f \circ x)b \\
&= p\ (afb) \circ x = p\ ((a \otimes b) \circ f) \circ x \\
&= (p((a \otimes b) \circ f)) \circ x
\end{aligned}$$

Согласно теореме 3.6.2, отображение (3.6.2) является линейным отображением. Согласно определению [7]-2.1.4 отображение (3.6.2) является представлением модуля $A_2 \otimes A_2$ в модуле $\mathcal{L}(A_1; A_2)$. $\square$

**Теорема 3.6.4.** *Пусть $A$ - алгебра над коммутативным кольцом $D$. Алгебра $A \otimes A$, в которой произведение определено согласно правилу*

(3.6.3) $\qquad (a \otimes b) \circ (c \otimes d) = (ac) \otimes (db)$

*порождает представление в модуле $\mathcal{L}(A; A)$. Это представление позволяет отождествить тензор $d \in A \otimes A$ с преобразованием $d \circ \delta$, где $\delta$ - тождественное преобразование.*

Доказательство. Согласно теореме 3.6.2, отображение $f \in \mathcal{L}(A; A)$ и тензор $d \in A \otimes A$ порождают отображение

(3.6.4) $\qquad\qquad\qquad x \to (d \circ f) \circ x$

Если мы положим $f = \delta$, $d = a \otimes b$, то равенство (3.6.4) приобретает вид

(3.6.5) $\qquad ((a \otimes b) \circ \delta) \circ x = (a\delta b) \circ x = a\ (\delta \circ x)\ b = axb$

Если мы положим

(3.6.6) $\qquad ((a \otimes b) \circ \delta) \circ x = (a \otimes b) \circ (\delta \circ x) = (a \otimes b) \circ x$

то сравнение равенств (3.6.5) и (3.6.6) даёт основание отождествить действие тензора $a \otimes b$ с преобразованием $(a \otimes b) \circ \delta$. Следовательно, отображение

(3.6.7) $\qquad\qquad d \in A \otimes A \to d \circ \delta \in \mathcal{L}(A; A)$

является гомоморфизмом модуля $A \otimes A$ в модуль $\mathcal{L}(A; A)$.

---

[3.10]Определение представления $\Omega$-алгебры дано в определении [7]-2.1.4.

Отображение (3.6.7) является также гомоморфизмом алгебр, так как произведение преобразований $a \otimes b$ и $c \otimes d$ имеет вид

$$\begin{aligned}((a \otimes b) \circ (c \otimes d)) \circ x &= ((ac) \otimes (db)) \circ x \\ &= (ac)x(db) \\ &= a(cxd)b \\ &= (a \otimes b) \circ (cxd) \\ &= (a \otimes b) \circ ((c \otimes d) \circ x)\end{aligned}$$

□

Из теоремы 3.6.4 следует, что отображение (3.6.2) можно рассматривать как произведение отображений $a \otimes b$ и $f$. Это позволяет вместо представления модуля $A_2 \otimes A_2$ в модуле $\mathcal{L}(A_1; A_2)$ рассматривать представление алгебры $A_2 \otimes A_2$ в модуле $\mathcal{L}(A_1; A_2)$.

Тензор $a \in A_2 \otimes A_2$ **невырожден**, если существует тензор $b \in A_2 \otimes A_2$ такой, что $a \circ b = 1 \otimes 1$.

**Определение 3.6.5.** Рассмотрим представление алгебры $A_2 \otimes A_2$ в модуле $\mathcal{L}(A_1; A_2)$.[3.11] **Орбитой линейного отображения** $f \in \mathcal{L}(A_1; A_2)$ называется множество

$$(A_2 \otimes A_2) \circ f = \{g = d \circ f : d \in A_2 \otimes A_2\}$$

□

**Теорема 3.6.6.** *Рассмотрим $D$-алгебру $A_1$ и ассоциативную $D$-алгебру $A_2$. Рассмотрим представление алгебры $A_2 \otimes A_2$ в модуле $\mathcal{L}(A_1; A_2)$. Отображение*

$$h : A_1 \to A_2$$

*порождённое отображением*

$$f : A_1 \to A_2$$

*имеет вид*

$$h = (a_{s \cdot 0} \otimes a_{s \cdot 1}) \circ f = a_{s \cdot 0} f a_{s \cdot 1}$$

Доказательство. Произвольный тензор $a \in A_2 \otimes A_2$ можно представить в виде

$$a = a_{s \cdot 0} \otimes a_{s \cdot 1}$$

Согласно теореме 3.6.3, отображение (3.6.2) линейно. Это доказывает утверждение теоремы. □

**Теорема 3.6.7.** *Пусть $A_2$ - алгебра с единицей $e$. Пусть $a \in A_2 \otimes A_2$ - невырожденный тензор. Орбиты линейных отображений $f \in \mathcal{L}(A_1; A_2)$ и $g = a \circ f$ совпадают*

(3.6.8) $$(A_2 \otimes A_2) \circ f = (A_2 \otimes A_2) \circ g$$

---

[3.11]Определение дано по аналогии с определением [7]-2.4.12.

Доказательство. Если $h \in (A_2 \otimes A_2) \circ g$, то существует $b \in A_2 \otimes A_2$ такое, что $h = b \circ g$. Тогда

(3.6.9) $$h = b \circ (a \circ f) = (b \circ a) \circ f$$

Следовательно, $h \in (A_2 \otimes A_2) \circ f$,

(3.6.10) $$(A_2 \otimes A_2) \circ g \subset (A_2 \otimes A_2) \circ f$$

Так как $a$ - невырожденный тензор, то

(3.6.11) $$f = a^{-1} \circ g$$

Если $h \in (A_2 \otimes A_2) \circ f$, то существует $b \in A_2 \otimes A_2$ такое, что

(3.6.12) $$h = b \circ f$$

Из равенств (3.6.11), (3.6.12), следует, что

$$h = b \circ (a^{-1} \circ g) = (b \circ a^{-1}) \circ g$$

Следовательно, $h \in (A_2 \otimes A_2) \circ g$,

(3.6.13) $$(A_2 \otimes A_2) \circ f \subset (A_2 \otimes A_2) \circ g$$

(3.6.8) следует из равенств (3.6.10), (3.6.13). □

Из теоремы 3.6.7 также следует, что, если $g = a \circ f$ и $a \in A_2 \otimes A_2$ - вырожденный тензор, то отношение (3.6.10) верно. Однако основной результат теоремы 3.6.7 состоит в том, что представления алгебры $A_2 \otimes A_2$ в модуле $\mathcal{L}(A_1; A_2)$ порождает отношение эквивалентности в модуле $\mathcal{L}(A_1; A_2)$. Если удачно выбрать представители каждого класса эквивалентности, то полученное множество будет множеством образующих рассматриваемого представления.[3.12]

## 3.7. Линейное отображение в свободную конечно мерную ассоциативную алгебру

**Теорема 3.7.1.** *Пусть $A_1$ - алгебра над кольцом $D$. Пусть $A_2$ - свободная конечно мерная ассоциативная алгебра над кольцом $D$. Пусть $\overline{\overline{e}}$ - базис алгебры $A_2$ над кольцом $D$. Отображение*

(3.7.1) $$g = a \circ f$$

*порождённое отображением $f \in (A_1; A_2)$ посредством тензора $a \in A_2 \otimes A_2$, имеет стандартное представление*

(3.7.2) $$g = a^{ij}(\overline{e}_i \otimes \overline{e}_j) \circ f = a^{ij}\overline{e}_i f \overline{e}_j$$

Доказательство. Согласно теореме 3.5.7, стандартное представление тензора $a$ имеет вид

(3.7.3) $$a = a^{ij}\overline{e}_i \otimes \overline{e}_j$$

Равенство (3.7.2) следует из равенств (3.7.1), (3.7.3). □

---

[3.12]Множество образующих представления определено в определении [7]-2.6.5.

**Теорема 3.7.2.** *Пусть $\overline{\overline{e}}_1$ - базис свободной конечно мерной $D$-алгебры $A_1$. Пусть $\overline{\overline{e}}_2$ - базис свободной конечно мерной ассоциативной $D$-алгебры $A_2$. Пусть $B_2{}_{\cdot kl}^{\,p}$ - структурные константы алгебры $A_2$. Координаты отображения*

$$g = a \circ f$$

*порождённого отображением $f \in (A_1; A_2)$ посредством тензора $a \in A_2 \otimes A_2$, и его стандартные компоненты связаны равенством*

(3.7.4) $$g_l^k = f_l^m g^{ij} B_2{}_{\cdot im}^{\,p} B_2{}_{\cdot pj}^{\,k}$$

Доказательство. Относительно базисов $\overline{\overline{e}}_1$ и $\overline{\overline{e}}_2$, линейные отображения $f$ и $g$ имеют вид

(3.7.5) $$f \circ x = f_j^i x^j \overline{e}_{2\cdot i}$$

(3.7.6) $$g \circ x = g_j^i x^j \overline{e}_{2\cdot i}$$

Из равенств (3.7.5), (3.7.6), (3.7.2) следует

(3.7.7) $$\begin{aligned} g_l^k x^l \overline{e}_{2\cdot k} &= a^{ij} \overline{e}_{2\cdot i} f_l^m x^l \overline{e}_{2\cdot m} \overline{e}_{2\cdot j} \\ &= a^{ij} f_l^m x^l B_2{}_{\cdot im}^{\,p} B_2{}_{\cdot pj}^{\,k} \overline{e}_{2\cdot k} \end{aligned}$$

Так как векторы $\overline{e}_{2\cdot k}$ линейно независимы и $x^i$ произвольны, то равенство (3.7.4) следует из равенства (3.7.7). □

**Теорема 3.7.3.** *Пусть $\overline{\overline{e}}_1$ - базис свободной конечно мерной $D$-алгебры $A_1$. Пусть $\overline{\overline{e}}_2$ - базис свободной конечно мерной ассоциативной $D$-алгебры $A_2$. Пусть $B_2{}_{\cdot kl}^{\,p}$ - структурные константы алгебры $A_2$. Рассмотрим матрицу*

(3.7.8) $$\mathcal{B} = \left(\mathcal{B}_{\,m\,\cdot ij}^{\cdot k}\right) = \left(B_2{}_{\cdot im}^{\,p} B_2{}_{\cdot pj}^{\,k}\right)$$

*строки которой пронумерованы индексом $_{\,m}^{\cdot k}$ и столбцы пронумерованы индексом $_{\cdot ij}$. Если матрица $\mathcal{B}$ невырождена, то для заданных координат линейного преобразования $g_k^l$ и для отображения $f = \delta$, система линейных уравнений (3.7.4) относительно стандартных компонент этого преобразования $g^{kr}$ имеет единственное решение.*

*Если матрица $\mathcal{B}$ вырождена, то условием существования решения системы линейных уравнений (3.7.4) является равенство*

(3.7.9) $$\mathrm{rank}\left(\mathcal{B}_{\,m\,\cdot ij}^{\cdot k} \quad g_m^k\right) = \mathrm{rank}\,\mathcal{B}$$

*В этом случае система линейных уравнений (3.7.4) имеет бесконечно много решений и существует линейная зависимость между величинами $g_m^k$.*

Доказательство. Утверждение теоремы является следствием теории линейных уравнений над кольцом. □

**Теорема 3.7.4.** *Пусть $A$ - свободная конечно мерная ассоциативная алгебра над кольцом $D$. Пусть $\overline{\overline{e}}$ - базис алгебры $A$ над кольцом $D$. Пусть $B_{kl}^p$ - структурные константы алгебры $A$. Пусть матрица (3.7.8) вырождена. Пусть линейное отображение $f \in \mathcal{L}(A; A)$ невырождено. Если координаты линейных преобразований $f$ и $g$ удовлетворяют равенству*

(3.7.10) $$\mathrm{rank}\left(\mathcal{B}_{\,m\,\cdot ij}^{\cdot k} \quad g_m^k \quad f_m^k\right) = \mathrm{rank}\,\mathcal{B}$$

то система линейных уравнений

(3.7.11) $$g_l^k = f_l^m g^{ij} B_{im}^p B_{pj}^k$$

имеет бесконечно много решений.

Доказательство. Согласно равенству (3.7.10) и теореме 3.7.3, система линейных уравнений

(3.7.12) $$f_l^k = f^{ij} B_{il}^p B_{pj}^k$$

имеет бесконечно много решений, соответствующих линейному отображению

(3.7.13) $$f = f^{ij} \overline{e}_i \otimes \overline{e}_j$$

Согласно равенству (3.7.10) и теореме 3.7.3, система линейных уравнений

(3.7.14) $$g_l^k = g^{ij} B_{il}^p B_{pj}^k$$

имеет бесконечно много решений, соответствующих линейному отображению

(3.7.15) $$g = g^{ij} \overline{e}_i \otimes \overline{e}_j$$

Отображения $f$ и $g$ порождены отображением $\delta$. Согласно теореме 3.6.7, отображение $f$ порождает отображение $g$. Это доказывает утверждение теоремы. $\square$

**Теорема 3.7.5.** *Пусть $A$ - свободная конечно мерная ассоциативная алгебра над кольцом $D$. Представление алгебры $A \otimes A$ в алгебре $\mathcal{L}(A; A)$ имеет конечный базис $\overline{\overline{I}}$.*

*(1) Линейное отображение $f \in \mathcal{L}(A; A)$ имеет вид*

(3.7.16) $$f = \sum_k (a_{s_k \cdot 0} \otimes a_{s_k \cdot 1}) \circ I_k = \sum_k a_{s_k \cdot 0} I_k a_{s_k \cdot 1}$$

*(2) Его стандартное представление имеет вид*

(3.7.17) $$f = a^{k \cdot ij}(\overline{e}_i \otimes \overline{e}_j) \circ I_k = a^{k \cdot ij} \overline{e}_i I_k \overline{e}_j$$

Доказательство. Из теоремы 3.7.4 следует, что если матрица $\mathcal{B}$ вырождена и отображение $f$ удовлетворяет равенству

(3.7.18) $$\operatorname{rank} \begin{pmatrix} \mathcal{B}_{m \cdot ij}^{\cdot k} & f_m^k \end{pmatrix} = \operatorname{rank} \mathcal{B}$$

то отображение $f$ порождает то же самое множество отображений, что порождено отображением $\delta$. Следовательно, для того, чтобы построить базис представления алгебры $A \otimes A$ в модуле $\mathcal{L}(A; A)$ мы должны выполнить следующее построение.

Множество решений системы уравнений (3.7.11) порождает свободный подмодуль $\mathcal{L}$ модуля $\mathcal{L}(A; A)$. Мы строим базис $(\overline{h}_1, ..., \overline{h}_k)$ подмодуля $\mathcal{L}$. Затем дополняем этот базис линейно независимыми векторами $\overline{h}_{k+1}, ..., \overline{h}_m$, которые не принадлежат подмодулю $\mathcal{L}$, таким образом, что множество векторов $\overline{h}_1, ..., \overline{h}_m$ является базисом модуля $\mathcal{L}(A; A)$. Множество орбит $(A \otimes A) \circ \delta$, $(A \otimes A) \circ \overline{h}_{k+1}$, ..., $(A \otimes A) \circ \overline{h}_m$ порождает модуль $\mathcal{L}(A; A)$. Поскольку множество орбит конечно, мы можем выбрать орбиты так, чтобы они не пересекались. Для каждой орбиты мы можем выбрать представитель, порождающий эту орбиту. $\square$

**Пример 3.7.6.** Для поля комплексных чисел алгебра $\mathcal{L}(C;C)$ имеет базис

$$I_{\mathbf{0}} \circ z = z$$
$$I_{\mathbf{1}} \circ z = \overline{z}$$

Для алгебры кватернионов алгебра $\mathcal{L}(H;H)$ имеет базис

$$I_{\mathbf{0}} \circ z = z$$

□

### 3.8. Линейное отображение в неассоциативную алгебру

Так как произведение неассоциативно, мы можем предположить, что действие $a, b \in A$ на отображение $f$ может быть представленно либо в виде $a(fb)$, либо в виде $(af)b$. Однако это предположение приводит нас к довольно сложной структуре линейного отображения. Чтобы лучше представить насколько сложна структура линейного отображения, мы начнём с рассмотрения левого и правого сдвигов в неассоциативной алгебре.

**Теорема 3.8.1.** *Пусть*

(3.8.1) $$l(a) \circ x = ax$$

*отображение левого сдвига. Тогда*

(3.8.2) $$l(a) \circ l(b) = l(ab) - (a,b)_1$$

*где мы определили линейное отображение*

$$(a,b)_1 \circ x = (a,b,x)$$

Доказательство. Из равенств (3.2.2), (3.8.1) следует

(3.8.3) $$\begin{aligned}(l(a) \circ l(b)) \circ x &= l(a) \circ (l(b) \circ x) \\ &= a(bx) = (ab)x - (a,b,x) \\ &= l(ab) \circ x - (a,b)_1 \circ x\end{aligned}$$

Равенство (3.8.2) следует из равенства (3.8.3). □

**Теорема 3.8.2.** *Пусть*

(3.8.4) $$r(a) \circ x = xa$$

*отображение правого сдвига. Тогда*

(3.8.5) $$r(a) \circ r(b) = r(ba) + (b,a)_2$$

*где мы определили линейное отображение*

$$(b,a)_2 \circ x = (x,b,a)$$

Доказательство. Из равенств (3.2.2), (3.8.4) следует

(3.8.6) $$\begin{aligned}(r(a) \circ r(b)) \circ x &= r(a) \circ (r(b) \circ x) \\ &= (xb)a = x(ba) + (x,b,a) \\ &= r(ba) \circ x + (x,b,a)\end{aligned}$$

Равенство (3.8.5) следует из равенства (3.8.6). □

Пусть
$$f : A \to A \quad f = (ax)b$$
линейное отображение алгебры $A$. Согласно теореме 3.3.4, отображение
$$g : A \to A \quad g = (cf)d$$
также линейное отображение. Однако неочевидно, можем ли мы записать отображение $g$ в виде суммы слагаемых вида $(ax)b$ и $a(xb)$.

Если $A$ - свободная конечно мерная алгебра, то мы можем предположить, что линейное отображение имеет стандартное представление в виде[3.13]
$$(3.8.8) \qquad f \circ x = f^{ij} \, (\overline{e}_i x) \overline{e}_j$$
В этом случае мы можем применить теорему 3.7.5 для отображений в неассоциативную алгебру.

**Теорема 3.8.3.** *Пусть $\overline{\overline{e}}_1$ - базис свободной конечно мерной $D$-алгебры $A_1$. Пусть $\overline{\overline{e}}_2$ - базис свободной конечно мерной неассоциативной $D$-алгебры $A_2$. Пусть $B_2 \cdot {}_{kl}^{p}$ - структурные константы алгебры $A_2$. Пусть отображение*
$$(3.8.9) \qquad g = a \circ f$$
*порождённое отображением $f \in (A_1; A_2)$ посредством тензора $a \in A_2 \otimes A_2$, имеет стандартное представление*
$$(3.8.10) \qquad g = a^{ij}(\overline{e}_i \otimes \overline{e}_j) \circ f = a^{ij}(\overline{e}_i f)\overline{e}_j$$
*Координаты отображения (3.8.9) и его стандартные компоненты связаны равенством*
$$(3.8.11) \qquad g_l^k = f_l^m g^{ij} B_2 \cdot {}_{im}^{p} B_2 \cdot {}_{pj}^{k}$$

Доказательство. Относительно базисов $\overline{\overline{e}}_1$ и $\overline{\overline{e}}_2$, линейные отображения $f$ и $g$ имеют вид
$$(3.8.12) \qquad f \circ x = f_j^i x^j \overline{e}_{2 \cdot i}$$
$$(3.8.13) \qquad g \circ x = g_j^i x^j \overline{e}_{2 \cdot i}$$
Из равенств (3.8.12), (3.8.13), (3.8.10) следует
$$(3.8.14) \qquad \begin{aligned} g_l^k x^l \overline{e}_{2 \cdot k} &= a^{ij}(\overline{e}_{2 \cdot i}(f_l^m x^l \overline{e}_{2 \cdot m}))\overline{e}_{2 \cdot j} \\ &= a^{ij} f_l^m x^l B_2 \cdot {}_{im}^{p} B_2 \cdot {}_{pj}^{k} \overline{e}_{2 \cdot k} \end{aligned}$$

Так как векторы $\overline{e}_{2 \cdot k}$ линейно независимы и $x^i$ произвольны, то равенство (3.8.11) следует из равенства (3.8.14). □

**Теорема 3.8.4.** *Пусть $A$ - свободная конечно мерная неассоциативная алгебра над кольцом $D$. Представление алгебры $A \otimes A$ в алгебре $\mathcal{L}(A; A)$ имеет конечный базис $\overline{\overline{I}}$.*

---

[3.13]Выбор произволен. Мы можем рассмотреть стандартное представление в виде
$$f \circ x = f^{ij} \overline{e}_i (x \overline{e}_j)$$
Тогда равенство (3.8.11) имеет вид
$$(3.8.7) \qquad g_l^k = f_l^m g^{ij} B_2 \cdot {}_{ip}^{k} B_2 \cdot {}_{mj}^{p}$$
Я выбрал выражение (3.8.8) так как порядок сомножителей соответствует порядку, выбранному в теореме 3.7.5.

(1) *Линейное отображение* $f \in \mathcal{L}(A; A)$ *имеет вид*

(3.8.15) $$f = \sum_k (a_{s_k \cdot 0} \otimes a_{s_k \cdot 1}) \circ I_k = \sum_k (a_{s_k \cdot 0} I_k) a_{s_k \cdot 1}$$

(2) *Его стандартное представление имеет вид*

(3.8.16) $$f = a^{k \cdot ij}(\overline{e}_i \otimes \overline{e}_j) \circ I_k = a^{k \cdot ij}(\overline{e}_i I_k)\overline{e}_j$$

Доказательство. Рассмотрим матрицу (3.7.8). Если матрица $\mathcal{B}$ невырождена, то для заданных координат линейного преобразования $g_k^l$ и для отображения $f = \delta$, система линейных уравнений (3.8.11) относительно стандартных компонент этого преобразования $g^{kr}$ имеет единственное решение. Если матрица $\mathcal{B}$ вырождена, то согласно теореме 3.7.5 существует конечный базис $\overline{\overline{I}}$, порождающий множество линейных отображений. $\square$

В отличие от случая ассоциативной алгебры множество генераторов $I$ в теореме 3.8.4 не является минимальным. Из равенства (3.8.2) следует, что неверно равенство (3.6.9). Следовательно, орбиты отображений $I_k$ не порождают отношения эквивалентности в алгебре $L(A; A)$. Так как мы рассматриваем только отображения вида $(aI_k)b$, то возможно, что при $k \neq l$ отображение $I_k$ порождает отображение $I_l$, если рассмотреть все возможные операции в алгебре $A$. Поэтому множество образующих $I_k$ неассоциативной алгебры $A$ не играет такой критической роли как отображение сопряжения в поле комплексных чисел. Ответ на вопрос насколько важно отображение $I_k$ в неассоциативной алгебре требует дополнительного исследования.

Глава 4

# Алгебра с делением

### 4.1. Линейная функция комплексного поля

**Теорема 4.1.1** (Уравнения Коши-Римана). *Рассмотрим поле комплексных чисел $C$ как двумерную алгебру над полем действительных чисел. Положим*

$$(4.1.1) \qquad \overline{e}_{C\cdot\mathbf{0}} = 1 \quad \overline{e}_{C\cdot\mathbf{1}} = i$$

*базис алгебры $C$. Тогда в этом базисе произведение имеет вид*

$$(4.1.2) \qquad \overline{e}_{C\cdot\mathbf{1}}^2 = -\overline{e}_{C\cdot\mathbf{0}}$$

*и структурные константы имеют вид*

$$(4.1.3) \qquad \begin{aligned} B_{C\cdot\mathbf{00}}^{\mathbf{0}} &= 1 & B_{C\cdot\mathbf{01}}^{\mathbf{1}} &= 1 \\ B_{C\cdot\mathbf{10}}^{\mathbf{1}} &= 1 & B_{C\cdot\mathbf{11}}^{\mathbf{0}} &= -1 \end{aligned}$$

*Матрица линейной функции*

$$y^i = x^j f^i_j$$

*поля комплексных чисел над полем действительных чисел удовлетворяет соотношению*

$$(4.1.4) \qquad f^{\mathbf{0}}_{\mathbf{0}} = f^{\mathbf{1}}_{\mathbf{1}}$$

$$(4.1.5) \qquad f^{\mathbf{1}}_{\mathbf{0}} = -f^{\mathbf{0}}_{\mathbf{1}}$$

Доказательство. Равенства (4.1.2) и (4.1.3) следуют из равенства $i^2 = -1$. Пользуясь равенством [6]-(3.1.17) получаем соотношения

$$(4.1.6) \quad f^{\mathbf{0}}_{\mathbf{0}} = f^{kr} B_{C\cdot\mathbf{k0}}^{\mathbf{p}} B_{C\cdot\mathbf{pr}}^{\mathbf{0}} = f^{\mathbf{0}r} B_{C\cdot\mathbf{00}}^{\mathbf{0}} B_{C\cdot\mathbf{0}r}^{\mathbf{0}} + f^{\mathbf{1}r} B_{C\cdot\mathbf{10}}^{\mathbf{1}} B_{C\cdot\mathbf{1}r}^{\mathbf{0}} = f^{\mathbf{00}} - f^{\mathbf{11}}$$

$$(4.1.7) \quad f^{\mathbf{1}}_{\mathbf{0}} = f^{kr} B_{C\cdot\mathbf{k0}}^{\mathbf{p}} B_{C\cdot\mathbf{pr}}^{\mathbf{1}} = f^{\mathbf{0}r} B_{C\cdot\mathbf{00}}^{\mathbf{0}} B_{C\cdot\mathbf{0}r}^{\mathbf{1}} + f^{\mathbf{1}r} B_{C\cdot\mathbf{10}}^{\mathbf{1}} B_{C\cdot\mathbf{1}r}^{\mathbf{1}} = f^{\mathbf{01}} + f^{\mathbf{10}}$$

$$(4.1.8) \quad f^{\mathbf{0}}_{\mathbf{1}} = f^{kr} B_{C\cdot\mathbf{k1}}^{\mathbf{p}} B_{C\cdot\mathbf{pr}}^{\mathbf{0}} = f^{\mathbf{0}r} B_{C\cdot\mathbf{01}}^{\mathbf{1}} B_{C\cdot\mathbf{1}r}^{\mathbf{0}} + f^{\mathbf{1}r} B_{C\cdot\mathbf{11}}^{\mathbf{0}} B_{C\cdot\mathbf{0}r}^{\mathbf{0}} = -f^{\mathbf{01}} - f^{\mathbf{10}}$$

$$(4.1.9) \quad f^{\mathbf{1}}_{\mathbf{1}} = f^{kr} B_{C\cdot\mathbf{k1}}^{\mathbf{p}} B_{C\cdot\mathbf{pr}}^{\mathbf{1}} = f^{\mathbf{0}r} B_{C\cdot\mathbf{01}}^{\mathbf{1}} B_{C\cdot\mathbf{1}r}^{\mathbf{1}} + f^{\mathbf{1}r} B_{C\cdot\mathbf{11}}^{\mathbf{0}} B_{C\cdot\mathbf{0}r}^{\mathbf{1}} = f^{\mathbf{00}} - f^{\mathbf{11}}$$

Из равенств (4.1.6) и (4.1.9) следует (4.1.4). Из равенств (4.1.7) и (4.1.8) следует (4.1.5). $\square$

### 4.2. Алгебра кватернионов

В этой статье я рассматриваю множество кватернионных алгебр, определённых в [13].

**Определение 4.2.1.** Пусть $F$ - поле. Расширение $F(i, j, k)$ поля $F$ называется **алгеброй** $E(F, a, b)$ **кватернионов над полем** $F$[4.1], если произведение в алгебре $E$ определено согласно правилам

(4.2.1)

|   | $i$ | $j$ | $k$ |
|---|---|---|---|
| $i$ | $a$ | $k$ | $aj$ |
| $j$ | $-k$ | $b$ | $-bi$ |
| $k$ | $-aj$ | $bi$ | $-ab$ |

где $a, b \in F$, $ab \neq 0$.

Элементы алгебры $E(F, a, b)$ имеют вид

$$x = x^0 + x^1 i + x^2 j + x^3 k$$

где $x^i \in F$, $i = 0, 1, 2, 3$. Кватернион

$$\overline{x} = x^0 - x^1 i - x^2 j - x^3 k$$

называется сопряжённым кватерниону $x$. Мы определим **норму кватерниона** $x$ равенством

(4.2.2) $$|x|^2 = x\overline{x} = (x^0)^2 - a(x^1)^2 - b(x^2)^2 + ab(x^3)^2$$

Из равенства (4.2.2) следует, что $E(F, a, b)$ является алгеброй с делением только когда $a < 0$, $b < 0$. Тогда мы можем пронормировать базис так, что $a = -1$, $b = -1$.

Мы будем обозначать символом $E(F)$ алгебру $E(F, -1, -1)$ кватернионов с делением над полем $F$. Произведение в алгебре $E(F)$ определено согласно правилам

(4.2.3)

|   | $i$ | $j$ | $k$ |
|---|---|---|---|
| $i$ | $-1$ | $k$ | $-j$ |
| $j$ | $-k$ | $-1$ | $i$ |
| $k$ | $j$ | $-i$ | $-1$ |

$\square$

В алгебре $E(F)$ норма кватерниона имеет вид

(4.2.4) $$|x|^2 = x\overline{x} = (x^0)^2 + (x^1)^2 + (x^2)^2 + (x^3)^2$$

При этом обратный элемент имеет вид

(4.2.5) $$x^{-1} = |x|^{-2}\overline{x}$$

Мы будем полагать $H = E(R, -1, -1)$.

Внутренний автоморфизм алгебры кватернионов $H$[4.2]

(4.2.6) $$p \to qpq^{-1}$$
$$q(ix + jy + kz)q^{-1} = ix' + jy' + kz'$$

---

[4.1]Я буду следовать определению из [13].
[4.2]См. [15], с. 643.

описывает вращение вектора с координатами $x$, $y$, $z$. Если $q$ записан в виде суммы скаляра и вектора

$$q = \cos\alpha + (ia + jb + kc)\sin\alpha \quad a^2 + b^2 + c^2 = 1$$

то (4.2.6) описывает вращение вектора $(x, y, z)$ вокруг вектора $(a, b, c)$ на угол $2\alpha$.

### 4.3. Линейная функция алгебры кватернионов

**Теорема 4.3.1.** *Положим*

(4.3.1) $$\overline{e}_0 = 1 \quad \overline{e}_1 = i \quad \overline{e}_2 = j \quad \overline{e}_3 = k$$

*базис алгебры кватернионов $H$. Тогда в базисе* (4.3.1) *структурные константы имеют вид*

$$\begin{aligned}
B_{00}^0 &= 1 & B_{01}^1 &= 1 & B_{02}^2 &= 1 & B_{03}^3 &= 1 \\
B_{10}^1 &= 1 & B_{11}^0 &= -1 & B_{12}^3 &= 1 & B_{13}^2 &= -1 \\
B_{20}^2 &= 1 & B_{21}^3 &= -1 & B_{22}^0 &= -1 & B_{23}^1 &= 1 \\
B_{30}^3 &= 1 & B_{31}^2 &= 1 & B_{32}^1 &= -1 & B_{33}^0 &= -1
\end{aligned}$$

Доказательство. Значение структурных констант следует из таблицы умножения (4.2.3). □

Так как вычисления в этом разделе занимают много места, я собрал в одном месте ссылки на теоремы в этом разделе.

- **Теорема 4.3.2:** определение координат линейного отображения алгебры кватернионов $H$ через стандартные компоненты этого отображения.
- **Равенство** (4.3.22): матричная форма зависимости координат линейного отображения алгебры кватернионов $H$ от стандартных компонент этого отображения.
- **Равенство** (4.3.23): матричная форма зависимости стандартных компонент линейного отображения алгебры кватернионов $H$ от координат этого отображения.
- **Theorem 4.3.4:** зависимость стандартных компонент линейного отображения алгебры кватернионов $H$ от координат этого отображения.

**Теорема 4.3.2.** *Стандартные компоненты линейной функции алгебры кватернионов $H$ относительно базиса* (4.3.1) *и координаты соответствующего линейного преобразования удовлетворяют соотношениям*

(4.3.2) $$\begin{cases} f_0^0 = f^{00} - f^{11} - f^{22} - f^{33} \\ f_1^1 = f^{00} - f^{11} + f^{22} + f^{33} \\ f_2^2 = f^{00} + f^{11} - f^{22} + f^{33} \\ f_3^3 = f^{00} + f^{11} + f^{22} - f^{33} \end{cases}$$

$$(4.3.3) \quad \begin{cases} f_0^1 = f^{01} + f^{10} + f^{23} - f^{32} \\ f_1^0 = -f^{01} - f^{10} + f^{23} - f^{32} \\ f_2^3 = -f^{01} + f^{10} - f^{23} - f^{32} \\ f_3^2 = f^{01} - f^{10} - f^{23} - f^{32} \end{cases}$$

$$(4.3.4) \quad \begin{cases} f_0^2 = f^{02} - f^{13} + f^{20} + f^{31} \\ f_1^3 = f^{02} - f^{13} - f^{20} - f^{31} \\ f_2^0 = -f^{02} - f^{13} - f^{20} + f^{31} \\ f_3^1 = -f^{02} - f^{13} + f^{20} - f^{31} \end{cases}$$

$$(4.3.5) \quad \begin{cases} f_0^3 = f^{03} + f^{12} - f^{21} + f^{30} \\ f_1^2 = -f^{03} - f^{12} - f^{21} + f^{30} \\ f_2^1 = f^{03} - f^{12} - f^{21} - f^{30} \\ f_3^0 = -f^{03} + f^{12} - f^{21} - f^{30} \end{cases}$$

Доказательство. Пользуясь равенством (3.7.11) получаем соотношения

$$(4.3.6) \quad \begin{aligned} f_0^0 &= f^{kr} B_{k0}^p B_{pr}^0 \\ &= f^{00} B_{00}^0 B_{00}^0 + f^{11} B_{10}^1 B_{11}^0 + f^{22} B_{20}^2 B_{22}^0 + f^{33} B_{30}^3 B_{33}^0 \\ &= f^{00} - f^{11} - f^{22} - f^{33} \end{aligned}$$

$$(4.3.7) \quad \begin{aligned} f_0^1 &= f^{kr} B_{k0}^p B_{pr}^1 \\ &= f^{01} B_{00}^0 B_{01}^1 + f^{10} B_{10}^1 B_{10}^1 + f^{23} B_{20}^2 B_{23}^1 + f^{32} B_{30}^3 B_{32}^1 \\ &= f^{01} + f^{10} + f^{23} - f^{32} \end{aligned}$$

$$(4.3.8) \quad \begin{aligned} f_0^2 &= f^{kr} B_{k0}^p B_{pr}^2 \\ &= f^{02} B_{00}^0 B_{02}^2 + f^{13} B_{10}^1 B_{13}^2 + f^{20} B_{20}^2 B_{20}^2 + f^{31} B_{30}^3 B_{31}^2 \\ &= f^{02} - f^{13} + f^{20} + f^{31} \end{aligned}$$

$$(4.3.9) \quad \begin{aligned} f_0^3 &= f^{kr} B_{k0}^p B_{pr}^3 \\ &= f^{03} B_{00}^0 B_{03}^3 + f^{12} B_{10}^1 B_{12}^3 + f^{21} B_{20}^2 B_{21}^3 + f^{30} B_{30}^3 B_{30}^3 \\ &= f^{03} + f^{12} - f^{21} + f^{30} \end{aligned}$$

$$(4.3.10) \quad \begin{aligned} f_1^0 &= f^{kr} B_{k1}^p B_{pr}^0 \\ &= f^{01} B_{01}^1 B_{11}^0 + f^{10} B_{11}^0 B_{00}^0 + f^{23} B_{21}^3 B_{33}^0 + f^{32} B_{31}^2 B_{22}^0 \\ &= -f^{01} - f^{10} + f^{23} - f^{32} \end{aligned}$$

$$f_1^1 = f^{kr} B_{k1}^p B_{pr}^1$$
(4.3.11)
$$= f^{00} B_{01}^1 B_{10}^1 + f^{11} B_{11}^0 B_{01}^1 + f^{22} B_{21}^3 B_{32}^1 + f^{33} B_{31}^2 B_{23}^1$$
$$= f^{00} - f^{11} + f^{22} + f^{33}$$

$$f_1^2 = f^{kr} B_{k1}^p B_{pr}^2$$
(4.3.12)
$$= f^{03} B_{01}^1 B_{13}^2 + f^{12} B_{11}^0 B_{02}^2 + f^{21} B_{21}^3 B_{31}^2 + f^{30} B_{31}^2 B_{20}^2$$
$$= -f^{03} - f^{12} - f^{21} + f^{30}$$

$$f_1^3 = f^{kr} B_{k1}^p B_{pr}^3$$
(4.3.13)
$$= f^{02} B_{01}^1 B_{12}^3 + f^{13} B_{11}^0 B_{03}^3 + f^{20} B_{21}^3 B_{30}^3 + f^{31} B_{31}^2 B_{21}^3$$
$$= f^{02} - f^{13} - f^{20} - f^{31}$$

$$f_2^0 = f^{kr} B_{k2}^p B_{pr}^0$$
(4.3.14)
$$= f^{02} B_{02}^2 B_{22}^0 + f^{13} B_{12}^3 B_{33}^0 + f^{20} B_{22}^0 B_{00}^0 + f^{31} B_{32}^1 B_{11}^0$$
$$= -f^{02} - f^{13} - f^{20} + f^{31}$$

$$f_2^1 = f^{kr} B_{k2}^p B_{pr}^1$$
(4.3.15)
$$= f^{03} B_{02}^2 B_{23}^1 + f^{12} B_{12}^3 B_{32}^1 + f^{21} B_{22}^0 B_{01}^1 + f^{30} B_{32}^1 B_{10}^1$$
$$= f^{03} - f^{12} - f^{21} - f^{30}$$

$$f_2^2 = f^{kr} B_{k2}^p B_{pr}^2$$
(4.3.16)
$$= f^{00} B_{02}^2 B_{20}^2 + f^{11} B_{12}^3 B_{31}^2 + f^{22} B_{22}^0 B_{02}^2 + f^{33} B_{32}^1 B_{13}^2$$
$$= f^{00} + f^{11} - f^{22} + f^{33}$$

$$f_2^3 = f^{kr} B_{k2}^p B_{pr}^3$$
(4.3.17)
$$= f^{01} B_{02}^2 B_{21}^3 + f^{10} B_{12}^3 B_{30}^3 + f^{23} B_{22}^0 B_{03}^3 + f^{32} B_{32}^1 B_{12}^3$$
$$= -f^{01} + f^{10} - f^{23} - f^{32}$$

$$f_3^0 = f^{kr} B_{k3}^p B_{pr}^0$$
(4.3.18)
$$= f^{03} B_{03}^3 B_{33}^0 + f^{12} B_{13}^2 B_{22}^0 + f^{21} B_{23}^1 B_{11}^0 + f^{30} B_{33}^0 B_{00}^0$$
$$= -f^{03} + f^{12} - f^{21} - f^{30}$$

$$f_3^1 = f^{kr} B_{k3}^p B_{pr}^1$$
(4.3.19)
$$= f^{02} B_{03}^3 B_{32}^1 + f^{13} B_{13}^2 B_{23}^1 + f^{20} B_{23}^1 B_{10}^1 + f^{31} B_{33}^0 B_{01}^1$$
$$= -f^{02} - f^{13} + f^{20} - f^{31}$$

$$f^2_3 = f^{kr}B^p_{k3}B^2_{pr}$$
(4.3.20)
$$= f^{01}B^3_{03}B^2_{31} + f^{10}B^2_{13}B^2_{20} + f^{23}B^1_{23}B^2_{13} + f^{32}B^0_{33}B^2_{02}$$
$$= f^{01} - f^{10} - f^{23} - f^{32}$$

$$f^3_3 = f^{kr}B^p_{k3}B^3_{pr}$$
(4.3.21)
$$= f^{00}B^3_{03}B^3_{30} + f^{11}B^2_{13}B^3_{21} + f^{22}B^1_{23}B^3_{12} + f^{33}B^0_{33}B^3_{03}$$
$$= f^{00} + f^{11} + f^{22} - f^{33}$$

Уравнения (4.3.6), (4.3.11), (4.3.16), (4.3.21) формируют систему линейных уравнений (4.3.2).

Уравнения (4.3.7), (4.3.10), (4.3.17), (4.3.20) формируют систему линейных уравнений (4.3.3).

Уравнения (4.3.8), (4.3.13), (4.3.14), (4.3.19) формируют систему линейных уравнений (4.3.4).

Уравнения (4.3.9), (4.3.12), (4.3.15), (4.3.18) формируют систему линейных уравнений (4.3.5). $\square$

**Теорема 4.3.3.** *Рассмотрим алгебру кватернионов $H$ с базисом* (4.3.1). *Стандартные компоненты аддитивной функции над полем $F$ и координаты этой функции над полем $F$ удовлетворяют соотношениям*

(4.3.22)
$$\begin{pmatrix} f^0_0 & f^0_1 & f^0_2 & f^0_3 \\ f^1_1 & -f^1_0 & f^1_3 & -f^1_2 \\ f^2_2 & -f^2_3 & -f^2_0 & f^2_1 \\ f^3_3 & f^3_2 & -f^3_1 & -f^3_0 \end{pmatrix}$$
$$= \begin{pmatrix} 1 & -1 & -1 & -1 \\ 1 & -1 & 1 & 1 \\ 1 & 1 & -1 & 1 \\ 1 & 1 & 1 & -1 \end{pmatrix} \begin{pmatrix} f^{00} & -f^{01} & -f^{02} & -f^{03} \\ f^{11} & f^{10} & f^{13} & -f^{12} \\ f^{22} & -f^{23} & f^{20} & f^{21} \\ f^{33} & f^{32} & -f^{31} & f^{30} \end{pmatrix}$$

(4.3.23)
$$\begin{pmatrix} f^{00} & -f^{01} & -f^{02} & -f^{03} \\ f^{11} & f^{10} & f^{13} & -f^{12} \\ f^{22} & -f^{23} & f^{20} & f^{21} \\ f^{33} & f^{32} & -f^{31} & f^{30} \end{pmatrix}$$
$$= \frac{1}{4} \begin{pmatrix} 1 & 1 & 1 & 1 \\ -1 & -1 & 1 & 1 \\ -1 & 1 & -1 & 1 \\ -1 & 1 & 1 & -1 \end{pmatrix} \begin{pmatrix} f^0_0 & f^0_1 & f^0_2 & f^0_3 \\ f^1_1 & -f^1_0 & f^1_3 & -f^1_2 \\ f^2_2 & -f^2_3 & -f^2_0 & f^2_1 \\ f^3_3 & f^3_2 & -f^3_1 & -f^3_0 \end{pmatrix}$$

*где*

$$\begin{pmatrix} 1 & -1 & -1 & -1 \\ 1 & -1 & 1 & 1 \\ 1 & 1 & -1 & 1 \\ 1 & 1 & 1 & -1 \end{pmatrix}^{-1} = \frac{1}{4} \begin{pmatrix} 1 & 1 & 1 & 1 \\ -1 & -1 & 1 & 1 \\ -1 & 1 & -1 & 1 \\ -1 & 1 & 1 & -1 \end{pmatrix}$$

Доказательство. Запишем систему линейных уравнений (4.3.2) в виде произведения матриц

$$(4.3.24) \quad \begin{pmatrix} f_0^0 \\ f_1^1 \\ f_2^2 \\ f_3^3 \end{pmatrix} = \begin{pmatrix} 1 & -1 & -1 & -1 \\ 1 & -1 & 1 & 1 \\ 1 & 1 & -1 & 1 \\ 1 & 1 & 1 & -1 \end{pmatrix} \begin{pmatrix} f^{00} \\ f^{11} \\ f^{22} \\ f^{33} \end{pmatrix}$$

Запишем систему линейных уравнений (4.3.3) в виде произведения матриц

$$(4.3.25) \quad \begin{pmatrix} f_0^1 \\ f_1^0 \\ f_2^3 \\ f_3^2 \end{pmatrix} = \begin{pmatrix} 1 & 1 & 1 & -1 \\ -1 & -1 & 1 & -1 \\ -1 & 1 & -1 & -1 \\ 1 & -1 & -1 & -1 \end{pmatrix} \begin{pmatrix} f^{01} \\ f^{10} \\ f^{23} \\ f^{32} \end{pmatrix}$$

Из равенства (4.3.25) следует

$$\begin{pmatrix} -f_0^1 \\ f_1^0 \\ f_2^3 \\ -f_3^2 \end{pmatrix} = \begin{pmatrix} -1 & -1 & -1 & 1 \\ -1 & -1 & 1 & -1 \\ -1 & 1 & -1 & -1 \\ -1 & 1 & 1 & 1 \end{pmatrix} \begin{pmatrix} f^{01} \\ f^{10} \\ f^{23} \\ f^{32} \end{pmatrix}$$

$$= \begin{pmatrix} 1 & -1 & 1 & 1 \\ 1 & -1 & -1 & -1 \\ 1 & 1 & 1 & -1 \\ 1 & 1 & -1 & 1 \end{pmatrix} \begin{pmatrix} -f^{01} \\ f^{10} \\ -f^{23} \\ f^{32} \end{pmatrix}$$

$$(4.3.26) \quad \begin{pmatrix} f_1^0 \\ -f_0^1 \\ -f_3^2 \\ f_2^3 \end{pmatrix} = \begin{pmatrix} 1 & -1 & -1 & -1 \\ 1 & -1 & 1 & 1 \\ 1 & 1 & -1 & 1 \\ 1 & 1 & 1 & -1 \end{pmatrix} \begin{pmatrix} -f^{01} \\ f^{10} \\ -f^{23} \\ f^{32} \end{pmatrix}$$

Запишем систему линейных уравнений (4.3.4) в виде произведения матриц

$$(4.3.27) \quad \begin{pmatrix} f_0^2 \\ f_1^3 \\ f_2^0 \\ f_3^1 \end{pmatrix} = \begin{pmatrix} 1 & -1 & 1 & 1 \\ 1 & -1 & -1 & -1 \\ -1 & -1 & -1 & 1 \\ -1 & -1 & 1 & -1 \end{pmatrix} \begin{pmatrix} f^{02} \\ f^{13} \\ f^{20} \\ f^{31} \end{pmatrix}$$

Из равенства (4.3.27) следует

$$\begin{pmatrix} -f_0^2 \\ -f_1^3 \\ f_2^0 \\ f_3^1 \end{pmatrix} = \begin{pmatrix} -1 & 1 & -1 & -1 \\ -1 & 1 & 1 & 1 \\ -1 & -1 & -1 & 1 \\ -1 & -1 & 1 & -1 \end{pmatrix} \begin{pmatrix} f^{02} \\ f^{13} \\ f^{20} \\ f^{31} \end{pmatrix}$$

$$= \begin{pmatrix} 1 & 1 & -1 & 1 \\ 1 & 1 & 1 & -1 \\ 1 & -1 & -1 & -1 \\ 1 & -1 & 1 & 1 \end{pmatrix} \begin{pmatrix} -f^{02} \\ f^{13} \\ f^{20} \\ -f^{31} \end{pmatrix}$$

(4.3.28) $$\begin{pmatrix} f_2^0 \\ f_3^1 \\ -f_0^2 \\ -f_1^3 \end{pmatrix} = \begin{pmatrix} 1 & -1 & -1 & -1 \\ 1 & -1 & 1 & 1 \\ 1 & 1 & -1 & 1 \\ 1 & 1 & 1 & -1 \end{pmatrix} \begin{pmatrix} -f^{02} \\ f^{13} \\ f^{20} \\ -f^{31} \end{pmatrix}$$

Запишем систему линейных уравнений (4.3.5) в виде произведения матриц

(4.3.29) $$\begin{pmatrix} f_0^3 \\ f_1^2 \\ f_2^1 \\ f_3^0 \end{pmatrix} = \begin{pmatrix} 1 & 1 & -1 & 1 \\ -1 & -1 & -1 & 1 \\ 1 & -1 & -1 & -1 \\ -1 & 1 & -1 & -1 \end{pmatrix} \begin{pmatrix} f^{03} \\ f^{12} \\ f^{21} \\ f^{30} \end{pmatrix}$$

Из равенства (4.3.29) следует

$$\begin{pmatrix} -f_0^3 \\ f_1^2 \\ -f_2^1 \\ f_3^0 \end{pmatrix} = \begin{pmatrix} -1 & -1 & 1 & -1 \\ -1 & -1 & -1 & 1 \\ -1 & 1 & 1 & 1 \\ -1 & 1 & -1 & -1 \end{pmatrix} \begin{pmatrix} f^{03} \\ f^{12} \\ f^{21} \\ f^{30} \end{pmatrix}$$

$$= \begin{pmatrix} 1 & 1 & 1 & -1 \\ 1 & 1 & -1 & 1 \\ 1 & -1 & 1 & 1 \\ 1 & -1 & -1 & -1 \end{pmatrix} \begin{pmatrix} -f^{03} \\ -f^{12} \\ f^{21} \\ f^{30} \end{pmatrix}$$

(4.3.30) $$\begin{pmatrix} f_3^0 \\ -f_2^1 \\ f_1^2 \\ -f_0^3 \end{pmatrix} = \begin{pmatrix} 1 & -1 & -1 & -1 \\ 1 & -1 & 1 & 1 \\ 1 & 1 & -1 & 1 \\ 1 & 1 & 1 & -1 \end{pmatrix} \begin{pmatrix} -f^{03} \\ -f^{12} \\ f^{21} \\ f^{30} \end{pmatrix}$$

Мы объединяем равенства (4.3.24), (4.3.26), (4.3.28), (4.3.30) в равенстве (4.3.22).
$\square$

**Теорема 4.3.4.** *Стандартные компоненты линейной функции алгебры кватернионов H относительно базиса* (4.3.1) *и координаты соответствующего* *линейного преобразования удовлетворяют соотношениям*

$$(4.3.31) \quad \begin{cases} 4f^{00} = f_0^0 + f_1^1 + f_2^2 + f_3^3 \\ 4f^{11} = -f_0^0 - f_1^1 + f_2^2 + f_3^3 \\ 4f^{22} = -f_0^0 + f_1^1 - f_2^2 + f_3^3 \\ 4f^{33} = -f_0^0 + f_1^1 + f_2^2 - f_3^3 \end{cases}$$

$$(4.3.32) \quad \begin{cases} 4f^{10} = -f_1^0 + f_0^1 - f_3^2 + f_2^3 \\ 4f^{01} = -f_1^0 + f_0^1 + f_3^2 - f_2^3 \\ 4f^{32} = -f_1^0 - f_0^1 - f_3^2 - f_2^3 \\ 4f^{23} = f_1^0 + f_0^1 - f_3^2 - f_2^3 \end{cases}$$

$$(4.3.33) \quad \begin{cases} 4f^{20} = -f_2^0 + f_3^1 + f_0^2 - f_1^3 \\ 4f^{31} = f_2^0 - f_3^1 + f_0^2 - f_1^3 \\ 4f^{02} = -f_2^0 - f_3^1 + f_0^2 + f_1^3 \\ 4f^{13} = -f_2^0 - f_3^1 - f_0^2 - f_1^3 \end{cases}$$

$$(4.3.34) \quad \begin{cases} 4f^{30} = -f_3^0 - f_2^1 + f_1^2 + f_0^3 \\ 4f^{21} = -f_3^0 - f_2^1 - f_1^2 - f_0^3 \\ 4f^{12} = f_3^0 - f_2^1 - f_1^2 + f_0^3 \\ 4f^{03} = -f_3^0 + f_2^1 - f_1^2 + f_0^3 \end{cases}$$

Доказательство. Системы линейных уравнений (4.3.31), (4.3.32), (4.3.33), (4.3.34) получены в результате перемножения матриц в равенстве (4.3.23). □

### 4.4. Алгебра октонионов

**Определение 4.4.1.** Алгебра $O$ называется **алгеброй октонионов** если алгебра имеет базис

$$(4.4.1) \quad \begin{array}{cccc} \overline{e}_0 = 1 & \overline{e}_1 = i & \overline{e}_2 = j & \overline{e}_3 = k \\ \overline{e}_4 = -l & \overline{e}_5 = il & \overline{e}_6 = jl & \overline{e}_7 = kl \end{array}$$

и произведение в алгебре $O$ определено согласно правилам

(4.4.2)

| | $\overline{e}_1$ | $\overline{e}_2$ | $\overline{e}_3$ | $\overline{e}_4$ | $\overline{e}_5$ | $\overline{e}_6$ | $\overline{e}_7$ |
|---|---|---|---|---|---|---|---|
| $\overline{e}_1$ | $-\overline{e}_0$ | $\overline{e}_3$ | $-\overline{e}_2$ | $\overline{e}_5$ | $-\overline{e}_4$ | $-\overline{e}_7$ | $\overline{e}_6$ |
| $\overline{e}_2$ | $-\overline{e}_3$ | $-\overline{e}_0$ | $\overline{e}_1$ | $\overline{e}_6$ | $\overline{e}_7$ | $-\overline{e}_4$ | $-\overline{e}_5$ |
| $\overline{e}_3$ | $\overline{e}_2$ | $-\overline{e}_1$ | $-\overline{e}_0$ | $\overline{e}_7$ | $-\overline{e}_6$ | $\overline{e}_5$ | $-\overline{e}_4$ |
| $\overline{e}_4$ | $-\overline{e}_5$ | $-\overline{e}_6$ | $-\overline{e}_7$ | $-\overline{e}_0$ | $\overline{e}_1$ | $\overline{e}_2$ | $\overline{e}_3$ |
| $\overline{e}_5$ | $\overline{e}_4$ | $-\overline{e}_7$ | $\overline{e}_6$ | $-\overline{e}_1$ | $-\overline{e}_0$ | $-\overline{e}_3$ | $\overline{e}_2$ |
| $\overline{e}_6$ | $\overline{e}_7$ | $\overline{e}_4$ | $-\overline{e}_5$ | $-\overline{e}_2$ | $\overline{e}_3$ | $-\overline{e}_0$ | $-\overline{e}_1$ |
| $\overline{e}_7$ | $-\overline{e}_6$ | $\overline{e}_5$ | $\overline{e}_4$ | $-\overline{e}_3$ | $-\overline{e}_2$ | $\overline{e}_1$ | $-\overline{e}_0$ |

□

**Теорема 4.4.2.** *Структурные константы алгебры октонионов $O$ относительно базиса* (4.4.1) *имеют вид*

$$
\begin{array}{llll}
B^0_{00}= 1 & B^1_{01}= 1 & B^2_{02}= 1 & B^3_{03}= 1 \\
B^4_{04}= 1 & B^5_{05}= 1 & B^6_{06}= 1 & B^7_{07}= 1 \\
B^1_{10}= 1 & B^0_{11}= -1 & B^3_{12}= 1 & B^2_{13}= -1 \\
B^5_{14}= 1 & B^4_{15}= -1 & B^7_{16}= -1 & B^6_{17}= 1 \\
B^2_{20}= 1 & B^3_{21}= -1 & B^0_{22}= -1 & B^1_{23}= 1 \\
B^6_{24}= 1 & B^7_{25}= 1 & B^4_{26}= -1 & B^5_{27}= -1 \\
B^3_{30}= 1 & B^2_{31}= 1 & B^1_{32}= -1 & B^0_{33}= -1 \\
B^7_{34}= 1 & B^6_{35}= -1 & B^5_{36}= 1 & B^4_{37}= -1 \\
B^4_{40}= 1 & B^5_{41}= -1 & B^6_{42}= -1 & B^7_{43}= -1 \\
B^0_{44}= -1 & B^1_{45}= 1 & B^2_{46}= 1 & B^3_{47}= 1 \\
B^5_{50}= 1 & B^4_{51}= 1 & B^7_{52}= -1 & B^6_{53}= 1 \\
B^1_{54}= -1 & B^0_{55}= -1 & B^3_{56}= -1 & B^2_{57}= 1 \\
B^6_{60}= 1 & B^7_{61}= 1 & B^4_{62}= 1 & B^5_{63}= -1 \\
B^2_{64}= -1 & B^3_{65}= 1 & B^0_{66}= -1 & B^1_{67}= -1 \\
B^7_{70}= 1 & B^6_{71}= -1 & B^5_{72}= 1 & B^4_{73}= 1 \\
B^3_{74}= -1 & B^2_{75}= -1 & B^1_{76}= 1 & B^0_{77}= -1
\end{array}
$$

Доказательство. Значение структурных констант следует из таблицы умножения (4.4.2). □

### 4.5. Линейная функция алгебры октонионов

Так как вычисления в этом разделе занимают много места, я собрал в одном месте ссылки на теоремы в этом разделе.

**Теорема 4.5.1:** определение координат линейного отображения алгебры октонионов $O$ через стандартные компоненты этого отображения.

**Равенство** (4.5.73): матричная форма зависимости координат линейного отображения алгебры октонионов $O$ от стандартных компонент этого отображения.

**Равенство** (4.5.74): матричная форма зависимости стандартных компонент линейного отображения алгебры октонионов $O$ от координат этого отображения.

**Theorem 4.5.3:** зависимость стандартных компонент линейного отображения алгебры октонионов $O$ от координат этого отображения.

**Теорема 4.5.1.** *Стандартные компоненты линейной функции алгебры октонионов $O$ относительно базиса* (4.4.1) *и координаты соответствующего* *линейного преобразования удовлетворяют соотношениям*

$$(4.5.1) \quad \begin{cases} f_0^0 = f^{00} - f^{11} - f^{22} - f^{33} - f^{44} - f^{55} - f^{66} - f^{77} \\ f_1^1 = f^{00} - f^{11} + f^{22} + f^{33} + f^{44} + f^{55} + f^{66} + f^{77} \\ f_2^2 = f^{00} + f^{11} - f^{22} + f^{33} + f^{44} + f^{55} + f^{66} + f^{77} \\ f_3^3 = f^{00} + f^{11} + f^{22} - f^{33} + f^{44} + f^{55} + f^{66} + f^{77} \\ f_4^4 = f^{00} + f^{11} + f^{22} + f^{33} - f^{44} + f^{55} + f^{66} + f^{77} \\ f_5^5 = f^{00} + f^{11} + f^{22} + f^{33} + f^{44} - f^{55} + f^{66} + f^{77} \\ f_6^6 = f^{00} + f^{11} + f^{22} + f^{33} + f^{44} + f^{55} - f^{66} + f^{77} \\ f_7^7 = f^{00} + f^{11} + f^{22} + f^{33} + f^{44} + f^{55} + f^{66} - f^{77} \end{cases}$$

$$(4.5.2) \quad \begin{cases} f_0^1 = \phantom{-}f^{01} + f^{10} + f^{23} - f^{32} + f^{45} - f^{54} - f^{67} + f^{76} \\ f_1^0 = -f^{01} - f^{10} + f^{23} - f^{32} + f^{45} - f^{54} - f^{67} + f^{76} \\ f_2^3 = -f^{01} + f^{10} - f^{23} - f^{32} - f^{45} + f^{54} + f^{67} - f^{76} \\ f_3^2 = \phantom{-}f^{01} - f^{10} - f^{23} - f^{32} + f^{45} - f^{54} - f^{67} + f^{76} \\ f_4^5 = -f^{01} + f^{10} - f^{23} + f^{32} - f^{45} - f^{54} + f^{67} - f^{76} \\ f_5^4 = \phantom{-}f^{01} - f^{10} + f^{23} - f^{32} - f^{45} - f^{54} - f^{67} + f^{76} \\ f_6^7 = \phantom{-}f^{01} - f^{10} + f^{23} - f^{32} + f^{45} - f^{54} - f^{67} - f^{76} \\ f_7^6 = -f^{01} + f^{10} - f^{23} + f^{32} - f^{45} + f^{54} - f^{67} - f^{76} \end{cases}$$

$$(4.5.3) \quad \begin{cases} f_0^2 = \phantom{-}f^{02} - f^{13} + f^{20} + f^{31} + f^{46} + f^{57} - f^{64} - f^{75} \\ f_1^3 = \phantom{-}f^{02} - f^{13} - f^{20} - f^{31} + f^{46} + f^{57} - f^{64} - f^{75} \\ f_2^0 = -f^{02} - f^{13} - f^{20} + f^{31} + f^{46} + f^{57} - f^{64} - f^{75} \\ f_3^1 = -f^{02} - f^{13} + f^{20} - f^{31} - f^{46} - f^{57} + f^{64} + f^{75} \\ f_4^6 = -f^{02} + f^{13} + f^{20} - f^{31} - f^{46} - f^{57} - f^{64} + f^{75} \\ f_5^7 = -f^{02} + f^{13} + f^{20} + f^{31} - f^{46} - f^{57} + f^{64} - f^{75} \\ f_6^4 = \phantom{-}f^{02} - f^{13} - f^{20} + f^{31} - f^{46} + f^{57} - f^{64} - f^{75} \\ f_7^5 = \phantom{-}f^{02} - f^{13} - f^{20} + f^{31} + f^{46} - f^{57} - f^{64} - f^{75} \end{cases}$$

$$(4.5.4)\begin{cases} f_0^3 = \phantom{-}f^{03} + f^{12} - f^{21} + f^{30} + f^{47} - f^{56} + f^{65} - f^{74} \\ f_1^2 = -f^{03} - f^{12} - f^{21} + f^{30} - f^{47} + f^{56} - f^{65} + f^{74} \\ f_2^1 = \phantom{-}f^{03} - f^{12} - f^{21} - f^{30} + f^{47} - f^{56} + f^{65} - f^{74} \\ f_3^0 = -f^{03} + f^{12} - f^{21} - f^{30} + f^{47} - f^{56} + f^{65} - f^{74} \\ f_4^7 = -f^{03} - f^{12} + f^{21} + f^{30} - f^{47} + f^{56} - f^{65} - f^{74} \\ f_5^6 = \phantom{-}f^{03} + f^{12} - f^{21} - f^{30} + f^{47} - f^{56} - f^{65} - f^{74} \\ f_6^5 = -f^{03} - f^{12} + f^{21} + f^{30} - f^{47} - f^{56} - f^{65} + f^{74} \\ f_7^4 = \phantom{-}f^{03} + f^{12} - f^{21} - f^{30} - f^{47} - f^{56} + f^{65} - f^{74} \end{cases}$$

$$(4.5.5)\begin{cases} f_0^4 = \phantom{-}f^{04} - f^{15} - f^{26} - f^{37} + f^{40} + f^{51} + f^{62} + f^{73} \\ f_1^5 = \phantom{-}f^{04} - f^{15} - f^{26} - f^{37} - f^{40} - f^{51} + f^{62} + f^{73} \\ f_2^6 = \phantom{-}f^{04} - f^{15} - f^{26} - f^{37} - f^{40} + f^{51} - f^{62} + f^{73} \\ f_3^7 = \phantom{-}f^{04} - f^{15} - f^{26} - f^{37} - f^{40} + f^{51} + f^{62} - f^{73} \\ f_4^0 = -f^{04} - f^{15} - f^{26} - f^{37} - f^{40} + f^{51} + f^{62} + f^{73} \\ f_5^1 = -f^{04} - f^{15} + f^{26} + f^{37} + f^{40} - f^{51} - f^{62} - f^{73} \\ f_6^2 = -f^{04} + f^{15} - f^{26} + f^{37} + f^{40} - f^{51} - f^{62} - f^{73} \\ f_7^3 = -f^{04} + f^{15} + f^{26} - f^{37} + f^{40} - f^{51} - f^{62} - f^{73} \end{cases}$$

$$(4.5.6)\begin{cases} f_0^5 = \phantom{-}f^{05} + f^{14} - f^{27} + f^{36} - f^{41} + f^{50} - f^{63} + f^{72} \\ f_1^4 = -f^{05} - f^{14} + f^{27} - f^{36} - f^{41} + f^{50} + f^{63} - f^{72} \\ f_2^7 = \phantom{-}f^{05} + f^{14} - f^{27} + f^{36} - f^{41} - f^{50} - f^{63} - f^{72} \\ f_3^6 = -f^{05} - f^{14} + f^{27} - f^{36} + f^{41} + f^{50} - f^{63} - f^{72} \\ f_4^1 = \phantom{-}f^{05} - f^{14} - f^{27} + f^{36} - f^{41} - f^{50} - f^{63} + f^{72} \\ f_5^0 = -f^{05} + f^{14} - f^{27} + f^{36} - f^{41} - f^{50} - f^{63} + f^{72} \\ f_6^3 = \phantom{-}f^{05} + f^{14} - f^{27} - f^{36} - f^{41} - f^{50} - f^{63} + f^{72} \\ f_7^2 = -f^{05} - f^{14} - f^{27} - f^{36} + f^{41} + f^{50} + f^{63} - f^{72} \end{cases}$$

$$(4.5.7)\begin{cases} f_0^6 = \phantom{-}f^{06} + f^{17} + f^{24} - f^{35} - f^{42} + f^{53} + f^{60} - f^{71} \\ f_1^7 = -f^{06} - f^{17} - f^{24} + f^{35} + f^{42} - f^{53} + f^{60} - f^{71} \\ f_2^4 = -f^{06} - f^{17} - f^{24} + f^{35} - f^{42} - f^{53} - f^{60} + f^{71} \\ f_3^5 = \phantom{-}f^{06} + f^{17} + f^{24} - f^{35} - f^{42} - f^{53} - f^{60} - f^{71} \\ f_4^2 = \phantom{-}f^{06} + f^{17} - f^{24} - f^{35} - f^{42} + f^{53} - f^{60} - f^{71} \\ f_5^3 = -f^{06} - f^{17} - f^{24} - f^{35} + f^{42} - f^{53} + f^{60} + f^{71} \\ f_6^0 = -f^{06} + f^{17} + f^{24} - f^{35} - f^{42} + f^{53} - f^{60} - f^{71} \\ f_7^1 = \phantom{-}f^{06} - f^{17} + f^{24} - f^{35} - f^{42} + f^{53} - f^{60} - f^{71} \end{cases}$$

$$(4.5.8)\begin{cases} f_0^7 = \phantom{-}f^{07} - f^{16} + f^{25} + f^{34} - f^{43} - f^{52} + f^{61} + f^{70} \\ f_1^6 = \phantom{-}f^{07} - f^{16} + f^{25} + f^{34} - f^{43} - f^{52} - f^{61} - f^{70} \\ f_2^5 = -f^{07} + f^{16} - f^{25} - f^{34} + f^{43} - f^{52} - f^{61} + f^{70} \\ f_3^4 = -f^{07} + f^{16} - f^{25} - f^{34} - f^{43} + f^{52} - f^{61} + f^{70} \\ f_4^3 = \phantom{-}f^{07} - f^{16} + f^{25} - f^{34} - f^{43} - f^{52} + f^{61} - f^{70} \\ f_5^2 = \phantom{-}f^{07} - f^{16} - f^{25} + f^{34} - f^{43} - f^{52} + f^{61} - f^{70} \\ f_6^1 = -f^{07} - f^{16} - f^{25} - f^{34} + f^{43} + f^{52} - f^{61} + f^{70} \\ f_7^0 = -f^{07} - f^{16} + f^{25} + f^{34} - f^{43} - f^{52} + f^{61} - f^{70} \end{cases}$$

Доказательство. Пользуясь равенством (3.7.11) получаем соотношения

$$(4.5.9)\begin{aligned} f_0^0 &= f^{kr} B_{k0}^p B_{pr}^0 \\ &= f^{00} B_{00}^0 B_{00}^0 + f^{11} B_{10}^1 B_{11}^0 + f^{22} B_{20}^2 B_{22}^0 + f^{33} B_{30}^3 B_{33}^0 \\ &\quad + f^{44} B_{40}^4 B_{44}^0 + f^{55} B_{50}^5 B_{55}^0 + f^{66} B_{60}^6 B_{66}^0 + f^{77} B_{70}^7 B_{77}^0 \\ &= f^{00} - f^{11} - f^{22} - f^{33} - f^{44} - f^{55} - f^{66} - f^{77} \end{aligned}$$

$$(4.5.10)\begin{aligned} f_0^1 &= f^{kr} B_{k0}^p B_{pr}^1 \\ &= f^{01} B_{00}^0 B_{01}^1 + f^{10} B_{10}^1 B_{10}^1 + f^{23} B_{20}^2 B_{23}^1 + f^{32} B_{30}^3 B_{32}^1 \\ &\quad + f^{45} B_{40}^4 B_{45}^1 + f^{54} B_{50}^5 B_{54}^1 + f^{67} B_{60}^6 B_{67}^1 + f^{76} B_{70}^7 B_{76}^1 \\ &= f^{01} + f^{10} + f^{23} - f^{32} + f^{45} - f^{54} - f^{67} + f^{76} \end{aligned}$$

$$(4.5.11)\begin{aligned} f_0^2 &= f^{kr} B_{k0}^p B_{pr}^2 \\ &= f^{02} B_{00}^0 B_{02}^2 + f^{13} B_{10}^1 B_{13}^2 + f^{20} B_{20}^2 B_{20}^2 + f^{31} B_{30}^3 B_{31}^2 \\ &\quad + f^{46} B_{40}^4 B_{46}^2 + f^{57} B_{50}^5 B_{57}^2 + f^{64} B_{60}^6 B_{64}^2 + f^{75} B_{70}^7 B_{75}^2 \\ &= f^{02} - f^{13} + f^{20} + f^{31} + f^{46} + f^{57} - f^{64} - f^{75} \end{aligned}$$

$$f_0^3 = f^{kr} B_{k0}^p B_{pr}^3$$
$$= f^{03} B_{00}^0 B_{03}^3 + f^{12} B_{10}^1 B_{12}^3 + f^{21} B_{20}^2 B_{21}^3 + f^{30} B_{30}^3 B_{30}^3$$
(4.5.12)
$$+ f^{47} B_{40}^4 B_{47}^3 + f^{56} B_{50}^5 B_{56}^3 + f^{65} B_{60}^6 B_{65}^3 + f^{74} B_{70}^7 B_{74}^3$$
$$= f^{03} + f^{12} - f^{21} + f^{30} + f^{47} - f^{56} + f^{65} - f^{74}$$

$$f_0^4 = f^{kr} B_{k0}^p B_{pr}^4$$
$$= f^{04} B_{00}^0 B_{04}^4 + f^{15} B_{10}^1 B_{15}^4 + f^{26} B_{20}^2 B_{26}^4 + f^{37} B_{30}^3 B_{37}^4$$
(4.5.13)
$$+ f^{40} B_{40}^4 B_{40}^4 + f^{51} B_{50}^5 B_{51}^4 + f^{62} B_{60}^6 B_{62}^4 + f^{73} B_{70}^7 B_{73}^4$$
$$= f^{04} - f^{15} - f^{26} - f^{37} + f^{40} + f^{51} + f^{62} + f^{73}$$

$$f_0^5 = f^{kr} B_{k0}^p B_{pr}^5$$
$$= f^{05} B_{00}^0 B_{05}^5 + f^{14} B_{10}^1 B_{14}^5 + f^{27} B_{20}^2 B_{27}^5 + f^{36} B_{30}^3 B_{36}^5$$
(4.5.14)
$$+ f^{41} B_{40}^4 B_{41}^5 + f^{50} B_{50}^5 B_{50}^5 + f^{63} B_{60}^6 B_{63}^5 + f^{72} B_{70}^7 B_{72}^5$$
$$= f^{05} + f^{14} - f^{27} + f^{36} - f^{41} + f^{50} - f^{63} + f^{72}$$

$$f_0^6 = f^{kr} B_{k0}^p B_{pr}^6$$
$$= f^{06} B_{00}^0 B_{06}^6 + f^{17} B_{10}^1 B_{17}^6 + f^{24} B_{20}^2 B_{24}^6 + f^{35} B_{30}^3 B_{35}^6$$
(4.5.15)
$$+ f^{42} B_{40}^4 B_{42}^6 + f^{53} B_{50}^5 B_{53}^6 + f^{60} B_{60}^6 B_{60}^6 + f^{71} B_{70}^7 B_{71}^6$$
$$= f^{06} + f^{17} + f^{24} - f^{35} - f^{42} + f^{53} + f^{60} - f^{71}$$

$$f_0^7 = f^{kr} B_{k0}^p B_{pr}^7$$
$$= f^{07} B_{00}^0 B_{07}^7 + f^{16} B_{10}^1 B_{16}^7 + f^{25} B_{20}^2 B_{25}^7 + f^{34} B_{30}^3 B_{34}^7$$
(4.5.16)
$$+ f^{43} B_{40}^4 B_{43}^7 + f^{52} B_{50}^5 B_{52}^7 + f^{61} B_{60}^6 B_{61}^7 + f^{70} B_{70}^7 B_{70}^7$$
$$= f^{07} - f^{16} + f^{25} + f^{34} - f^{43} - f^{52} + f^{61} + f^{70}$$

$$f_1^0 = f^{kr} B_{k1}^p B_{pr}^0$$
$$= f^{01} B_{01}^1 B_{11}^0 + f^{10} B_{11}^0 B_{00}^0 + f^{23} B_{21}^3 B_{33}^0 + f^{32} B_{31}^2 B_{22}^0$$
(4.5.17)
$$+ f^{45} B_{41}^5 B_{55}^0 + f^{54} B_{51}^4 B_{44}^0 + f^{67} B_{61}^7 B_{77}^0 + f^{76} B_{71}^6 B_{66}^0$$
$$= -f^{01} - f^{10} + f^{23} - f^{32} + f^{45} - f^{54} - f^{67} + f^{76}$$

$$f_1^1 = f^{kr} B_{k1}^p B_{pr}^1$$
$$= f^{00} B_{01}^1 B_{10}^1 + f^{11} B_{11}^0 B_{01}^1 + f^{22} B_{21}^3 B_{32}^1 + f^{33} B_{31}^2 B_{23}^1$$
(4.5.18)
$$+ f^{44} B_{41}^5 B_{54}^1 + f^{55} B_{51}^4 B_{45}^1 + f^{66} B_{61}^7 B_{76}^1 + f^{77} B_{71}^6 B_{67}^1$$
$$= f^{00} - f^{11} + f^{22} + f^{33} + f^{44} + f^{55} + f^{66} + f^{77}$$

$$(4.5.19) \quad \begin{aligned} f_1^2 &= f^{kr} B_{k1}^p B_{pr}^2 \\ &= f^{03} B_{01}^1 B_{13}^2 + f^{12} B_{11}^0 B_{02}^2 + f^{21} B_{21}^3 B_{31}^2 + f^{30} B_{31}^2 B_{20}^2 \\ &\quad + f^{47} B_{41}^5 B_{57}^2 + f^{56} B_{51}^4 B_{46}^2 + f^{65} B_{61}^7 B_{75}^2 + f^{74} B_{71}^6 B_{64}^2 \\ &= -f^{03} - f^{12} - f^{21} + f^{30} - f^{47} + f^{56} - f^{65} + f^{74} \end{aligned}$$

$$(4.5.20) \quad \begin{aligned} f_1^3 &= f^{kr} B_{k1}^p B_{pr}^3 \\ &= f^{02} B_{01}^1 B_{12}^3 + f^{13} B_{11}^0 B_{03}^3 + f^{20} B_{21}^3 B_{30}^3 + f^{31} B_{31}^2 B_{21}^3 \\ &\quad + f^{46} B_{41}^5 B_{56}^3 + f^{57} B_{51}^4 B_{47}^3 + f^{64} B_{61}^7 B_{74}^3 + f^{75} B_{71}^6 B_{65}^3 \\ &= f^{02} - f^{13} - f^{20} - f^{31} + f^{46} + f^{57} - f^{64} - f^{75} \end{aligned}$$

$$(4.5.21) \quad \begin{aligned} f_1^4 &= f^{kr} B_{k1}^p B_{pr}^4 \\ &= f^{05} B_{01}^1 B_{15}^4 + f^{14} B_{11}^0 B_{04}^4 + f^{27} B_{21}^3 B_{37}^4 + f^{36} B_{31}^2 B_{26}^4 \\ &\quad + f^{41} B_{41}^5 B_{51}^4 + f^{50} B_{51}^4 B_{40}^4 + f^{63} B_{61}^7 B_{73}^4 + f^{72} B_{71}^6 B_{62}^4 \\ &= -f^{05} - f^{14} + f^{27} - f^{36} - f^{41} + f^{50} + f^{63} - f^{72} \end{aligned}$$

$$(4.5.22) \quad \begin{aligned} f_1^5 &= f^{kr} B_{k1}^p B_{pr}^5 \\ &= f^{04} B_{01}^1 B_{14}^5 + f^{15} B_{11}^0 B_{05}^5 + f^{26} B_{21}^3 B_{36}^5 + f^{37} B_{31}^2 B_{27}^5 \\ &\quad + f^{40} B_{41}^5 B_{50}^5 + f^{51} B_{51}^4 B_{41}^5 + f^{62} B_{61}^7 B_{72}^5 + f^{73} B_{71}^6 B_{63}^5 \\ &= f^{04} - f^{15} - f^{26} - f^{37} - f^{40} - f^{51} + f^{62} + f^{73} \end{aligned}$$

$$(4.5.23) \quad \begin{aligned} f_1^6 &= f^{kr} B_{k1}^p B_{pr}^6 \\ &= f^{07} B_{01}^1 B_{17}^6 + f^{16} B_{11}^0 B_{06}^6 + f^{25} B_{21}^3 B_{35}^6 + f^{34} B_{31}^2 B_{24}^6 \\ &\quad + f^{43} B_{41}^5 B_{53}^6 + f^{52} B_{51}^4 B_{42}^6 + f^{61} B_{61}^7 B_{71}^6 + f^{70} B_{71}^6 B_{60}^6 \\ &= f^{07} - f^{16} + f^{25} + f^{34} - f^{43} - f^{52} - f^{61} - f^{70} \end{aligned}$$

$$(4.5.24) \quad \begin{aligned} f_1^7 &= f^{kr} B_{k1}^p B_{pr}^7 \\ &= f^{06} B_{01}^1 B_{16}^7 + f^{17} B_{11}^0 B_{07}^7 + f^{24} B_{21}^3 B_{34}^7 + f^{35} B_{31}^2 B_{25}^7 \\ &\quad + f^{42} B_{41}^5 B_{52}^7 + f^{53} B_{51}^4 B_{43}^7 + f^{60} B_{61}^7 B_{70}^7 + f^{71} B_{71}^6 B_{61}^7 \\ &= -f^{06} - f^{17} - f^{24} + f^{35} + f^{42} - f^{53} + f^{60} - f^{71} \end{aligned}$$

$$(4.5.25) \quad \begin{aligned} f_2^0 &= f^{kr} B_{k2}^p B_{pr}^0 \\ &= f^{02} B_{02}^2 B_{22}^0 + f^{13} B_{12}^3 B_{33}^0 + f^{20} B_{22}^0 B_{00}^0 + f^{31} B_{32}^1 B_{11}^0 \\ &\quad + f^{46} B_{42}^6 B_{66}^0 + f^{57} B_{52}^7 B_{77}^0 + f^{64} B_{62}^4 B_{44}^0 + f^{75} B_{72}^5 B_{55}^0 \\ &= -f^{02} - f^{13} - f^{20} + f^{31} + f^{46} + f^{57} - f^{64} - f^{75} \end{aligned}$$

$$f_2^1 = f^{kr} B_{k2}^p B_{pr}^1$$

(4.5.26)
$$= f^{03} B_{02}^2 B_{23}^1 + f^{12} B_{12}^3 B_{32}^1 + f^{21} B_{22}^0 B_{01}^1 + f^{30} B_{32}^1 B_{10}^1$$
$$+ f^{47} B_{42}^6 B_{67}^1 + f^{56} B_{52}^7 B_{76}^1 + f^{65} B_{62}^4 B_{45}^1 + f^{74} B_{72}^5 B_{54}^1$$
$$= f^{03} - f^{12} - f^{21} - f^{30} + f^{47} - f^{56} + f^{65} - f^{74}$$

$$f_2^2 = f^{kr} B_{k2}^p B_{pr}^2$$

(4.5.27)
$$= f^{00} B_{02}^2 B_{20}^2 + f^{11} B_{12}^3 B_{31}^2 + f^{22} B_{22}^0 B_{02}^2 + f^{33} B_{32}^1 B_{13}^2$$
$$+ f^{44} B_{42}^6 B_{64}^2 + f^{55} B_{52}^7 B_{75}^2 + f^{66} B_{62}^4 B_{46}^2 + f^{77} B_{72}^5 B_{57}^2$$
$$= f^{00} + f^{11} - f^{22} + f^{33} + f^{44} + f^{55} + f^{66} + f^{77}$$

$$f_2^3 = f^{kr} B_{k2}^p B_{pr}^3$$

(4.5.28)
$$= f^{01} B_{02}^2 B_{21}^3 + f^{10} B_{12}^3 B_{30}^3 + f^{23} B_{22}^0 B_{03}^3 + f^{32} B_{32}^1 B_{12}^3$$
$$+ f^{45} B_{42}^6 B_{65}^3 + f^{54} B_{52}^7 B_{74}^3 + f^{67} B_{62}^4 B_{47}^3 + f^{76} B_{72}^5 B_{56}^3$$
$$= -f^{01} + f^{10} - f^{23} - f^{32} - f^{45} + f^{54} + f^{67} - f^{76}$$

$$f_2^4 = f^{kr} B_{k2}^p B_{pr}^4$$

(4.5.29)
$$= f^{06} B_{02}^2 B_{26}^4 + f^{17} B_{12}^3 B_{37}^4 + f^{24} B_{22}^0 B_{04}^4 + f^{35} B_{32}^1 B_{15}^4$$
$$+ f^{42} B_{42}^6 B_{62}^4 + f^{53} B_{52}^7 B_{73}^4 + f^{60} B_{62}^4 B_{40}^4 + f^{71} B_{72}^5 B_{51}^4$$
$$= -f^{06} - f^{17} - f^{24} + f^{35} - f^{42} - f^{53} + f^{60} + f^{71}$$

$$f_2^5 = f^{kr} B_{k2}^p B_{pr}^5$$

(4.5.30)
$$= f^{07} B_{02}^2 B_{27}^5 + f^{16} B_{12}^3 B_{36}^5 + f^{25} B_{22}^0 B_{05}^5 + f^{34} B_{32}^1 B_{14}^5$$
$$+ f^{43} B_{42}^6 B_{63}^5 + f^{52} B_{52}^7 B_{72}^5 + f^{61} B_{62}^4 B_{41}^5 + f^{70} B_{72}^5 B_{50}^5$$
$$= -f^{07} + f^{16} - f^{25} - f^{34} + f^{43} - f^{52} - f^{61} + f^{70}$$

$$f_2^6 = f^{kr} B_{k2}^p B_{pr}^6$$

(4.5.31)
$$= f^{04} B_{02}^2 B_{24}^6 + f^{15} B_{12}^3 B_{35}^6 + f^{26} B_{22}^0 B_{06}^6 + f^{37} B_{32}^1 B_{17}^6$$
$$+ f^{40} B_{42}^6 B_{60}^6 + f^{51} B_{52}^7 B_{71}^6 + f^{62} B_{62}^4 B_{42}^6 + f^{73} B_{72}^5 B_{53}^6$$
$$= f^{04} - f^{15} - f^{26} - f^{37} - f^{40} + f^{51} - f^{62} + f^{73}$$

$$f_2^7 = f^{kr} B_{k2}^p B_{pr}^7$$

(4.5.32)
$$= f^{05} B_{02}^2 B_{25}^7 + f^{14} B_{12}^3 B_{34}^7 + f^{27} B_{22}^0 B_{07}^7 + f^{36} B_{32}^1 B_{16}^7$$
$$+ f^{41} B_{42}^6 B_{61}^7 + f^{50} B_{52}^7 B_{70}^7 + f^{63} B_{62}^4 B_{43}^7 + f^{72} B_{72}^5 B_{52}^7$$
$$= f^{05} + f^{14} - f^{27} + f^{36} - f^{41} - f^{50} - f^{63} - f^{72}$$

$$f_3^0 = f^{kr} B_{k3}^p B_{pr}^0$$
$$= f^{03} B_{03}^3 B_{33}^0 + f^{12} B_{13}^2 B_{22}^0 + f^{21} B_{23}^1 B_{11}^0 + f^{30} B_{33}^0 B_{00}^0$$
$$+ f^{47} B_{43}^7 B_{77}^0 + f^{56} B_{53}^6 B_{66}^0 + f^{65} B_{63}^5 B_{55}^0 + f^{74} B_{73}^4 B_{44}^0$$
$$= -f^{03} + f^{12} - f^{21} - f^{30} + f^{47} - f^{56} + f^{65} - f^{74} \quad (4.5.33)$$

$$f_3^1 = f^{kr} B_{k3}^p B_{pr}^1$$
$$= f^{02} B_{03}^3 B_{32}^1 + f^{13} B_{13}^2 B_{23}^1 + f^{20} B_{23}^1 B_{10}^1 + f^{31} B_{33}^0 B_{01}^1$$
$$+ f^{46} B_{43}^7 B_{76}^1 + f^{57} B_{53}^6 B_{67}^1 + f^{64} B_{63}^5 B_{54}^1 + f^{75} B_{73}^4 B_{45}^1$$
$$= -f^{02} - f^{13} + f^{20} - f^{31} - f^{46} - f^{57} + f^{64} + f^{75} \quad (4.5.34)$$

$$f_3^2 = f^{kr} B_{k3}^p B_{pr}^2$$
$$= f^{01} B_{03}^3 B_{31}^2 + f^{10} B_{13}^2 B_{20}^2 + f^{23} B_{23}^1 B_{13}^2 + f^{32} B_{33}^0 B_{02}^2$$
$$+ f^{45} B_{43}^7 B_{75}^2 + f^{54} B_{53}^6 B_{64}^2 + f^{67} B_{63}^5 B_{57}^2 + f^{76} B_{73}^4 B_{46}^2$$
$$= f^{01} - f^{10} - f^{23} - f^{32} + f^{45} - f^{54} - f^{67} + f^{76} \quad (4.5.35)$$

$$f_3^3 = f^{kr} B_{k3}^p B_{pr}^3$$
$$= f^{00} B_{03}^3 B_{30}^3 + f^{11} B_{13}^2 B_{21}^3 + f^{22} B_{23}^1 B_{12}^3 + f^{33} B_{33}^0 B_{03}^3$$
$$+ f^{44} B_{43}^7 B_{74}^3 + f^{55} B_{53}^6 B_{65}^3 + f^{66} B_{63}^5 B_{56}^3 + f^{77} B_{73}^4 B_{47}^3$$
$$= f^{00} + f^{11} + f^{22} - f^{33} + f^{44} + f^{55} + f^{66} + f^{77} \quad (4.5.36)$$

$$f_3^4 = f^{kr} B_{k3}^p B_{pr}^4$$
$$= f^{07} B_{03}^3 B_{37}^4 + f^{16} B_{13}^2 B_{26}^4 + f^{25} B_{23}^1 B_{15}^4 + f^{34} B_{33}^0 B_{04}^4$$
$$+ f^{43} B_{43}^7 B_{73}^4 + f^{52} B_{53}^6 B_{62}^4 + f^{61} B_{63}^5 B_{51}^4 + f^{70} B_{73}^4 B_{40}^4$$
$$= -f^{07} + f^{16} - f^{25} - f^{34} - f^{43} + f^{52} - f^{61} + f^{70} \quad (4.5.37)$$

$$f_3^5 = f^{kr} B_{k3}^p B_{pr}^5$$
$$= f^{06} B_{03}^3 B_{36}^5 + f^{17} B_{13}^2 B_{27}^5 + f^{24} B_{23}^1 B_{14}^5 + f^{35} B_{33}^0 B_{05}^5$$
$$+ f^{42} B_{43}^7 B_{72}^5 + f^{53} B_{53}^6 B_{63}^5 + f^{60} B_{63}^5 B_{50}^5 + f^{71} B_{73}^4 B_{41}^5$$
$$= f^{06} + f^{17} + f^{24} - f^{35} - f^{42} - f^{53} - f^{60} - f^{71} \quad (4.5.38)$$

$$f_3^6 = f^{kr} B_{k3}^p B_{pr}^6$$
$$= f^{05} B_{03}^3 B_{35}^6 + f^{14} B_{13}^2 B_{24}^6 + f^{27} B_{23}^1 B_{17}^6 + f^{36} B_{33}^0 B_{06}^6$$
$$+ f^{41} B_{43}^7 B_{71}^6 + f^{50} B_{53}^6 B_{60}^6 + f^{63} B_{63}^5 B_{53}^6 + f^{72} B_{73}^4 B_{42}^6$$
$$= -f^{05} - f^{14} + f^{27} - f^{36} + f^{41} + f^{50} - f^{63} - f^{72} \quad (4.5.39)$$

$$(4.5.40)\quad\begin{aligned}f^7_3&=f^{kr}B^p_{k3}B^7_{pr}\\&=f^{04}B^3_{03}B^7_{34}+f^{15}B^2_{13}B^7_{25}+f^{26}B^1_{23}B^7_{16}+f^{37}B^0_{33}B^7_{07}\\&\quad+f^{40}B^7_{43}B^7_{70}+f^{51}B^6_{53}B^7_{61}+f^{62}B^5_{63}B^7_{52}+f^{73}B^4_{73}B^7_{43}\\&=f^{04}-f^{15}-f^{26}-f^{37}-f^{40}+f^{51}+f^{62}-f^{73}\end{aligned}$$

$$(4.5.41)\quad\begin{aligned}f^0_4&=f^{kr}B^p_{k4}B^0_{pr}\\&=f^{04}B^4_{04}B^0_{44}+f^{15}B^5_{14}B^0_{55}+f^{26}B^6_{24}B^0_{66}+f^{37}B^7_{34}B^0_{77}\\&\quad+f^{40}B^0_{44}B^0_{00}+f^{51}B^1_{54}B^0_{11}+f^{62}B^2_{64}B^0_{22}+f^{73}B^3_{74}B^0_{33}\\&=-f^{04}-f^{15}-f^{26}-f^{37}-f^{40}+f^{51}+f^{62}+f^{73}\end{aligned}$$

$$(4.5.42)\quad\begin{aligned}f^1_4&=f^{kr}B^p_{k4}B^1_{pr}\\&=f^{05}B^4_{04}B^1_{45}+f^{14}B^5_{14}B^1_{54}+f^{27}B^6_{24}B^1_{67}+f^{36}B^7_{34}B^1_{76}\\&\quad+f^{41}B^0_{44}B^1_{01}+f^{50}B^1_{54}B^1_{10}+f^{63}B^2_{64}B^1_{23}+f^{72}B^3_{74}B^1_{32}\\&=f^{05}-f^{14}-f^{27}+f^{36}-f^{41}-f^{50}-f^{63}+f^{72}\end{aligned}$$

$$(4.5.43)\quad\begin{aligned}f^2_4&=f^{kr}B^p_{k4}B^2_{pr}\\&=f^{06}B^4_{04}B^2_{46}+f^{17}B^5_{14}B^2_{57}+f^{24}B^6_{24}B^2_{64}+f^{35}B^7_{34}B^2_{75}\\&\quad+f^{42}B^0_{44}B^2_{02}+f^{53}B^1_{54}B^2_{13}+f^{60}B^2_{64}B^2_{20}+f^{71}B^3_{74}B^2_{31}\\&=f^{06}+f^{17}-f^{24}-f^{35}-f^{42}+f^{53}-f^{60}-f^{71}\end{aligned}$$

$$(4.5.44)\quad\begin{aligned}f^3_4&=f^{kr}B^p_{k4}B^3_{pr}\\&=f^{07}B^4_{04}B^3_{47}+f^{16}B^5_{14}B^3_{56}+f^{25}B^6_{24}B^3_{65}+f^{34}B^7_{34}B^3_{74}\\&\quad+f^{43}B^0_{44}B^3_{03}+f^{52}B^1_{54}B^3_{12}+f^{61}B^2_{64}B^3_{21}+f^{70}B^3_{74}B^3_{30}\\&=f^{07}-f^{16}+f^{25}-f^{34}-f^{43}-f^{52}+f^{61}-f^{70}\end{aligned}$$

$$(4.5.45)\quad\begin{aligned}f^4_4&=f^{kr}B^p_{k4}B^4_{pr}\\&=f^{00}B^4_{04}B^4_{40}+f^{11}B^5_{14}B^4_{51}+f^{22}B^6_{24}B^4_{62}+f^{33}B^7_{34}B^4_{73}\\&\quad+f^{44}B^0_{44}B^4_{04}+f^{55}B^1_{54}B^4_{15}+f^{66}B^2_{64}B^4_{26}+f^{77}B^3_{74}B^4_{37}\\&=f^{00}+f^{11}+f^{22}+f^{33}-f^{44}+f^{55}+f^{66}+f^{77}\end{aligned}$$

$$(4.5.46)\quad\begin{aligned}f^5_4&=f^{kr}B^p_{k4}B^5_{pr}\\&=f^{01}B^4_{04}B^5_{41}+f^{10}B^5_{14}B^5_{50}+f^{23}B^6_{24}B^5_{63}+f^{32}B^7_{34}B^5_{72}\\&\quad+f^{45}B^0_{44}B^5_{05}+f^{54}B^1_{54}B^5_{14}+f^{67}B^2_{64}B^5_{27}+f^{76}B^3_{74}B^5_{36}\\&=-f^{01}+f^{10}-f^{23}+f^{32}-f^{45}-f^{54}+f^{67}-f^{76}\end{aligned}$$

$$f_4^6 = f^{kr} B_{k4}^p B_{pr}^6$$
$$= f^{02} B_{04}^4 B_{42}^6 + f^{13} B_{14}^5 B_{53}^6 + f^{20} B_{24}^6 B_{60}^6 + f^{31} B_{34}^7 B_{71}^6$$
$$+ f^{46} B_{44}^0 B_{06}^6 + f^{57} B_{54}^1 B_{17}^6 + f^{64} B_{64}^2 B_{24}^6 + f^{75} B_{74}^3 B_{35}^6$$
$$= -f^{02} + f^{13} + f^{20} - f^{31} - f^{46} - f^{57} - f^{64} + f^{75}$$
(4.5.47)

$$f_4^7 = f^{kr} B_{k4}^p B_{pr}^7$$
$$= f^{03} B_{04}^4 B_{43}^7 + f^{12} B_{14}^5 B_{52}^7 + f^{21} B_{24}^6 B_{61}^7 + f^{30} B_{34}^7 B_{70}^7$$
$$+ f^{47} B_{44}^0 B_{07}^7 + f^{56} B_{54}^1 B_{16}^7 + f^{65} B_{64}^2 B_{25}^7 + f^{74} B_{74}^3 B_{34}^7$$
$$= -f^{03} - f^{12} + f^{21} + f^{30} - f^{47} + f^{56} - f^{65} - f^{74}$$
(4.5.48)

$$f_5^0 = f^{kr} B_{k5}^p B_{pr}^0$$
$$= f^{05} B_{05}^5 B_{55}^0 + f^{14} B_{15}^4 B_{44}^0 + f^{27} B_{25}^7 B_{77}^0 + f^{36} B_{35}^6 B_{66}^0$$
$$+ f^{41} B_{45}^1 B_{11}^0 + f^{50} B_{55}^0 B_{00}^0 + f^{63} B_{65}^3 B_{33}^0 + f^{72} B_{75}^2 B_{22}^0$$
$$= -f^{05} + f^{14} - f^{27} + f^{36} - f^{41} - f^{50} - f^{63} + f^{72}$$
(4.5.49)

$$f_5^1 = f^{kr} B_{k5}^p B_{pr}^1$$
$$= f^{04} B_{05}^5 B_{54}^1 + f^{15} B_{15}^4 B_{45}^1 + f^{26} B_{25}^7 B_{76}^1 + f^{37} B_{35}^6 B_{67}^1$$
$$+ f^{40} B_{45}^1 B_{10}^1 + f^{51} B_{55}^0 B_{01}^1 + f^{62} B_{65}^3 B_{32}^1 + f^{73} B_{75}^2 B_{23}^1$$
$$= -f^{04} - f^{15} + f^{26} + f^{37} + f^{40} - f^{51} - f^{62} - f^{73}$$
(4.5.50)

$$f_5^2 = f^{kr} B_{k5}^p B_{pr}^2$$
$$= f^{07} B_{05}^5 B_{57}^2 + f^{16} B_{15}^4 B_{46}^2 + f^{25} B_{25}^7 B_{75}^2 + f^{34} B_{35}^6 B_{64}^2$$
$$+ f^{43} B_{45}^1 B_{13}^2 + f^{52} B_{55}^0 B_{02}^2 + f^{61} B_{65}^3 B_{31}^2 + f^{70} B_{75}^2 B_{20}^2$$
$$= f^{07} - f^{16} - f^{25} + f^{34} - f^{43} - f^{52} + f^{61} - f^{70}$$
(4.5.51)

$$f_5^3 = f^{kr} B_{k5}^p B_{pr}^3$$
$$= f^{06} B_{05}^5 B_{56}^3 + f^{17} B_{15}^4 B_{47}^3 + f^{24} B_{25}^7 B_{74}^3 + f^{35} B_{35}^6 B_{65}^3$$
$$+ f^{42} B_{45}^1 B_{12}^3 + f^{53} B_{55}^0 B_{03}^3 + f^{60} B_{65}^3 B_{30}^3 + f^{71} B_{75}^2 B_{21}^3$$
$$= -f^{06} - f^{17} - f^{24} - f^{35} + f^{42} - f^{53} + f^{60} + f^{71}$$
(4.5.52)

$$f_5^4 = f^{kr} B_{k5}^p B_{pr}^4$$
$$= f^{01} B_{05}^5 B_{51}^4 + f^{10} B_{15}^4 B_{40}^4 + f^{23} B_{25}^7 B_{73}^4 + f^{32} B_{35}^6 B_{62}^4$$
$$+ f^{45} B_{45}^1 B_{15}^4 + f^{54} B_{55}^0 B_{04}^4 + f^{67} B_{65}^3 B_{37}^4 + f^{76} B_{75}^2 B_{26}^4$$
$$= f^{01} - f^{10} + f^{23} - f^{32} - f^{45} - f^{54} - f^{67} + f^{76}$$
(4.5.53)

$$(4.5.54) \quad \begin{aligned} f_5^5 &= f^{kr} B_{k5}^p B_{pr}^5 \\ &= f^{00} B_{05}^5 B_{50}^5 + f^{11} B_{15}^4 B_{41}^5 + f^{22} B_{25}^7 B_{72}^5 + f^{33} B_{35}^6 B_{63}^5 \\ &\quad + f^{44} B_{45}^1 B_{14}^5 + f^{55} B_{55}^0 B_{05}^5 + f^{66} B_{65}^3 B_{36}^5 + f^{77} B_{75}^2 B_{27}^5 \\ &= f^{00} + f^{11} + f^{22} + f^{33} + f^{44} - f^{55} + f^{66} + f^{77} \end{aligned}$$

$$(4.5.55) \quad \begin{aligned} f_5^6 &= f^{kr} B_{k5}^p B_{pr}^6 \\ &= f^{03} B_{05}^5 B_{53}^6 + f^{12} B_{15}^4 B_{42}^6 + f^{21} B_{25}^7 B_{71}^6 + f^{30} B_{35}^6 B_{60}^6 \\ &\quad + f^{47} B_{45}^1 B_{17}^6 + f^{56} B_{55}^0 B_{06}^6 + f^{65} B_{65}^3 B_{35}^6 + f^{74} B_{75}^2 B_{24}^6 \\ &= f^{03} + f^{12} - f^{21} - f^{30} + f^{47} - f^{56} - f^{65} - f^{74} \end{aligned}$$

$$(4.5.56) \quad \begin{aligned} f_5^7 &= f^{kr} B_{k5}^p B_{pr}^7 \\ &= f^{02} B_{05}^5 B_{52}^7 + f^{13} B_{15}^4 B_{43}^7 + f^{20} B_{25}^7 B_{70}^7 + f^{31} B_{35}^6 B_{61}^7 \\ &\quad + f^{46} B_{45}^1 B_{16}^7 + f^{57} B_{55}^0 B_{07}^7 + f^{64} B_{65}^3 B_{34}^7 + f^{75} B_{75}^2 B_{25}^7 \\ &= -f^{02} + f^{13} + f^{20} - f^{31} - f^{46} - f^{57} + f^{64} - f^{75} \end{aligned}$$

$$(4.5.57) \quad \begin{aligned} f_6^0 &= f^{kr} B_{k6}^p B_{pr}^0 \\ &= f^{06} B_{06}^6 B_{66}^0 + f^{17} B_{16}^7 B_{77}^0 + f^{24} B_{26}^4 B_{44}^0 + f^{35} B_{36}^5 B_{55}^0 \\ &\quad + f^{42} B_{46}^2 B_{22}^0 + f^{53} B_{56}^3 B_{33}^0 + f^{60} B_{66}^0 B_{00}^0 + f^{71} B_{76}^1 B_{11}^0 \\ &= -f^{06} + f^{17} + f^{24} - f^{35} - f^{42} + f^{53} - f^{60} - f^{71} \end{aligned}$$

$$(4.5.58) \quad \begin{aligned} f_6^1 &= f^{kr} B_{k6}^p B_{pr}^1 \\ &= f^{07} B_{06}^6 B_{67}^1 + f^{16} B_{16}^7 B_{76}^1 + f^{25} B_{26}^4 B_{45}^1 + f^{34} B_{36}^5 B_{54}^1 \\ &\quad + f^{43} B_{46}^2 B_{23}^1 + f^{52} B_{56}^3 B_{32}^1 + f^{61} B_{66}^0 B_{01}^1 + f^{70} B_{76}^1 B_{10}^1 \\ &= -f^{07} - f^{16} - f^{25} - f^{34} + f^{43} + f^{52} - f^{61} + f^{70} \end{aligned}$$

$$(4.5.59) \quad \begin{aligned} f_6^2 &= f^{kr} B_{k6}^p B_{pr}^2 \\ &= f^{04} B_{06}^6 B_{64}^2 + f^{15} B_{16}^7 B_{75}^2 + f^{26} B_{26}^4 B_{46}^2 + f^{37} B_{36}^5 B_{57}^2 \\ &\quad + f^{40} B_{46}^2 B_{20}^2 + f^{51} B_{56}^3 B_{31}^2 + f^{62} B_{66}^0 B_{02}^2 + f^{73} B_{76}^1 B_{13}^2 \\ &= -f^{04} + f^{15} - f^{26} + f^{37} + f^{40} - f^{51} - f^{62} - f^{73} \end{aligned}$$

$$(4.5.60) \quad \begin{aligned} f_6^3 &= f^{kr} B_{k6}^p B_{pr}^3 \\ &= f^{05} B_{06}^6 B_{65}^3 + f^{14} B_{16}^7 B_{74}^3 + f^{27} B_{26}^4 B_{47}^3 + f^{36} B_{36}^5 B_{56}^3 \\ &\quad + f^{41} B_{46}^2 B_{21}^3 + f^{50} B_{56}^3 B_{30}^3 + f^{63} B_{66}^0 B_{03}^3 + f^{72} B_{76}^1 B_{12}^3 \\ &= f^{05} + f^{14} - f^{27} - f^{36} - f^{41} - f^{50} - f^{63} + f^{72} \end{aligned}$$

$$(4.5.61) \quad \begin{aligned} f_6^4 &= f^{kr} B_{k6}^p B_{pr}^4 \\ &= f^{02} B_{06}^6 B_{62}^4 + f^{13} B_{16}^7 B_{73}^4 + f^{20} B_{26}^4 B_{40}^4 + f^{31} B_{36}^5 B_{51}^4 \\ &\quad + f^{46} B_{46}^2 B_{26}^4 + f^{57} B_{56}^3 B_{37}^4 + f^{64} B_{66}^0 B_{04}^4 + f^{75} B_{76}^1 B_{15}^4 \\ &= f^{02} - f^{13} - f^{20} + f^{31} - f^{46} + f^{57} - f^{64} - f^{75} \end{aligned}$$

$$(4.5.62) \quad \begin{aligned} f_6^5 &= f^{kr} B_{k6}^p B_{pr}^5 \\ &= f^{03} B_{06}^6 B_{63}^5 + f^{12} B_{16}^7 B_{72}^5 + f^{21} B_{26}^4 B_{41}^5 + f^{30} B_{36}^5 B_{50}^5 \\ &\quad + f^{47} B_{46}^2 B_{27}^5 + f^{56} B_{56}^3 B_{36}^5 + f^{65} B_{66}^0 B_{05}^5 + f^{74} B_{76}^1 B_{14}^5 \\ &= -f^{03} - f^{12} + f^{21} + f^{30} - f^{47} - f^{56} - f^{65} + f^{74} \end{aligned}$$

$$(4.5.63) \quad \begin{aligned} f_6^6 &= f^{kr} B_{k6}^p B_{pr}^6 \\ &= f^{00} B_{06}^6 B_{60}^6 + f^{11} B_{16}^7 B_{71}^6 + f^{22} B_{26}^4 B_{42}^6 + f^{33} B_{36}^5 B_{53}^6 \\ &\quad + f^{44} B_{46}^2 B_{24}^6 + f^{55} B_{56}^3 B_{35}^6 + f^{66} B_{66}^0 B_{06}^6 + f^{77} B_{76}^1 B_{17}^6 \\ &= f^{00} + f^{11} + f^{22} + f^{33} + f^{44} + f^{55} - f^{66} + f^{77} \end{aligned}$$

$$(4.5.64) \quad \begin{aligned} f_6^7 &= f^{kr} B_{k6}^p B_{pr}^7 \\ &= f^{01} B_{06}^6 B_{61}^7 + f^{10} B_{16}^7 B_{70}^7 + f^{23} B_{26}^4 B_{43}^7 + f^{32} B_{36}^5 B_{52}^7 \\ &\quad + f^{45} B_{46}^2 B_{25}^7 + f^{54} B_{56}^3 B_{34}^7 + f^{67} B_{66}^0 B_{07}^7 + f^{76} B_{76}^1 B_{16}^7 \\ &= f^{01} - f^{10} + f^{23} - f^{32} + f^{45} - f^{54} - f^{67} - f^{76} \end{aligned}$$

$$(4.5.65) \quad \begin{aligned} f_7^0 &= f^{kr} B_{k7}^p B_{pr}^0 \\ &= f^{07} B_{07}^7 B_{77}^0 + f^{16} B_{17}^6 B_{66}^0 + f^{25} B_{27}^5 B_{55}^0 + f^{34} B_{37}^4 B_{44}^0 \\ &\quad + f^{43} B_{47}^3 B_{33}^0 + f^{52} B_{57}^2 B_{22}^0 + f^{61} B_{67}^1 B_{11}^0 + f^{70} B_{77}^0 B_{00}^0 \\ &= -f^{07} - f^{16} + f^{25} + f^{34} - f^{43} - f^{52} + f^{61} - f^{70} \end{aligned}$$

$$(4.5.66) \quad \begin{aligned} f_7^1 &= f^{kr} B_{k7}^p B_{pr}^1 \\ &= f^{06} B_{07}^7 B_{76}^1 + f^{17} B_{17}^6 B_{67}^1 + f^{24} B_{27}^5 B_{54}^1 + f^{35} B_{37}^4 B_{45}^1 \\ &\quad + f^{42} B_{47}^3 B_{32}^1 + f^{53} B_{57}^2 B_{23}^1 + f^{60} B_{67}^1 B_{10}^1 + f^{71} B_{77}^0 B_{01}^1 \\ &= f^{06} - f^{17} + f^{24} - f^{35} - f^{42} + f^{53} - f^{60} - f^{71} \end{aligned}$$

$$(4.5.67) \quad \begin{aligned} f_7^2 &= f^{kr} B_{k7}^p B_{pr}^2 \\ &= f^{05} B_{07}^7 B_{75}^2 + f^{14} B_{17}^6 B_{64}^2 + f^{27} B_{27}^5 B_{57}^2 + f^{36} B_{37}^4 B_{46}^2 \\ &\quad + f^{41} B_{47}^3 B_{31}^2 + f^{50} B_{57}^2 B_{20}^2 + f^{63} B_{67}^1 B_{13}^2 + f^{72} B_{77}^0 B_{02}^2 \\ &= -f^{05} - f^{14} - f^{27} - f^{36} + f^{41} + f^{50} + f^{63} - f^{72} \end{aligned}$$

$$f_7^3 = f^{kr} B_{k7}^p B_{pr}^3$$
$$= f^{04} B_{07}^7 B_{74}^3 + f^{15} B_{17}^6 B_{65}^3 + f^{26} B_{27}^5 B_{56}^3 + f^{37} B_{37}^4 B_{47}^3$$
(4.5.68)
$$+ f^{40} B_{47}^3 B_{30}^3 + f^{51} B_{57}^2 B_{21}^3 + f^{62} B_{67}^1 B_{12}^3 + f^{73} B_{77}^0 B_{03}^3$$
$$= -f^{04} + f^{15} + f^{26} - f^{37} + f^{40} - f^{51} - f^{62} - f^{73}$$

$$f_7^4 = f^{kr} B_{k7}^p B_{pr}^4$$
$$= f^{03} B_{07}^7 B_{73}^4 + f^{12} B_{17}^6 B_{62}^4 + f^{21} B_{27}^5 B_{51}^4 + f^{30} B_{37}^4 B_{40}^4$$
(4.5.69)
$$+ f^{47} B_{47}^3 B_{37}^4 + f^{56} B_{57}^2 B_{26}^4 + f^{65} B_{67}^1 B_{15}^4 + f^{74} B_{77}^0 B_{04}^4$$
$$= f^{03} + f^{12} - f^{21} - f^{30} - f^{47} - f^{56} + f^{65} - f^{74}$$

$$f_7^5 = f^{kr} B_{k7}^p B_{pr}^5$$
$$= f^{02} B_{07}^7 B_{72}^5 + f^{13} B_{17}^6 B_{63}^5 + f^{20} B_{27}^5 B_{50}^5 + f^{31} B_{37}^4 B_{41}^5$$
(4.5.70)
$$+ f^{46} B_{47}^3 B_{36}^5 + f^{57} B_{57}^2 B_{27}^5 + f^{64} B_{67}^1 B_{14}^5 + f^{75} B_{77}^0 B_{05}^5$$
$$= f^{02} - f^{13} - f^{20} + f^{31} + f^{46} - f^{57} - f^{64} - f^{75}$$

$$f_7^6 = f^{kr} B_{k7}^p B_{pr}^6$$
$$= f^{01} B_{07}^7 B_{71}^6 + f^{10} B_{17}^6 B_{60}^6 + f^{23} B_{27}^5 B_{53}^6 + f^{32} B_{37}^4 B_{42}^6$$
(4.5.71)
$$+ f^{45} B_{47}^3 B_{35}^6 + f^{54} B_{57}^2 B_{24}^6 + f^{67} B_{67}^1 B_{17}^6 + f^{76} B_{77}^0 B_{06}^6$$
$$= -f^{01} + f^{10} - f^{23} + f^{32} - f^{45} + f^{54} - f^{67} - f^{76}$$

$$f_7^7 = f^{kr} B_{k7}^p B_{pr}^7$$
$$= f^{00} B_{07}^7 B_{70}^7 + f^{11} B_{17}^6 B_{61}^7 + f^{22} B_{27}^5 B_{52}^7 + f^{33} B_{37}^4 B_{43}^7$$
(4.5.72)
$$+ f^{44} B_{47}^3 B_{34}^7 + f^{55} B_{57}^2 B_{25}^7 + f^{66} B_{67}^1 B_{16}^7 + f^{77} B_{77}^0 B_{07}^7$$
$$= f^{00} + f^{11} + f^{22} + f^{33} + f^{44} + f^{55} + f^{66} - f^{77}$$

Уравнения (4.5.9), (4.5.18), (4.5.27), (4.5.36), (4.5.45), (4.5.54), (4.5.63), (4.5.72) формируют систему линейных уравнений (4.5.1).

Уравнения (4.5.10), (4.5.17), (4.5.28), (4.5.35), (4.5.46), (4.5.53), (4.5.64), (4.5.71) формируют систему линейных уравнений (4.5.2).

Уравнения (4.5.11), (4.5.20), (4.5.25), (4.5.34), (4.5.47), (4.5.56), (4.5.61), (4.5.70) формируют систему линейных уравнений (4.5.3).

Уравнения (4.5.12), (4.5.19), (4.5.26), (4.5.33), (4.5.48), (4.5.55), (4.5.62), (4.5.69) формируют систему линейных уравнений (4.5.4).

Уравнения (4.5.13), (4.5.22), (4.5.31), (4.5.40), (4.5.41), (4.5.50), (4.5.59), (4.5.68) формируют систему линейных уравнений (4.5.5).

Уравнения (4.5.14), (4.5.21), (4.5.32), (4.5.39), (4.5.42), (4.5.49), (4.5.60), (4.5.67) формируют систему линейных уравнений (4.5.6).

Уравнения (4.5.15), (4.5.24), (4.5.29), (4.5.38), (4.5.43), (4.5.52), (4.5.57), (4.5.66) формируют систему линейных уравнений (4.5.7).

Уравнения (4.5.16), (4.5.23), (4.5.30), (4.5.37), (4.5.44), (4.5.51), (4.5.58), (4.5.65) формируют систему линейных уравнений (4.5.8).

($4.5.90$) - это решение системы линейных уравнений ($4.5.1$). $\square$

**Теорема 4.5.2.** *Рассмотрим алгебру октонионов $O$ с базисом ($4.4.1$). Стандартные компоненты линейной функции и координаты этой функции удовлетворяют соотношениям*

(4.5.73) $$A = FB$$
(4.5.74) $$B = F^{-1}A$$

*где*

$$A = \begin{pmatrix} f_0^0 & f_1^0 & f_2^0 & f_3^0 & f_4^0 & f_5^0 & f_6^0 & f_7^0 \\ f_1^1 & -f_0^1 & f_3^1 & -f_2^1 & f_5^1 & -f_4^1 & -f_7^1 & f_6^1 \\ f_2^2 & -f_3^2 & -f_0^2 & f_1^2 & f_6^2 & f_7^2 & -f_4^2 & -f_5^2 \\ f_3^3 & f_2^3 & -f_1^3 & -f_0^3 & f_7^3 & -f_6^3 & f_5^3 & -f_4^3 \\ f_4^4 & -f_5^4 & -f_6^4 & -f_7^4 & -f_0^4 & f_1^4 & f_2^4 & f_3^4 \\ f_5^5 & f_4^5 & -f_7^5 & f_6^5 & -f_1^5 & -f_0^5 & -f_3^5 & f_2^5 \\ f_6^6 & f_7^6 & f_4^6 & -f_5^6 & -f_2^6 & f_3^6 & -f_0^6 & -f_1^6 \\ f_7^7 & -f_6^7 & f_5^7 & f_4^7 & -f_3^7 & -f_2^7 & f_1^7 & -f_0^7 \end{pmatrix}$$

$$B = \begin{pmatrix} f^{00} & -f^{01} & -f^{02} & -f^{03} & -f^{04} & -f^{05} & -f^{06} & -f^{07} \\ f^{11} & f^{10} & f^{13} & -f^{12} & f^{15} & -f^{14} & -f^{17} & f^{16} \\ f^{22} & -f^{23} & f^{20} & f^{21} & f^{26} & f^{27} & -f^{24} & -f^{25} \\ f^{33} & f^{32} & -f^{31} & f^{30} & f^{37} & -f^{36} & f^{35} & -f^{34} \\ f^{44} & -f^{45} & -f^{46} & -f^{47} & f^{40} & f^{41} & f^{42} & f^{43} \\ f^{55} & f^{54} & -f^{57} & f^{56} & -f^{51} & f^{50} & -f^{53} & f^{52} \\ f^{66} & f^{67} & f^{64} & -f^{65} & -f^{62} & f^{63} & f^{60} & -f^{61} \\ f^{77} & -f^{76} & f^{75} & f^{74} & -f^{73} & -f^{72} & f^{71} & f^{70} \end{pmatrix}$$

$$F = \begin{pmatrix} 1 & -1 & -1 & -1 & -1 & -1 & -1 & -1 \\ 1 & -1 & 1 & 1 & 1 & 1 & 1 & 1 \\ 1 & 1 & -1 & 1 & 1 & 1 & 1 & 1 \\ 1 & 1 & 1 & -1 & 1 & 1 & 1 & 1 \\ 1 & 1 & 1 & 1 & -1 & 1 & 1 & 1 \\ 1 & 1 & 1 & 1 & 1 & -1 & 1 & 1 \\ 1 & 1 & 1 & 1 & 1 & 1 & -1 & 1 \\ 1 & 1 & 1 & 1 & 1 & 1 & 1 & -1 \end{pmatrix}$$

$$F^{-1} = \frac{1}{12}\begin{pmatrix} 5 & 1 & 1 & 1 & 1 & 1 & 1 & 1 \\ -1 & -5 & 1 & 1 & 1 & 1 & 1 & 1 \\ -1 & 1 & -5 & 1 & 1 & 1 & 1 & 1 \\ -1 & 1 & 1 & -5 & 1 & 1 & 1 & 1 \\ -1 & 1 & 1 & 1 & -5 & 1 & 1 & 1 \\ -1 & 1 & 1 & 1 & 1 & -5 & 1 & 1 \\ -1 & 1 & 1 & 1 & 1 & 1 & -5 & 1 \\ -1 & 1 & 1 & 1 & 1 & 1 & 1 & -5 \end{pmatrix}$$

Доказательство. Запишем систему линейных уравнений (4.5.1) в виде произведения матриц

$$(4.5.75) \quad \begin{pmatrix} f_0^0 \\ f_1^1 \\ f_2^2 \\ f_3^3 \\ f_4^4 \\ f_5^5 \\ f_6^6 \\ f_7^7 \end{pmatrix} = \begin{pmatrix} 1 & -1 & -1 & -1 & -1 & -1 & -1 & -1 \\ 1 & -1 & 1 & 1 & 1 & 1 & 1 & 1 \\ 1 & 1 & -1 & 1 & 1 & 1 & 1 & 1 \\ 1 & 1 & 1 & -1 & 1 & 1 & 1 & 1 \\ 1 & 1 & 1 & 1 & -1 & 1 & 1 & 1 \\ 1 & 1 & 1 & 1 & 1 & -1 & 1 & 1 \\ 1 & 1 & 1 & 1 & 1 & 1 & -1 & 1 \\ 1 & 1 & 1 & 1 & 1 & 1 & 1 & -1 \end{pmatrix} \begin{pmatrix} f^{00} \\ f^{11} \\ f^{22} \\ f^{33} \\ f^{44} \\ f^{55} \\ f^{66} \\ f^{77} \end{pmatrix}$$

Запишем систему линейных уравнений (4.5.2) в виде произведения матриц

$$(4.5.76) \quad \begin{pmatrix} f_0^1 \\ f_1^0 \\ f_2^3 \\ f_3^2 \\ f_4^5 \\ f_5^4 \\ f_6^7 \\ f_7^6 \end{pmatrix} = \begin{pmatrix} 1 & 1 & 1 & -1 & 1 & -1 & -1 & 1 \\ -1 & -1 & 1 & -1 & 1 & -1 & -1 & 1 \\ -1 & 1 & -1 & -1 & -1 & 1 & 1 & -1 \\ 1 & -1 & -1 & -1 & 1 & -1 & -1 & 1 \\ -1 & 1 & -1 & 1 & -1 & -1 & 1 & -1 \\ 1 & -1 & 1 & -1 & -1 & -1 & -1 & 1 \\ 1 & -1 & 1 & -1 & 1 & -1 & -1 & -1 \\ -1 & 1 & -1 & 1 & -1 & 1 & -1 & -1 \end{pmatrix} \begin{pmatrix} f^{01} \\ f^{10} \\ f^{23} \\ f^{32} \\ f^{45} \\ f^{54} \\ f^{67} \\ f^{76} \end{pmatrix}$$

Из равенства (4.5.76) следует

$$\begin{pmatrix} -f_0^1 \\ f_1^0 \\ f_2^3 \\ -f_3^2 \\ f_4^5 \\ -f_5^4 \\ -f_6^7 \\ f_7^6 \end{pmatrix} = \begin{pmatrix} -1 & -1 & -1 & 1 & -1 & 1 & 1 & -1 \\ -1 & -1 & 1 & -1 & 1 & -1 & -1 & 1 \\ -1 & 1 & -1 & -1 & -1 & 1 & 1 & -1 \\ -1 & 1 & 1 & 1 & -1 & 1 & 1 & -1 \\ -1 & 1 & -1 & 1 & -1 & -1 & 1 & -1 \\ -1 & 1 & -1 & 1 & 1 & 1 & 1 & -1 \\ -1 & 1 & -1 & 1 & -1 & 1 & 1 & 1 \\ -1 & 1 & -1 & 1 & -1 & 1 & -1 & -1 \end{pmatrix} \begin{pmatrix} f^{01} \\ f^{10} \\ f^{23} \\ f^{32} \\ f^{45} \\ f^{54} \\ f^{67} \\ f^{76} \end{pmatrix}$$

$$= \begin{pmatrix} 1 & -1 & 1 & 1 & 1 & 1 & 1 & 1 \\ 1 & -1 & -1 & -1 & -1 & -1 & -1 & -1 \\ 1 & 1 & 1 & -1 & 1 & 1 & 1 & 1 \\ 1 & 1 & -1 & 1 & 1 & 1 & 1 & 1 \\ 1 & 1 & 1 & 1 & 1 & -1 & 1 & 1 \\ 1 & 1 & 1 & 1 & -1 & 1 & 1 & 1 \\ 1 & 1 & 1 & 1 & 1 & 1 & 1 & -1 \\ 1 & 1 & 1 & 1 & 1 & 1 & -1 & 1 \end{pmatrix} \begin{pmatrix} -f^{01} \\ f^{10} \\ -f^{23} \\ f^{32} \\ -f^{45} \\ f^{54} \\ f^{67} \\ -f^{76} \end{pmatrix}$$

$$(4.5.77)\quad \begin{pmatrix} f_1^0 \\ -f_0^1 \\ -f_3^2 \\ f_2^3 \\ -f_5^4 \\ f_4^5 \\ f_7^6 \\ -f_6^7 \end{pmatrix} = \begin{pmatrix} 1 & -1 & -1 & -1 & -1 & -1 & -1 & -1 \\ 1 & -1 & 1 & 1 & 1 & 1 & 1 & 1 \\ 1 & 1 & -1 & 1 & 1 & 1 & 1 & 1 \\ 1 & 1 & 1 & -1 & 1 & 1 & 1 & 1 \\ 1 & 1 & 1 & 1 & -1 & 1 & 1 & 1 \\ 1 & 1 & 1 & 1 & 1 & -1 & 1 & 1 \\ 1 & 1 & 1 & 1 & 1 & 1 & -1 & 1 \\ 1 & 1 & 1 & 1 & 1 & 1 & 1 & -1 \end{pmatrix} \begin{pmatrix} -f^{01} \\ f^{10} \\ -f^{23} \\ f^{32} \\ -f^{45} \\ f^{54} \\ f^{67} \\ -f^{76} \end{pmatrix}$$

Запишем систему линейных уравнений (4.5.3) в виде произведения матриц

$$(4.5.78)\quad \begin{pmatrix} f_0^2 \\ f_1^3 \\ f_2^0 \\ f_3^1 \\ f_4^6 \\ f_5^7 \\ f_6^4 \\ f_7^5 \end{pmatrix} = \begin{pmatrix} 1 & -1 & 1 & 1 & 1 & 1 & -1 & -1 \\ 1 & -1 & -1 & -1 & 1 & 1 & -1 & -1 \\ -1 & -1 & -1 & 1 & 1 & 1 & -1 & -1 \\ -1 & -1 & 1 & -1 & -1 & -1 & 1 & 1 \\ -1 & 1 & 1 & -1 & -1 & -1 & -1 & 1 \\ -1 & 1 & 1 & -1 & -1 & -1 & 1 & -1 \\ 1 & -1 & -1 & 1 & -1 & 1 & -1 & -1 \\ 1 & -1 & -1 & 1 & 1 & -1 & -1 & -1 \end{pmatrix} \begin{pmatrix} f^{02} \\ f^{13} \\ f^{20} \\ f^{31} \\ f^{46} \\ f^{57} \\ f^{64} \\ f^{75} \end{pmatrix}$$

Из равенства (4.5.78) следует

$$\begin{pmatrix} -f_0^2 \\ -f_1^3 \\ f_2^0 \\ f_3^1 \\ f_4^6 \\ f_5^7 \\ -f_6^4 \\ -f_7^5 \end{pmatrix} = \begin{pmatrix} -1 & 1 & -1 & -1 & -1 & -1 & 1 & 1 \\ -1 & 1 & 1 & 1 & -1 & -1 & 1 & 1 \\ -1 & -1 & -1 & 1 & 1 & 1 & -1 & -1 \\ -1 & -1 & 1 & -1 & -1 & -1 & 1 & 1 \\ -1 & 1 & 1 & -1 & -1 & -1 & -1 & 1 \\ -1 & 1 & 1 & -1 & -1 & -1 & 1 & -1 \\ -1 & 1 & 1 & -1 & 1 & -1 & 1 & 1 \\ -1 & 1 & 1 & -1 & -1 & 1 & 1 & 1 \end{pmatrix} \begin{pmatrix} f^{02} \\ f^{13} \\ f^{20} \\ f^{31} \\ f^{46} \\ f^{57} \\ f^{64} \\ f^{75} \end{pmatrix}$$

$$= \begin{pmatrix} 1 & 1 & -1 & 1 & 1 & 1 & 1 & 1 \\ 1 & 1 & 1 & -1 & 1 & 1 & 1 & 1 \\ 1 & -1 & -1 & -1 & -1 & -1 & -1 & -1 \\ 1 & -1 & 1 & 1 & 1 & 1 & 1 & 1 \\ 1 & 1 & 1 & 1 & 1 & 1 & -1 & 1 \\ 1 & 1 & 1 & 1 & 1 & 1 & 1 & -1 \\ 1 & 1 & 1 & 1 & -1 & 1 & 1 & 1 \\ 1 & 1 & 1 & 1 & 1 & -1 & 1 & 1 \end{pmatrix} \begin{pmatrix} -f^{02} \\ f^{13} \\ f^{20} \\ -f^{31} \\ -f^{46} \\ -f^{57} \\ f^{64} \\ f^{75} \end{pmatrix}$$

$$(4.5.79)\quad \begin{pmatrix} f_2^0 \\ f_3^1 \\ -f_0^2 \\ -f_1^3 \\ -f_6^4 \\ -f_7^5 \\ f_4^6 \\ f_5^7 \end{pmatrix} = \begin{pmatrix} 1 & -1 & -1 & -1 & -1 & -1 & -1 & -1 \\ 1 & -1 & 1 & 1 & 1 & 1 & 1 & 1 \\ 1 & 1 & -1 & 1 & 1 & 1 & 1 & 1 \\ 1 & 1 & 1 & -1 & 1 & 1 & 1 & 1 \\ 1 & 1 & 1 & 1 & -1 & 1 & 1 & 1 \\ 1 & 1 & 1 & 1 & 1 & -1 & 1 & 1 \\ 1 & 1 & 1 & 1 & 1 & 1 & -1 & 1 \\ 1 & 1 & 1 & 1 & 1 & 1 & 1 & -1 \end{pmatrix} \begin{pmatrix} -f^{02} \\ f^{13} \\ f^{20} \\ -f^{31} \\ -f^{46} \\ -f^{57} \\ f^{64} \\ f^{75} \end{pmatrix}$$

Запишем систему линейных уравнений (4.5.4) в виде произведения матриц

$$(4.5.80)\quad \begin{pmatrix} f_0^3 \\ f_1^2 \\ f_2^1 \\ f_3^0 \\ f_4^7 \\ f_5^6 \\ f_6^5 \\ f_7^4 \end{pmatrix} = \begin{pmatrix} 1 & 1 & -1 & 1 & 1 & -1 & 1 & -1 \\ -1 & -1 & -1 & 1 & -1 & 1 & -1 & 1 \\ 1 & -1 & -1 & -1 & 1 & -1 & 1 & -1 \\ -1 & 1 & -1 & -1 & 1 & -1 & 1 & -1 \\ -1 & -1 & 1 & 1 & -1 & 1 & -1 & -1 \\ 1 & 1 & -1 & -1 & 1 & -1 & -1 & -1 \\ -1 & -1 & 1 & 1 & -1 & -1 & -1 & 1 \\ 1 & 1 & -1 & -1 & -1 & -1 & 1 & -1 \end{pmatrix} \begin{pmatrix} f^{03} \\ f^{12} \\ f^{21} \\ f^{30} \\ f^{47} \\ f^{56} \\ f^{65} \\ f^{74} \end{pmatrix}$$

Из равенства (4.5.80) следует

$$\begin{pmatrix} -f_0^3 \\ f_1^2 \\ -f_2^1 \\ f_3^0 \\ f_4^7 \\ -f_5^6 \\ f_6^5 \\ -f_7^4 \end{pmatrix} = \begin{pmatrix} -1 & -1 & 1 & -1 & -1 & 1 & -1 & 1 \\ -1 & -1 & -1 & 1 & -1 & 1 & -1 & 1 \\ -1 & 1 & 1 & 1 & -1 & 1 & -1 & 1 \\ -1 & 1 & -1 & -1 & 1 & -1 & 1 & -1 \\ -1 & -1 & 1 & 1 & -1 & 1 & -1 & -1 \\ -1 & -1 & 1 & 1 & -1 & 1 & 1 & 1 \\ -1 & -1 & 1 & 1 & -1 & -1 & -1 & 1 \\ -1 & -1 & 1 & 1 & 1 & 1 & -1 & 1 \end{pmatrix} \begin{pmatrix} f^{03} \\ f^{12} \\ f^{21} \\ f^{30} \\ f^{47} \\ f^{56} \\ f^{65} \\ f^{74} \end{pmatrix}$$

$$= \begin{pmatrix} 1 & 1 & 1 & -1 & 1 & 1 & 1 & 1 \\ 1 & 1 & -1 & 1 & 1 & 1 & 1 & 1 \\ 1 & -1 & 1 & 1 & 1 & 1 & 1 & 1 \\ 1 & -1 & -1 & -1 & -1 & -1 & -1 & -1 \\ 1 & 1 & 1 & 1 & 1 & 1 & 1 & -1 \\ 1 & 1 & 1 & 1 & 1 & 1 & -1 & 1 \\ 1 & 1 & 1 & 1 & 1 & -1 & 1 & 1 \\ 1 & 1 & 1 & 1 & -1 & 1 & 1 & 1 \end{pmatrix} \begin{pmatrix} -f^{03} \\ -f^{12} \\ f^{21} \\ f^{30} \\ -f^{47} \\ f^{56} \\ -f^{65} \\ f^{74} \end{pmatrix}$$

$$(4.5.81) \begin{pmatrix} f_3^0 \\ -f_2^1 \\ f_1^2 \\ -f_0^3 \\ -f_7^4 \\ f_6^5 \\ -f_5^6 \\ f_4^7 \end{pmatrix} = \begin{pmatrix} 1 & -1 & -1 & -1 & -1 & -1 & -1 & -1 \\ 1 & -1 & 1 & 1 & 1 & 1 & 1 & 1 \\ 1 & 1 & -1 & 1 & 1 & 1 & 1 & 1 \\ 1 & 1 & 1 & -1 & 1 & 1 & 1 & 1 \\ 1 & 1 & 1 & 1 & -1 & 1 & 1 & 1 \\ 1 & 1 & 1 & 1 & 1 & -1 & 1 & 1 \\ 1 & 1 & 1 & 1 & 1 & 1 & -1 & 1 \\ 1 & 1 & 1 & 1 & 1 & 1 & 1 & -1 \end{pmatrix} \begin{pmatrix} -f^{03} \\ -f^{12} \\ f^{21} \\ f^{30} \\ -f^{47} \\ f^{56} \\ -f^{65} \\ f^{74} \end{pmatrix}$$

Запишем систему линейных уравнений (4.5.5) в виде произведения матриц

$$(4.5.82) \begin{pmatrix} f_0^4 \\ f_1^5 \\ f_2^6 \\ f_3^7 \\ f_4^0 \\ f_5^1 \\ f_6^2 \\ f_7^3 \end{pmatrix} = \begin{pmatrix} 1 & -1 & -1 & -1 & 1 & 1 & 1 & 1 \\ 1 & -1 & -1 & -1 & -1 & -1 & 1 & 1 \\ 1 & -1 & -1 & -1 & -1 & 1 & -1 & 1 \\ 1 & -1 & -1 & -1 & -1 & 1 & 1 & -1 \\ -1 & -1 & -1 & -1 & -1 & 1 & 1 & 1 \\ -1 & -1 & 1 & 1 & 1 & -1 & -1 & -1 \\ -1 & 1 & -1 & 1 & 1 & -1 & -1 & -1 \\ -1 & 1 & 1 & -1 & 1 & -1 & -1 & -1 \end{pmatrix} \begin{pmatrix} f^{04} \\ f^{15} \\ f^{26} \\ f^{37} \\ f^{40} \\ f^{51} \\ f^{62} \\ f^{73} \end{pmatrix}$$

Из равенства (4.5.82) следует

$$\begin{pmatrix} -f_0^4 \\ -f_1^5 \\ -f_2^6 \\ -f_3^7 \\ f_4^0 \\ f_5^1 \\ f_6^2 \\ f_7^3 \end{pmatrix} = \begin{pmatrix} -1 & 1 & 1 & 1 & -1 & -1 & -1 & -1 \\ -1 & 1 & 1 & 1 & 1 & 1 & -1 & -1 \\ -1 & 1 & 1 & 1 & 1 & -1 & 1 & -1 \\ -1 & 1 & 1 & 1 & 1 & -1 & -1 & 1 \\ -1 & -1 & -1 & -1 & -1 & 1 & 1 & 1 \\ -1 & -1 & 1 & 1 & 1 & -1 & -1 & -1 \\ -1 & 1 & -1 & 1 & 1 & -1 & -1 & -1 \\ -1 & 1 & 1 & -1 & 1 & -1 & -1 & -1 \end{pmatrix} \begin{pmatrix} f^{04} \\ f^{15} \\ f^{26} \\ f^{37} \\ f^{40} \\ f^{51} \\ f^{62} \\ f^{73} \end{pmatrix}$$

$$= \begin{pmatrix} 1 & 1 & 1 & 1 & -1 & 1 & 1 & 1 \\ 1 & 1 & 1 & 1 & 1 & -1 & 1 & 1 \\ 1 & 1 & 1 & 1 & 1 & 1 & -1 & 1 \\ 1 & 1 & 1 & 1 & 1 & 1 & 1 & -1 \\ 1 & -1 & -1 & -1 & -1 & -1 & -1 & -1 \\ 1 & -1 & 1 & 1 & 1 & 1 & 1 & 1 \\ 1 & 1 & -1 & 1 & 1 & 1 & 1 & 1 \\ 1 & 1 & 1 & -1 & 1 & 1 & 1 & 1 \end{pmatrix} \begin{pmatrix} -f^{04} \\ f^{15} \\ f^{26} \\ f^{37} \\ f^{40} \\ -f^{51} \\ -f^{62} \\ -f^{73} \end{pmatrix}$$

$$(4.5.83)\begin{pmatrix}f_4^0\\f_5^1\\f_6^2\\f_7^3\\-f_0^4\\-f_1^5\\-f_2^6\\-f_3^7\end{pmatrix}=\begin{pmatrix}1&-1&-1&-1&-1&-1&-1&-1\\1&-1&1&1&1&1&1&1\\1&1&-1&1&1&1&1&1\\1&1&1&-1&1&1&1&1\\1&1&1&1&-1&1&1&1\\1&1&1&1&1&-1&1&1\\1&1&1&1&1&1&-1&1\\1&1&1&1&1&1&1&-1\end{pmatrix}\begin{pmatrix}-f^{04}\\f^{15}\\f^{26}\\f^{37}\\f^{40}\\-f^{51}\\-f^{62}\\-f^{73}\end{pmatrix}$$

Запишем систему линейных уравнений (4.5.6) в виде произведения матриц

$$(4.5.84)\begin{pmatrix}f_0^5\\f_1^4\\f_2^7\\f_3^6\\f_4^1\\f_5^0\\f_6^3\\f_7^2\end{pmatrix}=\begin{pmatrix}1&1&-1&1&-1&1&-1&1\\-1&-1&1&-1&-1&1&1&-1\\1&1&-1&1&-1&-1&-1&-1\\-1&-1&1&-1&1&1&-1&-1\\1&-1&-1&1&-1&-1&-1&1\\-1&1&-1&1&-1&-1&-1&1\\1&1&-1&-1&-1&-1&-1&1\\-1&-1&-1&-1&1&1&1&-1\end{pmatrix}\begin{pmatrix}f^{05}\\f^{14}\\f^{27}\\f^{36}\\f^{41}\\f^{50}\\f^{63}\\f^{72}\end{pmatrix}$$

Из равенства (4.5.84) следует

$$\begin{pmatrix}-f_0^5\\f_1^4\\-f_2^7\\f_3^6\\-f_4^1\\f_5^0\\-f_6^3\\f_7^2\end{pmatrix}=\begin{pmatrix}-1&-1&1&-1&1&-1&1&-1\\-1&-1&1&-1&-1&1&1&-1\\-1&-1&1&-1&1&1&1&1\\-1&-1&1&-1&1&1&-1&-1\\-1&1&1&-1&1&1&1&-1\\-1&1&-1&1&-1&-1&-1&1\\-1&-1&1&1&1&1&1&-1\\-1&-1&-1&-1&1&1&1&-1\end{pmatrix}\begin{pmatrix}f^{05}\\f^{14}\\f^{27}\\f^{36}\\f^{41}\\f^{50}\\f^{63}\\f^{72}\end{pmatrix}$$
$$=\begin{pmatrix}1&1&1&1&1&-1&1&1\\1&1&1&1&-1&1&1&1\\1&1&1&1&1&1&1&-1\\1&1&1&1&1&1&-1&1\\1&-1&1&1&1&1&1&1\\1&-1&-1&-1&-1&-1&-1&-1\\1&1&1&-1&1&1&1&1\\1&1&-1&1&1&1&1&1\end{pmatrix}\begin{pmatrix}-f^{05}\\-f^{14}\\f^{27}\\-f^{36}\\f^{41}\\f^{50}\\f^{63}\\-f^{72}\end{pmatrix}$$

$$(4.5.85) \quad \begin{pmatrix} f_5^0 \\ -f_4^1 \\ f_7^2 \\ -f_6^3 \\ f_1^4 \\ -f_0^5 \\ f_3^6 \\ -f_2^7 \end{pmatrix} = \begin{pmatrix} 1 & -1 & -1 & -1 & -1 & -1 & -1 & -1 \\ 1 & -1 & 1 & 1 & 1 & 1 & 1 & 1 \\ 1 & 1 & -1 & 1 & 1 & 1 & 1 & 1 \\ 1 & 1 & 1 & -1 & 1 & 1 & 1 & 1 \\ 1 & 1 & 1 & 1 & -1 & 1 & 1 & 1 \\ 1 & 1 & 1 & 1 & 1 & -1 & 1 & 1 \\ 1 & 1 & 1 & 1 & 1 & 1 & -1 & 1 \\ 1 & 1 & 1 & 1 & 1 & 1 & 1 & -1 \end{pmatrix} \begin{pmatrix} -f^{05} \\ -f^{14} \\ f^{27} \\ -f^{36} \\ f^{41} \\ f^{50} \\ f^{63} \\ -f^{72} \end{pmatrix}$$

Запишем систему линейных уравнений (4.5.7) в виде произведения матриц

$$(4.5.86) \quad \begin{pmatrix} f_0^6 \\ f_1^7 \\ f_2^4 \\ f_3^5 \\ f_4^2 \\ f_5^3 \\ f_6^0 \\ f_7^1 \end{pmatrix} = \begin{pmatrix} 1 & 1 & 1 & -1 & -1 & 1 & 1 & -1 \\ -1 & -1 & -1 & 1 & 1 & -1 & 1 & -1 \\ -1 & -1 & -1 & 1 & -1 & -1 & 1 & 1 \\ 1 & 1 & 1 & -1 & -1 & -1 & -1 & -1 \\ 1 & 1 & -1 & -1 & -1 & 1 & -1 & -1 \\ -1 & -1 & -1 & -1 & 1 & -1 & 1 & 1 \\ -1 & 1 & 1 & -1 & -1 & 1 & -1 & -1 \\ 1 & -1 & 1 & -1 & -1 & 1 & -1 & -1 \end{pmatrix} \begin{pmatrix} f^{06} \\ f^{17} \\ f^{24} \\ f^{35} \\ f^{42} \\ f^{53} \\ f^{60} \\ f^{71} \end{pmatrix}$$

Из равенства (4.5.86) следует

$$\begin{pmatrix} -f_0^6 \\ f_1^7 \\ f_2^4 \\ -f_3^5 \\ -f_4^2 \\ f_5^3 \\ f_6^0 \\ -f_7^1 \end{pmatrix} = \begin{pmatrix} -1 & -1 & -1 & 1 & 1 & -1 & -1 & 1 \\ -1 & -1 & -1 & 1 & 1 & -1 & 1 & -1 \\ -1 & -1 & -1 & 1 & -1 & -1 & 1 & 1 \\ -1 & -1 & -1 & 1 & 1 & 1 & 1 & 1 \\ -1 & -1 & 1 & 1 & 1 & -1 & 1 & 1 \\ -1 & -1 & -1 & -1 & 1 & -1 & 1 & 1 \\ -1 & 1 & 1 & -1 & -1 & 1 & -1 & -1 \\ -1 & 1 & -1 & 1 & 1 & -1 & 1 & 1 \end{pmatrix} \begin{pmatrix} f^{06} \\ f^{17} \\ f^{24} \\ f^{35} \\ f^{42} \\ f^{53} \\ f^{60} \\ f^{71} \end{pmatrix}$$

$$= \begin{pmatrix} 1 & 1 & 1 & 1 & 1 & 1 & -1 & 1 \\ 1 & 1 & 1 & 1 & 1 & 1 & 1 & -1 \\ 1 & 1 & 1 & 1 & -1 & 1 & 1 & 1 \\ 1 & 1 & 1 & 1 & 1 & -1 & 1 & 1 \\ 1 & 1 & -1 & 1 & 1 & 1 & 1 & 1 \\ 1 & 1 & 1 & -1 & 1 & 1 & 1 & 1 \\ 1 & -1 & -1 & -1 & -1 & -1 & -1 & -1 \\ 1 & -1 & 1 & 1 & 1 & 1 & 1 & 1 \end{pmatrix} \begin{pmatrix} -f^{06} \\ -f^{17} \\ -f^{24} \\ f^{35} \\ f^{42} \\ -f^{53} \\ f^{60} \\ f^{71} \end{pmatrix}$$

$$(4.5.87) \quad \begin{pmatrix} f_6^0 \\ -f_7^1 \\ -f_4^2 \\ f_5^3 \\ f_2^4 \\ -f_3^5 \\ -f_0^6 \\ f_1^7 \end{pmatrix} = \begin{pmatrix} 1 & -1 & -1 & -1 & -1 & -1 & -1 & -1 \\ 1 & -1 & 1 & 1 & 1 & 1 & 1 & 1 \\ 1 & 1 & -1 & 1 & 1 & 1 & 1 & 1 \\ 1 & 1 & 1 & -1 & 1 & 1 & 1 & 1 \\ 1 & 1 & 1 & 1 & -1 & 1 & 1 & 1 \\ 1 & 1 & 1 & 1 & 1 & -1 & 1 & 1 \\ 1 & 1 & 1 & 1 & 1 & 1 & -1 & 1 \\ 1 & 1 & 1 & 1 & 1 & 1 & 1 & -1 \end{pmatrix} \begin{pmatrix} -f^{06} \\ -f^{17} \\ -f^{24} \\ f^{35} \\ f^{42} \\ -f^{53} \\ f^{60} \\ f^{71} \end{pmatrix}$$

Запишем систему линейных уравнений (4.5.8) в виде произведения матриц

$$(4.5.88) \quad \begin{pmatrix} f_0^7 \\ f_1^6 \\ f_2^5 \\ f_3^4 \\ f_4^3 \\ f_5^2 \\ f_6^1 \\ f_7^0 \end{pmatrix} = \begin{pmatrix} 1 & -1 & 1 & 1 & -1 & -1 & 1 & 1 \\ 1 & -1 & 1 & 1 & -1 & -1 & -1 & -1 \\ -1 & 1 & -1 & -1 & 1 & -1 & -1 & 1 \\ -1 & 1 & -1 & -1 & -1 & 1 & -1 & 1 \\ 1 & -1 & 1 & -1 & -1 & -1 & 1 & -1 \\ 1 & -1 & -1 & 1 & -1 & -1 & 1 & -1 \\ -1 & -1 & -1 & -1 & 1 & 1 & -1 & 1 \\ -1 & -1 & 1 & 1 & -1 & -1 & 1 & -1 \end{pmatrix} \begin{pmatrix} f^{07} \\ f^{16} \\ f^{25} \\ f^{34} \\ f^{43} \\ f^{52} \\ f^{61} \\ f^{70} \end{pmatrix}$$

Из равенства (4.5.88) следует

$$\begin{pmatrix} -f_0^7 \\ -f_1^6 \\ f_2^5 \\ f_3^4 \\ -f_4^3 \\ -f_5^2 \\ f_6^1 \\ f_7^0 \end{pmatrix} = \begin{pmatrix} -1 & 1 & -1 & -1 & 1 & 1 & -1 & -1 \\ -1 & 1 & -1 & -1 & 1 & 1 & 1 & 1 \\ -1 & 1 & -1 & -1 & 1 & -1 & -1 & 1 \\ -1 & 1 & -1 & -1 & -1 & 1 & -1 & 1 \\ -1 & 1 & -1 & 1 & 1 & 1 & -1 & 1 \\ -1 & 1 & 1 & -1 & 1 & 1 & -1 & 1 \\ -1 & -1 & -1 & -1 & 1 & 1 & -1 & 1 \\ -1 & -1 & 1 & 1 & -1 & -1 & 1 & -1 \end{pmatrix} \begin{pmatrix} f^{07} \\ f^{16} \\ f^{25} \\ f^{34} \\ f^{43} \\ f^{52} \\ f^{61} \\ f^{70} \end{pmatrix}$$

$$= \begin{pmatrix} 1 & 1 & 1 & 1 & 1 & 1 & 1 & -1 \\ 1 & 1 & 1 & 1 & 1 & 1 & -1 & 1 \\ 1 & 1 & 1 & 1 & 1 & -1 & 1 & 1 \\ 1 & 1 & 1 & 1 & -1 & 1 & 1 & 1 \\ 1 & 1 & 1 & -1 & 1 & 1 & 1 & 1 \\ 1 & 1 & -1 & 1 & 1 & 1 & 1 & 1 \\ 1 & -1 & 1 & 1 & 1 & 1 & 1 & 1 \\ 1 & -1 & -1 & -1 & -1 & -1 & -1 & -1 \end{pmatrix} \begin{pmatrix} -f^{07} \\ f^{16} \\ -f^{25} \\ -f^{34} \\ f^{43} \\ f^{52} \\ -f^{61} \\ f^{70} \end{pmatrix}$$

$$(4.5.89)\begin{pmatrix} f_7^0 \\ f_6^1 \\ -f_5^2 \\ -f_4^3 \\ f_3^4 \\ f_2^5 \\ -f_1^6 \\ -f_0^7 \end{pmatrix} = \begin{pmatrix} 1 & -1 & -1 & -1 & -1 & -1 & -1 & -1 \\ 1 & -1 & 1 & 1 & 1 & 1 & 1 & 1 \\ 1 & 1 & -1 & 1 & 1 & 1 & 1 & 1 \\ 1 & 1 & 1 & -1 & 1 & 1 & 1 & 1 \\ 1 & 1 & 1 & 1 & -1 & 1 & 1 & 1 \\ 1 & 1 & 1 & 1 & 1 & -1 & 1 & 1 \\ 1 & 1 & 1 & 1 & 1 & 1 & -1 & 1 \\ 1 & 1 & 1 & 1 & 1 & 1 & 1 & -1 \end{pmatrix} \begin{pmatrix} -f^{07} \\ f^{16} \\ -f^{25} \\ -f^{34} \\ f^{43} \\ f^{52} \\ -f^{61} \\ f^{70} \end{pmatrix}$$

Мы объединяем равенства (4.5.75), (4.5.77), (4.5.79), (4.5.81), (4.5.83), (4.5.85), (4.5.87), (4.5.89) в равенстве (4.5.73). □

**Теорема 4.5.3.** *Стандартные компоненты линейной функции алгебры октонионов $O$ относительно базиса (4.4.1) и координаты соответствующего линейного преобразования удовлетворяют соотношениям*

$$(4.5.90)\begin{cases} 12f^{00} = 5f_0^0 + f_1^1 + f_2^2 + f_3^3 + f_4^4 + f_5^5 + f_6^6 + f_7^7 \\ 12f^{11} = -f_0^0 - 5f_1^1 + f_2^2 + f_3^3 + f_4^4 + f_5^5 + f_6^6 + f_7^7 \\ 12f^{22} = -f_0^0 + f_1^1 - 5f_2^2 + f_3^3 + f_4^4 + f_5^5 + f_6^6 + f_7^7 \\ 12f^{33} = -f_0^0 + f_1^1 + f_2^2 - 5f_3^3 + f_4^4 + f_5^5 + f_6^6 + f_7^7 \\ 12f^{44} = -f_0^0 + f_1^1 + f_2^2 + f_3^3 - 5f_4^4 + f_5^5 + f_6^6 + f_7^7 \\ 12f^{55} = -f_0^0 + f_1^1 + f_2^2 + f_3^3 + f_4^4 - 5f_5^5 + f_6^6 + f_7^7 \\ 12f^{66} = -f_0^0 + f_1^1 + f_2^2 + f_3^3 + f_4^4 + f_5^5 - 5f_6^6 + f_7^7 \\ 12f^{77} = -f_0^0 + f_1^1 + f_2^2 + f_3^3 + f_4^4 + f_5^5 + f_6^6 - 5f_7^7 \end{cases}$$

$$(4.5.91)\begin{cases} -12f^{01} = 5f_1^0 - f_0^1 - f_3^2 + f_2^3 - f_5^4 + f_4^5 + f_7^6 - f_6^7 \\ 12f^{10} = -f_1^0 + 5f_0^1 - f_3^2 + f_2^3 - f_5^4 + f_4^5 + f_7^6 - f_6^7 \\ -12f^{23} = -f_1^0 - f_0^1 + 5f_3^2 + f_2^3 - f_5^4 + f_4^5 + f_7^6 - f_6^7 \\ 12f^{32} = -f_1^0 - f_0^1 - f_3^2 - 5f_2^3 - f_5^4 + f_4^5 + f_7^6 - f_6^7 \\ -12f^{45} = -f_1^0 - f_0^1 - f_3^2 + f_2^3 + 5f_5^4 + f_4^5 + f_7^6 - f_6^7 \\ 12f^{54} = -f_1^0 - f_0^1 - f_3^2 + f_2^3 - f_5^4 - 5f_4^5 + f_7^6 - f_6^7 \\ 12f^{67} = -f_1^0 - f_0^1 - f_3^2 + f_2^3 - f_5^4 + f_4^5 - 5f_7^6 - f_6^7 \\ -12f^{76} = -f_1^0 - f_0^1 - f_3^2 + f_2^3 - f_5^4 + f_4^5 + f_7^6 + 5f_6^7 \end{cases}$$

$$(4.5.92)\begin{cases} -12f^{02} = 5f_2^0 + f_3^1 - f_0^2 - f_1^3 - f_6^4 - f_7^5 + f_4^6 + f_5^7 \\ 12f^{13} = -f_2^0 - 5f_3^1 - f_0^2 - f_1^3 - f_6^4 - f_7^5 + f_4^6 + f_5^7 \\ 12f^{20} = -f_2^0 + f_3^1 + 5f_0^2 - f_1^3 - f_6^4 - f_7^5 + f_4^6 + f_5^7 \\ -12f^{31} = -f_2^0 + f_3^1 - f_0^2 + 5f_1^3 - f_6^4 - f_7^5 + f_4^6 + f_5^7 \\ -12f^{46} = -f_2^0 + f_3^1 - f_0^2 - f_1^3 + 5f_6^4 - f_7^5 + f_4^6 + f_5^7 \\ -12f^{57} = -f_2^0 + f_3^1 - f_0^2 - f_1^3 - f_6^4 + 5f_7^5 + f_4^6 + f_5^7 \\ 12f^{64} = -f_2^0 + f_3^1 - f_0^2 - f_1^3 - f_6^4 - f_7^5 - 5f_4^6 + f_5^7 \\ 12f^{75} = -f_2^0 + f_3^1 - f_0^2 - f_1^3 - f_6^4 - f_7^5 + f_4^6 - 5f_5^7 \end{cases}$$

$$(4.5.93)\begin{cases}-12f^{03}=5f_3^0-f_2^1+f_1^2-f_0^3-f_7^4+f_6^5-f_5^6+f_4^7\\-12f^{12}=-f_3^0+5f_2^1+f_1^2-f_0^3-f_7^4+f_6^5-f_5^6+f_4^7\\12f^{21}=-f_3^0-f_2^1-5f_1^2-f_0^3-f_7^4+f_6^5-f_5^6+f_4^7\\12f^{30}=-f_3^0-f_2^1+f_1^2+5f_0^3-f_7^4+f_6^5-f_5^6+f_4^7\\-12f^{47}=-f_3^0-f_2^1+f_1^2-f_0^3+5f_7^4+f_6^5-f_5^6+f_4^7\\12f^{56}=-f_3^0-f_2^1+f_1^2-f_0^3-f_7^4-5f_6^5-f_5^6+f_4^7\\-12f^{65}=-f_3^0-f_2^1+f_1^2-f_0^3-f_7^4+f_6^5+5f_5^6+f_4^7\\12f^{74}=-f_3^0-f_2^1+f_1^2-f_0^3-f_7^4+f_6^5-f_5^6-5f_4^7\end{cases}$$

$$(4.5.94)\begin{cases}-12f^{04}=5f_4^0+f_5^1+f_6^2+f_7^3-f_0^4-f_1^5-f_2^6-f_3^7\\12f^{15}=-f_4^0-5f_5^1+f_6^2+f_7^3-f_0^4-f_1^5-f_2^6-f_3^7\\12f^{26}=-f_4^0+f_5^1-5f_6^2+f_7^3-f_0^4-f_1^5-f_2^6-f_3^7\\12f^{37}=-f_4^0+f_5^1+f_6^2-5f_7^3-f_0^4-f_1^5-f_2^6-f_3^7\\12f^{40}=-f_4^0+f_5^1+f_6^2+f_7^3+5f_0^4-f_1^5-f_2^6-f_3^7\\-12f^{51}=-f_4^0+f_5^1+f_6^2+f_7^3-f_0^4+5f_1^5-f_2^6-f_3^7\\-12f^{62}=-f_4^0+f_5^1+f_6^2+f_7^3-f_0^4-f_1^5+5f_2^6-f_3^7\\-12f^{73}=-f_4^0+f_5^1+f_6^2+f_7^3-f_0^4-f_1^5-f_2^6+5f_3^7\end{cases}$$

$$(4.5.95)\begin{cases}-12f^{05}=5f_5^0-f_4^1+f_7^2-f_6^3+f_1^4-f_0^5+f_3^6-f_2^7\\-12f^{14}=-f_5^0+5f_4^1+f_7^2-f_6^3+f_1^4-f_0^5+f_3^6-f_2^7\\12f^{27}=-f_5^0-f_4^1-5f_7^2-f_6^3+f_1^4-f_0^5+f_3^6-f_2^7\\-12f^{36}=-f_5^0-f_4^1+f_7^2+5f_6^3+f_1^4-f_0^5+f_3^6-f_2^7\\12f^{41}=-f_5^0-f_4^1+f_7^2-f_6^3-5f_1^4-f_0^5+f_3^6-f_2^7\\12f^{50}=-f_5^0-f_4^1+f_7^2-f_6^3+f_1^4+5f_0^5+f_3^6-f_2^7\\12f^{63}=-f_5^0-f_4^1+f_7^2-f_6^3+f_1^4-f_0^5-5f_3^6-f_2^7\\12f^{72}=-f_5^0-f_4^1+f_7^2-f_6^3+f_1^4-f_0^5+f_3^6+5f_2^7\end{cases}$$

$$(4.5.96)\begin{cases}-12f^{06}=5f_6^0-f_7^1-f_4^2+f_5^3+f_2^4-f_3^5-f_0^6+f_1^7\\-12f^{17}=-f_6^0+5f_7^1-f_4^2+f_5^3+f_2^4-f_3^5-f_0^6+f_1^7\\-12f^{24}=-f_6^0-f_7^1+5f_4^2+f_5^3+f_2^4-f_3^5-f_0^6+f_1^7\\12f^{35}=-f_6^0-f_7^1-f_4^2-5f_5^3+f_2^4-f_3^5-f_0^6+f_1^7\\12f^{42}=-f_6^0-f_7^1-f_4^2+f_5^3-5f_2^4-f_3^5-f_0^6+f_1^7\\-12f^{53}=-f_6^0-f_7^1-f_4^2+f_5^3+f_2^4+5f_3^5-f_0^6+f_1^7\\12f^{60}=-f_6^0-f_7^1-f_4^2+f_5^3+f_2^4-f_3^5+5f_0^6+f_1^7\\12f^{71}=-f_6^0-f_7^1-f_4^2+f_5^3+ff_2^4-f_3^5-f_0^6-5f_1^7\end{cases}$$

$$(4.5.97) \begin{cases} -12f^{07} = 5f_7^0 + f_6^1 - f_5^2 - f_4^3 + f_3^4 + f_2^5 - f_1^6 - f_0^7 \\ 12f^{16} = -f_7^0 - 5f_6^1 - f_5^2 - f_4^3 + f_3^4 + f_2^5 - f_1^6 - f_0^7 \\ -12f^{25} = -f_7^0 + f_6^1 + 5f_5^2 - f_4^3 + f_3^4 + f_2^5 - f_1^6 - f_0^7 \\ -12f^{34} = -f_7^0 + f_6^1 - f_5^2 + 5f_4^3 + f_3^4 + f_2^5 - f_1^6 - f_0^7 \\ 12f^{43} = -f_7^0 + f_6^1 - f_5^2 - f_4^3 - 5f_3^4 + f_2^5 - f_1^6 - f_0^7 \\ 12f^{52} = -f_7^0 + f_6^1 - f_5^2 - f_4^3 + f_3^4 - 5f_2^5 - f_1^6 - f_0^7 \\ -12f^{61} = -f_7^0 + f_6^1 - f_5^2 - f_4^3 + f_3^4 + f_2^5 + 5f_1^6 - f_0^7 \\ 12f^{70} - f_7^0 + f_6^1 - f_5^2 - f_4^3 + f_3^4 + f_2^5 - f_1^6 + 5f_0^7 \end{cases}$$

Доказательство. Системы линейных уравнений (4.5.90), (4.5.91), (4.5.92), (4.5.93), (4.5.94), (4.5.95), (4.5.96), (4.5.97) получены в результате перемножения матриц в равенстве (4.5.74). □

Для того, чтобы найти линейное отображение, соответствующее операции сопряжения, я положу

$$(4.5.98) \qquad f_0^0 = 1 \quad f_1^1 = f_2^2 = f_3^3 = f_4^4 = f_5^5 = f_6^6 = f_7^7 = -1$$

Подставив (4.5.98) в систему уравнений (4.5.90), мы получим

$$(4.5.99) \qquad f^{00} = f^{11} = f^{22} = f^{33} = f^{44} = f^{55} = f^{66} = f^{77} = -\frac{1}{6}$$

Следовательно,

$$(4.5.100) \quad \overline{z} = -\frac{1}{6}(z + (iz)i + (jz)j + (kz)k + ((il)z)(il) + ((jl)z)(jl) + ((kl)z)(kl))$$

Глава 5

# Список литературы

# Глава 6

# Предметный указатель

# Глава 7

# Специальные символы и обозначения